\documentclass{article}

\title{Mathematical Proofs of Two Conjectures:\\ The Four Color Problem and \\ The Uniquely 4-colorable Planar Graph}
\author{Jin Xu}
\date{Version 2.0, 2012/01/15}

\usepackage{amsmath,amsthm}
\usepackage{graphicx}
\usepackage{psfrag}
\usepackage{float}
\usepackage{amssymb}
\DeclareGraphicsExtensions{.pdf}

\chardef\bslash=`\\ 

 \newtheorem{theorem}{Conjecture}[section]
 \newtheorem{theorem2}{Theorem}[section]
 \newtheorem{corollary}[theorem2]{Corollary}
 \newtheorem{lemma}{Lemma}[section]
 \newtheorem{Prop}[theorem2]{Proposition}
 \newtheorem{problem}[theorem]{Problem}

\theoremstyle{definition}
\newtheorem{definition}{Definition}[section]

\theoremstyle{remark}

\newcommand{\eval}[2][\right]{\relax
  \ifx#1\right\relax \left.\fi#2#1\rvert}

\begin{document}
\maketitle \markboth{The Mathematical Proofs of Two Conjectures}
{The Mathematical Proofs of Two Conjectures}
\renewcommand{\sectionmark}[1]{}

 \begin{abstract}

The Four-Color Conjecture says that every planar graph is
4-colorable. One method essential to attacking this conjecture is
through finding reducible unavoidable sets, which  goes back to
Kempe's `proof' in 1879. In accordance with this idea, Appel and
Haken found a unavoidable set containing 1936 reducible
configurations by the help of computer programs, so that the
Four-Color Conjecture was proved for the first time. In 1997, this
proof was simplified by Robertson, Sanders, Seymour and Thomas. They
found a unavoidable set containing 633 reducible configurations.
However, this proof is still a computer-assisted proof, and up until
now, there is not a proof of the Four-Color Theorem that can be
completed by hand. We shall write a series of articles to try to
prove the Four-Color Conjecture using different mathematical
methods. This paper is the first in the series, in which we
introduced the so-called color-coordinate system theory. Taking
advantage of this theory, we prove not only the Four-Color
Conjecture, but also the uniquely 4-colorable planar graph
conjecture, both by mathematical method. Our work is organized into
four parts: the first part sets up the recursion formula for
contracting vertices of maximal planar graphs, which is the basis of
proving these two conjectures; the second part depicts the basic
structure of maximal planar graphs and
 builds the generating operation system for maximal planar graphs, from which we can know that one maximal planar graph can be constructed from another lower order maximal planar graph by extending wheel operation; the third part sets up the color-coordinate system theory of graphs, mainly for the maximal planar graphs. According to this theory, a $k$-chromatic graph $G$ is a $k$-colorable-coordinate if and only if $G$ is a uniquely $k$-chromatic graph, or a quasi-uniquely $k$-chromatic graph, or a pseudo-uniquely $k$-chromatic graph. Furthermore, we depict the basic characteristic of uniquely 4-chromatic maximal planar graphs, quasi-uniquely 4-chromatic maximal planar graphs and pseudo-uniquely  4-chromatic maximal planar graphs, respectively. On the basis of these works, we prove that a 4-chromatic maximal planar graph $G$ is uniquely 4-colorable if and only if $G$ is a recursive maximal planar graph, and thus prove
 the Frioini-Wilson-Fisk Conjecture and the Jensen-Toft Conjecture on the
 uniquely 4-colorable planar graphs. In addition, we give a necessary and sufficient condition for the quasi-uniquely 4-chromatic maximal planar graphs and obtain some basic properties of the pseudo-uniquely  4-chromatic maximal planar graphs. Based on the above three parts, the last part proves that for any maximal planar graph $G$, its chromatic polynomial $f(G,t)$ has
 $f(G,4)>0$ for $t\geq 4$.

 \end{abstract}

 \begin{center}
 \textbf{Contents}
 \vspace{3mm}
 \end{center}

 1.\quad Introduction

 2.\quad Notations

 \quad \quad  2.1.  Basic notations and definitions

 \quad \quad  2.2.  Graph coloring

 \quad \quad  2.3.  Maximal planar graphs

 3.\quad Chromatic polynomials of graphs

 \quad \quad  3.1.  Introduction

 \quad \quad  3.2.  Some related results on chromatic polynomial

 \quad \quad  3.3.  Chromatic polynomial of maximal planar graphs

4.\quad  Generating operation system of maximal planar graphs

 \quad \quad  4.1.  Chromatic isomorphism of graphs
 
 \quad \quad  4.2.  Basic generating operation system of maximal planar graphs

 \quad \quad  4.3.  Construction of maximal planar graphs

 \quad \quad  4.4.  Extending and contracting $k$-wheel operations based on

 \quad \quad \quad \quad coloring

 5.\quad Recursive maximal planar graphs

 \quad \quad 5.1. Basic properties

 \quad \quad 5.2. (2,2)-FWF graphs

 \quad \quad 5.3. Color sequence of a (2,2)-FWF graph

 \quad \quad 5.4. Chromaticity of induced graphs by extending 4-wheel operation

6.\quad Chromatic isomorphism of maximal planar graphs

 \quad \quad 6.1. Cycle-coloring and tree-coloring

 \quad \quad 6.2. Equivalency of coloring  between tricolored induced subgraphs and maximal planar 
 
 \quad \quad \quad \quad graphs

 \quad \quad 6.3. Union structure of two bicolored induced subgraphs

 \quad \quad 6.4. Construction of tricolored induced subgraphs

7.\quad Black-White coloring, and necessary and sufficient conditions for 2-colorable cycle

 \quad \quad 7.1. Characteristics and distribution of the even-cycles

 \quad \quad 7.2. Enumeration of even-cycles
 
 \quad \quad 7.3. Black-White coloring operation
 
 \quad \quad 7.4. The necessary and sufficient condition of the 2-colorable cycle based on \quad \quad \quad \quad petal-syndrome
 
 \quad \quad 7.5. Necessary and sufficient conditions of 2-colorable cycles based on structure
 
 \quad \quad 7.6. Construction of semi-maximal planar graphs with 2-colorable cycles
 
 \quad \quad 7.7. Summary

8.\quad Coloring operation system for maximal planar graphs(I)--protected-cycle operation

 \quad \quad 8.1. Definitions and properties

 \quad \quad 8.2. Eigen graph of protected-cycle coloring

 \quad \quad 8.3. Properties of protected-cycle operation














 Acknowledgements

 References

 Appendix

\section{Introduction}
In 1852, Francis Guthrie \cite{F1880} put forward the Four-Color
Conjecture \cite{O1967,H1969,FF1998,BLW1976,BW1978,TK1986,CZ2005}:
four colors are sufficient to color any maps drawn in the plane or a
sphere so that no two regions with a common boundary line are
colored with the same color, where two regions which have a point or
a finite number of points in common are permitted to have the same
color. It is also required that all countries should be connected
together. Four-Color Conjecture can be converted into
vertex-coloring problem of a planar graph. The specific method is
that every country in a map is considered as a vertex, and if two
countries in the map share a common boundary line, then we will
connect those two vertices representing two countries with a line.
In this way, Four-Color Map Problem can be equivalently converted
into $4$-vertex-coloring problem of a planar graph as follow:

\begin{theorem}[The Four-Color Conjecture]\label{th1}
Any planar graph is 4-colorable.
\end{theorem}

The question was raised to the London Mathematical Society in 1878
by Cayley, the most famous British mathematician at that time
\cite{C1878,C1879}. From then on, it became the question in the
spotlight in the mathematical world. Then, Kempe
\cite{K1879,K1879(2)} and Tait \cite{T1880} respectively showed
their two different ``proofs'' of the Four-Color Conjecture in 1879
and 1880. Although the two proofs were both incorrect, their works
had greatly promoted the development of Graph Theory. Particularly,
Kempe's work had played a fundamental role in computer-assisted
proof \cite{AH1977,AH1977(2),AHK1977} in 1976.

The first concept introduced in Kempe's proof was
``\textbf{configuration}'' \cite{K1879,K1879(2)}. He proved that in
each normal map, there was at least one country being adjacent to 2,
3, 4 or 5 countries. A normal map with each country being adjacent
to six or more countries included didn't exist. In other words, a
set of ``configuration'' consisting of two adjacent countries, three
adjacent countries, four adjacent countries, five adjacent countries
is unavoidable, so every map at least includes one of those four
configurations. Another concept of ``reducibility'' was indirectly
introduced by kempe in his proof \cite{K1879,K1879(2)}. He proved
that if one country is adjacent to four countries in a colored map
with 5 colors, then the number of countries in this map can be
reduced. Since the concepts of ``configuration'' and
``reducibility'' were introduced, several standard methods by
checking a configuration to determine whether it was reducible or
was not gradually developing. To seek unavoidable sets of reducible
configurations is the important method in proving Four-Color
Problem. However, to prove whether a larger configuration is
reducible, a lot of details is needed to be checked and this is very
complicated. It was on the basis of those concepts and methods that
the computer-assisted proof of Four-Color Conjecture was gradually
developing.

In 1890, Heawood \cite{H1890} pointed out a fatal error in Kempe's
solution of Four-Color Problem by constructing a counterexample. And
then, using Kempe's method, Heawood proved Five-Color Theorem: Any
planar graph can be properly colored with five colors. Kempe
admitted the error indicated by Heawood. Meanwhile, he claimed he
was not able to correct it. However, Heawood deeply devoted himself
to the research of the Four-Color Problem later in the 60 years.

Another incorrect proof of Four-Color problem \cite{T1880} was given
by Tait in 1880. His proof was based on the following assumption:
Each 3-connected cubic planar graph was hamiltonian. After 11 years,
Petersen \cite{P1898} pointed out the assumption of Tait's proof was
incorrect. A 3-edge-coloring of a 3-regular graph is called a Tait
coloring. Tait proved that every 3-regular hamiltonian graph
admitted a Tait coloring, then he believed that he had done the
proof of the Four-Color Problem. Although the error in his proof was
found by Petersen \cite{P1898}, the counterexample was not given
until 1946 \cite{T1946}. Then, in 1968, Grinberg obtained a
necessary condition for producing many 3-connected non-Hamilton
cubic planar graphs. Although the proof of Tait was incorrect, his
work had had a strong influence on the research of graph theory,
especially edge-coloring theory.

In the early 20th century, the process of solving the Four-Color
Problem seemed to be stagnant. Basically, the proof of Four-Color
Conjecture still followed Kempe's theory was continually penetrate,
detail and refine. The basic method of this proof was so-called
smallest counterexample methods and its basic thought was to seek
\textbf{unavoidable sets of reducible configurations}.

The research on unavoidable sets originated from Wernicke's work
\cite{W1904} in 1904. The concept of reducibility was introduced by
Birkhoff \cite{B1913} in 1913. Combining Kempe's thought with his
new idea, Birkhoff proved some larger configurations to be
reducible. In 1922, Franklin proved that maps included at most 25
countries could be colored with four colors. The number was
increased to 27 by Reynolds in 1926 \cite{R1926}, to 32 by Franklin
in 1937 \cite{F1938}, to 35 by Winn in 1940 \cite{W1940}, to 39 by
Ore and Stemple in 1970 \cite{OS1970}, and finally to 95 by Mayer in
1976 \cite{M1978}. Obviously, this kind of promotion was very slow
and make little sense for the final solution to Four-Color Problem.

The planar graph is called a maximal planar graph if  each of its
faces is a triangle. A \textbf{configuration} of a maximal planar
graph is defined as a boundary cycle and the part inside the cycle.
In the Four-Color Problem, a set $\mathcal{U}$ of configurations is
\textbf{unavoidable} if every maximal planar graph necessarily
contains at least one member of $\mathcal{U}$. A configuration is
\textbf{reducible} if it cannot be a configuration of a smallest
counterexample to the Four-Color Conjecture.  On the research of
finding unavoidable sets of reducible configurations, the great
contribution was made by German mathematician Heesch \cite{H1969},
who laid the foundation for the final solution to Four-Color Problem
with electronic computer in 1976\cite{AH1977,AH1977(2),AHK1977}.
Heesch put forward a more systematic method. He firmly believed that
this method could effectively solve the Four-Color Problem. He
published it at a seminar in Hamburg University and Keele
University, which Haken attended. Heesch estimated that it might
contain about ten thousands unavoidable reducible configurations and
it was impossible to check the reducibility of so many
configurations by hand. He proposed a new method called
``discharging'' to prove the reducibility, which was a big step up
in the research on Four-Color Problem. From 1960s to 1976, the
research focused on how to find unavoidable sets of reducible
configurations by means of electronic computer. The main
contributors in the 1960s were Heesch, Haken, Durre, Shimamoto, etc.
And in the 1970s, Heesch, Allaire, Swart, Bernhart, Haken,Tutte etc.
Among them, Shimamoto was frustrated by checking that one of his
configurations was not a so-called D-reducible with computer
\cite{CZ2005}.

The algorithm for checking ``reducible obstacle'' was proposed by
Appel and Haken, which could greatly reduce computing time. With the
help of John Koch, they successfully constructed the unavoidable set
of 1936 reducible configurations (later reduced to 1476 kinds) in
June, 1976 \cite{AH1977,AH1977(2),AHK1977}. Each configuration
needed checking with electronic computer in turn. Different programs
and computers were used to independently reexamine this work. Using
three computers, spending 1200 hours, making ten billion judgements,
they definitely proved Four-Color Conjecture by computer-assisted
methods.

In 1997, Four-Color Theorem was proved with the method similar to
Appel and Haken's by Robertson, Sanders, Seymour, Thomas, etc. They
simplified the proof \cite{RSST1996,RSST1997}, in which only 633
configurations in the unavoidable sets were needed to be checked by
computer. Additionally, a new computer-assisted proof \cite{G2005}
was given by Georges Gonthier by a computer, called \emph{Coq}, that
equipped the mathematical software in 2005.

However, the above proofs were all to depend on computers and, it is
hard to be checked one by one by hand. Therefore, a computer-free
method to concisely solve the Four-Color Problem is still concerned
by the whole mathematical community.

The \textbf{vertex-coloring problem} of a graph is to classify the
vertices of a graph. There are two constraints: One is the adjacent
vertices need to be colored by different colors; the other one is
the chromatic number, namely, the number of classes in the vertex
partition or the number of color classes. Obviously, each color
class corresponds to a vertex independent set of the graph.

In this paper, let $t$ be the number of colors required to properly
color the vertices of a graph $G$, $\chi(G)$ be the chromatic number
of the graph $G$, $f(G,t)$ be the number of all possible colorings
to the vertices of the labeled graph $G$ with $t$ colors. Obviously,
if $t<\chi(G)$, the graph $G$ can not be properly colored, so
$f(G,t)=0$. However, if $\chi(G)\leq t$, this coloring must exist,
so that $f(G,t)>0$. For every planar graph $G$, if $f(G,4)>0$ can be
proved, it is equivalent to the proof of the Four-Color Problem!
This is the method that Birkhoff had proposed for Four-Color Problem
in 1912. Later on, it was found that $f(G,t)$ is a polynomial in
terms of the color number $t$, called the \textbf{chromatic
polynomial} of the graph, which has become a fascinating branch of
graph theory at present \cite{DKT2005}. But it was a pity that
Birkhoff's aim had not been reached. The best result given by Tutte
\cite{T1970} was that if $t=$ $\tau \sqrt{5}=3.618\cdots$ and
$\tau=\frac{1}{2}(1+\sqrt{5})$, then $f(G, \tau\sqrt{5})>0$. The
result seemed to be a pity that it brushed past the Four-Color
Problem, for all we need to do is prove $f(G,4)>0$, and then
Four-Color Conjecture holds true.

The first important tool in this paper is \textbf{chromatic
polynomials of graphs}, which is one of the key tools to prove the
Four-Color Conjecture. Some related researches will be introduced
and discussed in Section 3 in order to make a systematic
description.

The second important tool used in this paper is the contracting and
extending operations on maximal planar graphs, from which we will
establish the \textbf{generating operation system} of maximal planar
graphs, find out the method of constructing maximal planar graphs
and, research the operations of contracting and extending $k$-wheel
on coloring. This is the fundamental work and will be discussed in
the fourth section.

In the process of accomplishing the proof of FWF Conjecture and
Jensen-Toft's Conjecture (JT Conjecture, \cite{JT1995}), namely
obtaining the necessary and sufficient conditions of uniquely
4-chromatic maximal planar graphs, we find this problem is closely
relative to two kinds of maximal planar graphs. One is the so-called
recursive maximal planar graph, essentially the $(2,2)$-recursive
maximal planar graph. The other is the $1$-$m$ maximal planar graph.
The related coloring properties of the two kinds of graphs are
researched in this paper, especially the problem of chromaticity of
the graphs gained from operating the contracting and extending
operations on them  will be discussed in Sections 5 and 6.

Based on Sections 3, 4, 5 and 6, the JT Conjecture, namely FWF
Conjecture has been definitely proved in Section 7, and we obtained
the necessary and sufficient conditions for a maximal planar graph
being uniquely 4-colorable is that the graph is a recursive maximal
planar graph. We will give two different methods to prove these two
equivalent conjectures.

The third important tool is the \textbf{color-coordinate system
theory} of graphs, which points out graph $G$ is
$k$-colorable-coordinate if and only if $G$ is uniquely
$k$-colorable, or quasi-uniquely $k$-colorable, or pseudo-uniquely
$k$-colorable. We will investigate the related basic properties of
uniquely 4-colorable maximal planar graphs, or quasi-uniquely
4-colorable maximal planar graphs, or pseudo-uniquely 4-colorable
maximal planar graphs, which will be discussed in Sections 7, 8and
9. All above are the foundation of proving the Four-Color
Conjecture. Based on the results of these three sections, the
Four-Color Conjecture will be proved in Section 10.

In this paper, we give a new idea to prove the Four-Color
Conjecture, which combines the contracting and extending operations
of maximal planar graphs, the color-coordinate system theory with
chromatic polynomials of graphs together. The basic method is the
induction. First, for a maximal planar graph, we deduce its
chromatic polynomial on contracting 4-wheel and 5-wheel operations,
respectively. Second, we introduce the concepts of quasi-uniquely
colorable planar graphs and pseudo-uniquely colorable planar graphs
and then figure their basic characteristics comprehensively. In
addition, the methods for judging quasi-uniquely 4-colorable and
pseudo-uniquely 4-colorable maximal planar graphs are researched. In
order to solve the conjecture of uniquely 4-colorable planar graph,
the contracting and extending operations of maximal planar graphs
are introduced. The basic characteristics of uniquely 4-colorable
planar graphs are studied in detail. The necessary and sufficient
conditions of a uniquely 4-colorable planar graph is worked out.
Naturally, the conjecture of Frioini-Wilson-Fisk, which is also
called the open conjecture of uniquely 4-colorable planar graphs
proposed by Jensen-Toft is solved. On the basis of these results, we
prove: For a maximal planar graph $G$, its chromatic polynomial
$f(G,4)>0$. In this way, the Four-Color Conjecture is proved by the
mathematical method.


\section{Definitions and Notations}
This section gives some basic terminologies and notations that are
used in this paper. Other special definitions will be defined in the
corresponding chapters.

\subsection{Basic notation}

All graphs in this paper are restricted to be finite, simple and
undirected, except the graphs that contain the 2-wheel subgraphs. In
a given graph $G$, $V(G), E(G), d_{G}(u)$ and $\Gamma_{G}(u)$ denote
the \textbf{vertex set}, the \textbf{edge set}, the vertex $u$'s
\textbf{degree} and the \textbf{set of neighbors} of $u$
respectively, which are written as $V,E,d(u)$ and $\Gamma(u)$ for
short. Very often, we call the vertex $u$ a \textbf{$k$-degree
vertex} if $d_G(u)=k$. The number $|V(G)|$ of vertices is called
\textbf{order} of $G$ and the number $|E(G)|$ of edges is called
\textbf{size} of $G$. An \textbf{independent set} in a graph is a
set of vertices no two of which are adjacent. For graph
$H=(V(H),E(H))$, if $V(H) \subseteq V(G), E(H)\subseteq E(G)$, then
$H$ is called a \textbf{subgraph} of $G$. And whenever $u, v \in
V(H)$ are adjacent in the graph $G$, they are also adjacent in the
graph $H$, then $H$ is called s \textbf{induced subgraph} of $G$. An
induced subgraph of $G$ with vertex set $V'$ is denoted by $G[V']$.
Let $u$, $v$ be two different vertices in $V(G)$, the
\textbf{distance} between $u$ and $v$ is the length of the shortest
path from $u$ to $v$, denoted as $d_{G}(u,v)$. Two graphs $G$ and
$H$ are disjoint if they have no vertices in common. By starting
with two disjoint graphs $G$ and $H$, by adding edges joining every
vertex of $G$ to every vertex of $H$, one obtains the \textbf{join}
of $G$ and $H$, denoted as $G\vee H$. The join $C_{n} \vee K_{1}$ of
a cycle and a single vertex is referred to as a \textbf{wheel} with
$n$ spokes, denoted as $W_{n}$ ( four wheels
$W_{2},W_{3},W_{4},W_{5}$ are shown in Figure 2.1), where $C_n$ is
called the \textbf{cycle} of this wheel and the vertex of $K_1$ is
called the \textbf{center} of the wheel. A graph is
\emph{k}-\textbf{regular} if all of its vertices have the same
degree \emph{k}. A 3-regular graph is usually called a \textbf{cubic
graph}.

\begin{center}

\includegraphics [width=300pt]{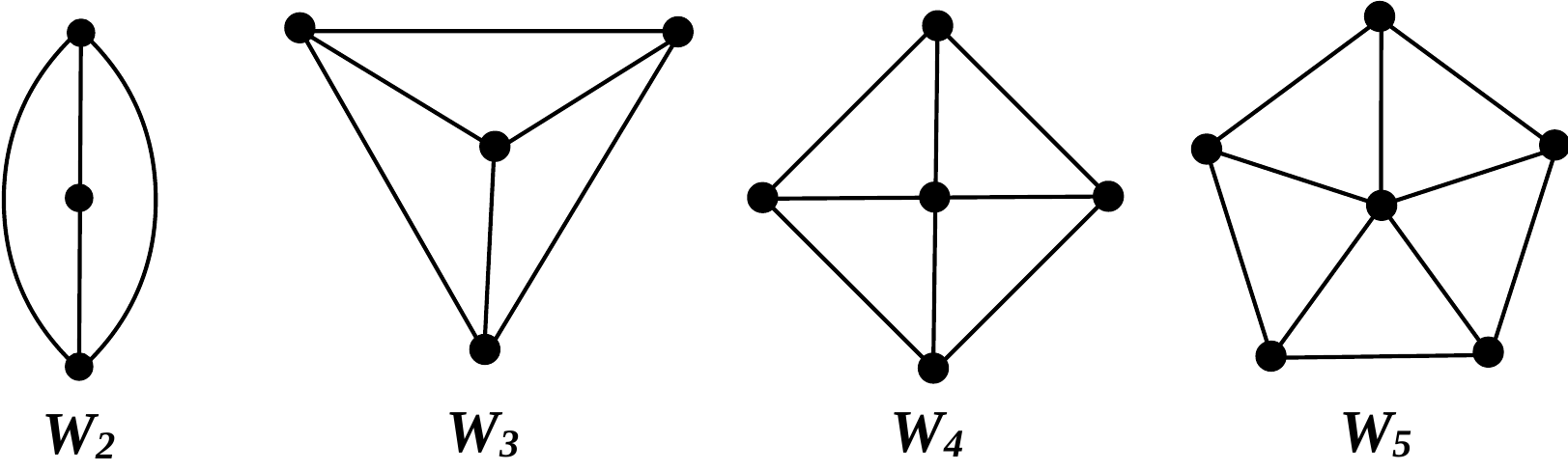}

\textbf{Figure 2.1.} Four wheels $W_{2},W_{3},W_{4},W_{5}$.
\end{center}

In order to \textbf{identify }the nonadjacent vertices $u$ and $v$
of a graph $G$, it is necessary to replace these two vertices by a
single vertex, and make it adjacent to all the edges which were
incident to either $u$ or $v$ in $G$. The resulting graph is denoted
as $G \circ \{u,v\}$. Contracting an edge $e=uv$ of a graph $G$
yields a new graph  $G \circ \{u,v\}$ by deleting the edge $e$ and
then identify its ends into one.

\subsection{Graph colorings}


A \textbf{$k$-vertex-coloring}, or simply a $k$-coloring, of a graph
\emph{G} is a mapping $f$ from the vertex set $V$ to the color sets
$C(k)=\{1,2,\ldots,k\}$ such that $f(x)\neq f(y)$ if vertex $x$ is
adjacent to vertex $y$.

A graph $G$ is \textbf{$k$-colorable} if it has a $k$-coloring. The
minimum number of $k$ colors required for which a graph $G$ is
$k$-colorable is called the \textbf{chromatic number}, denoted as
$\chi(G)$. If $\chi(G)=k$, then the graph $G$ can be colored with
$k$ colors, but not with $k-1$ colors. Alternatively, each
$k$-coloring $f$ of $G$ can be viewed as a partition
$\{V_{1},V_{2},\cdots,V_{k}\}$ of $V$, where $V_{i}$ denotes the set
of vertices assigned color $i$. So it can be written as
$f=(V_{1},V_{2},\cdots, V_{k})$. In other words, that is
$$V(G)=\bigcup_{i=1}^{k}V_{i},V_{i}\neq\emptyset,V_{i}\cap V_{j}= \emptyset, i \neq j,
 i,j=1,2,\ldots,k \eqno{(2.1)}$$
where $V_{i}$ is an independent set of $G$, $i=1,2,\ldots,k$. The
set of all $k$-colorings of a graph $G$ can be denoted by
$C_{k}(G)$. For a $k$-colorable graph $G$, the notation
$C^{0}_{k}(G)$ denotes the set consisting of all the partitions of
$k$-coloring class of $G$, simplified by the \textbf{partition set
of $k$-color class} of $G$. And define
$$\sigma_k^0(G)=|C_k^0(G)|
\eqno{(2.2)}$$

Suppose that $G$ is a $k(k\geq 3)$-chromatic graph. Let $f\in
C_k^0(G)$, and let $U=\{v_1,v_2,\cdots, v_t\}(t\leq |V(G)|)$ be a
subset of vertices of $G$. Now, define $f(U)=\{f(v_i)|i=1,2,\cdots,
t\}$, obviously, $f(U)\in C(k)$. Particularly, when $u\in V(G)$,
$f(\Gamma(u))$ denotes the set consisting of all the colors assigned
to the neighbors of the vertex $u$.

A $k$-colorable graph $G$ is called \textbf{uniquely $k$-colorable}
if each $k$-coloring of $G$ induces the same partition of $V(G)$
into $k$ independent sets.

Similarly, an \textbf{edge-coloring} of a graph
\cite{FW1977,FW1978} is an assignment from a color-set to its
edge-set such that no two distinct adjacent edges have the same
color.
An \textbf{$k$-edge-coloring} of a graph is an edge-coloring with
$k$ colors. A graph  is \textbf{$k$-edge-colorable} if it has an
$k$-edge-coloring. The chromatic index of $G$, $\chi^{\prime}(G)$,
is the minimum number of $k$ colors required for which $G$ is
$k$-edge-colorable. A graph is called \textbf{uniquely
$k$-edge-coloring} if there is a unique $k$-edge-coloring such that
any other colorings are equivalent to it. Alternatively, a graph $G$
is uniquely $k$-edge-colorable if there is exactly one partition of
the edge-set $E(G)$ into $k$ matchings.

In this paper, two isomorphic graphs $G$ and $H$ can be written as
$G \cong H$. A graph $G$ is labeled if each vertex is assigned with
a letter (or an integer). For a labeled graph $G$, two colorings are
different if there is at least one vertex receiving different
colors, and $f(G,t)$ is the number of $t$-colorings of $G$. $f(G,t)$
is called the \textbf{chromatic polynomial} of $G$ introduced by
Brikhoff for attacking the Four-Color Problem in 1912 \cite{B1912}.
More detail information can be found in
 \cite{B1946,R1968,BW1978,T1970,T1970(2),T1974,XL1995,X2004,DKT2005}.

\subsection{Maximal planar graphs}

A \textbf{maximal planar graph} is one planar graph to which no new
edges can be added without violating planarity. A
\textbf{triangulation} is a planar graph in which every face is
bounded by three edges (including its infinite face). It is not hard
to show that a maximal planar graph is equivalent to a
triangulation. Thus, we can say that each maximal planar graph is a
triangulation.

There exists a kind of uniquely 4-colorable planar graphs, called
the \textbf{recursive maximal planar graphs}, each of which can be
obtained from $K_{4}$ by embedding a 3-degree vertex in some
triangular face continuously. In this paper, $\Lambda$ denotes the
set consisting of all recursive maximal planar graphs and
$\Lambda_{n}$ the set of the graphs in $\Lambda$ with order $n$. Let
$\gamma_{n}=|\Lambda_{n}|$. Obviously,
$\gamma_{4}=\gamma_{5}=\gamma_{6}=1$, the corresponding recursive
maximal planar graphs are shown in Figure 2.2.

\begin{center}
\includegraphics [width=300pt]{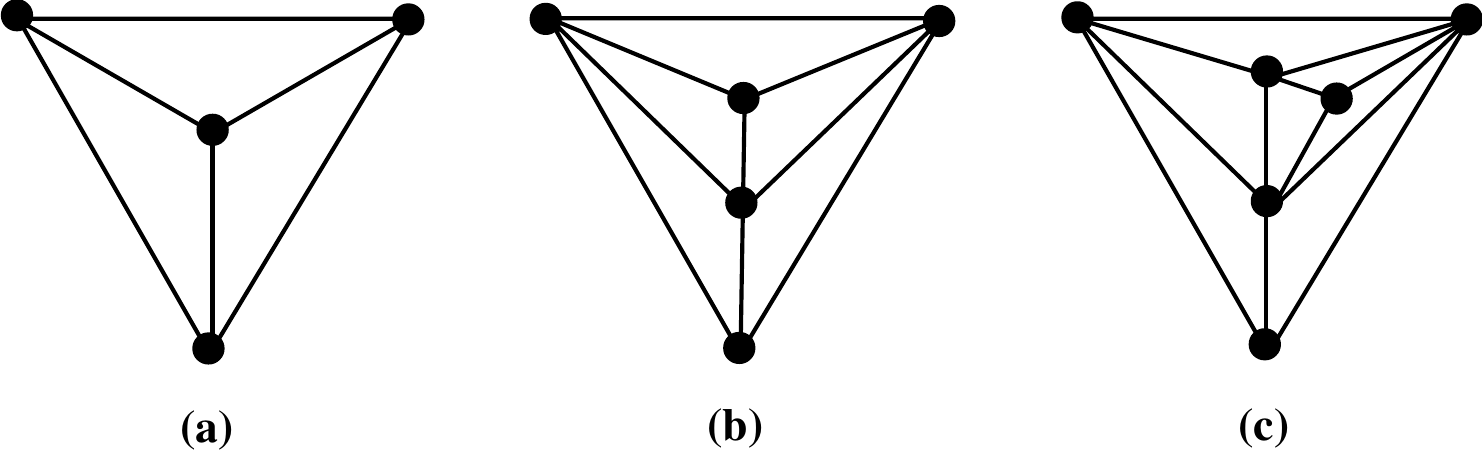}

\textbf{Figure 2.2.} Three recursive maximal planar graphs of
minimal order.
\end{center}

\noindent

In a maximal planar graph $G$, we mark a triangular face $a$-$b$-$c$
if the vertices on its boundary are labelled with $a,b,c$
respectively. The \textbf{ vertex addition} on a triangular face
$a$-$b$-$c$ is to add a new vertex $u$, and to join $u$ with $a,b,c$
in this face, denoted as $G+u$. Obviously, the resulting graph is
also a maximal planar graph. We refer to the vertex addition as the
\textbf{extending 3-wheel operation}. Another operation used in the
paper is the \textbf{vertex deletion}, which is the inverse
operation to the vertex addition. We also call the vertex deletion
the \textbf{contracting 3-wheel operation}. The operations of the
vertex addition and the vertex deletion are illustrated in Figure
2.3.

\begin{center}
\includegraphics [width=330pt]{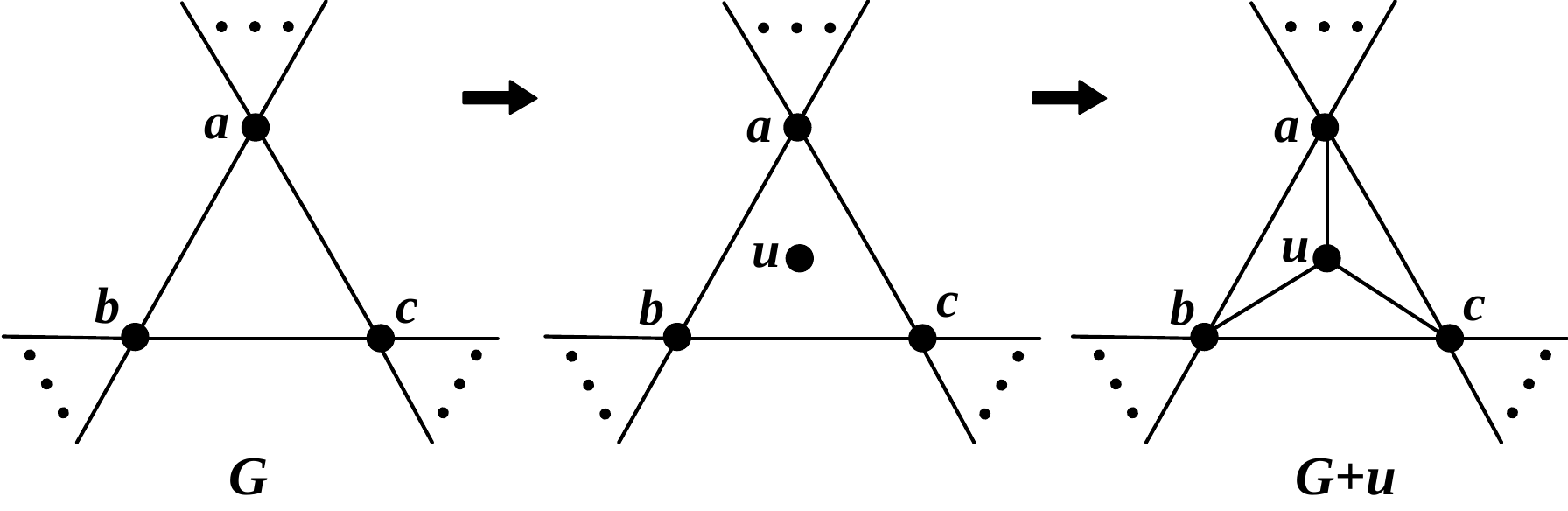}
\textbf{Figure 2.3.} Two operations of embedding and deleting a
3-degree vertex.
\end{center}

The definitions and notations not mentioned here can be found in
Boudy's book \cite{BM2008}.

\section{Chromatic polynomials of graphs}

In this section, we will introduce some correlative theories on
chromatic polynomials of graphs, and obtain two foundational
formulae of chromatic polynomials on contracting 4-wheel and 5-wheel
operations.

\subsection{Introduction}
Map coloring actually is a classification of all the countries in a
map such that no two adjacent countries are in the same class. We
can convert a map into a planar graph by dual transformation, then
the map coloring problem will be changed into the vertex-coloring
problem of a planar graph equivalently. Certainly, the
vertex-coloring problem of a graph is a type of vertex partitions,
in which adjacent vertices must receive different colors.
Accordingly, the basic scheme to attack the Four-Color Conjecture is
the partition of vertices of a graph. And this idea can be realized
by the chromatic polynomial, a much useful mathematical tool.
Although it was introduced for the labeled graphs, which would
generate a vast number of colorings, the chromatic polynomial
contains all the information on the partition of color classes.
Therefore, the chromatic polynomial may be a preferable tool to
prove the Four-Color Conjecture. Based on this tool, several
scholars had made some important contributions on attacking the
Four-Color Conjecture, mainly including Birkhoff
\cite{B1912,B1913,B1946}, Tutte \cite{T1970,T1970(2),T1974}, Read
\cite{R1968}, Whitney \cite{W1932}. Among them, the best result is
that for any planar graph $G$, the chromatic polynomial
$f(G,\tau\sqrt{5})>0$,$\tau\sqrt{5}=3.618\cdots$.

In order to prove the Four-Color Conjecture, the basic idea on the
chromatic polynomial proposed in this paper is to study the
recurrence relation for the chromatic polynomial of a maximal planar
graph directly when color number $t=4$.  Especially, the recurrence
formula when the minimal degree $\delta(G)=5$ points out the
direction for the mathematical proof of the Four-Color Problem.
These recurrence formulas are not only the foundation of the final
proof of the Four-Color Conjecture by mathematical deduction, but
also the basis of proving the uniquely 4-colorable planar graphs
conjecture. With the guide of the recurrence formulas, we need to
study the color-coordinate system theory of graphs. Specifically,
the basic characteristics of the colorable-coordinate graphs need to
be researched, since the colorable-coordinate graphs can be divided
into three types: uniquely colorable graphs, quasi-uniquely
colorable graphs and pseudo-uniquely  colorable graphs. We have to
accomplish three tasks: first, make the basic structure of the
uniquely 4-colorable graphs understood clearly; second, introduce
the concept of the quasi-uniquely $k$-colorable graphs and show the
basic characteristics of those graphs. That is to say, the necessary
and sufficient conditions will be given for quasi-uniquely
$k$-colorable graphs; third, introduce the concept of
pseudo-uniquely $k$-colorable graphs and make the basic
characteristics of those graphs clearly. However, the first task was
an unsolved conjecture proposed in 1977 \cite{FW1977,F1977}. Until
1998, a computer-assisted proof  was given by Thomas \cite{T1998}.

\subsection{Some results on chromatic polynomials of graphs}

$f(G,t)$ denotes the chromatic polynomial of a graph $G$ here. In
this subsection, the equivalent proposition of solving the
Four-Color Conjecture will be introduced again by means of the
method of chromatic polynomials.

\begin{theorem2}\label{th2}
A planar graph $G$ is 4-colorable if and only if $f(G,4)>0$.
\end{theorem2}

So far, one do not realize the dream of accomplishing Four-Color
Conjecture by chromatic polynomials purely, but the research on the
conjecture drew many scholars' interests. More detailed works on
this field can be found in the studies of Birkhoff, Lewis, Read,
Tutte and  Dong \cite{B1912,B1913,W1932,R1968,DKT2005}. In order to
calculate the chromatic polynomial of a given graph, the basic tool
is the \textbf{Deletion-Contract Edge Formula}.

For an edge $e$ of a graph $G$, the resulting graphs obtained by
deleting and contracting edge $e$ are denoted by $G-e$ and $G \circ
 e$, respectively.

\begin{lemma}\label{th3}\emph{[\textbf{The Deletion-Contract Edge Formula}]}
For a given graph $G$ and an edge $e\in E(G)$, we have
$$f(G,t)=f(G-e,t)-f(G\circ e,t)\eqno{(3.1)}$$
\end{lemma}

\begin{lemma}\label{th4}
Let $G$ be the union of two graphs $G_{1}$ and $G_{2}$, whose
intersection is a complete graph of order $k$, then
$$f(G,t)=\frac{f(G_{1},t)\times f(G_{2},t)}{t(t-1)\ldots(t-k+1)} \eqno{(3.2)}$$
\end{lemma}

Moreover, some recurrence formulas by vertex deletion \cite{X2004}
and the chromatic polynomial between graph and its complement were
given several years ago \cite{XL1995}.

Birkhoff introduced the chromatic polynomial in order to solve the
Four-Color Problem, while what we need to do is to prove that this
conjecture holds for maximal planar graphs(triangulations). So it is
meaningful to study the chromatic polynomials of  maximal planar
graphs. Another beautiful work had been made by Tutte as follows:

\begin{theorem2}\label{th5}\emph{[\textbf{the Vertex-Elimination Formula}]}
\cite{T1970} Let $G$ be a planar graph with a wheel $W_{m}$ as one
its subgraph. Then
$$f(G,\tau^{2})=(-1)^{m}\tau^{1-m}f(G-\nu,\tau^{2}) \eqno{(3.3)}$$
where $\tau^{2}=\frac{3+\sqrt{5}}{2}.$
\end{theorem2}

\begin{theorem2}\label{th6}
If  $G$ is a maximal planar graph on $n$ vertices, then
$$|f(G,t)|\leq\tau^{5-n} \eqno{(3.4)}$$
\end{theorem2}

Three operations for a plana r graph $G$ were introduced by Tutte
\cite{T1970}. Let $G$ be a planar graph with a 4-cycle $C=xyzlx$.
There does not exist any vertices or edges inside the cycle besides
a chord $e=xz$. $\theta_{e}$ means  replacing the edge $xz$ with the
edge $yl$; $\varphi_{e}$ stands for contracting the edge $xz$ to a
single vertex $x'$, and deleting duplicated edges; $\psi_{e}$,
similar as $\varphi_{e}$, indicates contracting the edge $yl$ in
$\theta_{e}(G)$. Where, we denote by $\theta_{e}(G)$,
$\varphi_{e}(G)$ and $\psi_{e}(G)$ the resulting graphs,
respectively (see Figure 3.1).

\begin{center}
\includegraphics [width=360pt]{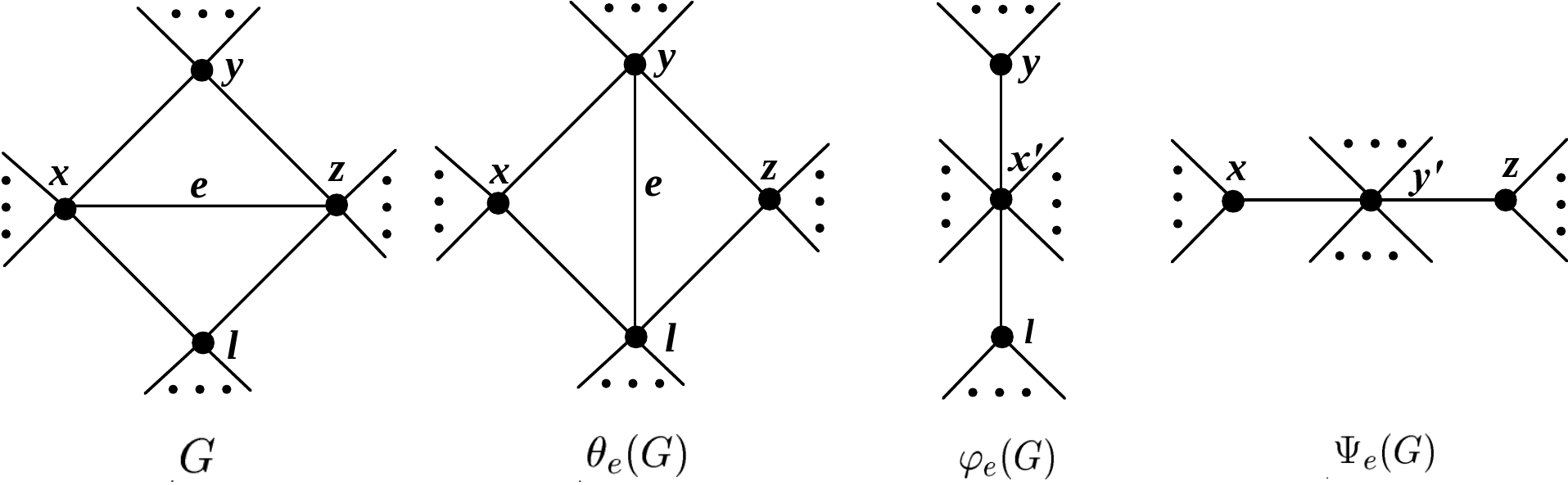}
\textbf{Figure 3.1.} Three operations introduced on a planar graph.
\end{center}

According to Lemma 3.1, an apparent result is easy to obtain as
follows:
$$f(G,t)-f(\theta_{e}(G),t)=f(\psi_{e}(G),t)-f(\varphi_{e}(G),t)\eqno{(3.5)}$$
It holds for any positive real number $t$. We have the following
result.

\begin{theorem2}\label{th7}
For $t=\tau^{2}$ we have
$$f(G,\tau^{2})+f(\theta_{e}(G),\tau^{2})=\tau^{-3}\{f(\psi_{e}(G),\tau^{2})+f(\varphi_{e}(G),\tau^{2})\}\eqno{(3.6)}$$
\end{theorem2}

Notice that $\tau=\frac{1+\sqrt{5}}{2}$, so $\tau
\sqrt{5}=\frac{5+\sqrt{5}}{2}=\tau+2 $. Then, the following theorem
can be obtained.

\begin{theorem2}\label{th8}\emph{[\textbf{Golden Identity}]}
Let $G$ be a maximal planar graph of order n, then
$$f(G,\tau\sqrt{5})=\sqrt{5}\times\tau^{3(K-3)}f^{2}(G,\tau^{2}) \eqno{(3.7)}$$
\end{theorem2}

\begin{theorem2}\label{th10}
Let $G$ be a connected graph with at least one edge. Then
$f(G,\tau+1)$ is non-zero.
\end{theorem2}

Based on  Theorems \ref{th8} and \ref{th10}, Tutte obtained an
interesting result
$$f(G,\tau\sqrt{5})>0 \eqno{(3.8)}$$

It is clear that this result is much closer to prove the Four-Color
Theorem, to which $f(G,4)>0$. However, $\tau\sqrt{5}=3.618\cdots$,
which is close to $4$ but 4. What a pity it is!

Perhaps for the perfection and excellence of Tutte's works and his
highly status in academia, once upon a time,  it was thought that to
attack the Four-Color Problem by chromatic polynomials is
impossible. Nevertheless, our works below will give a new way to
solve the Four-Color Problem by means of chromatic polynomials.

\subsection{Chromatic polynomials of maximal planar graphs}

We will give some useful results for the chromatic polynomials of
maximal planar graphs in this subsection.

\begin{theorem2}\label{th11}
\emph{(Chromatic polynomial on contracting 4-wheel operation)} Let
$G$ be a maximal planar graph, $v$ be a 4-degree vertex of $G$, and
$\Gamma(v)=\{v_{1},v_{2},v_{3},v_{4}\}$ \emph{(shown in Figure
3.2)}.

Then
$$f(G,4)=f(G_{1},4)+f(G_{2},4) \eqno{(3.9)}$$
where $G_{1}=(G-v) \circ \{v_{1},v_{3}\}$, and
$G_{2}=(G-v)\circ\{v_{2},v_{4}\}.$
\end{theorem2}

\begin{center}
\includegraphics [width=120pt]{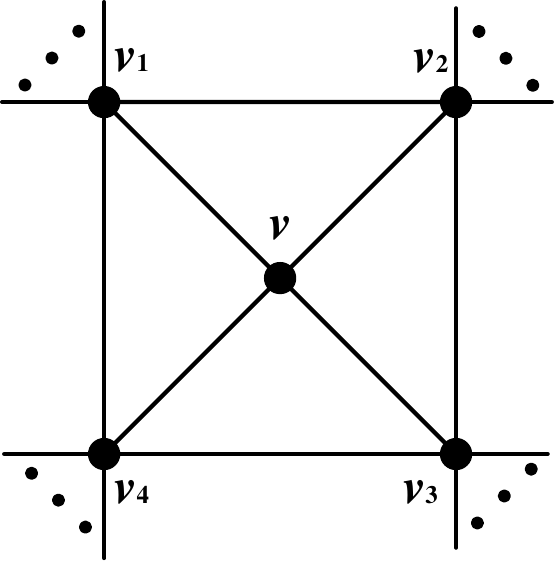}

\textbf{Figure 3.2.} A maximal planar graph having a $4$-degree
vertex $v$.
\end{center}
\begin{proof}
In the following diagram, we represent the whole $G$ by the picture
of $G[{\Gamma(v)\cup \{v\}}]$. Now we first compute the chromatic
polynomial of the graph $G$ by  Lemma \ref{th3}. For the sake of
understanding clearly, a method introduced by Zykov will be used
here \cite{Z1949}, where the chromatic polynomials are represented
by the corresponding graphical graphs without the color number $t$.
More details can be found in \cite{H1969,R1968}. Notice that if
there are at least two edges adjacent to two vertices, then only one
remains and others are deleted excluding the wheel $W_2$ on $2$
vertices.

\vspace{5mm}
\includegraphics[width=250pt]{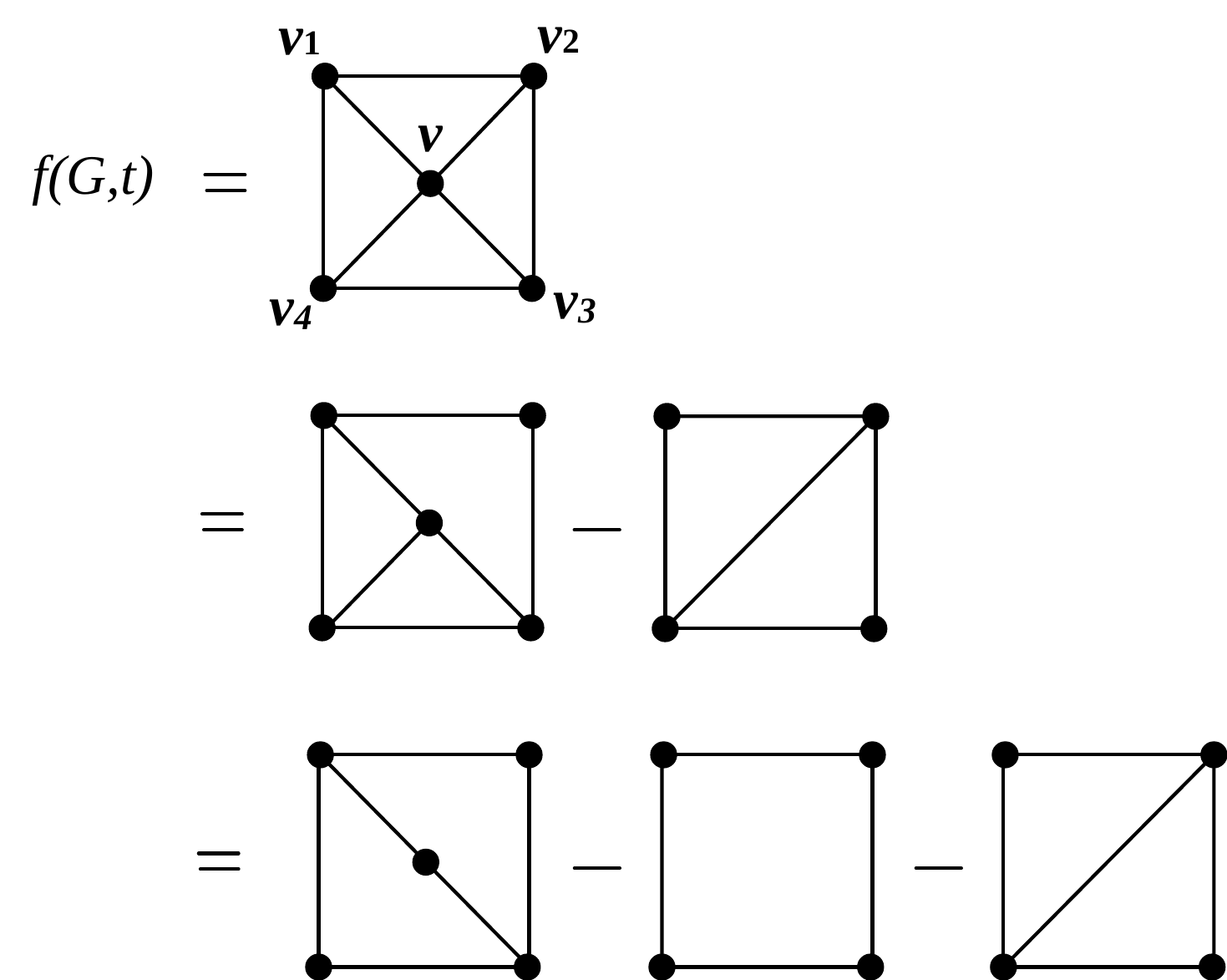}
\vspace{5mm}

\includegraphics[width=250pt]{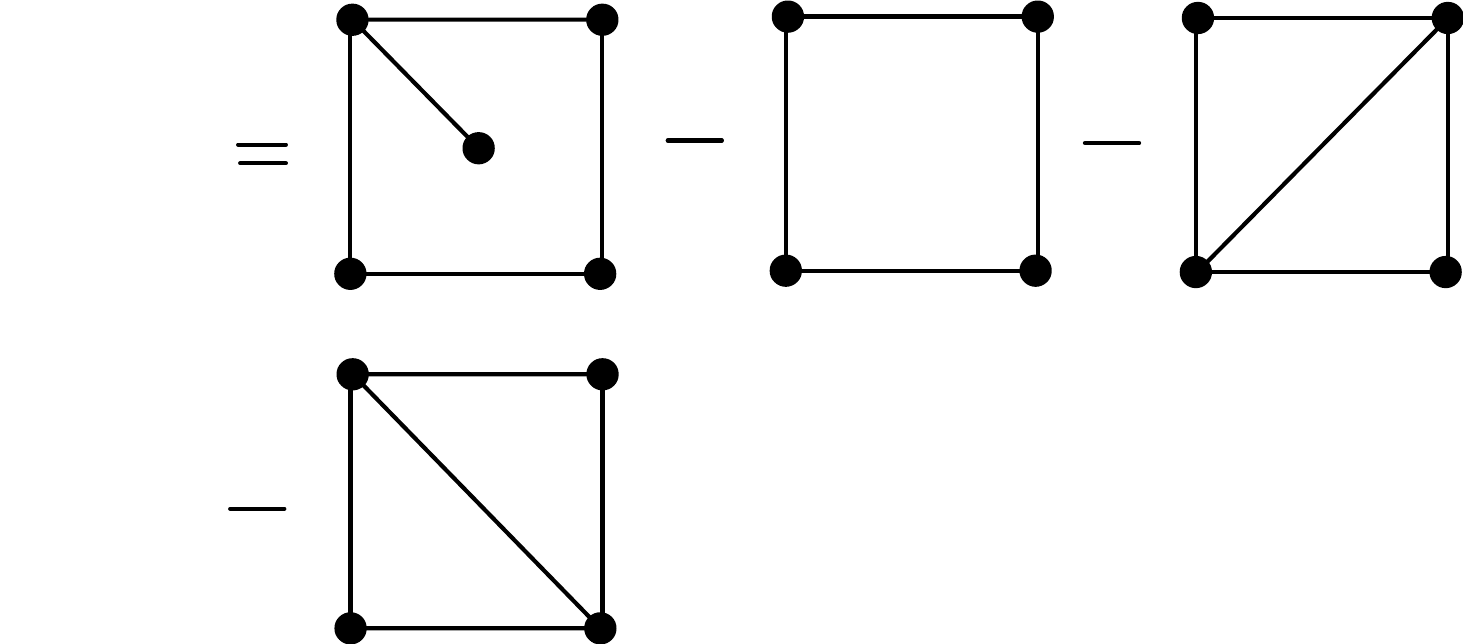}
\vspace{5mm}

\includegraphics[width=250pt]{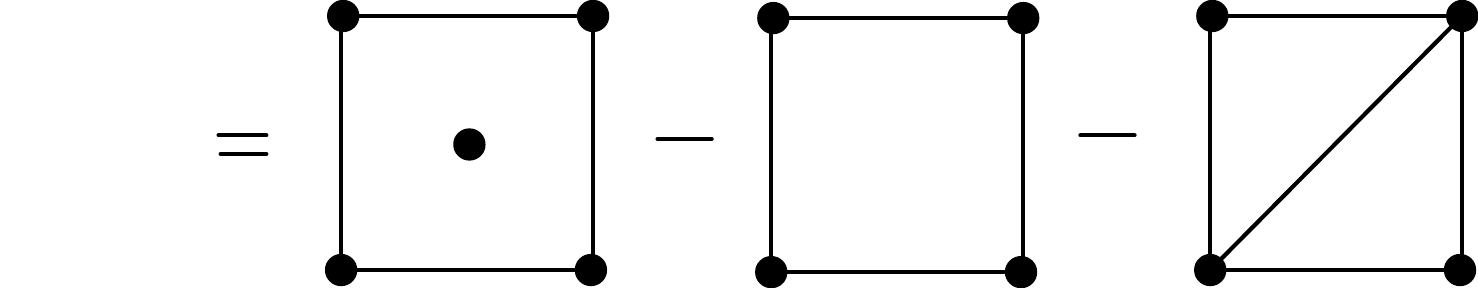}
\vspace{5mm}

\includegraphics[width=180pt]{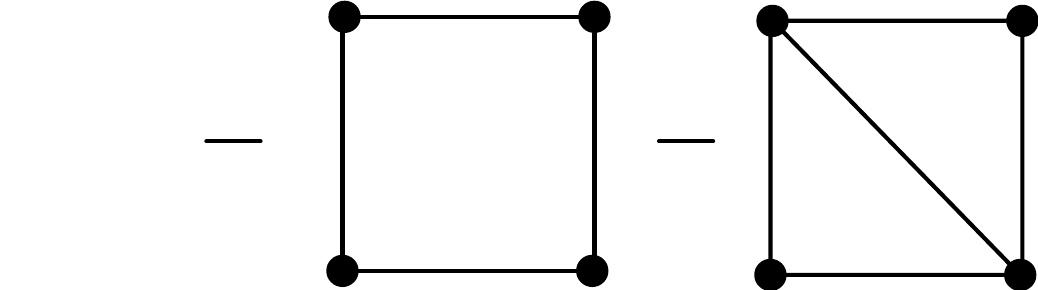}

By Lemma \ref{th4}, the chromatic polynomial of the first subgraph
is $tf(G-v,t)$. Therefore,

\vspace{1mm}
\includegraphics[width=280pt]{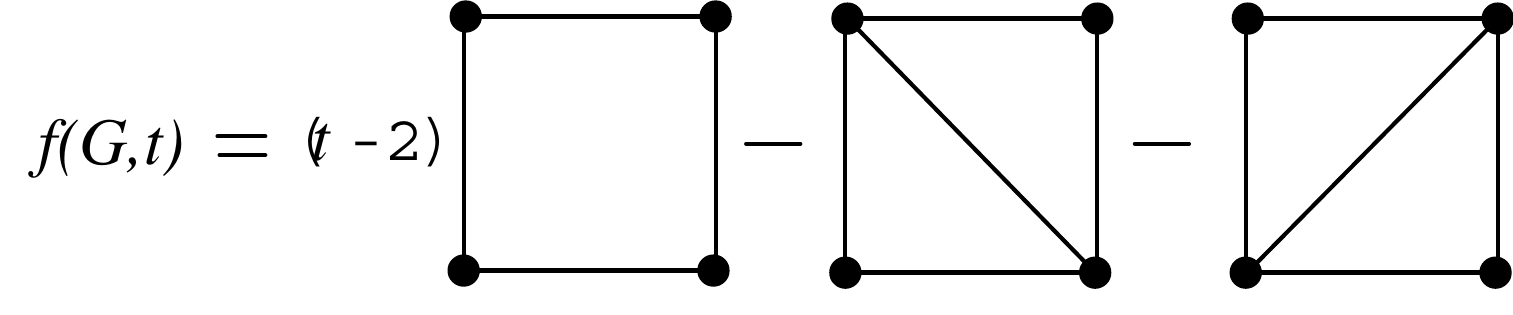}

For $t=4$, we can get

\vspace{5mm}
\includegraphics[width=200pt]{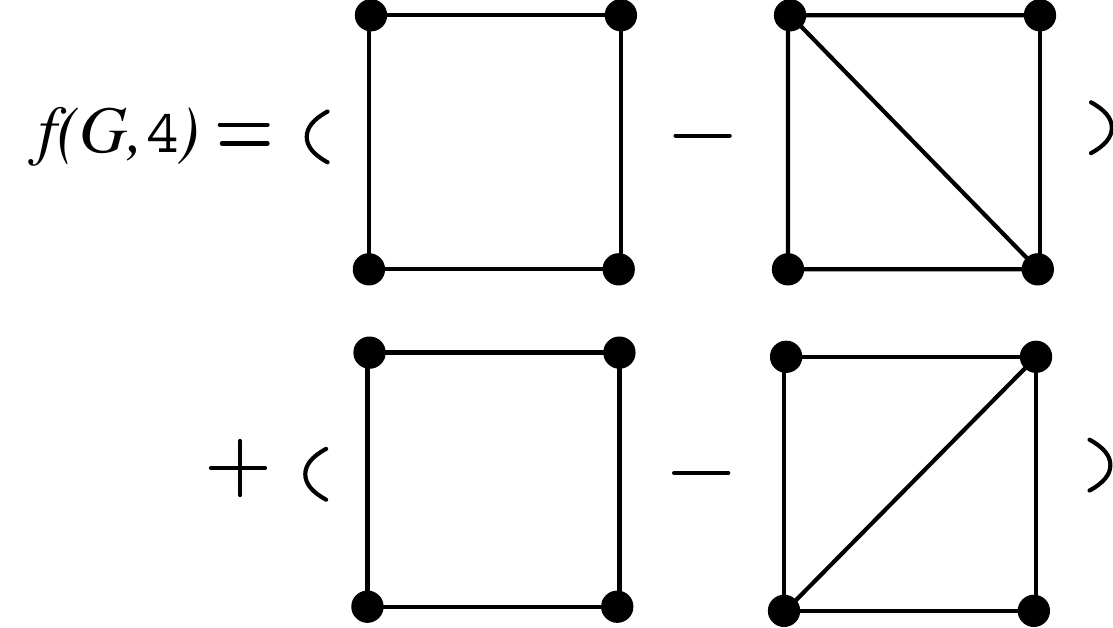}
\vspace{2mm}

\includegraphics[width=280pt]{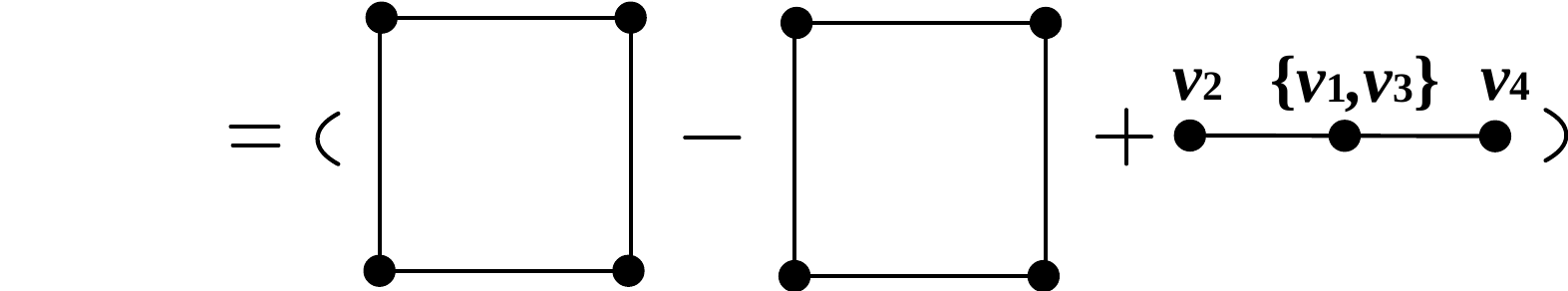}
\vspace{3mm}

\includegraphics[width=280pt]{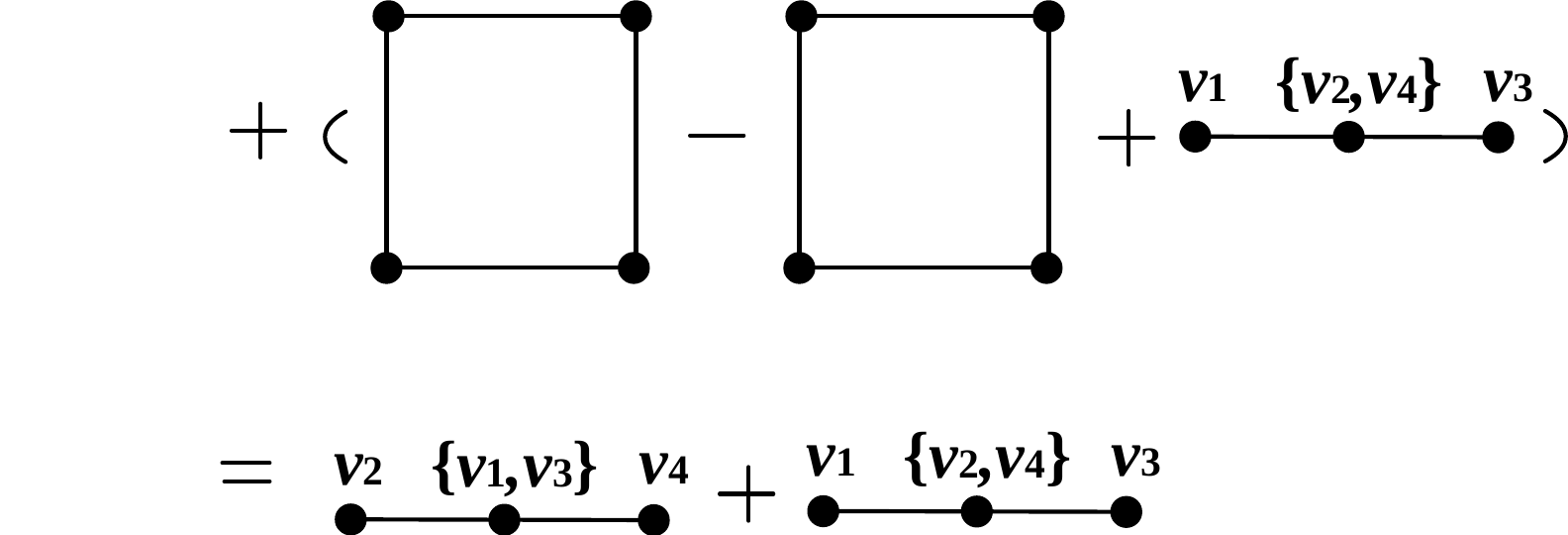}
\vspace{5mm}

Notice that two graphs at the last line above denote $(G-v)\circ
\{v_{1},v_{3}\}$ and $(G-v)\circ \{v_{2},v_{4}\}$, respectively, in
which ``$\circ$'' represents the operation of vertex contraction in
a graph. It is easily proved that they both are maximal planar
graphs of order $n-2$. Therefore, we obtain
$$f(G,4)=f((G-v)\circ\{v_{1},v_{3}\},4)+f((G-v)\circ\{v_{2},v_{4}\},4)=f(G_{1},4)+f(G_{2},4)$$
\end{proof}

\begin{center}
\includegraphics [width=100pt]{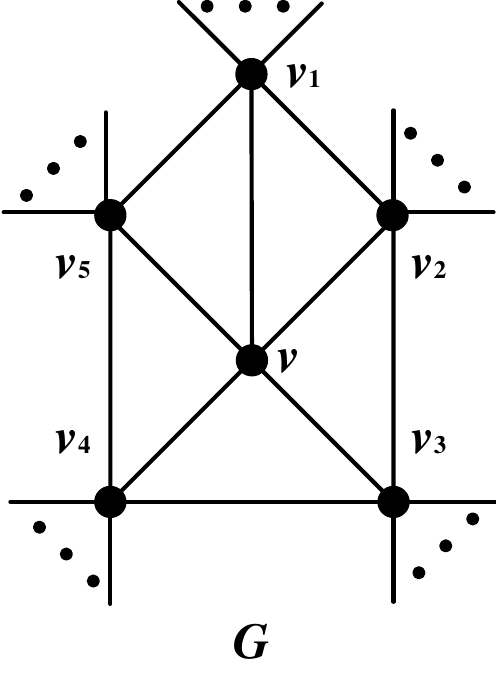}\\
\textbf{Figure 3.3.} A maximal planar graph having a 5-degree vertex
$v$.
\end{center}

\begin{theorem2}\label{th12}
\emph{(Chromatic polynomial on contracting 5-wheel operation)} Let
$G$ be a maximal planar graph, $v$ be a 5-degree vertex of $G$, and
$\Gamma(v)=\{v_{1},v_{2},v_{3},v_{4},v_{5}\}$ (shown in Figure 3.3).
Then
$$f(G,4)=[f(G_{1},4)-f(G_{1}\cup\{v_{1}v_{4},v_{1}v_{3}\},4)]+[f(G_{2},4)-f(G_{2}\cup\{v_{3}v_{1},v_{3}v_{5}\},4)]$$
$$+[f(G_{3},4)-f(G_{1}\cup\{v_{4}v_{1}\},4)]$$
where $G_{1}=(G-v)\circ \{v_{2},v_{5}\}$,
$G_{2}=(G-v)\circ\{v_{2},v_{4}\}$,$G_{3}=(G-v)\circ\{v_{3},v_{5}\}$.
\end{theorem2}

\begin{proof}
The maximal planar graph $G$ is represented by $G[{\Gamma(v)\cup
\{v\}}]$ in the following derivation. The chromatic polynomial of
$G$ can be calculated by applying Lemma 3.1 repeatedly. If parallel
edges appear in the process, reserve only one edge excluding $W_2$.
We use wheel $W_{5}$ to represent the chromatic polynomial of a
maximal planar graph. In this way, we can obtain

\includegraphics [width=130pt]{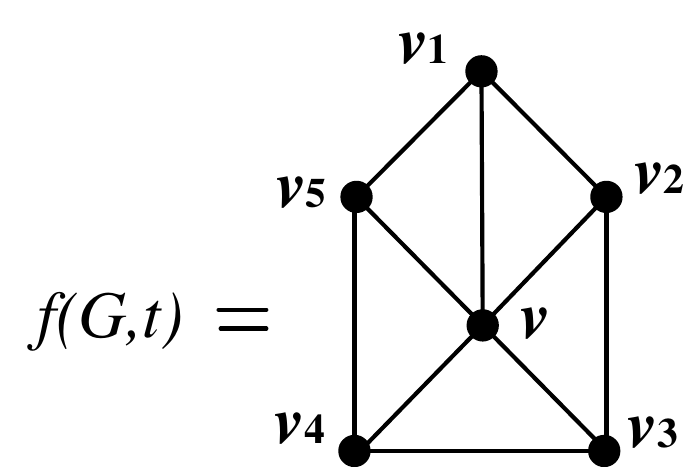}

\vspace{5mm}
\includegraphics [width=200pt]{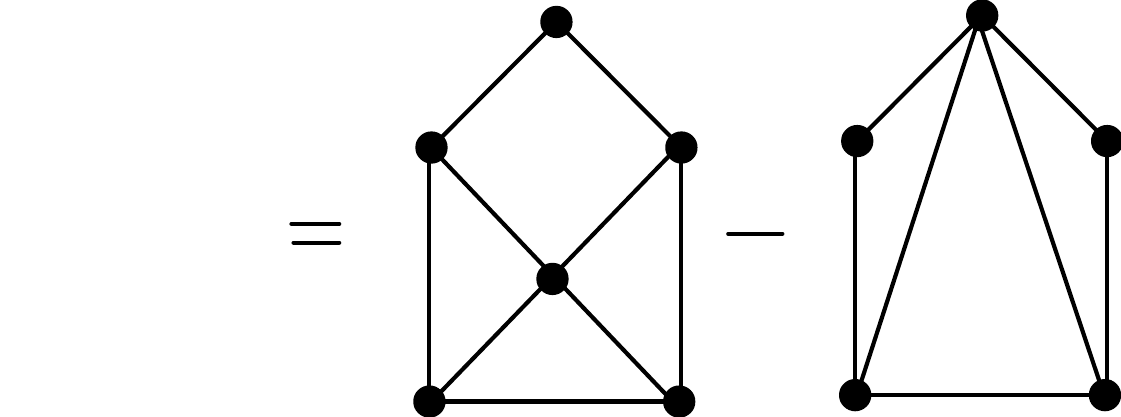}
\vspace{5mm}

\includegraphics [width=280pt]{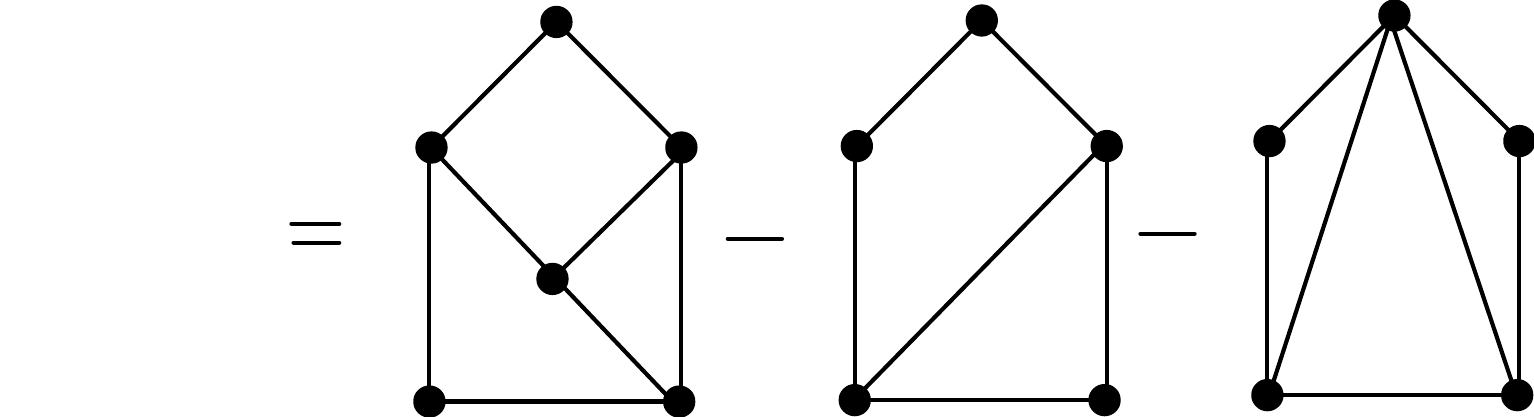}
\vspace{5mm}

\includegraphics [width=300pt]{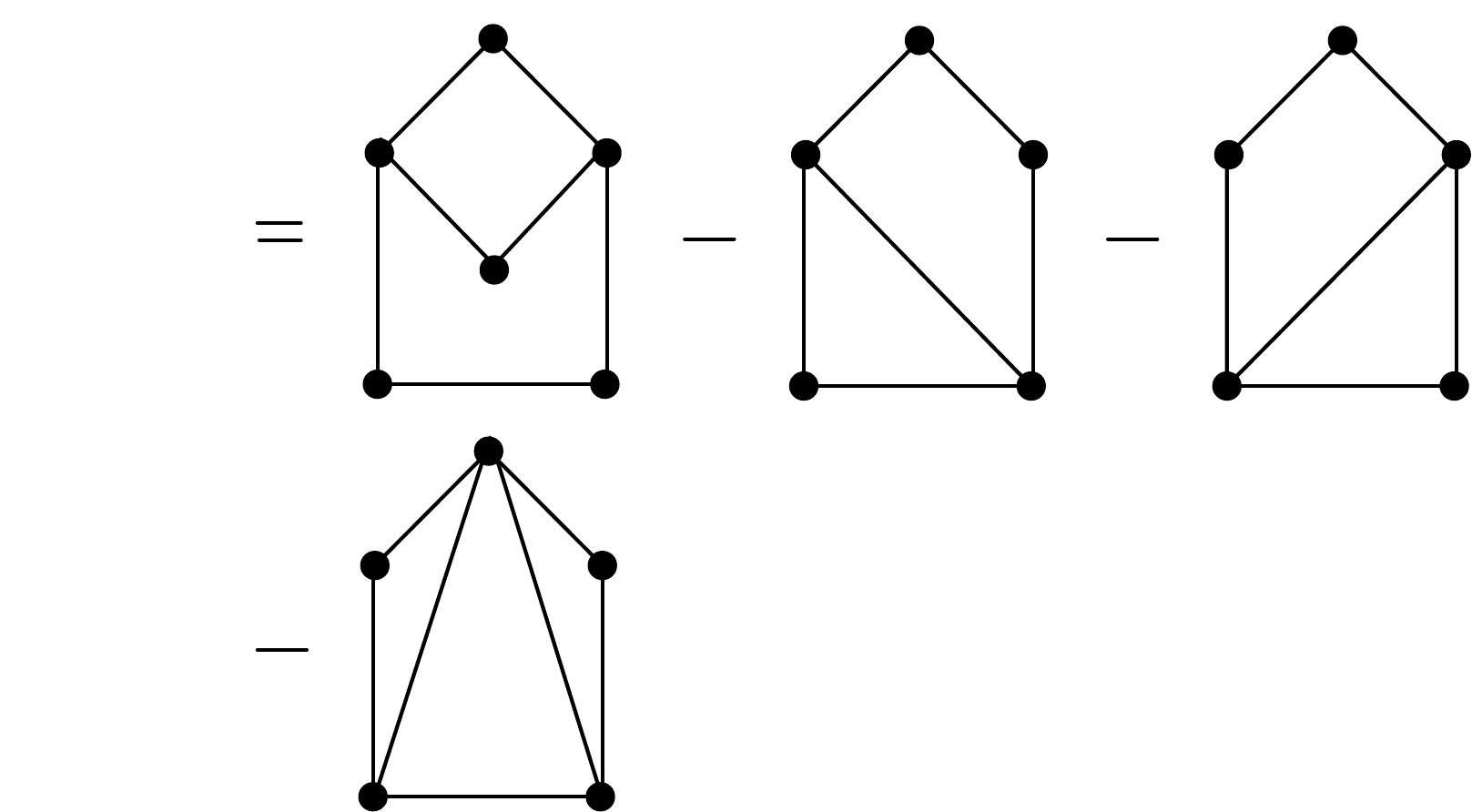}

\includegraphics [width=300pt]{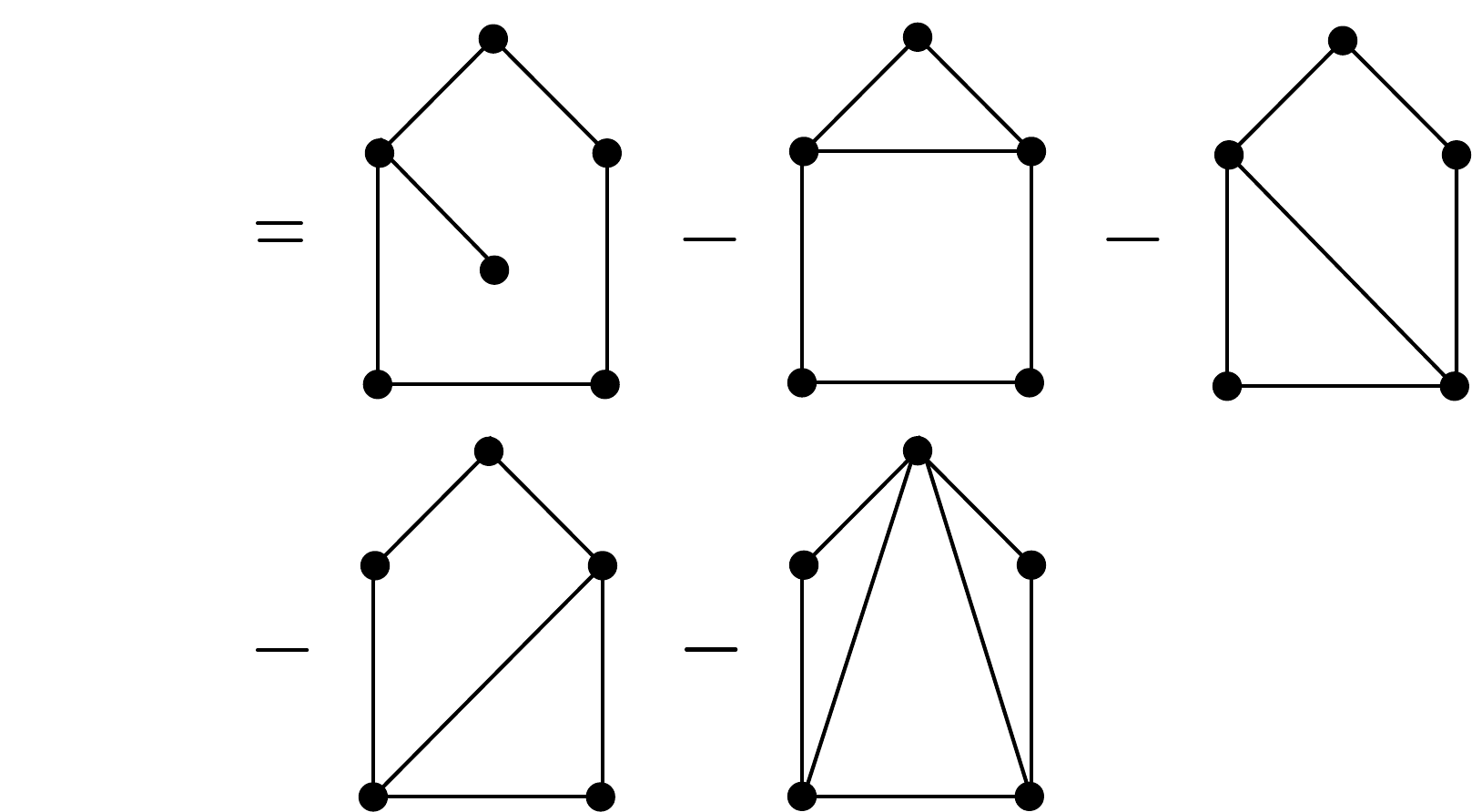}

\includegraphics [width=300pt]{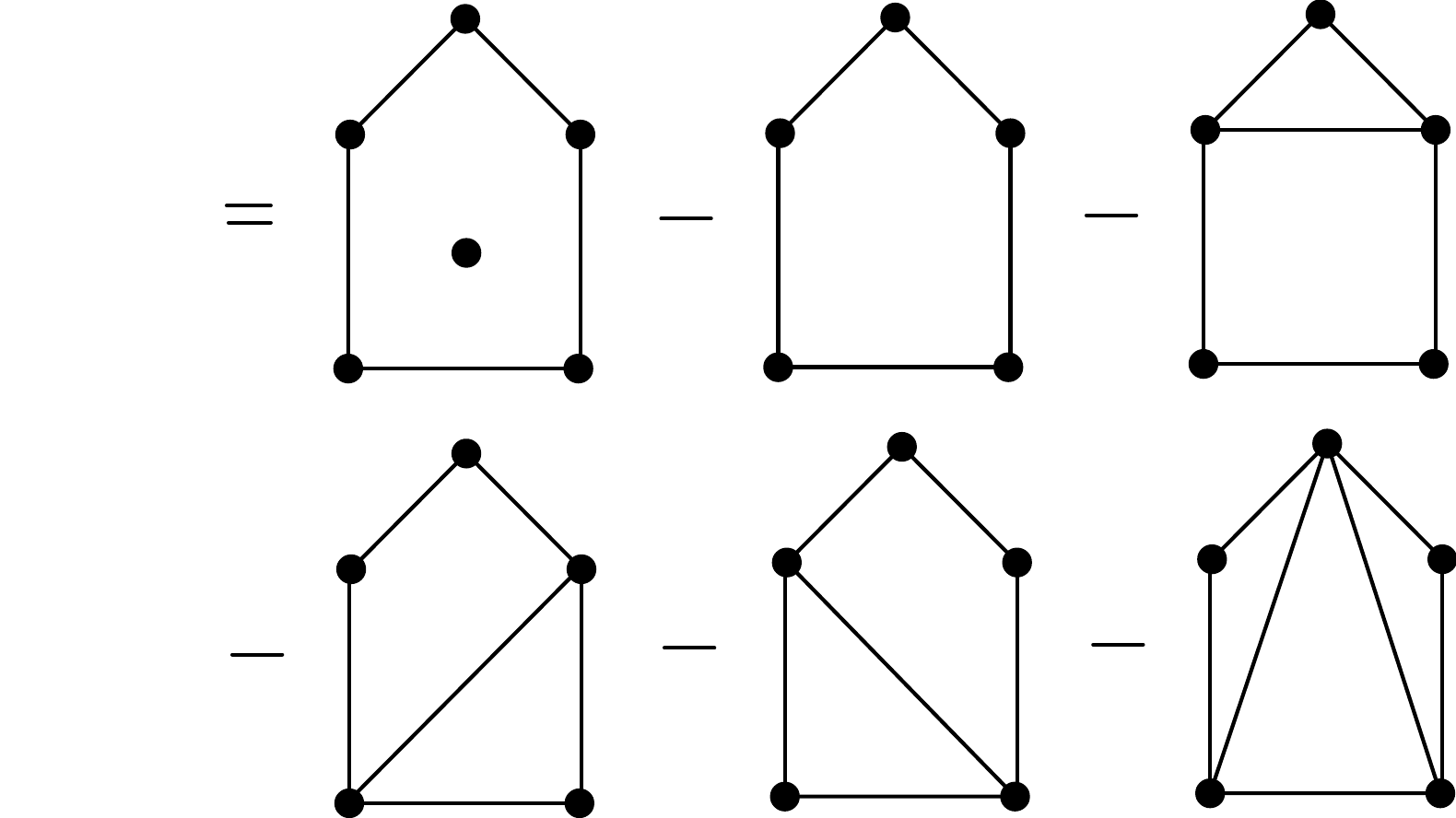}

By Lemma \ref{th4}, the chromatic polynomial of the first graph at
the righthand of the last equation is $tf(G-v,t)$. Therefore, we
have

\vspace{5mm}
\includegraphics [width=300pt]{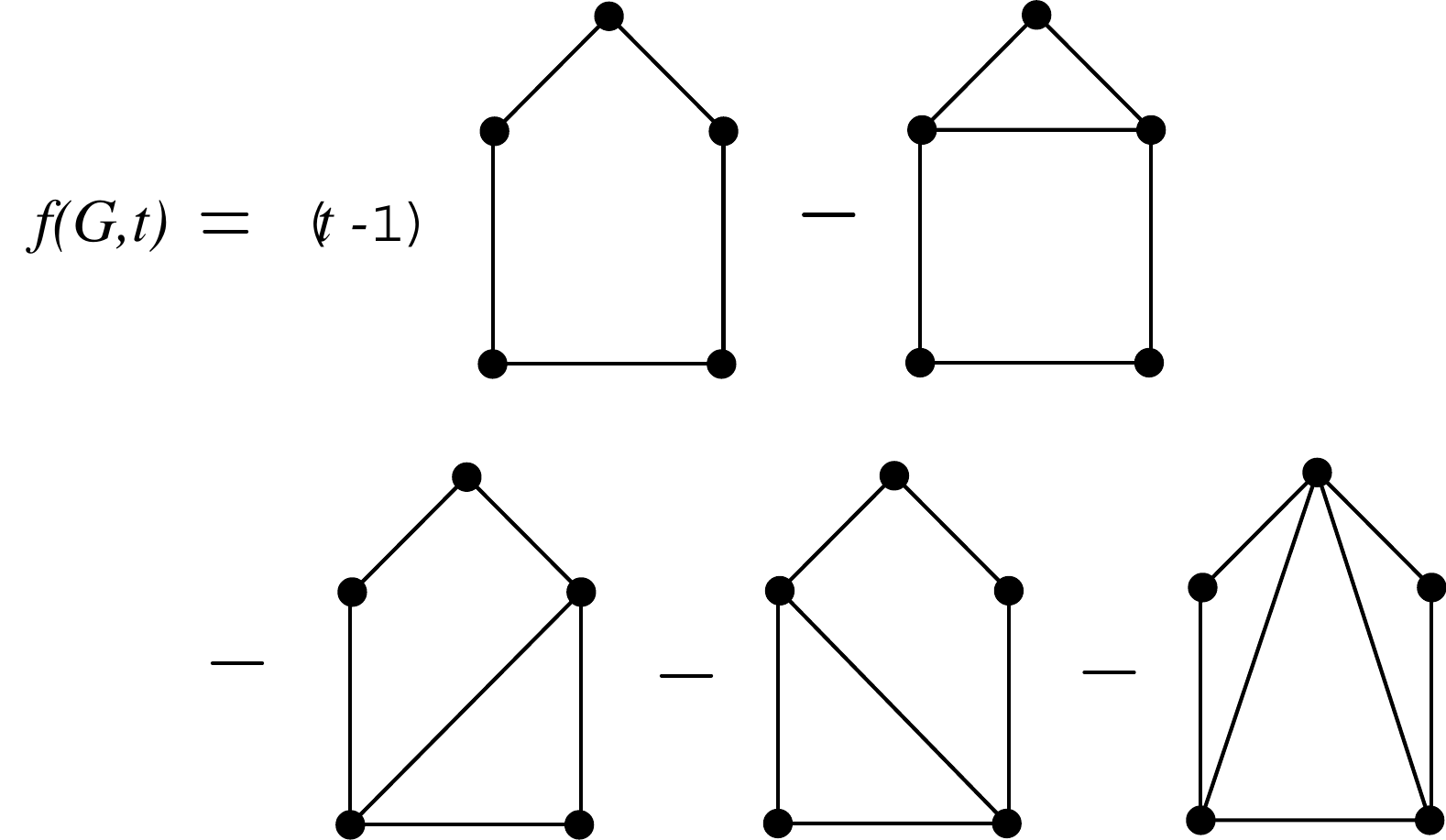}
\vspace{5mm}

For $t=4$, the following equation holds:

\vspace{5mm}
\includegraphics [width=300pt]{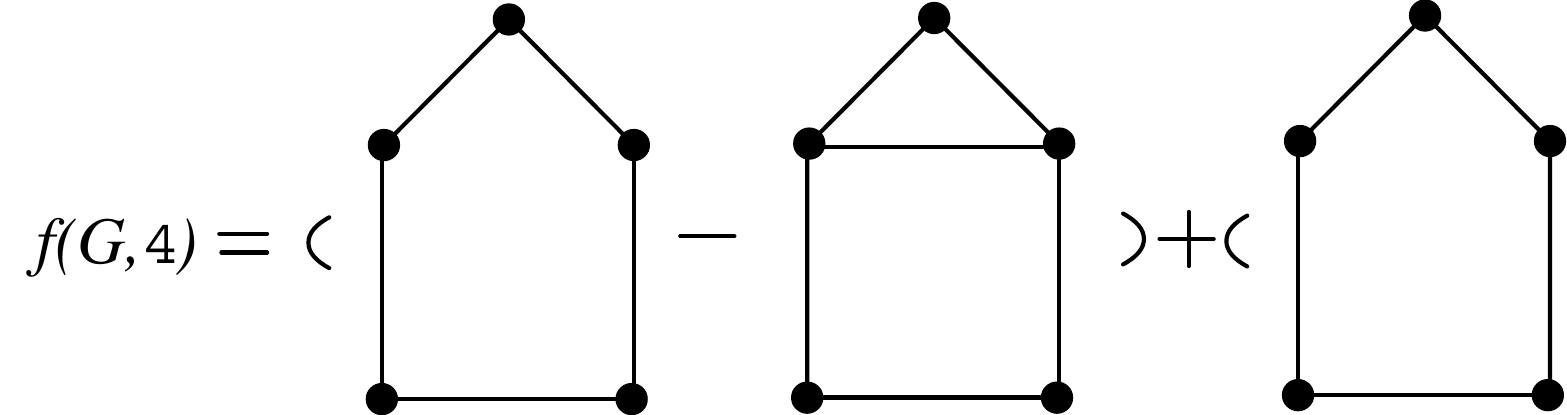}

\includegraphics [width=300pt]{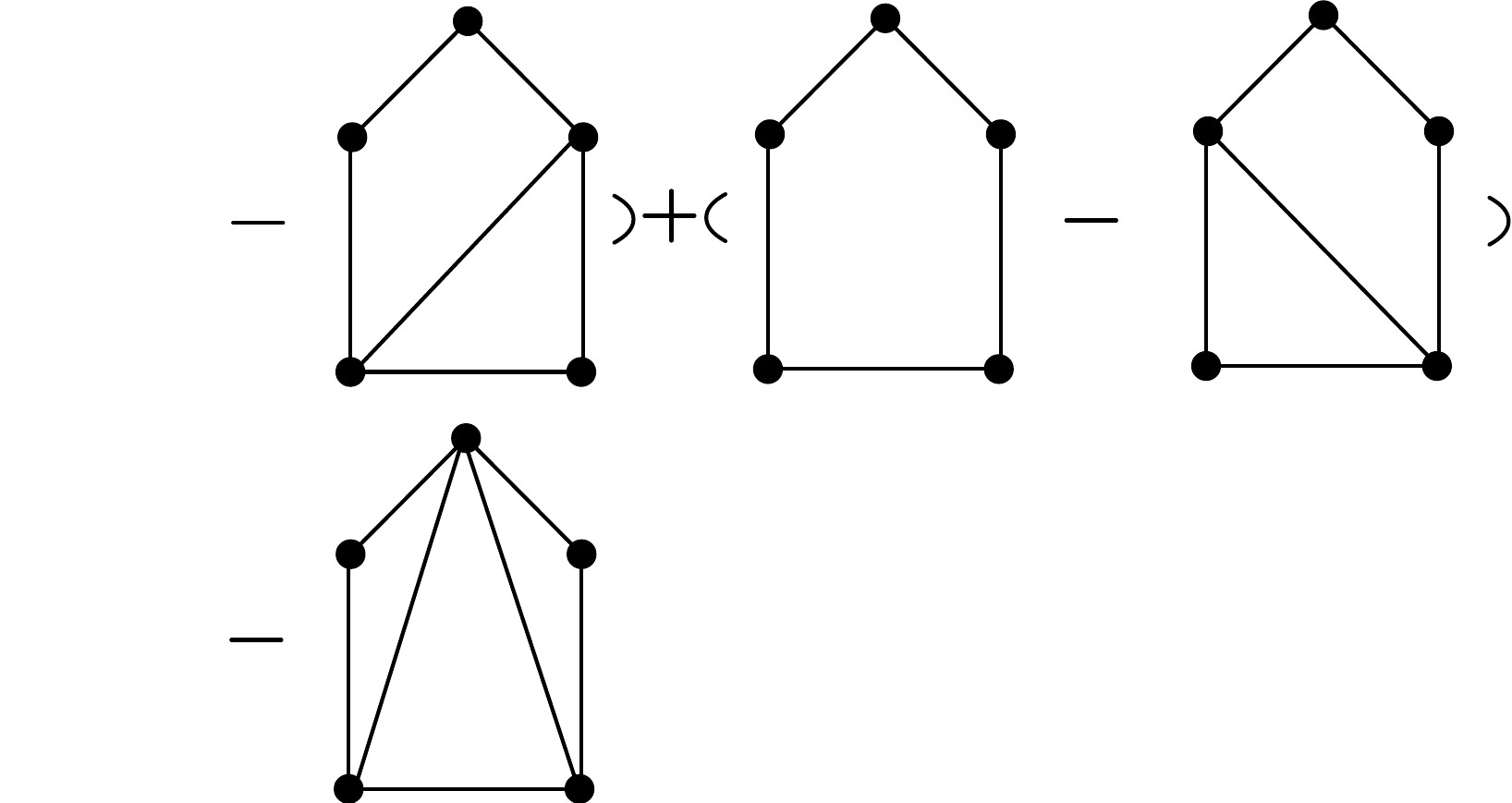}

\includegraphics [width=300pt]{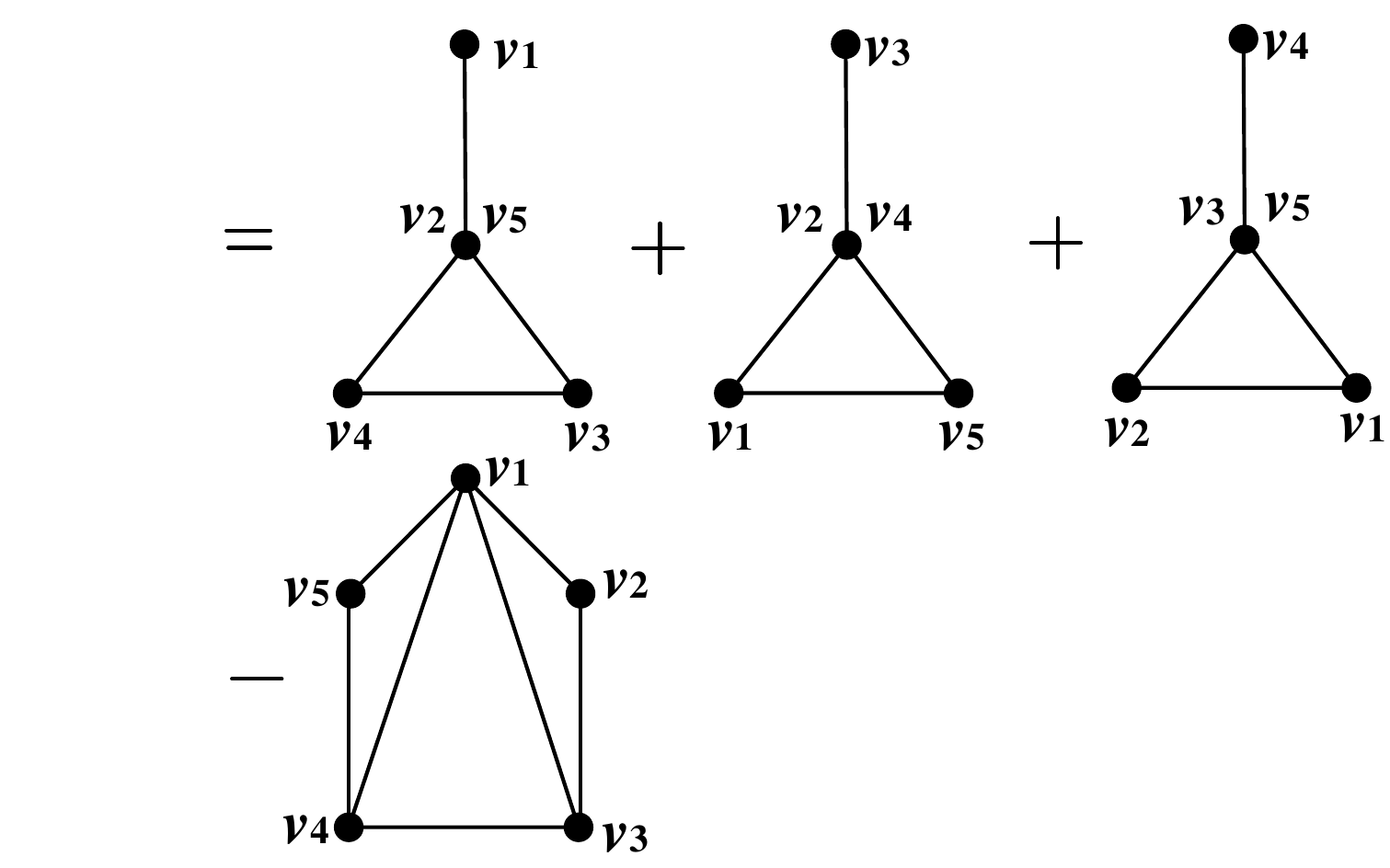}

\includegraphics [width=300pt]{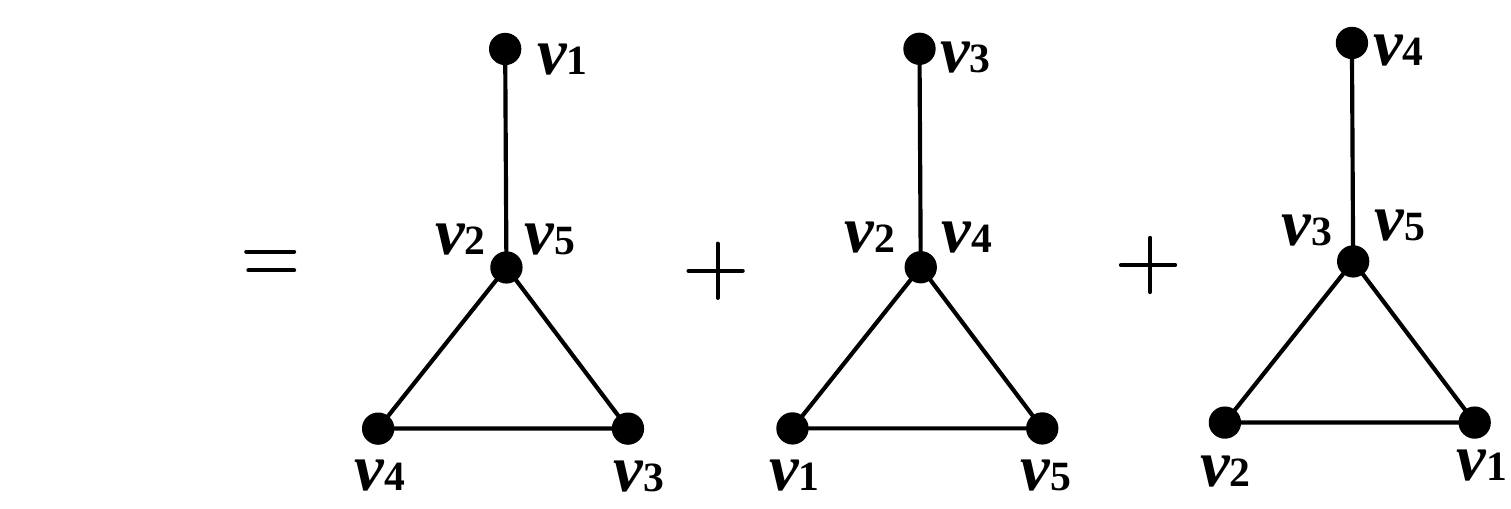}

\includegraphics [width=230pt]{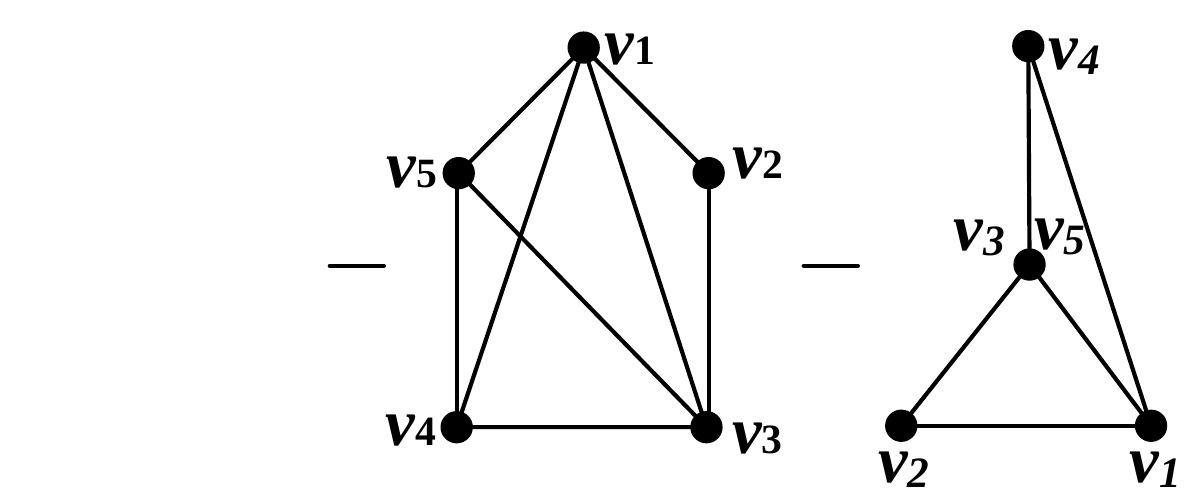}

\includegraphics [width=290pt]{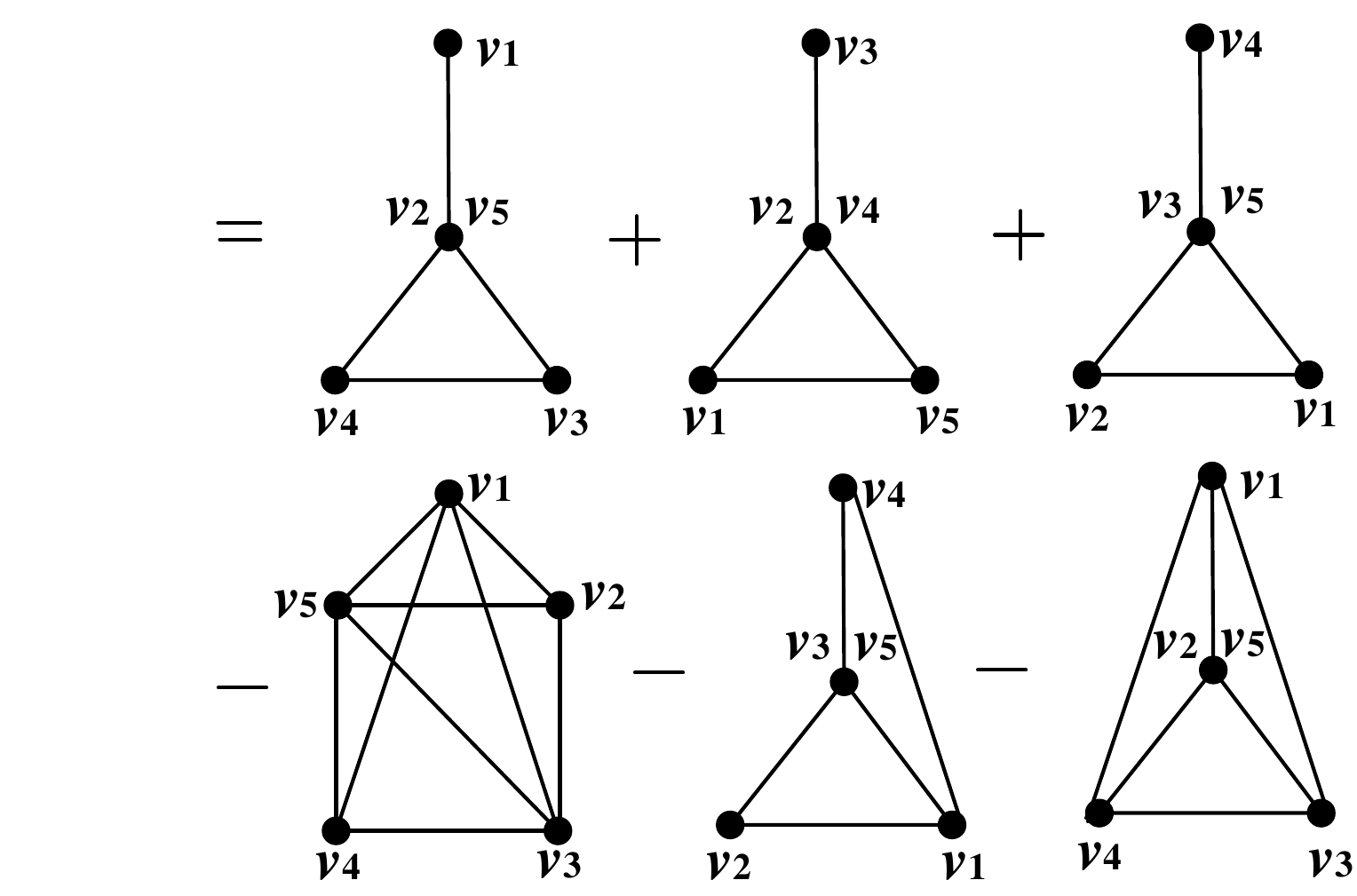}

\includegraphics [width=290pt]{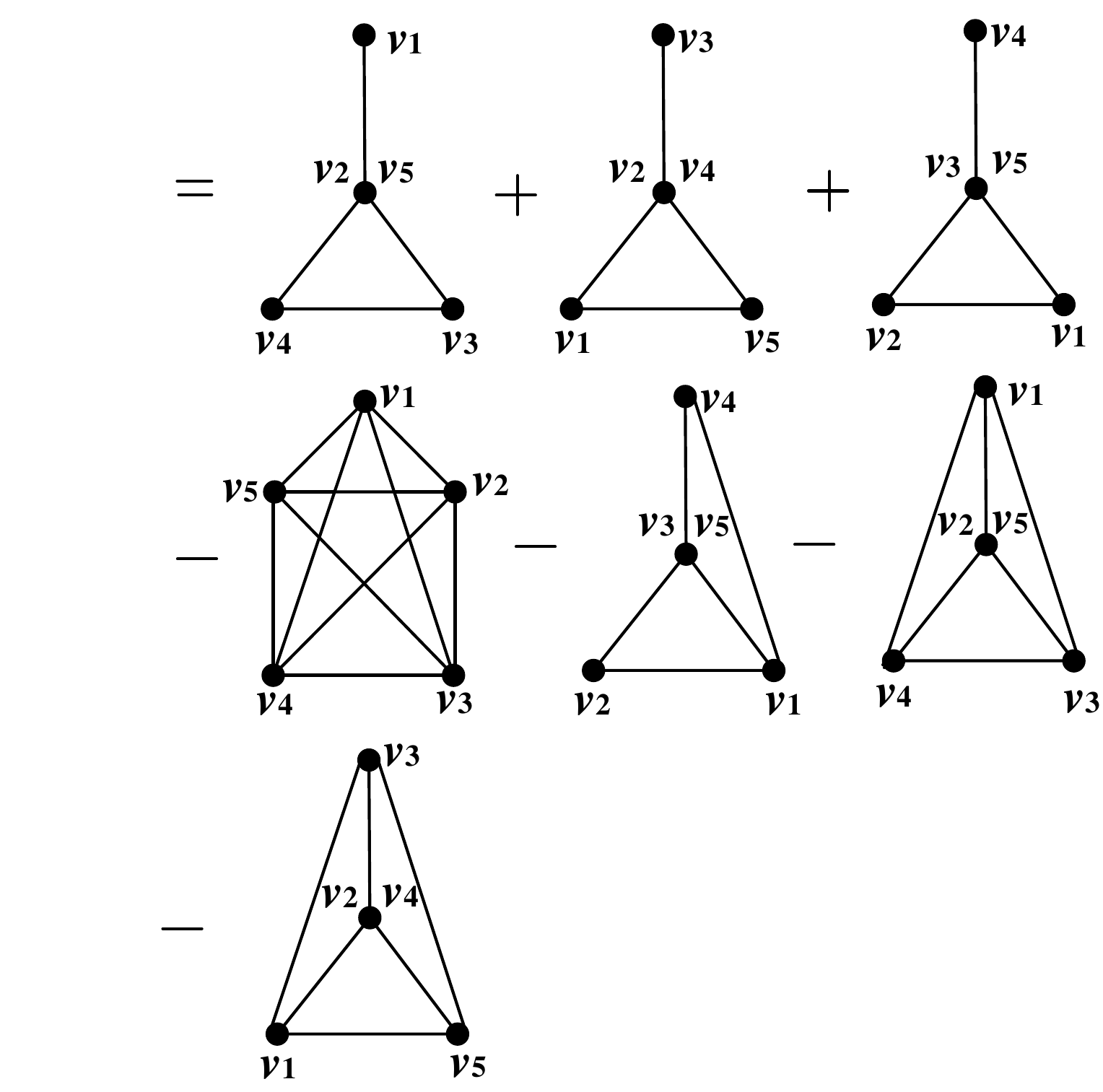}
\vspace{5mm}

Notice that the fourth graph at the righthand of the last equation,
denoted by $G'$, contains subgraph $K_{5}$ and so $f(G',4)=0$. Thus,
we can obtain that

\includegraphics [width=320pt]{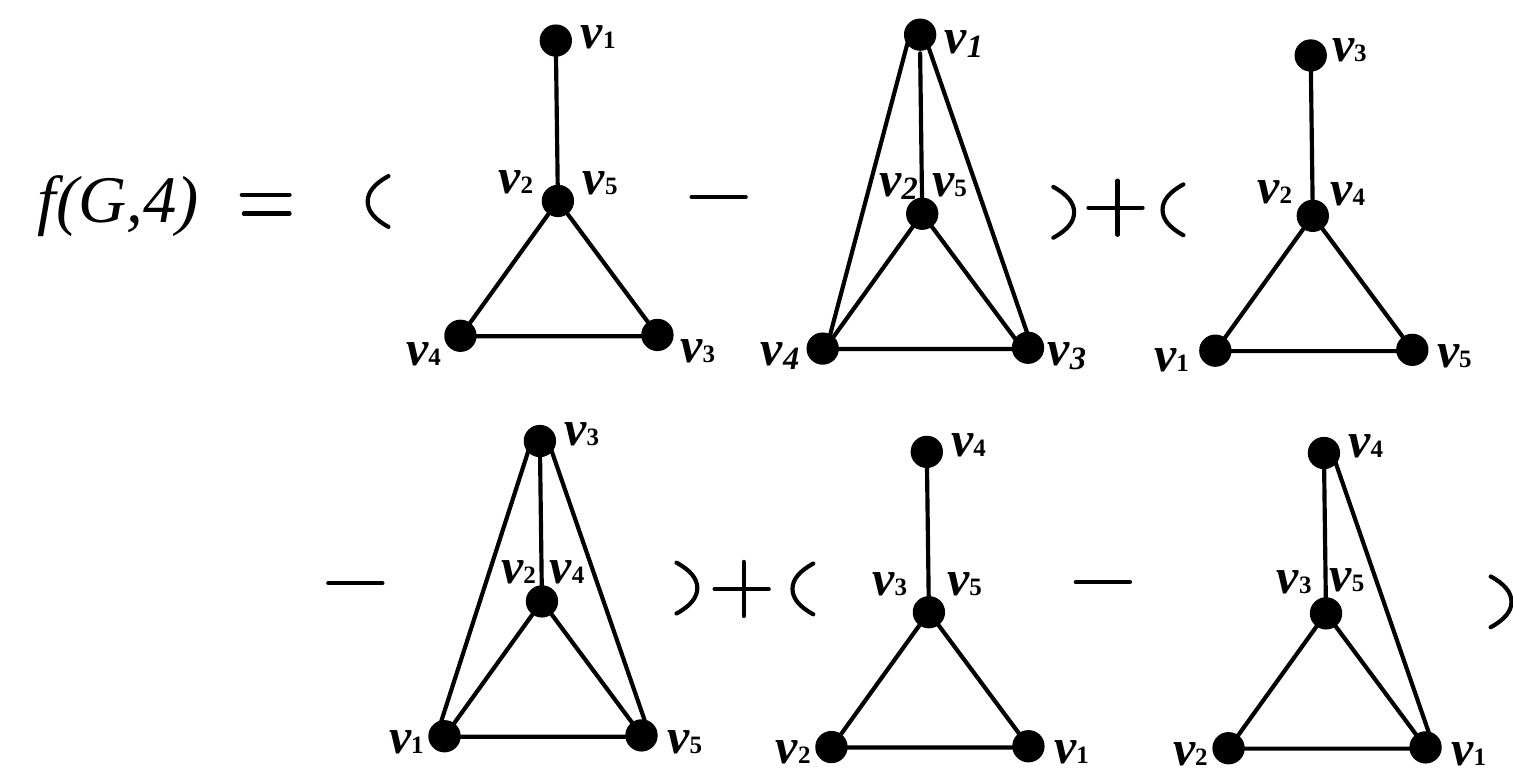}

Actually, the first graph in the first bracket of the equation is
$G_{1}=(G-v)\circ \{v_{2},v_{5}\}$; the first graph in the second
bracket is $G_{2}=(G-v)\circ\{v_{2},v_{4}\}$; and the first graph in
the third bracket is $G_{3}=(G-v)\circ \{v_{3},v_{5}\}$.
\end{proof}

This result is very important to the proof of the Four-Color
Theorem, since all resulted values of each bracket are no less than
zero. Obviously, The Four-Color Theorem can be proved if one
bracket's value is greater than zero. It may hold for every
bracket's value. Moreover, it is known that the graphs $G_{1},G_{2}$
and $G_{3}$ are 4-colorable, all of which are maximal planar graphs.
Take the graph in the first bracket for example. In the second graph
of the first bracket, the vertices $v_{1},v_{3},v_{4}$ and $v_{2}'$
(which is the new vertex from $v_{2}$ and $v_{5}$) can form a
complete subgraph of order 4, so they have to be colored with
different colors. Therefore, the Four-Color Theorem holds if there
exists one coloring in $C_{4}(G_{1})$ such that $v_{1}$ and $v_{3}$,
or $v_{1}$ and $v_{4}$ receive the same color. However, it is not an
easy task for proving it, which will depend on the works of Sections
4, 5, 6, 7, 8 and 9.

\section{The generating operation system on maximal planar graphs}
With respect to research on many problems about coloring of maximal
planar graphs, naturally, it is  important to understand clearly the
structure of maximal planar graphs. In this field, some relevant
results have been given in studying on methods of computer-assisted
proof of the Four-Color Conjecture [4-5,61]. Some scholars have
designed some constructive algorithms to research special properties
of maximal planar graphs.

Here, a generating operation system on maximal planar graphs is
given. This system consists of the operating objects and the basic
operators, where the operating objects are maximal planar graphs,
and with totally four pairs of the basic operators, which are the
extending $k$-wheel operation and the contracting $k$-wheel
operation of inverse of the extending $k$-wheel operation for
$k=2,3,4,5$. The basic function of this system is using $K_3$ as a
basic operating object and apply the four pairs of basic operators
repeatedly to generate a maximal planar graph. By the help of this
method, we can construct all maximal planar graphs with minimum
degree $\delta\geq 4$ and orders from 6 to 12, which will be used in
this paper. Moreover, we further study the contracting and extending
operations under 4-colorings.


\subsection{Chromatic isomorphism of graphs}
In order to make the following discussion convenient, the concept of
\textbf{chromatic isomorphism} will be introduced in this
subsection.

Figure 4.1(x) is a 3-chromatic graph $G$ of order 6, in which the
subgraph induced by $V'=\{v_{1},v_{2},v_{3},v_{4}\}$ is uniquely
3-colorable and the unique partition of color class is
$\{v_{1}\},\{v_{2},v_{4}\},\{v_{3}\}$ shown in Figure 4.1(y). Here,
we color the vertices with colors 1,2,3, respectively. Notice that
$\{1, 2, 3\}$ denotes the color set (see Figure 4.1(z)). That is
 $$\{v_{1}\}\rightarrow 1, \{v_{2},v_{4}\}\rightarrow 2, \{v_{3}\}\rightarrow 3 \eqno{(4.1)}$$

So we can obtain 3 kinds of colorings of $G$, denote by $f_{a},
f_{b}, f_{c}$, respectively.
\[f_{a}= (\begin{array}{cccccc}
v_{1} & v_{2} & v_{3} & v_{4} & v_{5} & v_{6} \\
1 & 2 & 3 & 2 & 1 & 2
\end{array}
 )\]
 \[f_{b}= (\begin{array}{cccccc}
    v_{1} & v_{2} & v_{3} & v_{4} & v_{5} & v_{6} \\
    1 & 2 & 3 & 2 & 3 & 1
     \end{array}
 )\]
 \[f_{c}= (\begin{array}{cccccc}
    v_{1} & v_{2} & v_{3} & v_{4} & v_{5} & v_{6} \\
    1 & 2 & 3 & 2 & 3 & 2
     \end{array}
 )\]

These 3 kinds of colorings are shown in Figure 4.1(a),(b),(c). With
the permutation of colors in (4.1), the corresponding new colorings
are generated, similar to $f_{a},f_{b},f_{c}$. This graph has
eighteen colorings, shown in Figure 4.1. In fact, it only has three
partitions of color class of colorings (see Figure 4.1(a),(b),(c)).
Other colorings can be obtained by the color permutation on these
three colorings. Here we present the notion of \textbf{chromatic
isomorphism}.
 \begin{center}

  \includegraphics [width=330pt]{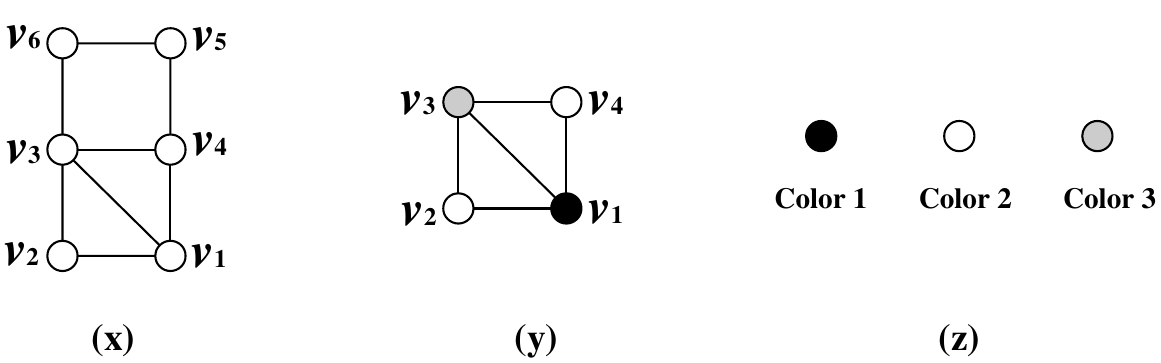}

  \vspace{5mm}
  \includegraphics [width=280pt]{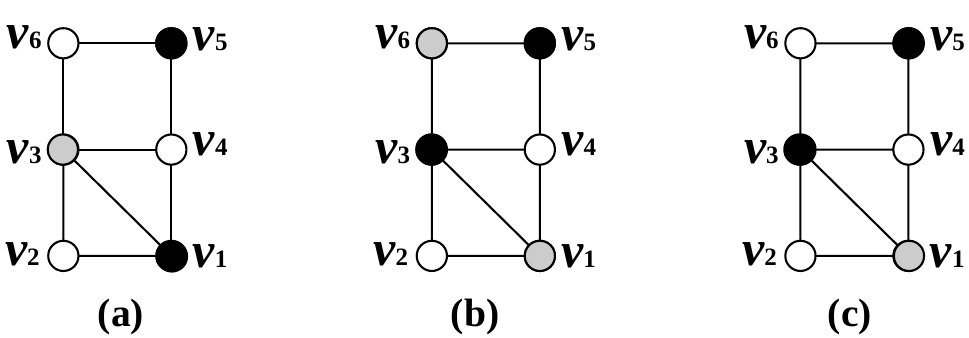}

  \includegraphics [width=280pt]{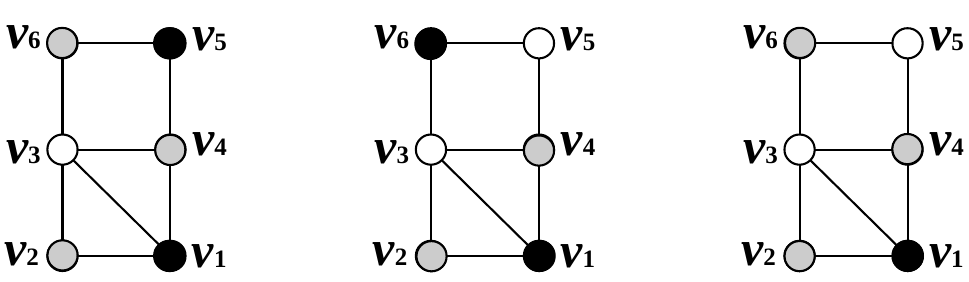}

  \includegraphics [width=280pt]{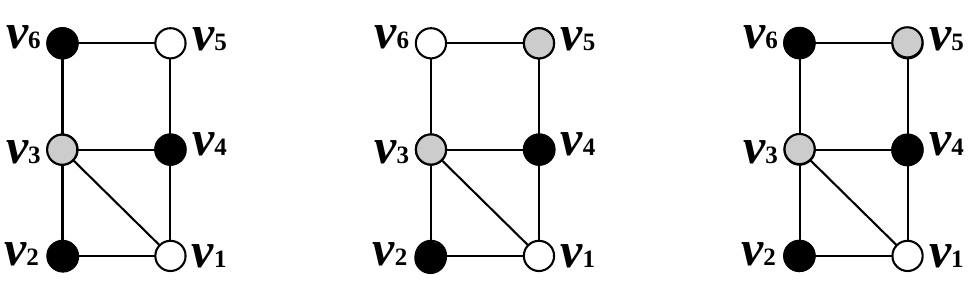}

  \includegraphics [width=280pt]{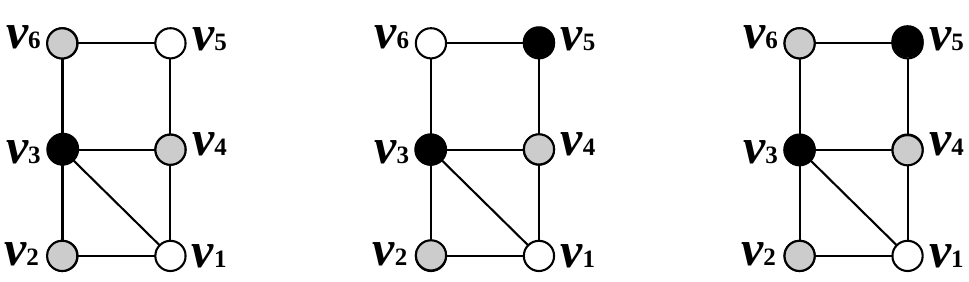}

  \includegraphics [width=280pt]{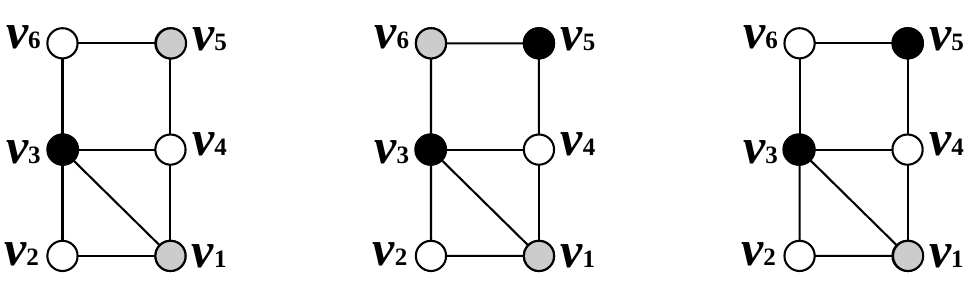}

  \includegraphics [width=280pt]{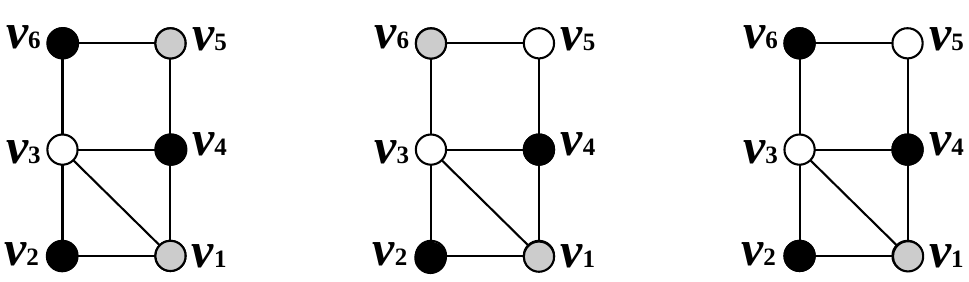}

  \vspace{5mm}
  \textbf{Figure 4.1.} The graph $G$ with 6 vertices and 18 colorings

 \end{center}

 \begin{definition}

Let $G$ be a $k$-chromatic graph, and let $f_{1},f_{2}\in C_{k}(G)$.
We say $f_{1},f_{2}$ to be \textbf{chromatic isomorphic} if the
partitions corresponding to the color classes of $f_{1}$ and $f_{2}$
are identical. The \textbf{$f_1$-chromatic isomorphic set} denotes
the coloring set containing such all of colorings of $G$ that are
chromatic isomorphic to $f_1$. If one representation is taken from
every chromatic isomorphic set, then the set containing such all of
representations is called the \textbf{chromatic isomorphic group},
denoted $C^{0}_{k}(G)$. Actually, $C^{0}_{k}(G)$ is the set formed
by the partitions that correspond to the color classes of all
$k$-colorings of $G$.
\end{definition}

If two colorings of a graph are chromatic isomorphic,  then they can
be easily inter-converted when the colors  are adjusted properly.
Therefore, these two colorings can be viewed as one coloring. We
just need to choose one of them. In Figure 4.1, it is easy to see
that six colorings in every column have the same color classes.
Therefore, all the colorings in every column are chromatic
isomorphic. So the graph $G$ in Figure 4.1  has only three kinds of
distinct colorings, represented by three colorings in the first row
(a), (b) and (c) of Figure 4.1. It is not hard to prove the
following result:

\begin{theorem2}\label{th20}
Let $G$ be a $k$-chromatic, simple and undirected graph. For the
chromatic isomorphic group $C^{0}_{k}(G)$ and the set $C_{k}(G)$ of
all the colorings of $G$, we then have
$$|C_{k}(G)|=k!|C^{0}_{k}(G)|\eqno{(4.2)}$$
\end{theorem2}

Therefore, when we discuss some properties of all $k$-colorings of a
graph $G$, we only consider the chromatic isomorphic group
$C^{0}_{k}(G)$.

\subsection{The basic generating operation system on maximal planar graphs}

In this section, we mainly give the definition of the basic
operators and some related properties which are used in the basic
generating operation system of maximal planar graphs without
colorings.

The extending $2$-wheel operation means: first, add a new edge
between two adjacent vertices, which shall generate a pair of
$2$-parallel edges, called a 2-cycle; second, add a new vertex in
the face bounded by the 2-cycle  and make the new vertex adjacent to
the two vertices on the 2-cycle. The contracting $2$-wheel operation
means: first, delete the center of this wheel and the two edges
incident to the center; second, delete one of the parallel edges in
the remainder.

In the second section of this paper, we have introduced the
extending $3$-wheel operation on maximal planar graphs: add a new
vertex in a certain face of the maximal planar graph, and then add
three edges linking the new vertex with the three vertices of the
face, respectively; meanwhile, we have also introduced the
contracting $3$-wheel operation: delete a certain $3$-degree vertex
and the edges incident with it.

The so-called contracting $4$-wheel operation is: in a maximal
planar graph, delete a certain $4$-degree vertex and the edges
incident with it, and then do the contracting operation to a pair of
the nonadjacent vertices in the neighborhood of the $4$-degree
vertex. The extending $4$-wheel operation is the inverse operation
of the contracting $4$-wheel operation. The following gives the
detailed definition.

Let $G$ be a maximal planar graph, and let $xuy$ be a path of length
2 in $G$. The so-called extending $4$-wheel operation on the path
$xuy$ is to replace the $xuy$ by a 4-cycle $xuu'yx$; that is, spilt
the vertex $u$ into two vertices $u$ and $u'$, and split the edges
$xu$ and $uy$ into two edges $xu,xu'$ and $uy,uy'$ respectively.
This process is shown from the first to the fourth graphs in Figure
4.2. Then add a new vertex $v$ on the face bounded by the $4$-cycle
$xu'yux$, and make $v$ adjacent to vertices $x,u',y,u$ respectively.
The resulting graph is referred to as a generated graph by the
extending $4$-wheel operation, denoted as $G*xuy$ (see the fifth
graph in Figure 4.2).

\begin{center}
\includegraphics [width=230pt]{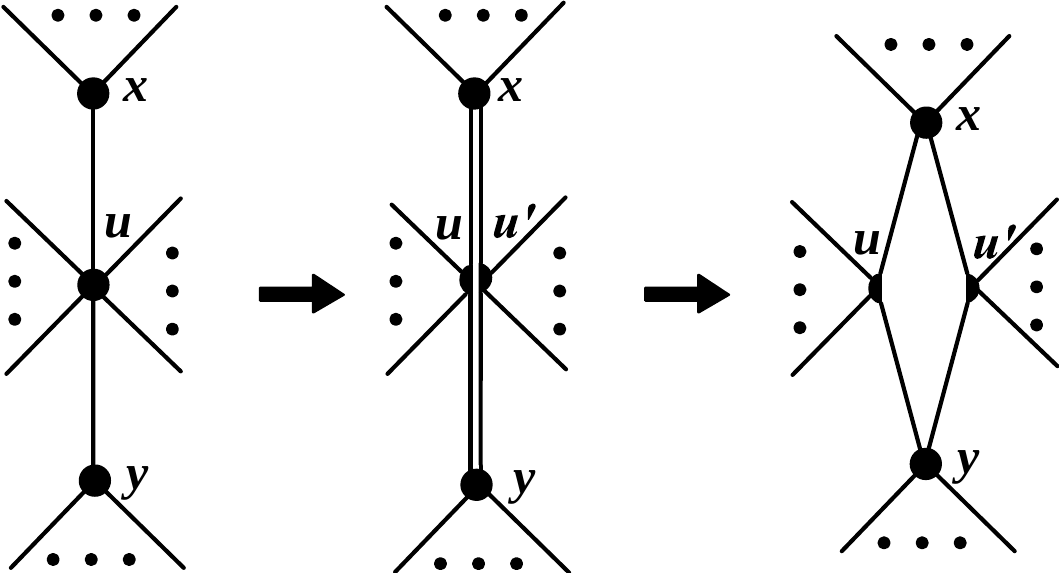}

\includegraphics [width=250pt]{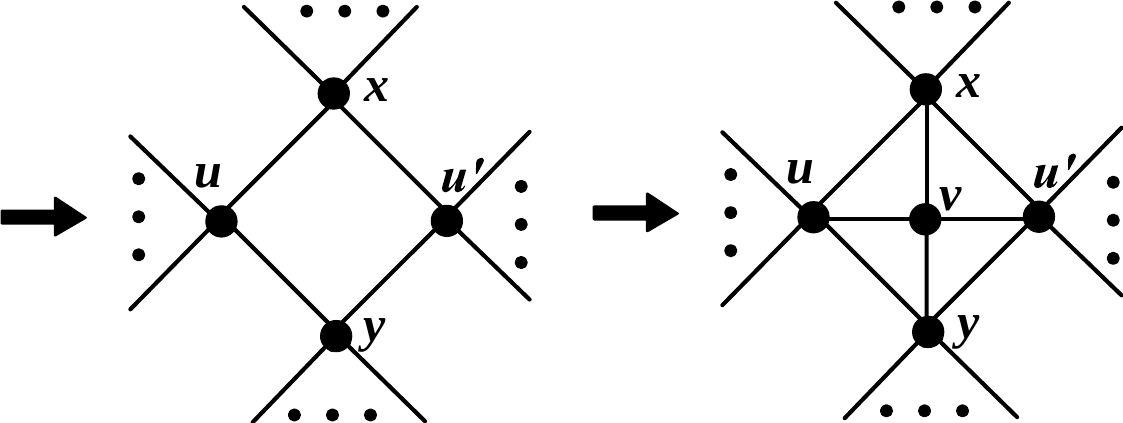}

\textbf{Figure 4.2.} The schematic diagram of the extending 4-wheel
operation
\end{center}

The graph shown in Figure 4.3  is called a \textbf{funnel}, where
the $1$-degree vertex is called the \textbf{top}, the $3$-degree
vertex is called the \textbf{stem} and two $2$-degree vertices are
called the \textbf{bottoms}.

\begin{center}
\includegraphics [width=300pt]{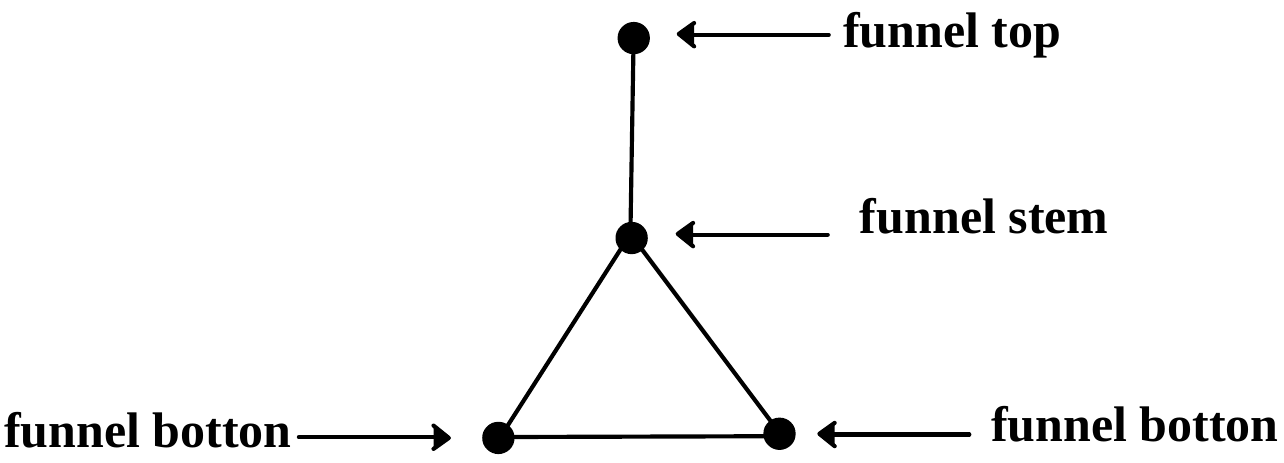}

\textbf{Figure 4.3.} A funnel.
\end{center}

For a maximal planar graph, the contracting 5-wheel operation and
the extending 5-wheel operation are similar to the contracting
4-wheel operation and the extending 4-wheel operation, except the
difference that the extending 5-wheel operation is on a funnel of
the maximal planar graph, while the extending 4-wheel operation is
on a $2$-path in the maximal planar graph. Here we only give a
graphical illustration shown in Figure 4.4 for  the contracting
5-wheel operation and the extending 5-wheel operation.
\begin{center}
  \includegraphics [width=370pt]{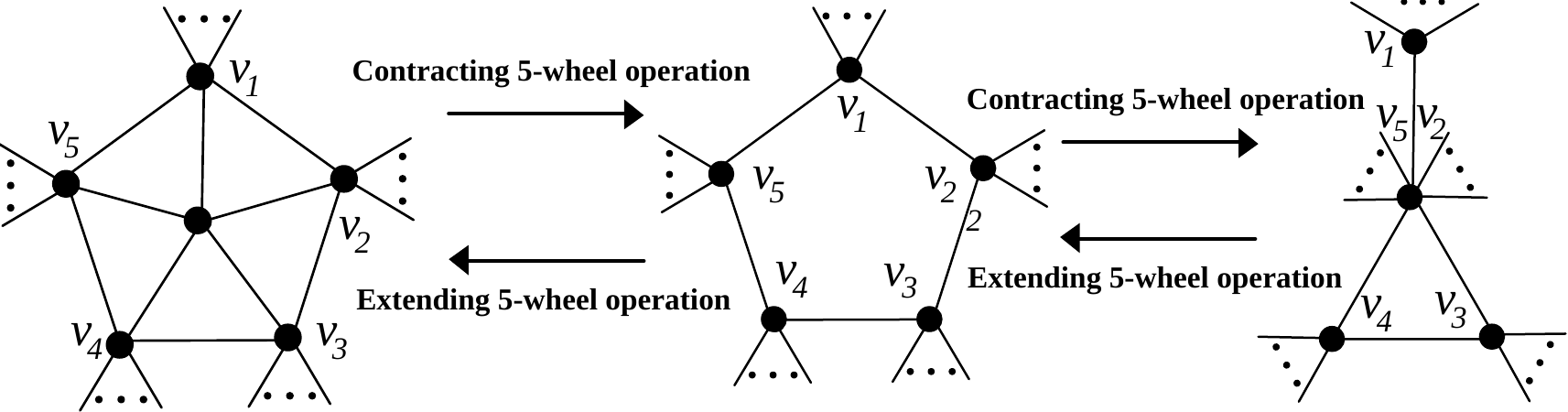}

 \textbf{Figure 4.4.} The schematic diagram of extending 5-wheel operation
\end{center}

In this paper, for $i=2,3,4,5$, we use $\zeta^{-}_{i}(G)$ and
$\zeta^{+}_{i}(G)$ to denote the resulted graphs obtained from the
contracting $i$-wheel operation  and  extending $i$-wheel operation
for the maximal planar graph $G$, respectively.

The following two propositions are easy to prove.

\begin{Prop}
Let $G$ be a maximal planar graph. Then $\zeta^{-}_{i}(G)$ and
$\zeta^{+}_{i}(G)$ ($i=2,3,4,5$) are still maximal planar graphs.
\end{Prop}

\begin{Prop}
Let $G$ be a maximal planar graph of order $n$. Then
$\zeta^{-}_{2}(G)$ and $\zeta^{-}_{3}(G)$  are  maximal planar
graphs with order $n-1$; $\zeta^{-}_{4}(G)$ and $\zeta^{-}_{5}(G)$
are maximal planar graphs with order $n-2$, namely

$$|\zeta^{-}_{2}(G)|=|\zeta^{-}_{3}(G)|=|V(G)|-1=n-1;$$
$$|\zeta^{-}_{4}(G)|=|\zeta^{-}_{5}(G)|=|V(G)|-2=n-2.$$
\end{Prop}

\begin{theorem2}\label{the4.4}
Let $G$ be a  maximal planar graph of order $n\geq 3$. Then by
properly implementing the contracting $i$-wheel operation for
$i=2,3,4,5$, then $G$ can be contracted to $K_3$.
\end{theorem2}
\begin{proof} When $n=4$, there is only one maximal planar graph $K_4$, so the
conclusion holds obviously. Suppose that the conclusion holds when
$n\leq p~(p\geq 4)$, which means that for any maximal planar graph
with order $p$ at most, it can be contracted to a complete graph
$K_3$  by applying properly the contracting $i$-wheel operation for
$i=2,3,4,5$.

Now we consider the case $n=p+1$. For any maximal planar graph $G$
of order $p+1$, if $G$ has any $2$-degree or $3$-degree vertex,
deleting the $2$-degree or $3$-degree vertex and the incident edges,
then we will get a maximal planar graph, $\zeta^{-}_{2}(G)$ or
$\zeta^{-}_{3}(G)$. By the induction hypothesis, the conclusion
holds. If $\delta(G)=4$ (or $\delta(G)=5$), doing the contracting
4-wheel operation (or the contracting 5-wheel operation) for some
$4$-degree (or $5$-degree) vertex, then we will get a graph
$\zeta^{-}_{4}(G)$ or $\zeta^{-}_{5}(G)$, which is a maximal planar
graph with order $(p-1)$. By the induction hypothesis, they can be
contracted to a complete graph $K_3$ by applying properly the
contracting $i$-wheel operations with respect to $2\leq i\leq 5$.
\end{proof}

Through Theorem \ref{the4.4}, we know that every maximal planar
graph of order $n$ can be contracted to $K_3$  by four contracting
operations. Of course, tracing back to the reverse directions of the
contracting $k$-wheel operations on a graph $G$, from $K_3$, doing
the corresponding extending $i$-wheel operations with respect to
$2\leq i\leq 5$, we can also get the original graph $G$. Hence,
\begin{corollary}
For any two maximal planar graphs $G$ and $G'$, we can get $G'$ from
$G$ by doing the four pairs of contracting and extending operations.
\end{corollary}

We use the notation
$\Psi=\{\zeta^{-}_{2},\zeta^{+}_{2},\zeta^{-}_{3},\zeta^{+}_{3},\zeta^{-}_{4},\zeta^{+}_{4},\zeta^{-}_{5},\zeta^{+}_{5}\}$
to denote the four pairs of contracting and extending operations,
and the symbol $S(G)=(K_3,\Psi)$  to denote the generating operation
system of maximal planar graphs, since every maximal planar graph
can be generated by this system.

\subsection{Construction of maximal planar graphs}

Based on the method mentioned in the last subsection, we will give
the construction method and steps of the entire maximal planar
graphs of order $n$ with $\delta(G)\geq 4$. Especially, we construct
all the maximal planar graphs with orders from 6 to 12 and
$\delta(G)\geq 4$.

Let $Aut(G)$ denote the automorphism group of graph $G$, $xuy$ and
$x'u'y'$ are two different paths of graph $G$. $xuy$ and $x'u'y'$
are called being \textbf{equivalent}, if there exists a $\sigma$  in
$Aut(G)$, which makes $\sigma(x)=x'$, $\sigma(u)=u'$,
$\sigma(z)=z'$. Otherwise, these two paths are
\textbf{nonequivalent}.

\subsubsection{Construction method and steps of maximal planar graph with order $n$}

\indent Step 1. Generated from the maximal planar graphs of order
$n-2$ by doing extending 4-wheel operation and extending 5-wheel
operation.

The detailed process is: for any maximal planar graph of order
$n-2$, first choose all of the 2-paths that are nonequivalent
mutually.  For the maximal planar graph of order 7 with
$\delta(G)=4$, $G_7$, shown in figure 4.5, there are 4 different
kinds of 2-paths: 444 type, 445 type, 454 type and 545 type
respectively, where 444 type means that the degree sequence of the
2-path is (444), and the other types are similar to this.  Then  for
each 2-path selected, do extending 4-wheel operation. For example,
in $G_7$, the two  maximal planar graphs of order 9 resulted of
doing extending 4-wheel operation on 444 type and 454 type 2-paths
respectively are isomorphic; on 545 type length-2 path, the degree
sequence of the maximal planar graph resulted of doing extending
4-wheel operation  is (444444477); the degree sequence of the
maximal planar graph resulted of doing extending 5-wheel operation
is (444455556). The above cases are shown in Figure 4.5,
respectively.

Step 2: Generated from the maximal planar graphs of order $n-3$ by
doing the combination operations of extending 2-wheel and extending
4-wheel; or the combination operations of extending 3-wheel and
extending 5-wheel.

For example,  the maximal planar graph of order 9, $G_9$, can be
only generated from  the maximal planar graph of order 6, $G_6$.
Because the maximal planar graph of order 6 with $\delta(G)=4$ is
only the regular octahedron(as the first graph in Figure 4.5), it
only generates a maximal planar graph $G$ of order 9 by doing the
combination operations of extending 2-wheel and extending 4-wheel.
Its degree sequence is (444444666) and $\delta(G)=4$. Through the
combination operations of extending 3-wheel and extending 5-wheel,
$G_6$ can also generate a maximal planar graph $G$ with
$\delta(G)=4$ and degree sequence (444555555).

We have in fact constructed all five maximal planar graphs of order
9 with $\delta(G)=4$  by above examples. Note that 445 type can not
generate any maximal planar graphs of order 9 with $\delta(G)=4$;
all other graphs resulted of doing extending 5-wheel operation are
isomorphic to one of the five graphs.

\begin{center}
  \includegraphics [width=270pt]{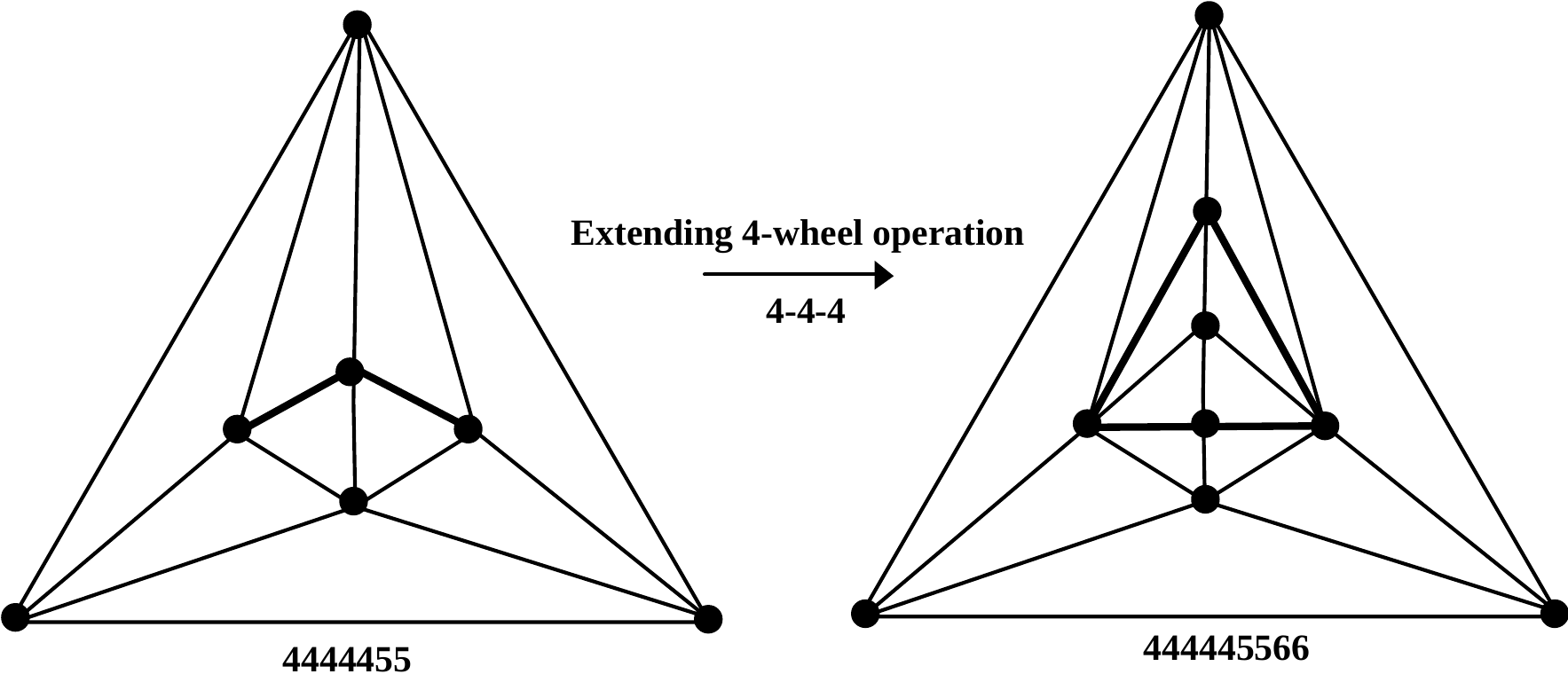}

  \includegraphics [width=270pt]{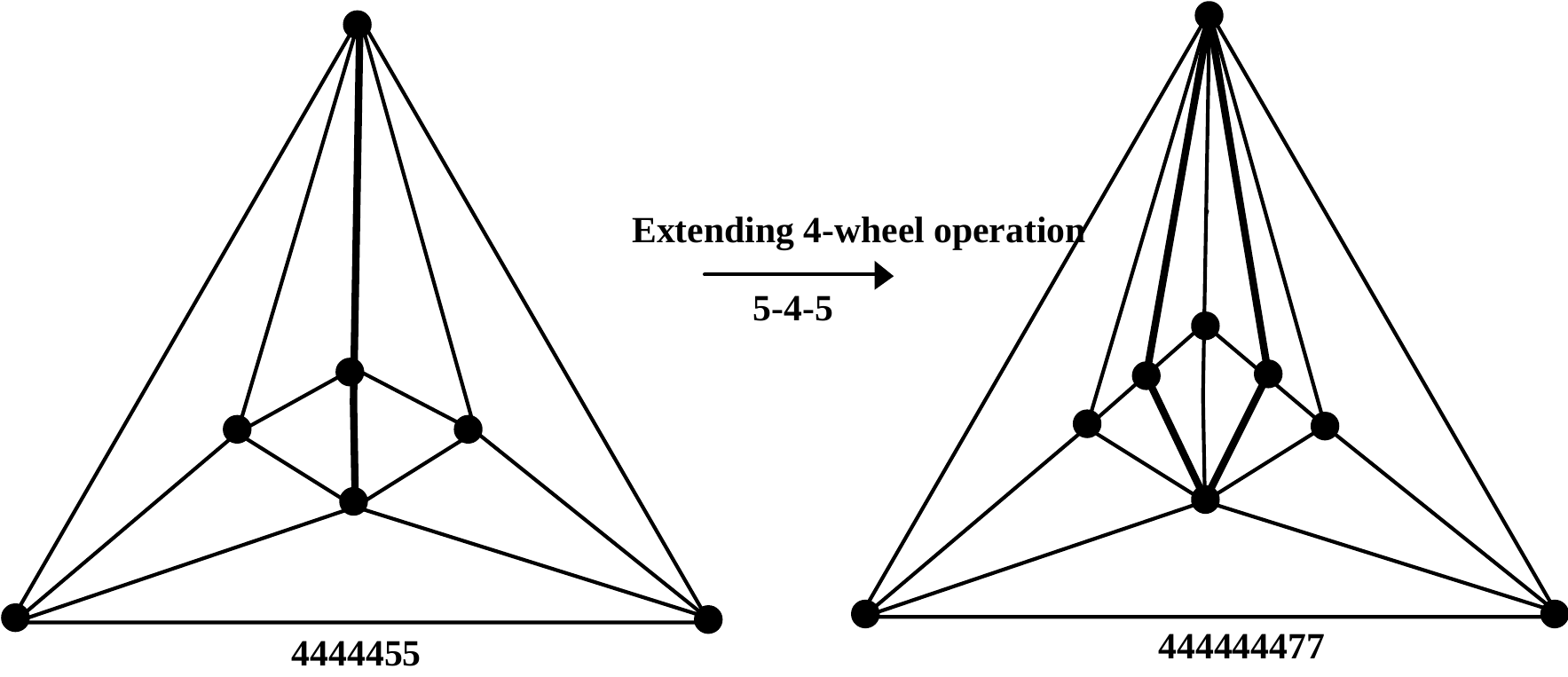}

  \includegraphics [width=280pt]{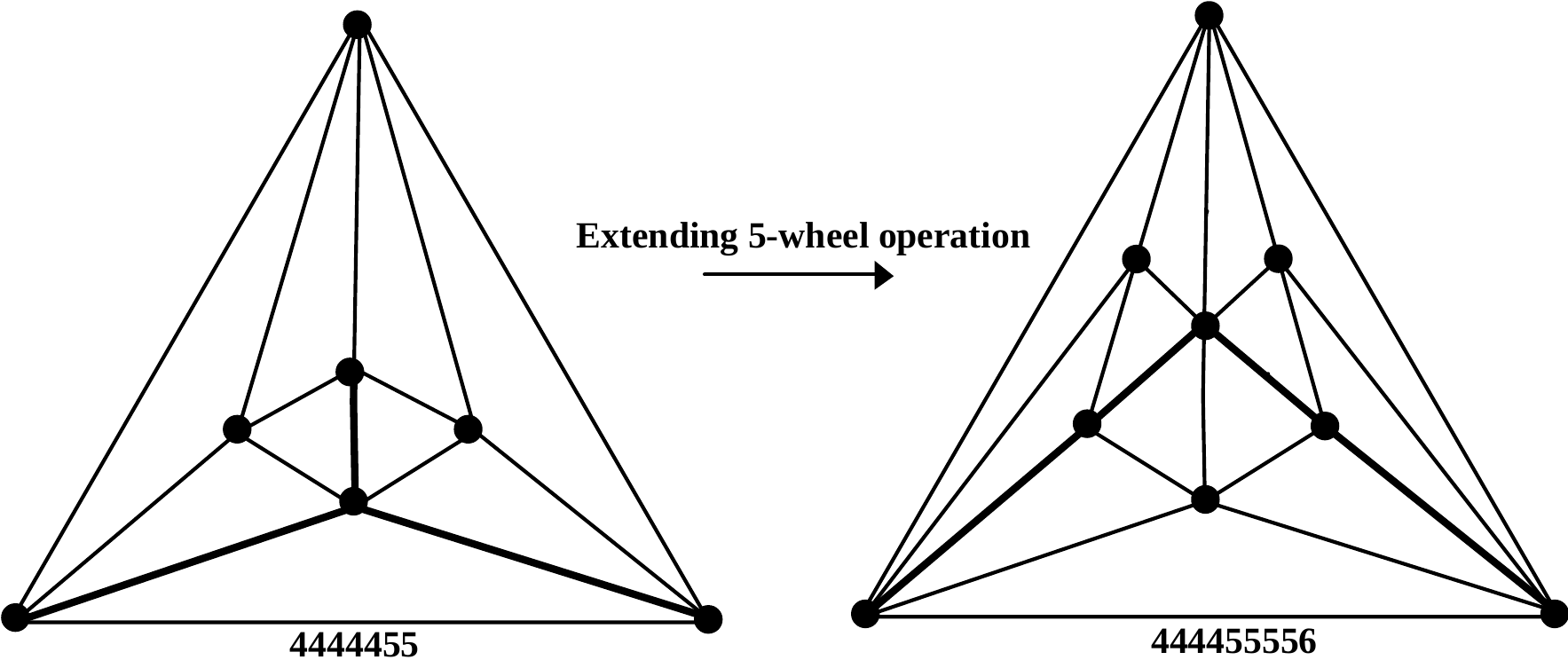}

  \includegraphics [width=370pt]{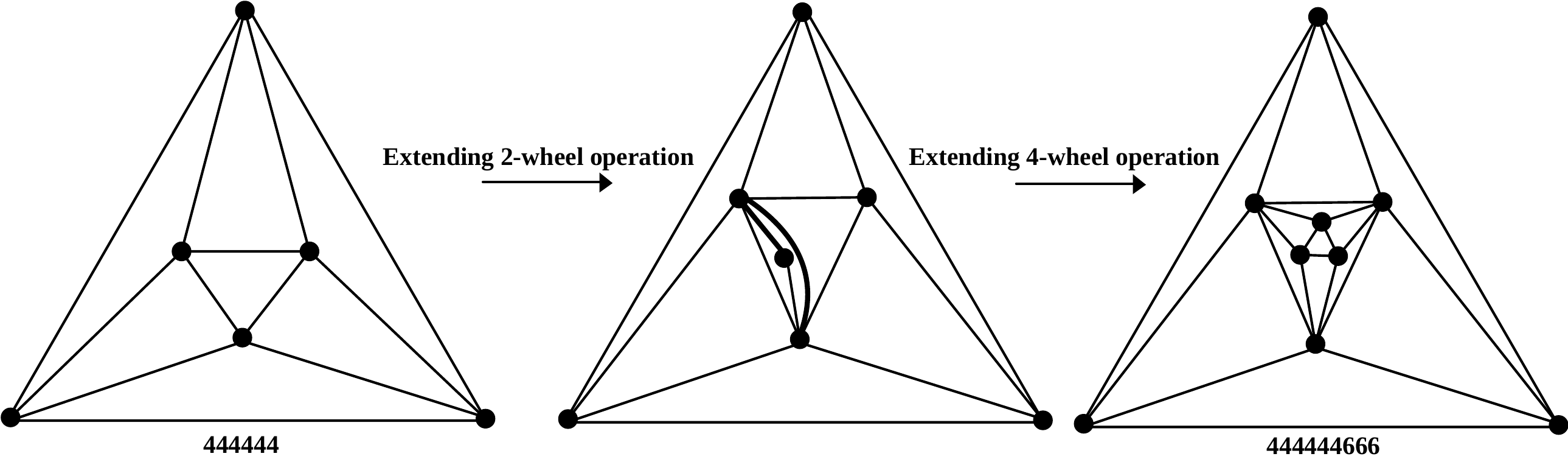}

  \includegraphics [width=370pt]{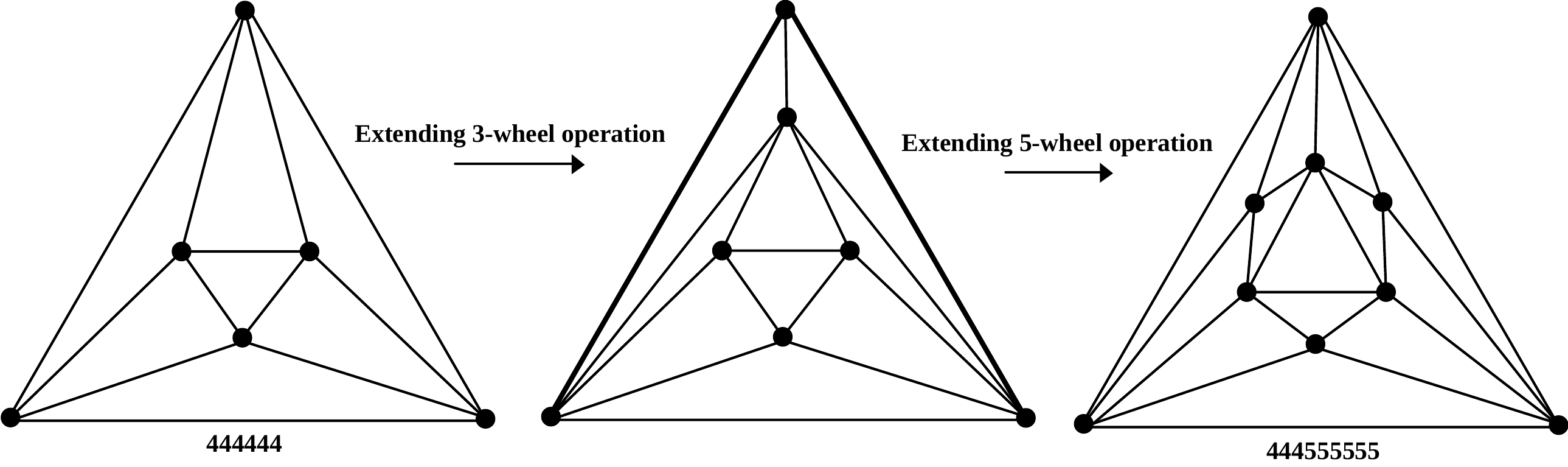}

  \textbf{Figure 4.5.} The schematic diagram of the generation procedure of maximal planar graph of order form 7 to 9 and $\delta(G)=4$
\end{center}

\subsubsection{All the maximal planar graphs with $\delta(G)\geq 4$ and orders from 6 to 12}
In order to prove the main result  of this section, we need to
investigate all maximal planar graphs with $\delta(G)=4$ and orders
from 6 to 12. Table 4.1 gives the number of maximal planar graphs
with $\delta(G)\geq 4$ in different orders. According to the above
generation methods, we construct these  maximal planar graphs, shown
in Figure 4.6 $\sim$ Figure 4.10 respectively.

\begin{center}
\scriptsize{ Table 4.1.the count chart of the maximal planar graphs
satisfied and the order of the graph is 6 to 12}\label{tab1}

\hspace{0.5cm}

 \large{ \begin{tabular}{ ccccccccccccccccc }

 \hline
 Order  &&6  && 7   &&8  && 9 &&  10&&  11 && 12\\
 \hline
 Graphs  &&1 &&  1  && 2 &&  5 &&  13 && 34 && 130\\
 \hline
  \end{tabular}}
\end{center}

\begin{center}
  \includegraphics [width=240pt]{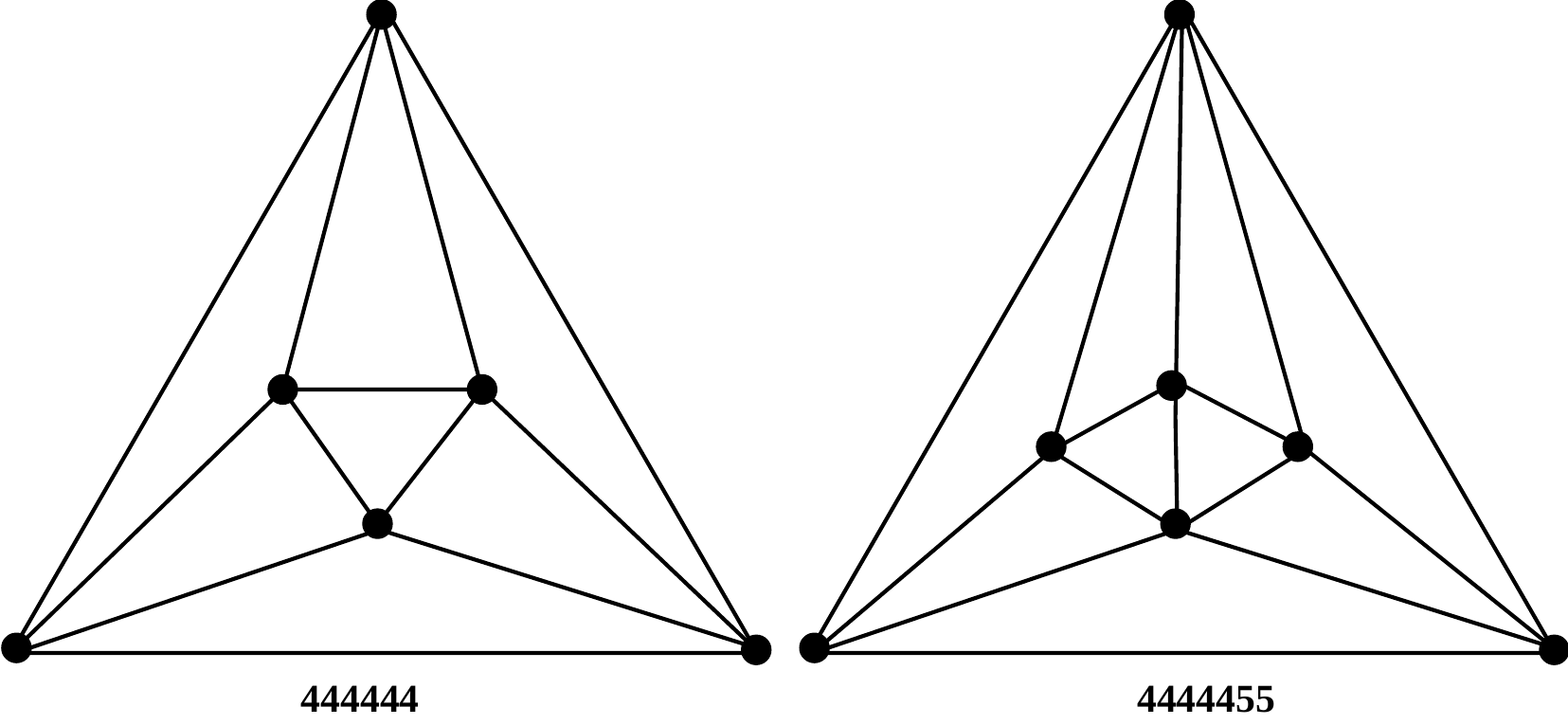}

  \textbf{Figure 4.6.} The maximal planar graph whose order is 6 and 7, $\delta(G)=4$
\end{center}
\begin{center}
  \includegraphics [width=240pt]{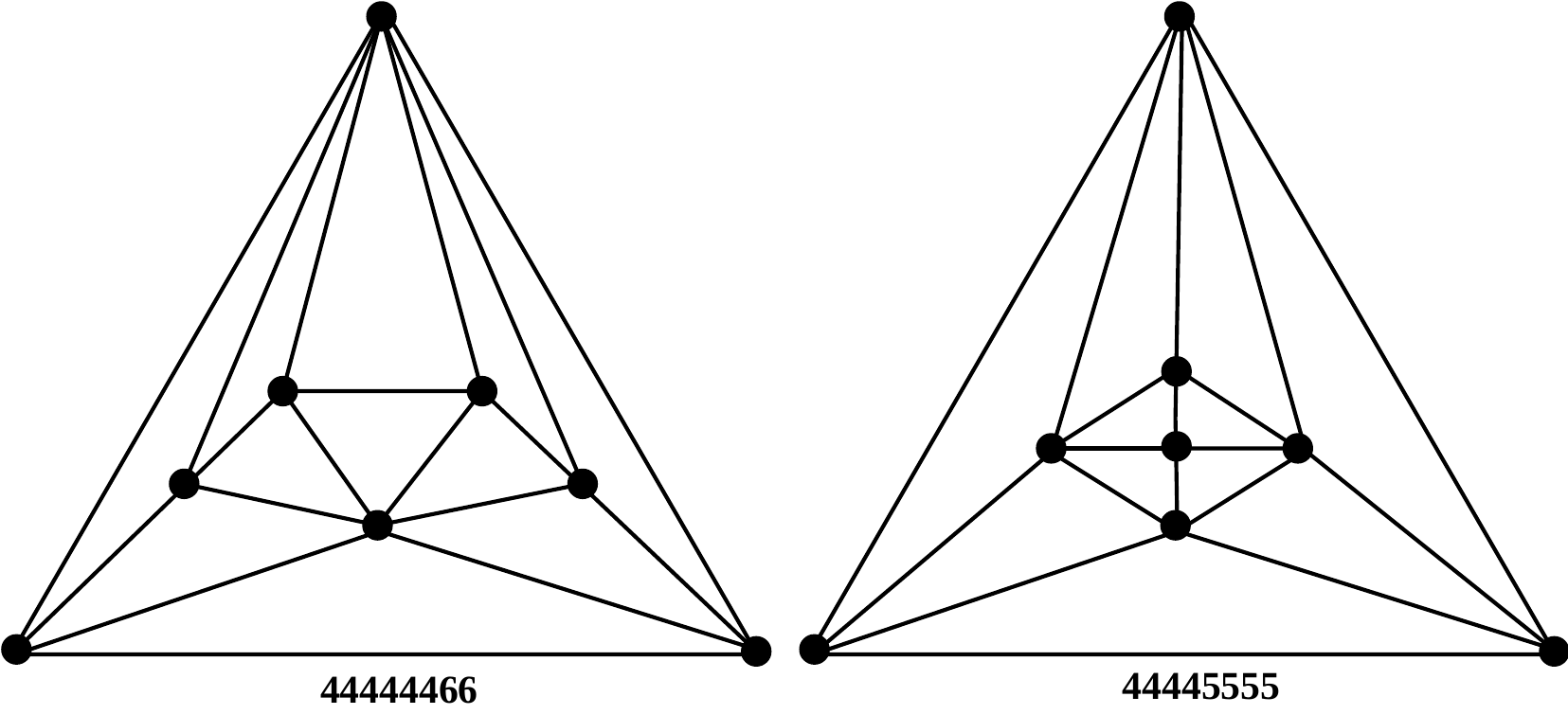}

  \textbf{Figure 4.7.} All two maximal planar graphs whose order is 8 and $\delta(G)=4$
\end{center}
\begin{center}
  \includegraphics [width=340pt]{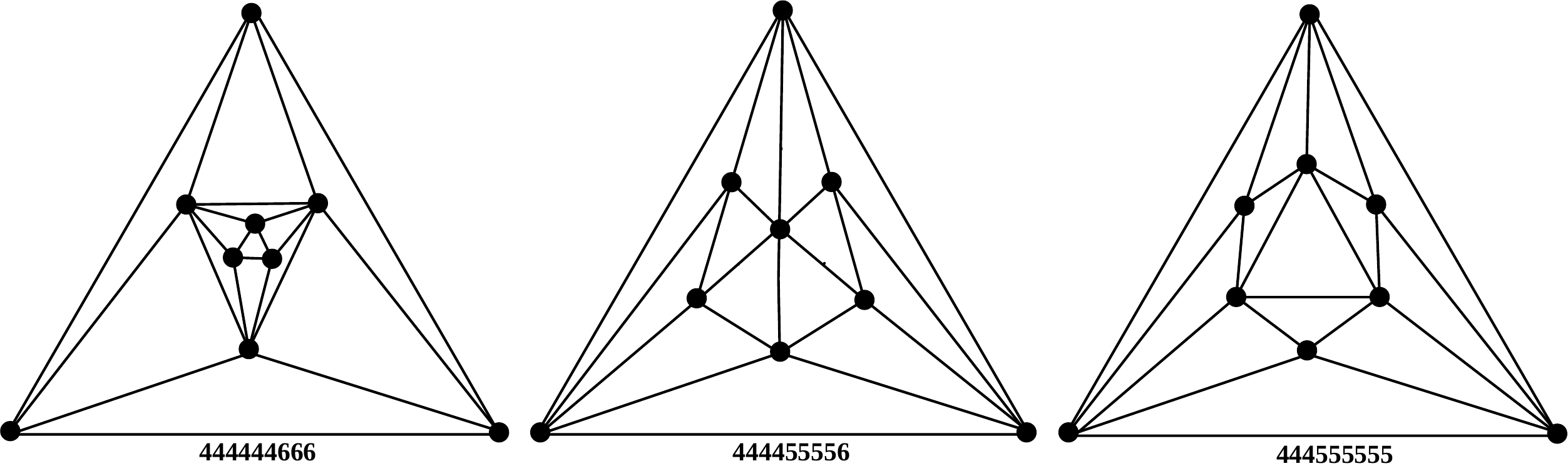}

  \includegraphics [width=240pt]{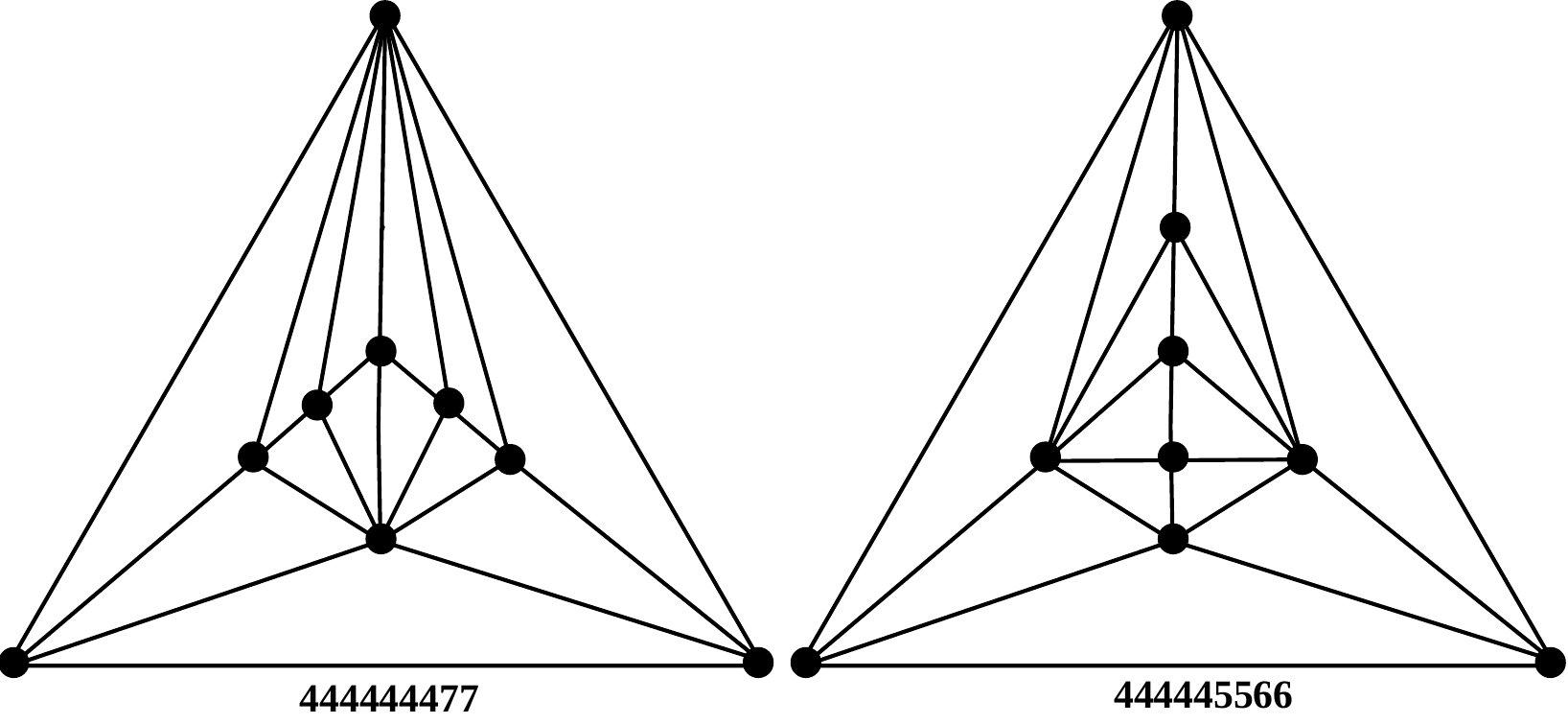}

  \textbf{Figure 4.8.} All five maximal planar graphs whose order is 9 and $\delta(G)=4$
\end{center}
\begin{center}
  \includegraphics [width=340pt]{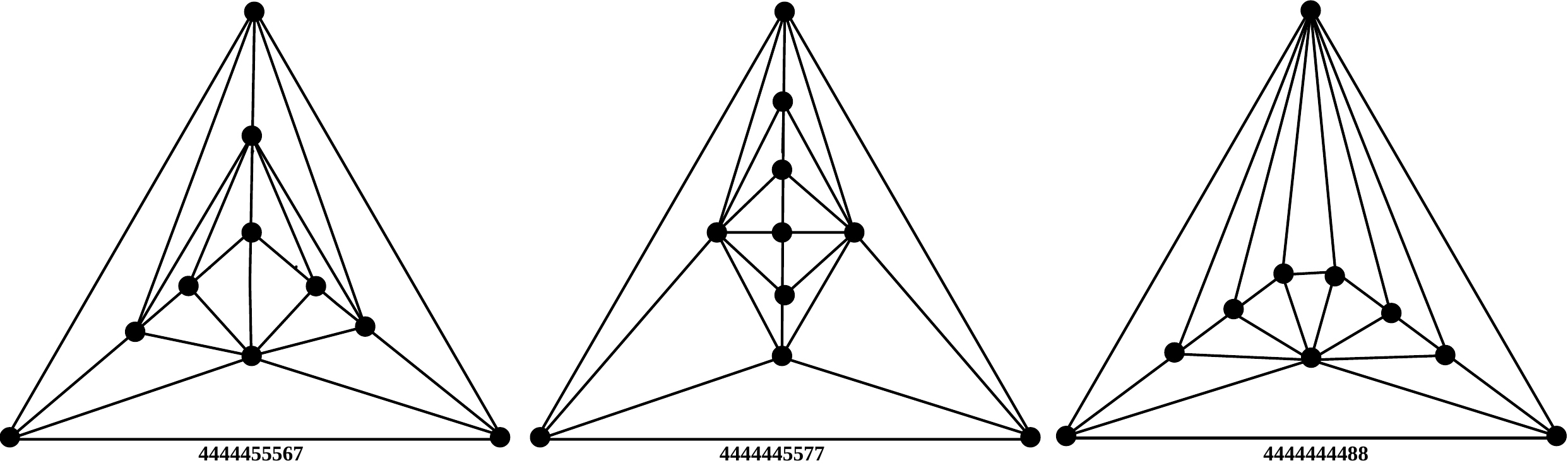}
  \includegraphics [width=340pt]{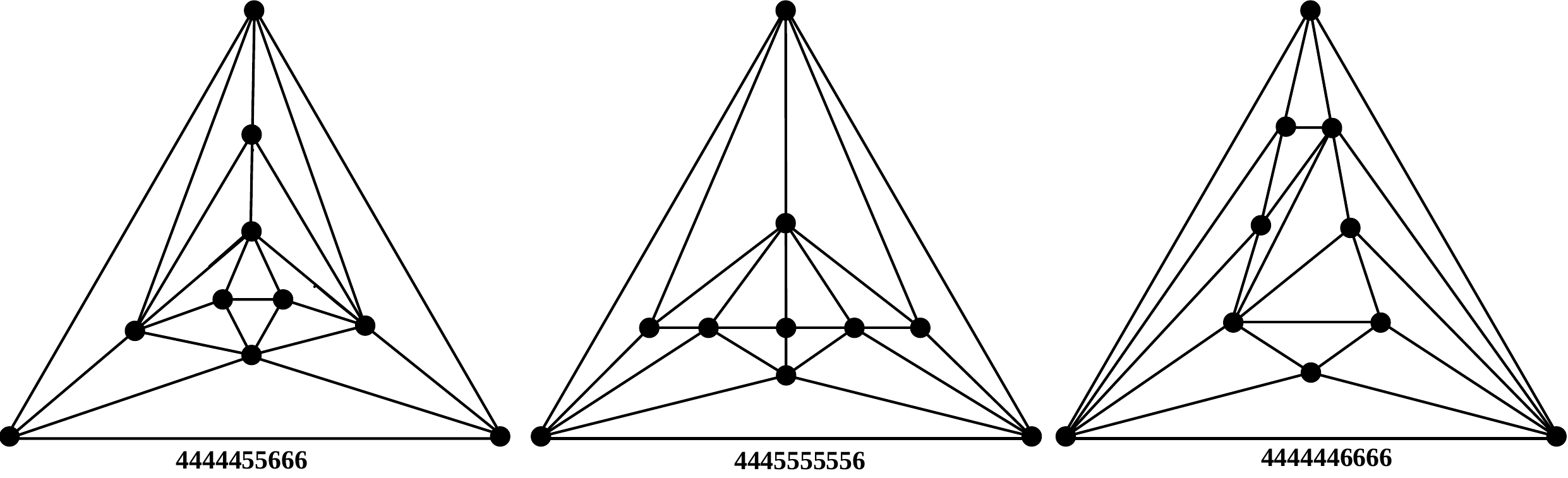}
  \includegraphics [width=340pt]{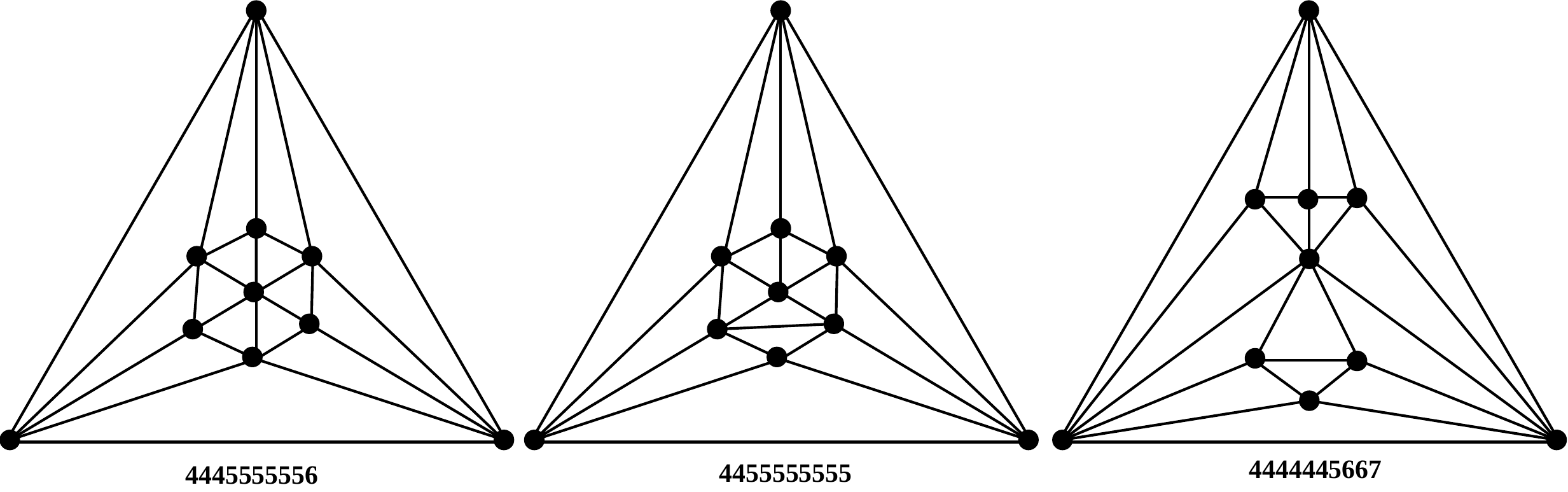}
  \includegraphics [width=340pt]{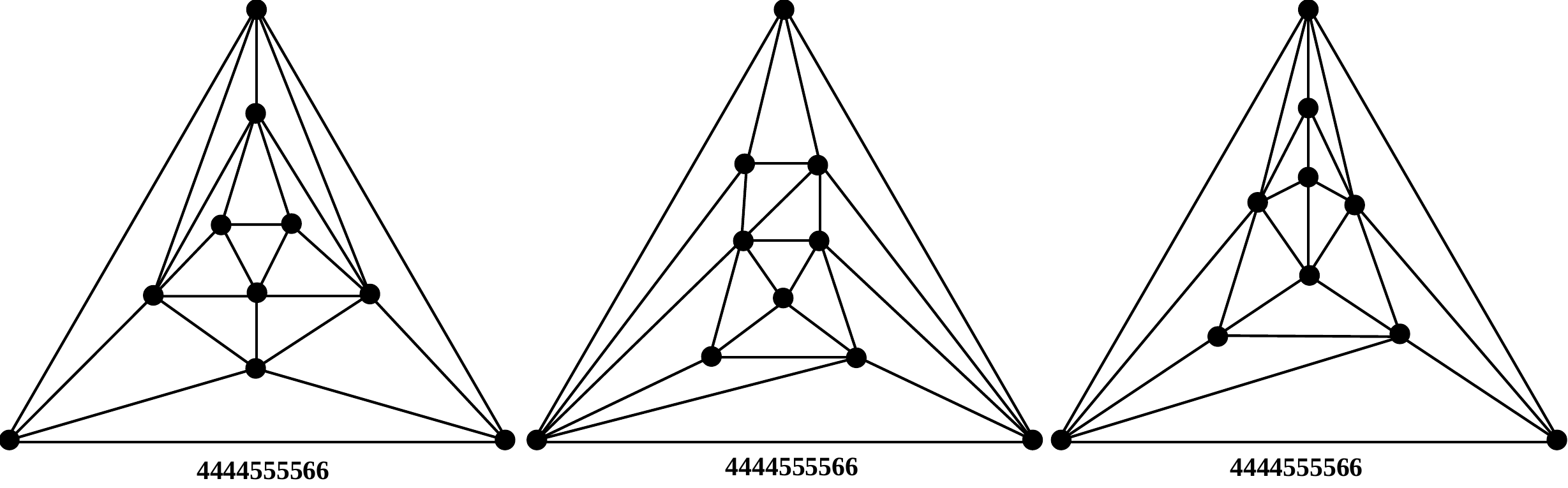}
  \includegraphics [width=120pt]{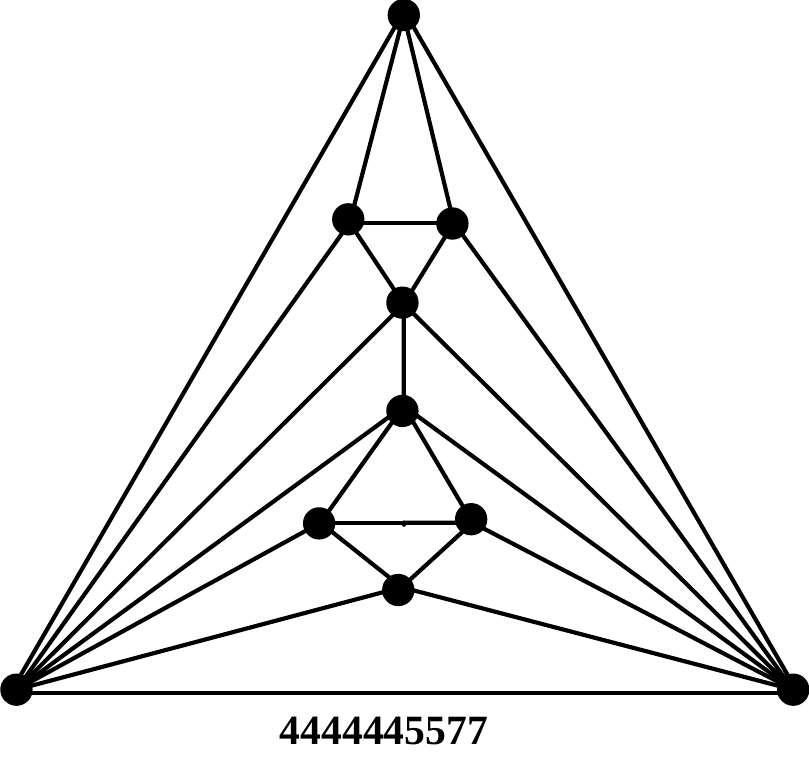}

  \textbf{Figure 4.9.} All 13 maximal planar graphs with $\delta(G)=4$ whose order is 10
\end{center}
\begin{center}
  \includegraphics [width=340pt]{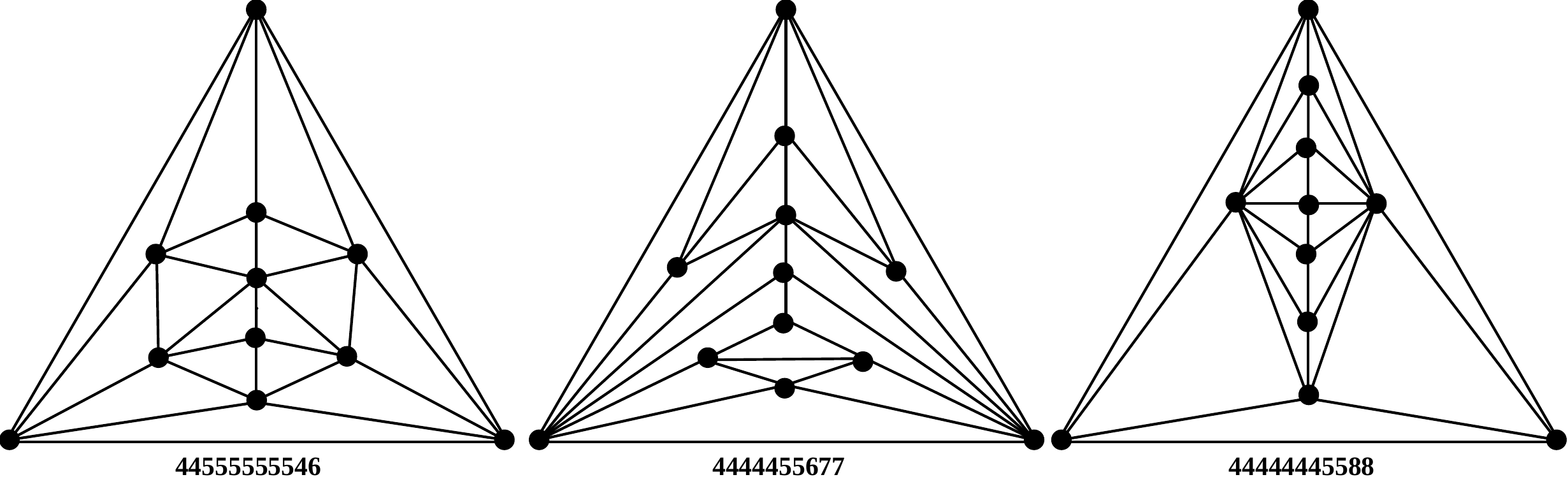}
  \includegraphics [width=340pt]{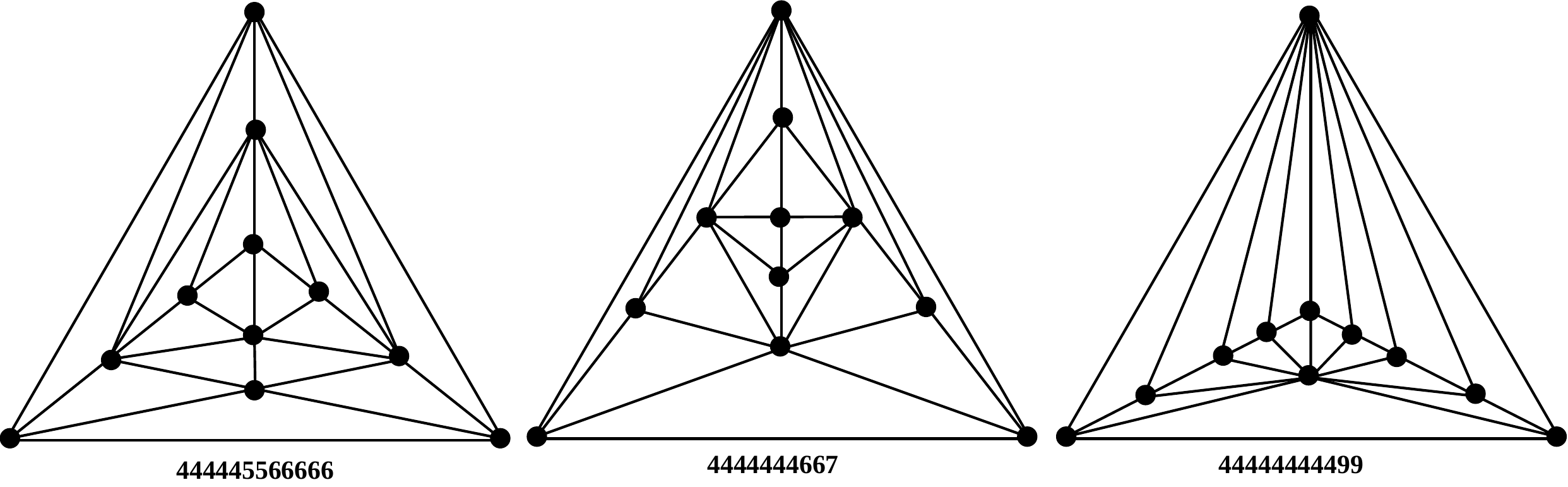}
  \includegraphics [width=340pt]{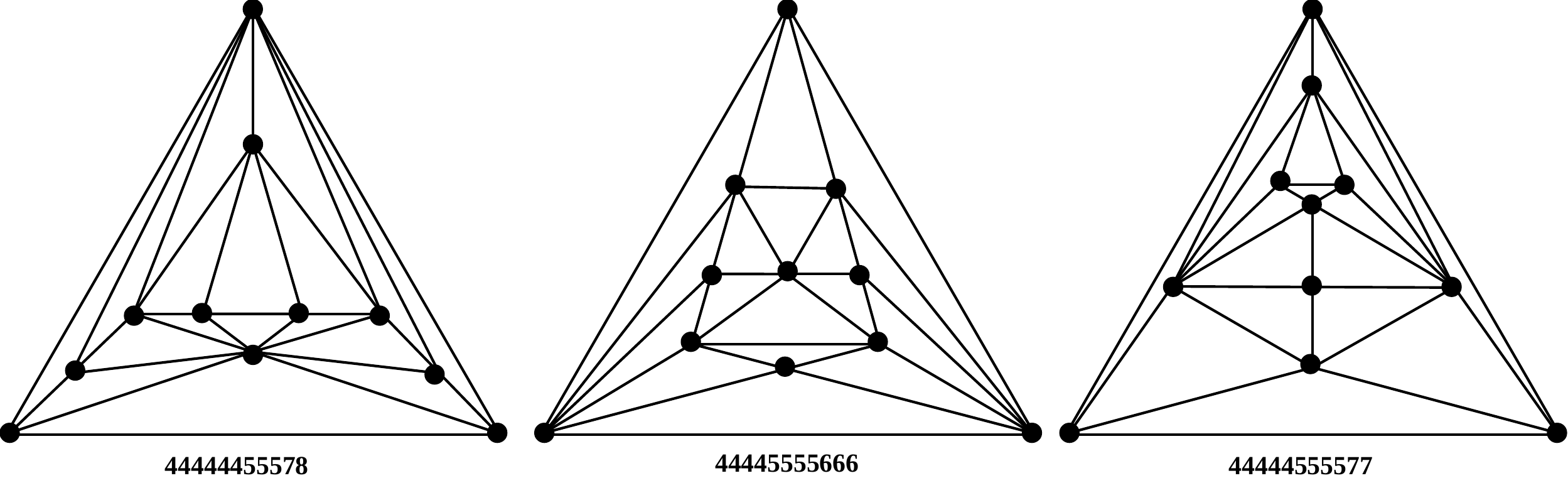}
  \includegraphics [width=340pt]{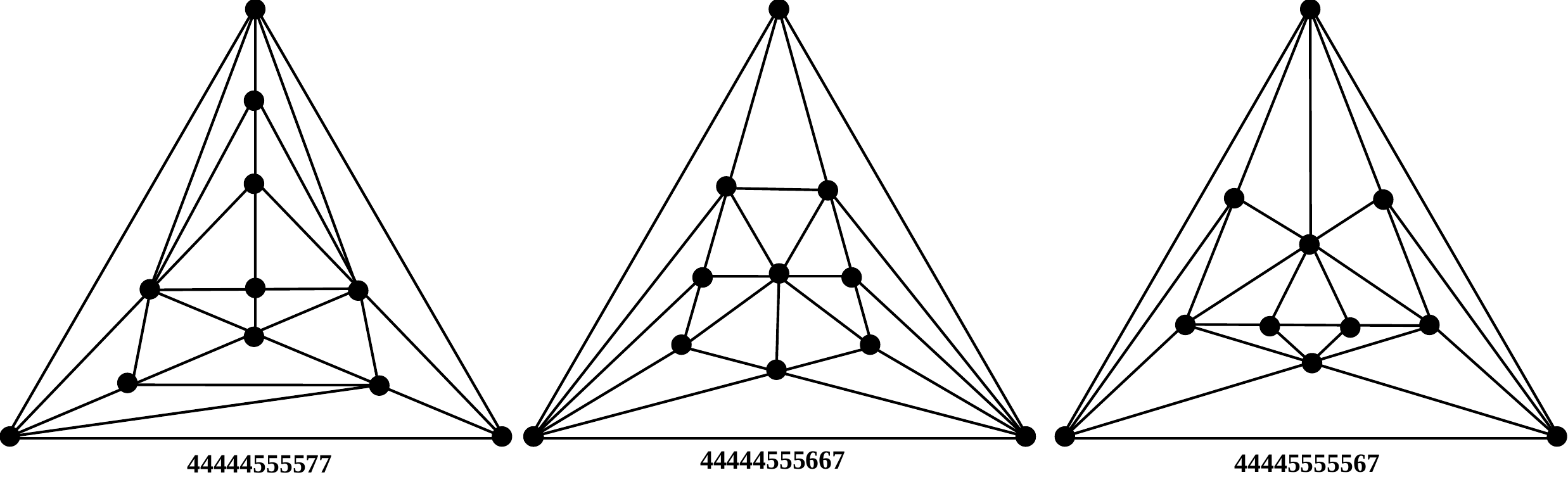}
  \includegraphics [width=340pt]{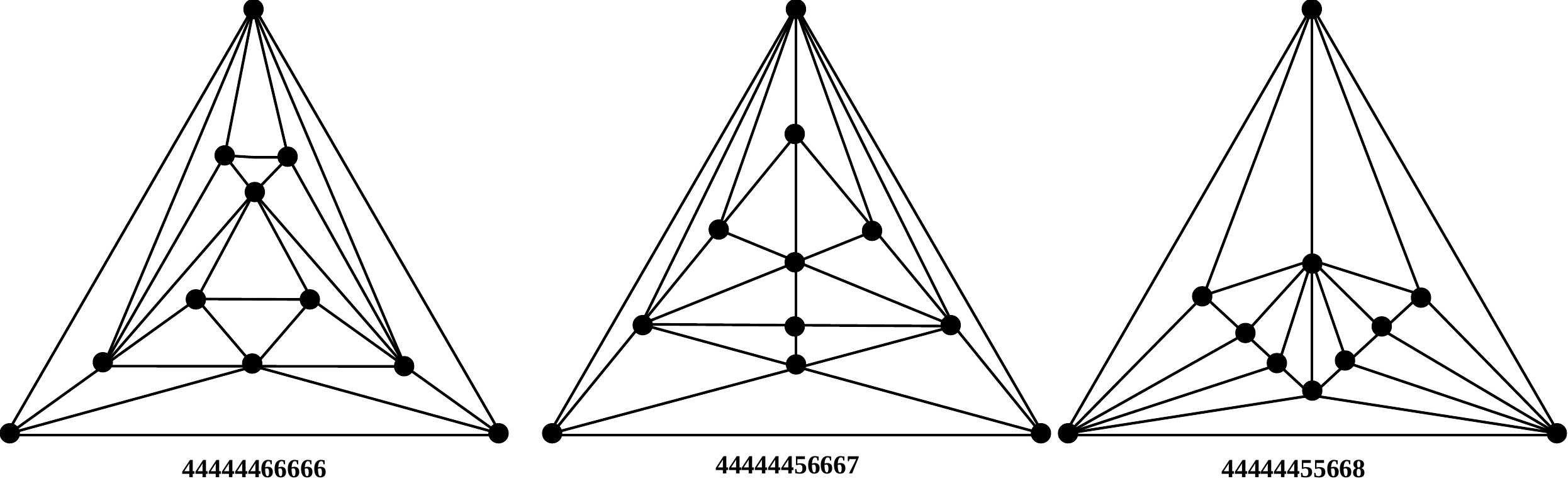}
  \includegraphics [width=340pt]{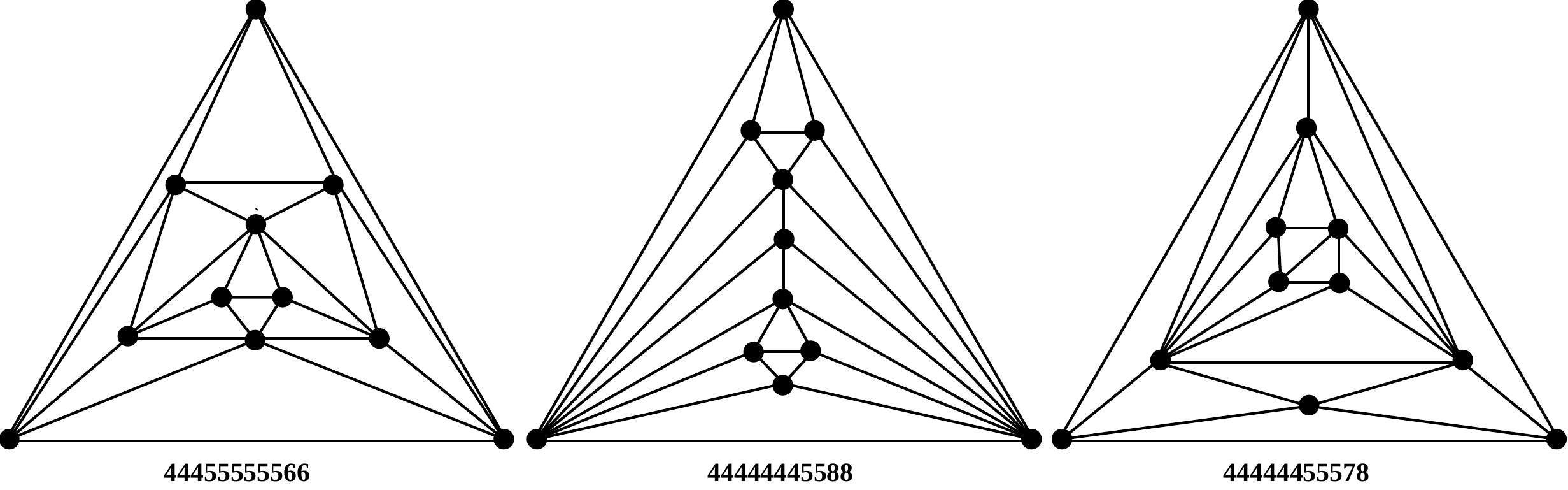}
  \includegraphics [width=340pt]{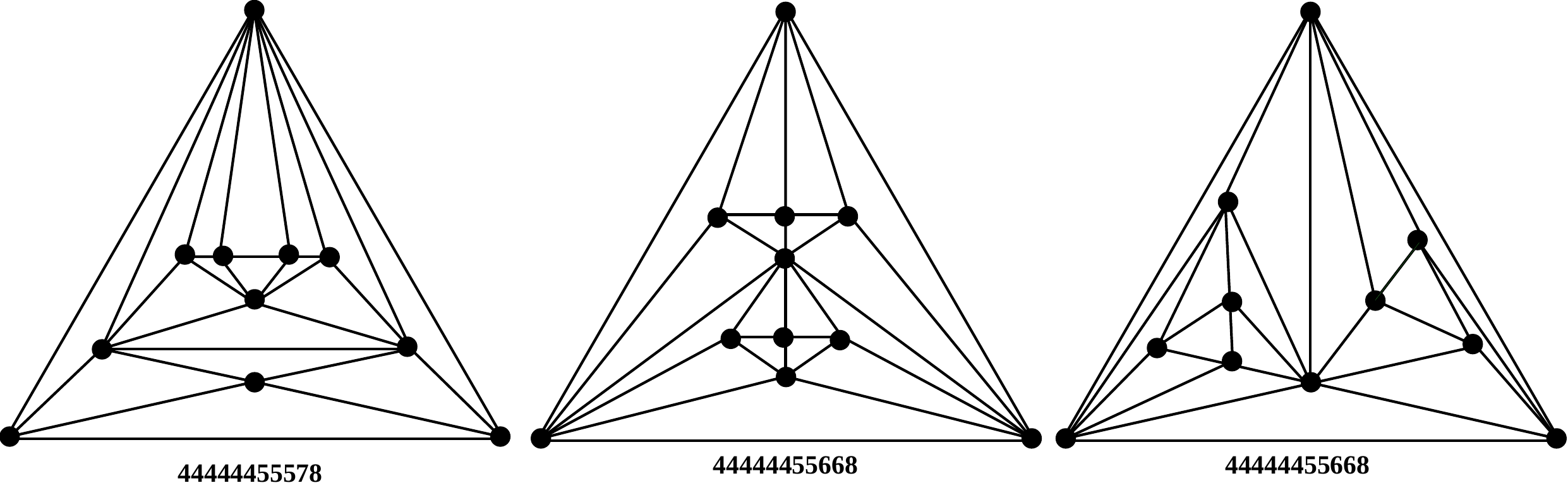}
  \includegraphics [width=340pt]{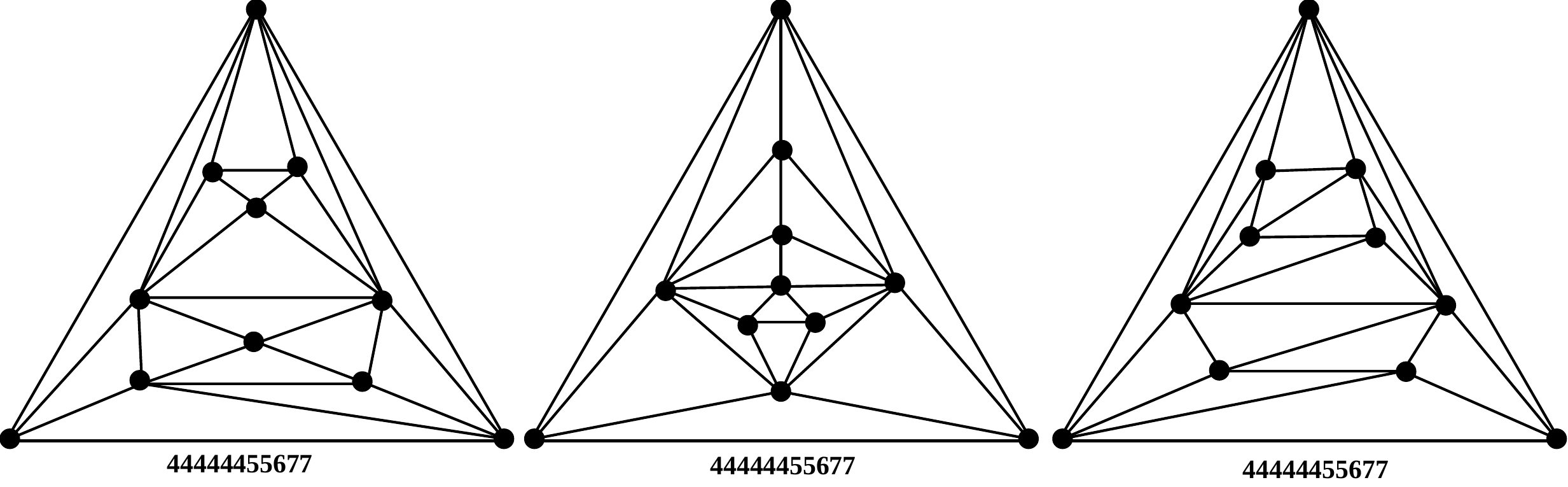}
  \includegraphics [width=340pt]{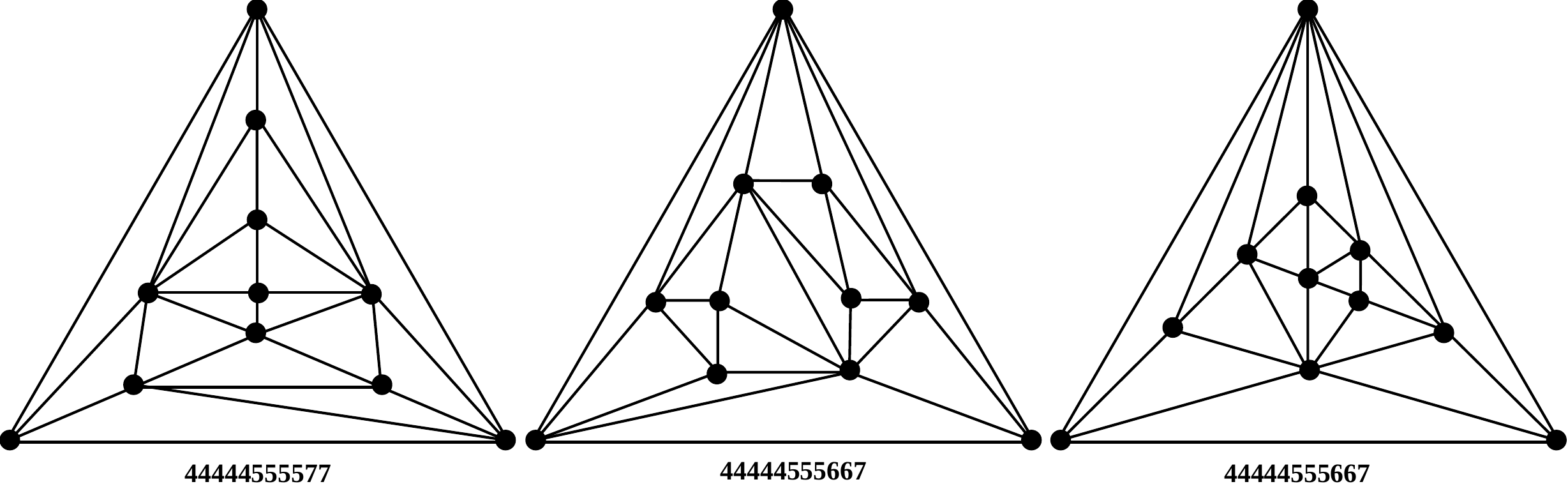}
  \includegraphics [width=340pt]{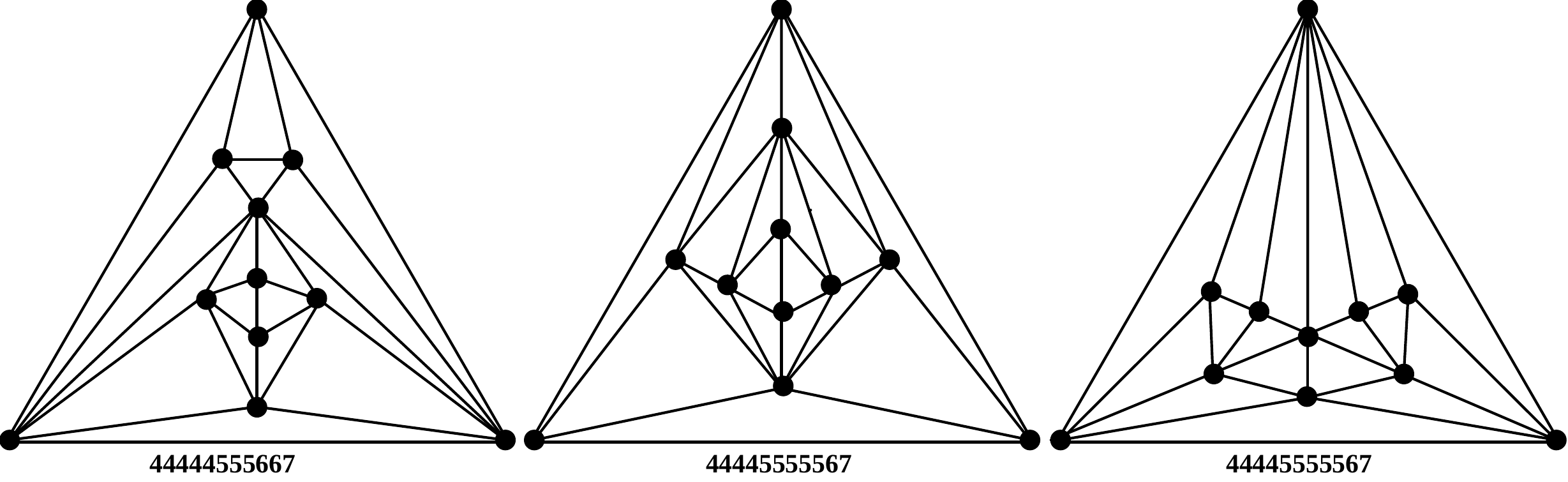}
  \includegraphics [width=340pt]{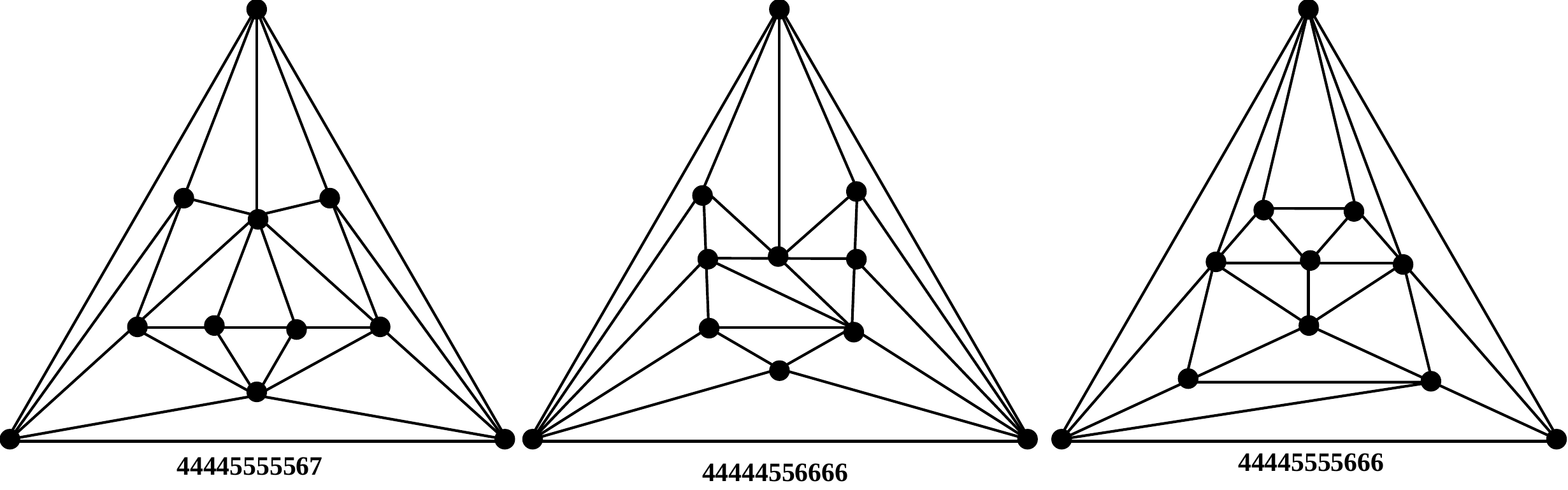}
  \includegraphics [width=340pt]{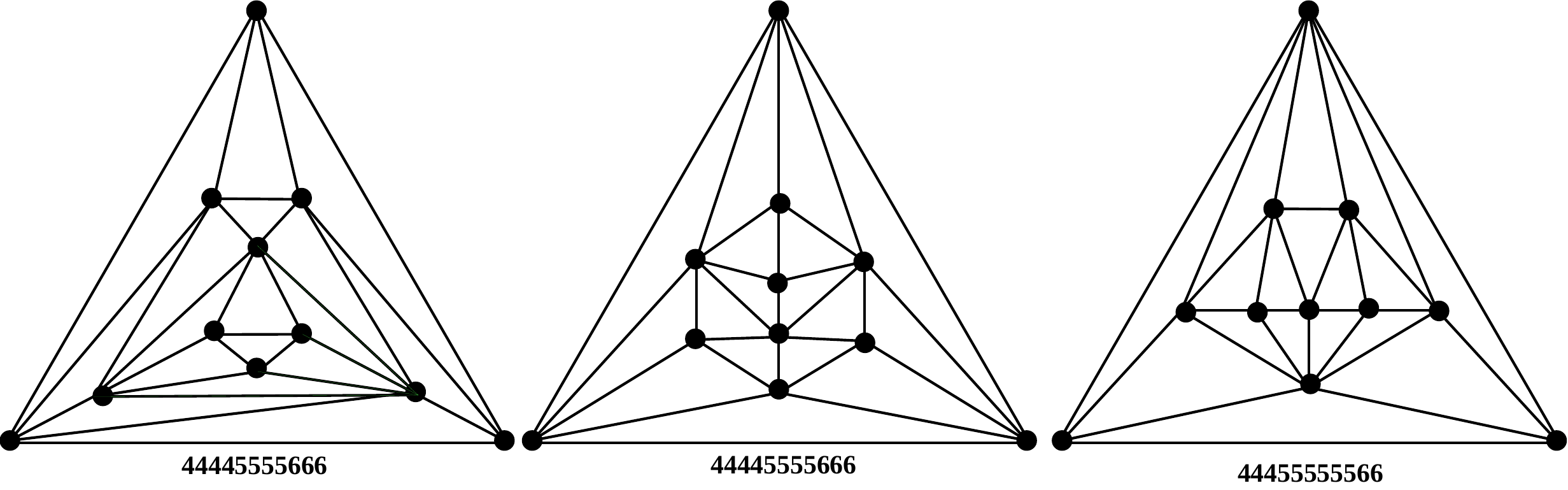}

  \textbf{Figure 4.10.} All 34 maximal planar graphs with $\delta(G)=4$ whose order is 11
\end{center}

¡¡

¡¡
\begin{center}
    \includegraphics [width=340pt]{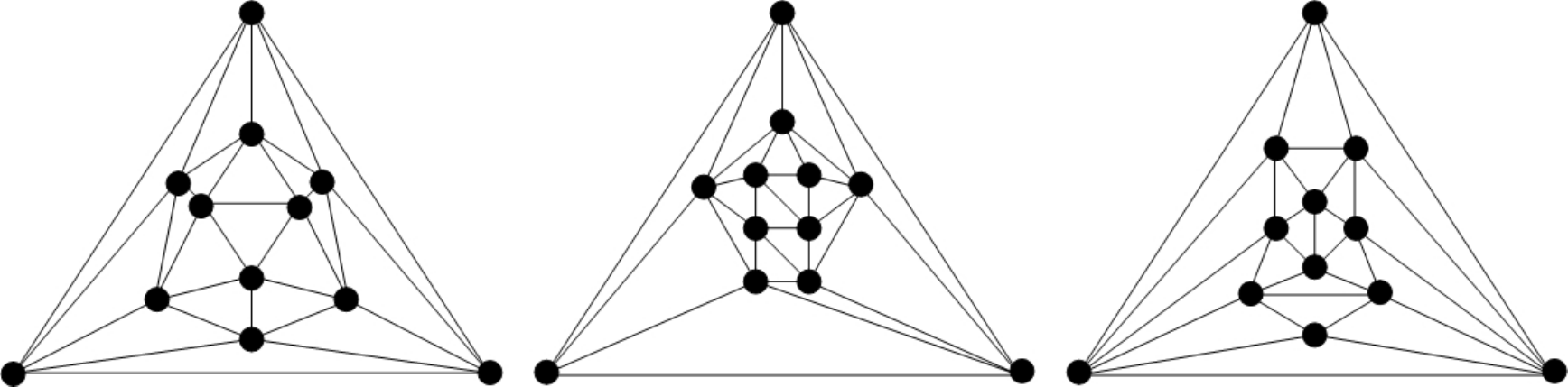}
    \includegraphics [width=340pt]{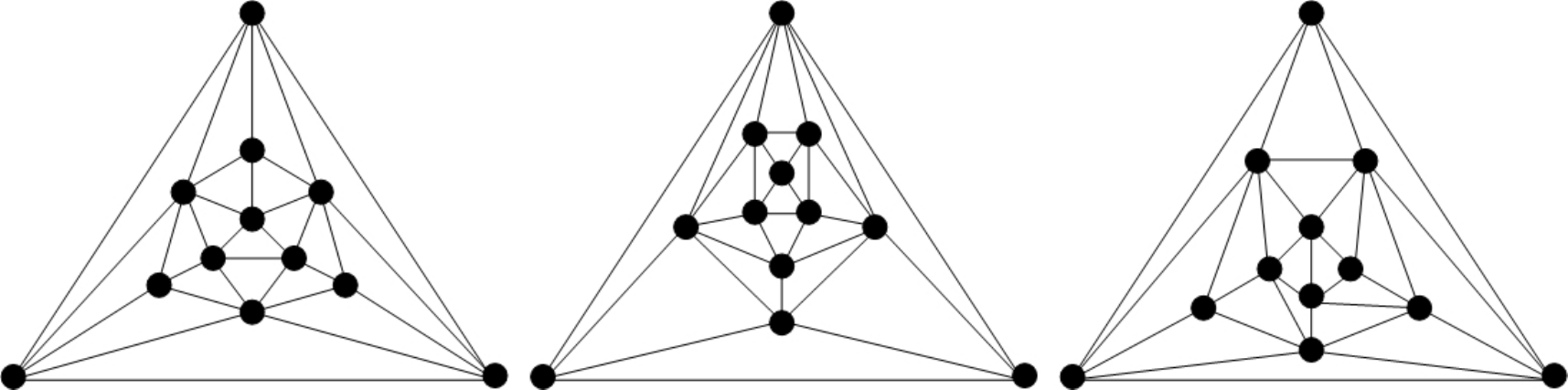}
    \includegraphics [width=340pt]{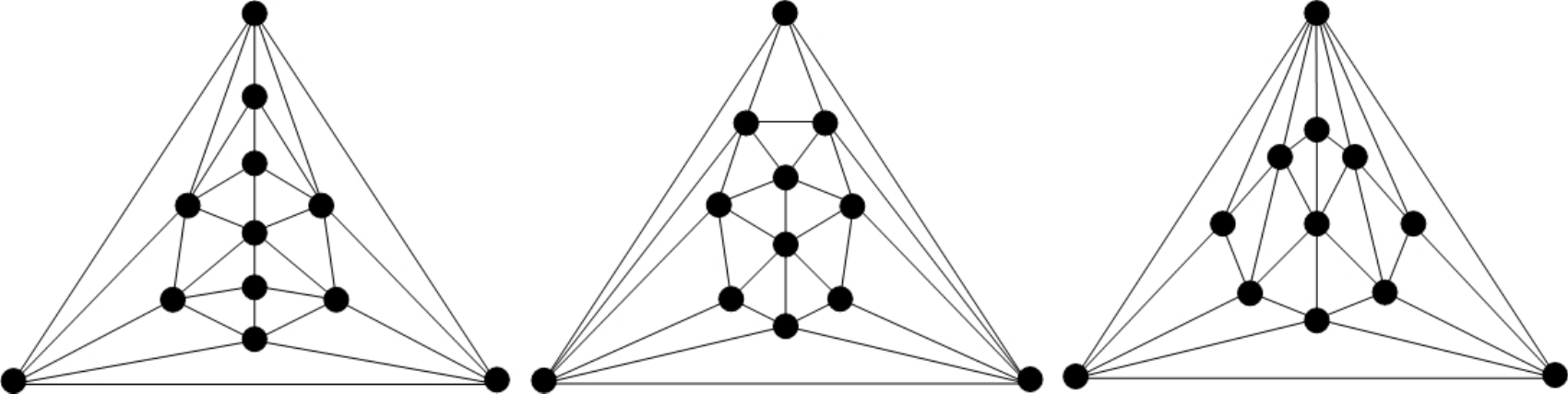}
    \includegraphics [width=340pt]{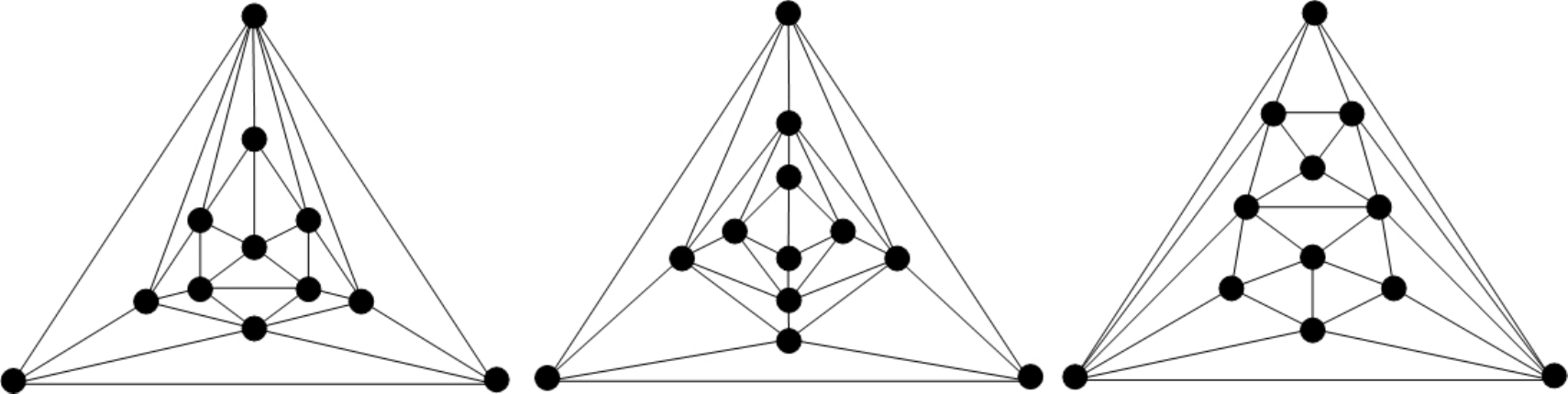}
    \includegraphics [width=340pt]{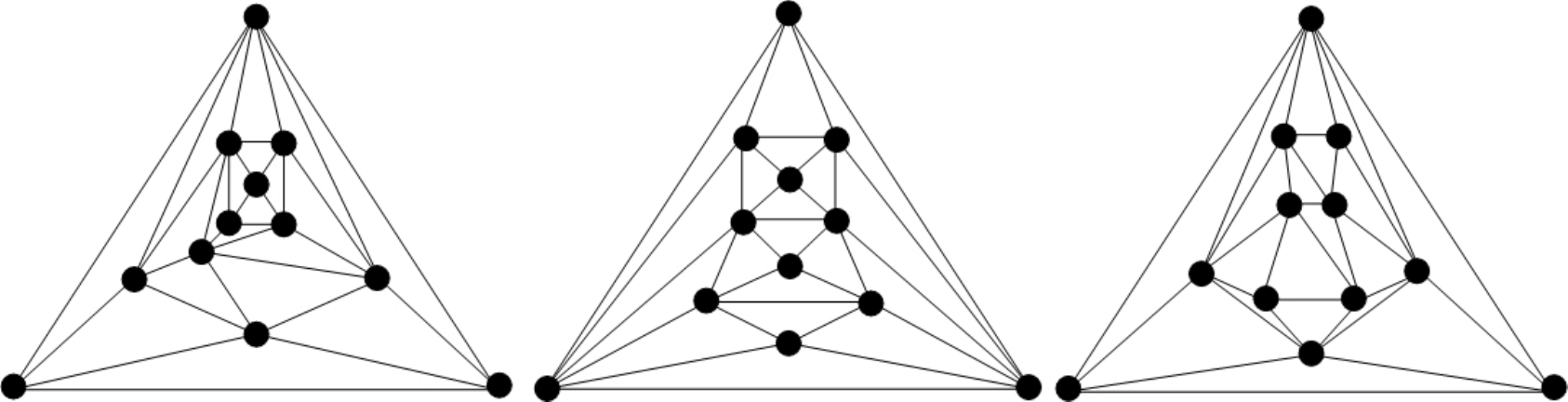}
    \includegraphics [width=340pt]{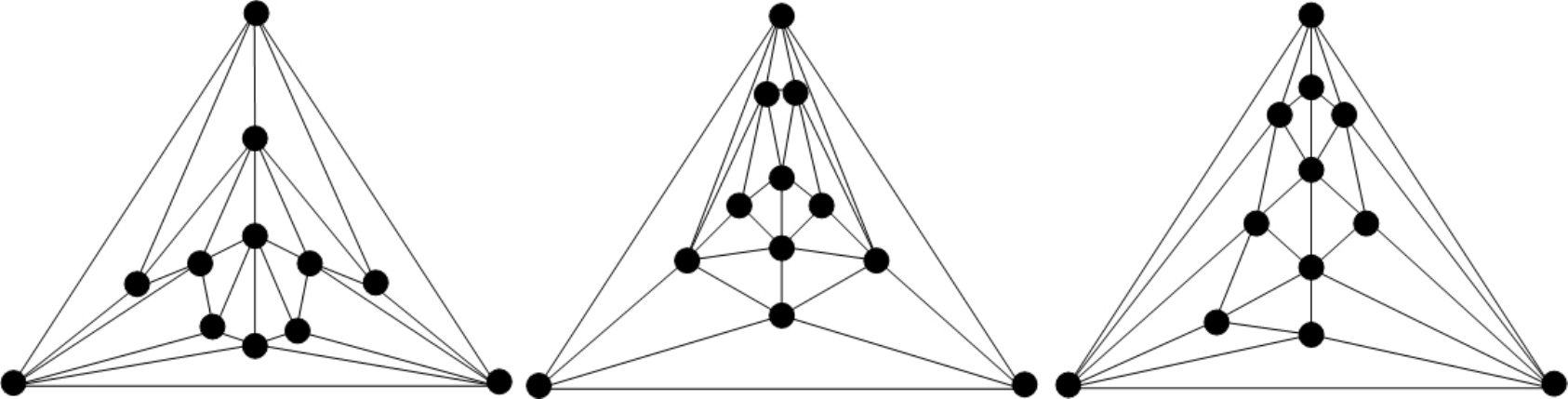}
    \includegraphics [width=340pt]{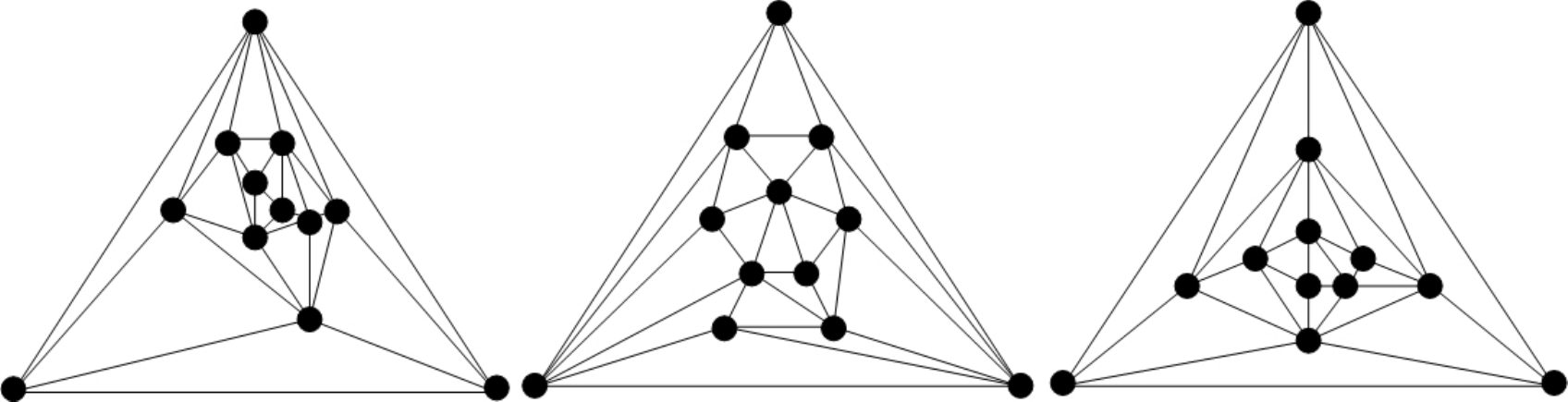}
    \includegraphics [width=340pt]{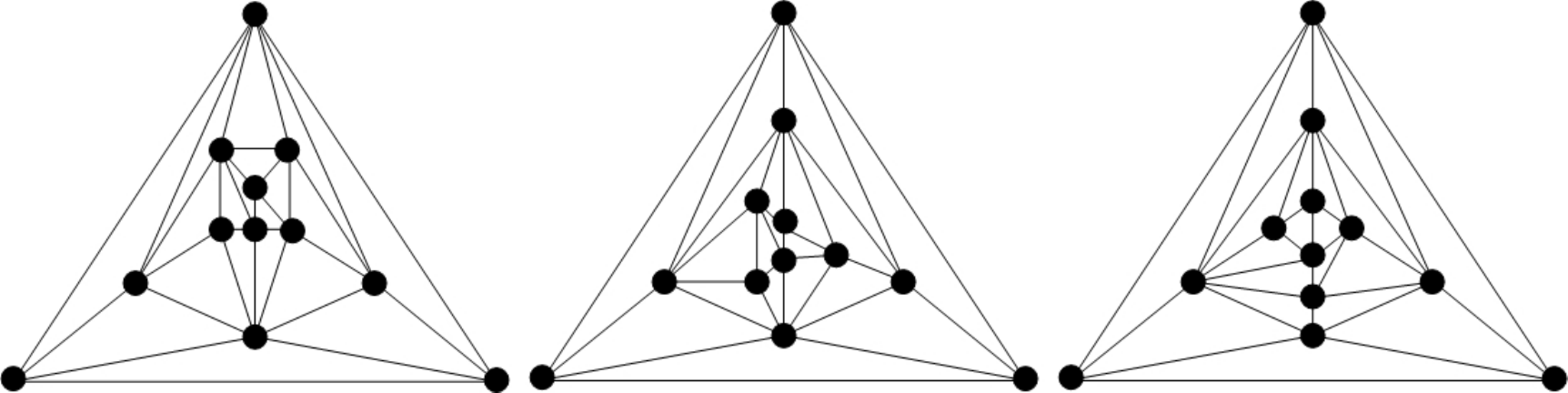}
    \includegraphics [width=340pt]{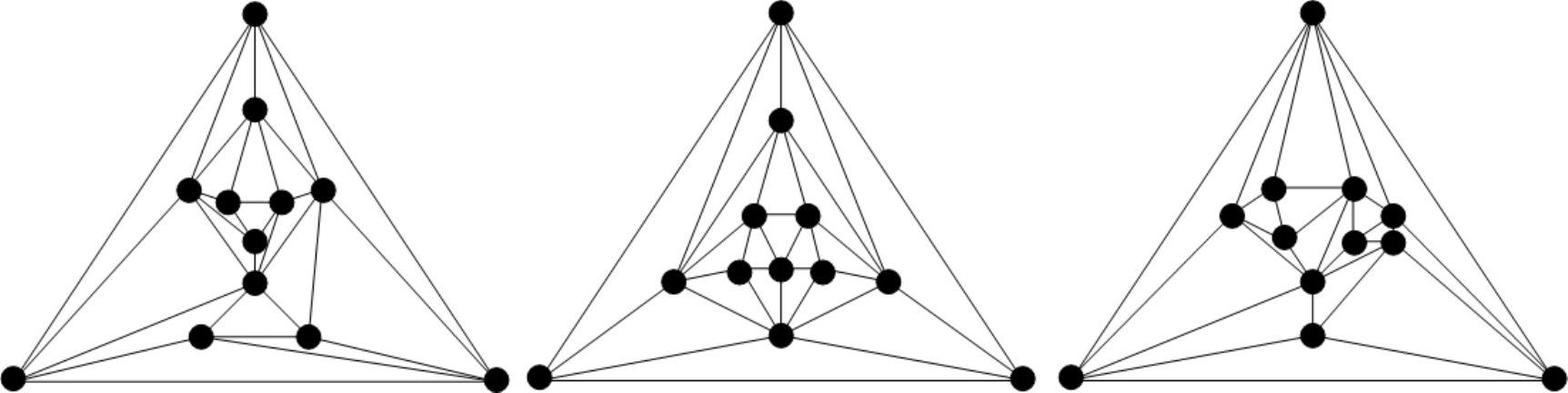}
    \includegraphics [width=340pt]{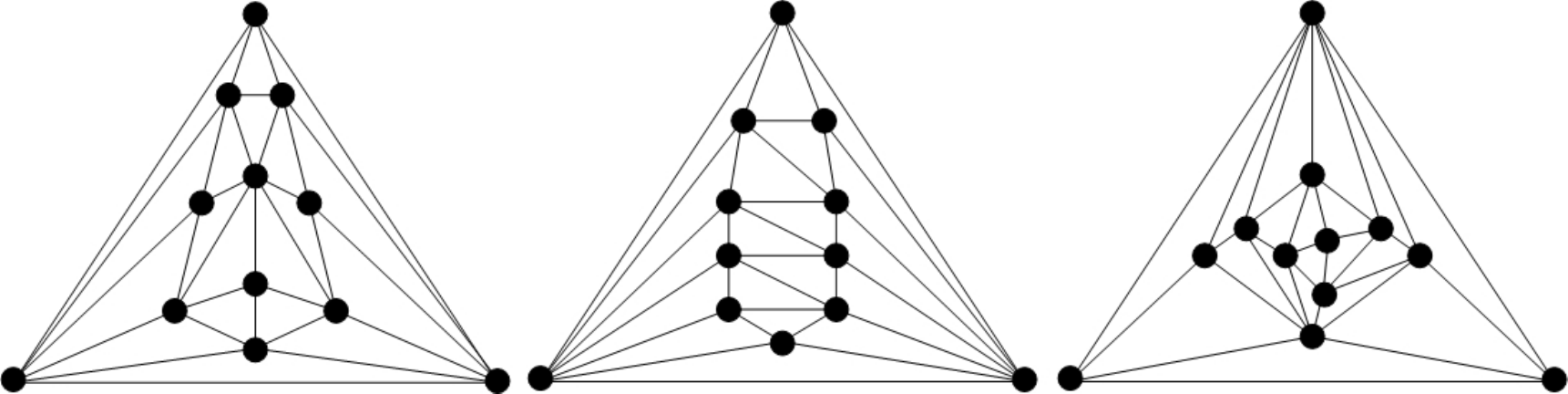}
    \includegraphics [width=340pt]{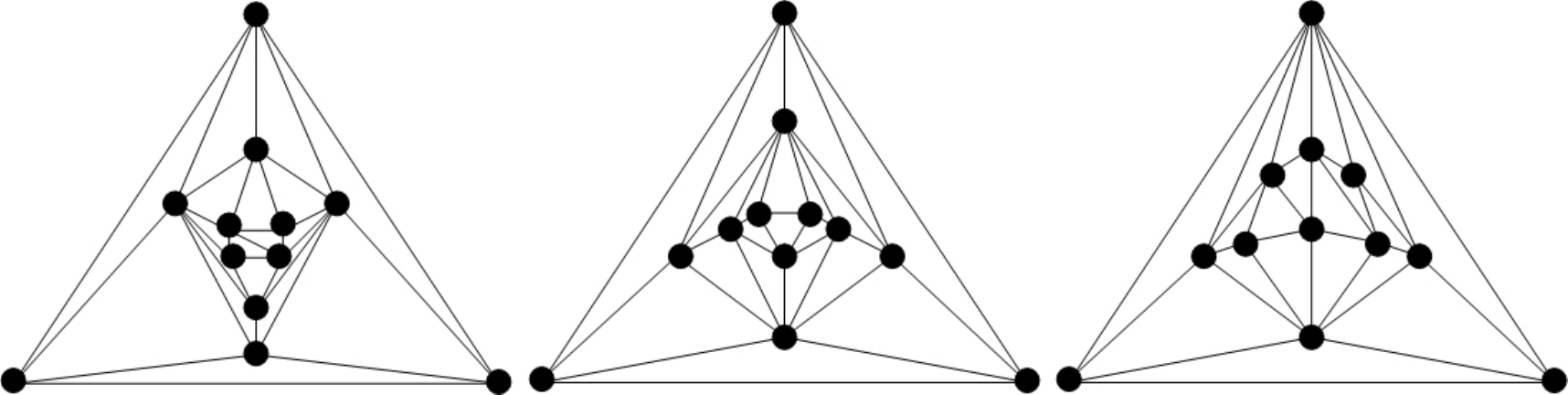}
    \includegraphics [width=340pt]{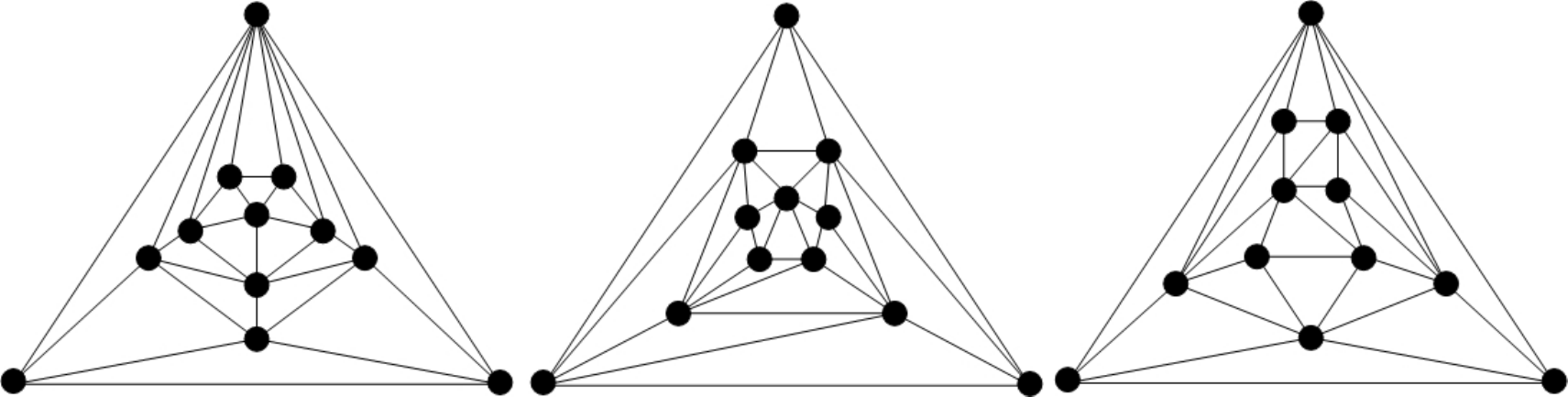}
    \includegraphics [width=340pt]{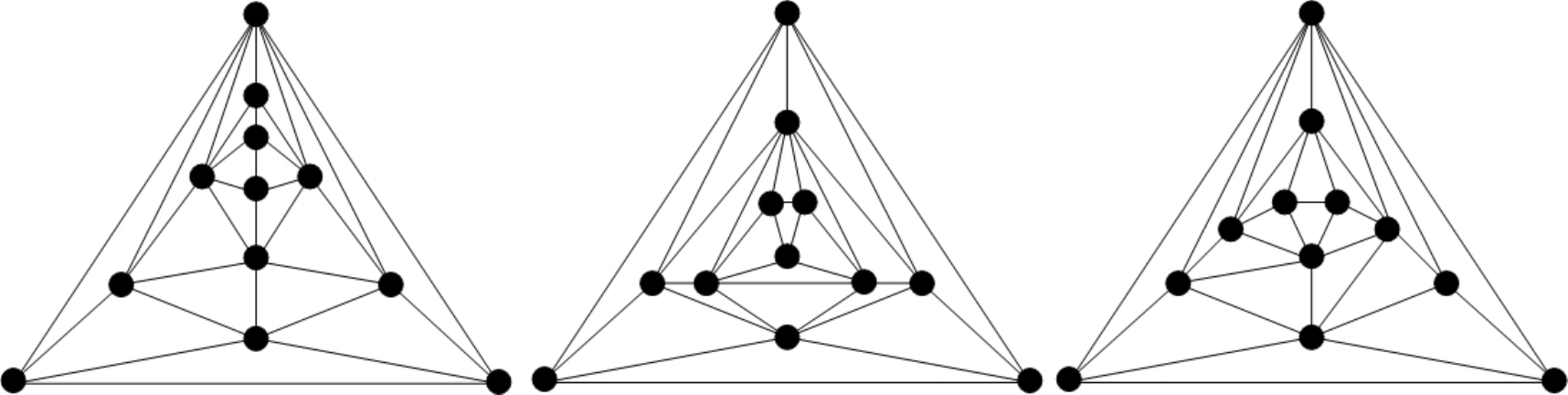}
    \includegraphics [width=340pt]{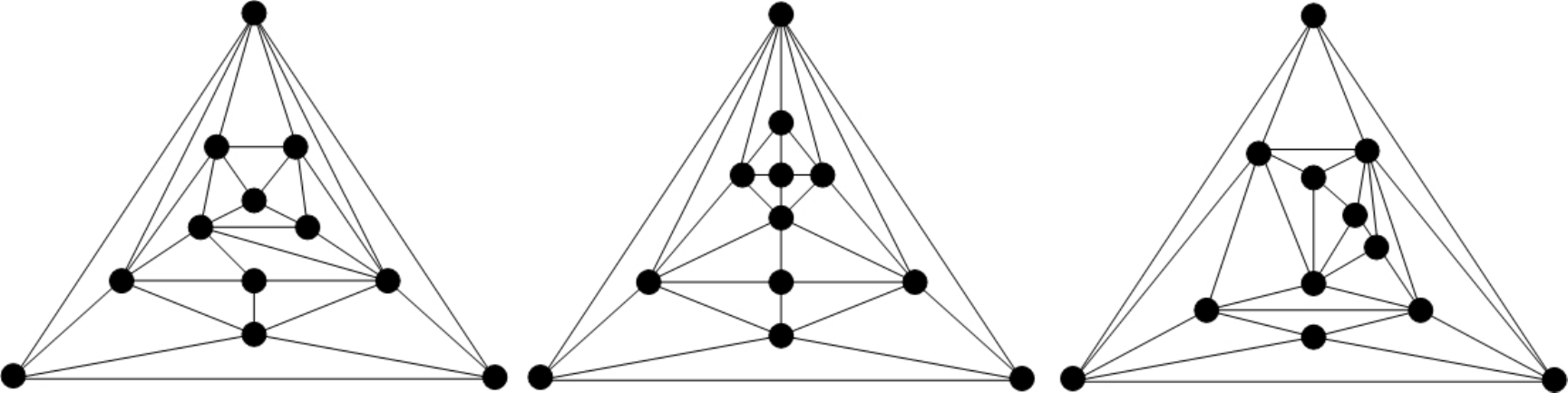}
    \includegraphics [width=340pt]{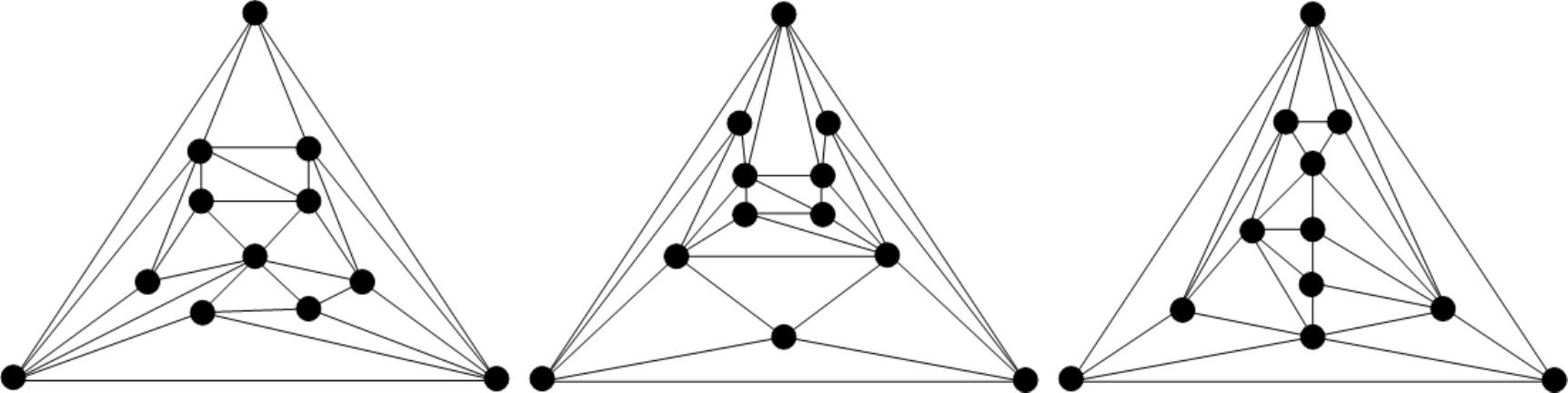}
    \includegraphics [width=340pt]{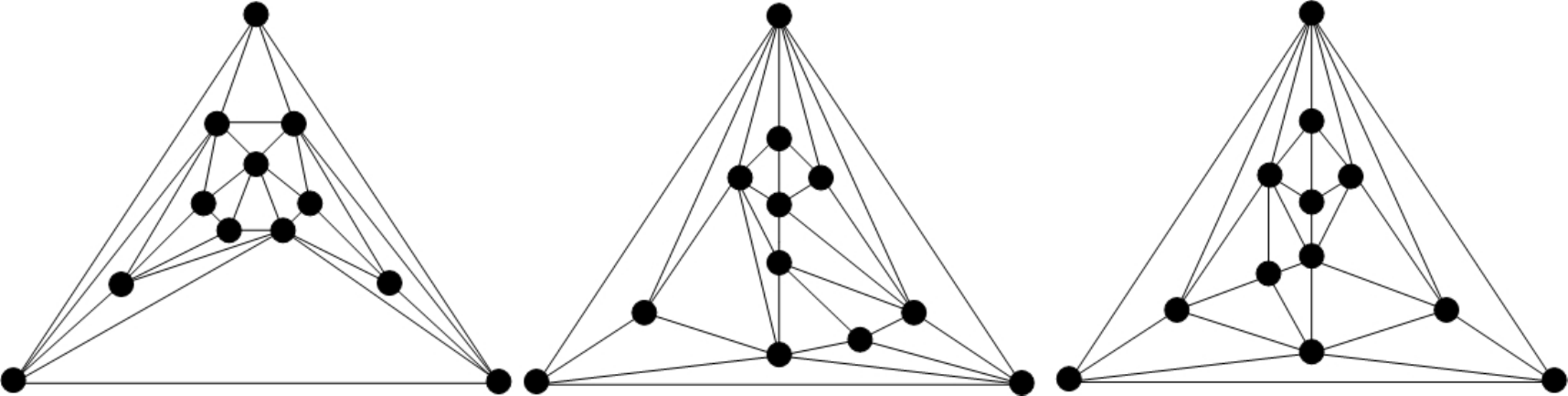}
    \includegraphics [width=340pt]{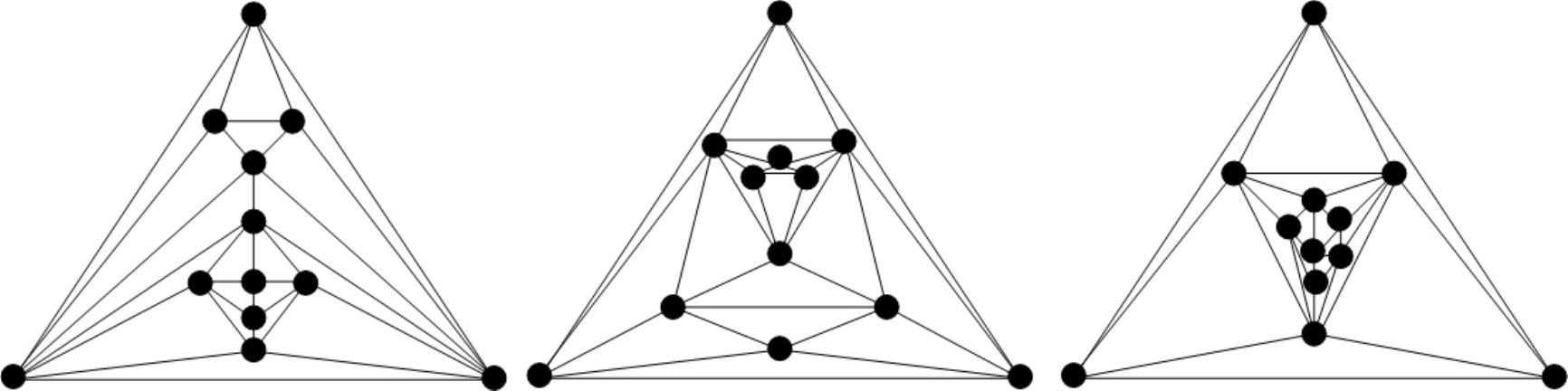}
    \includegraphics [width=340pt]{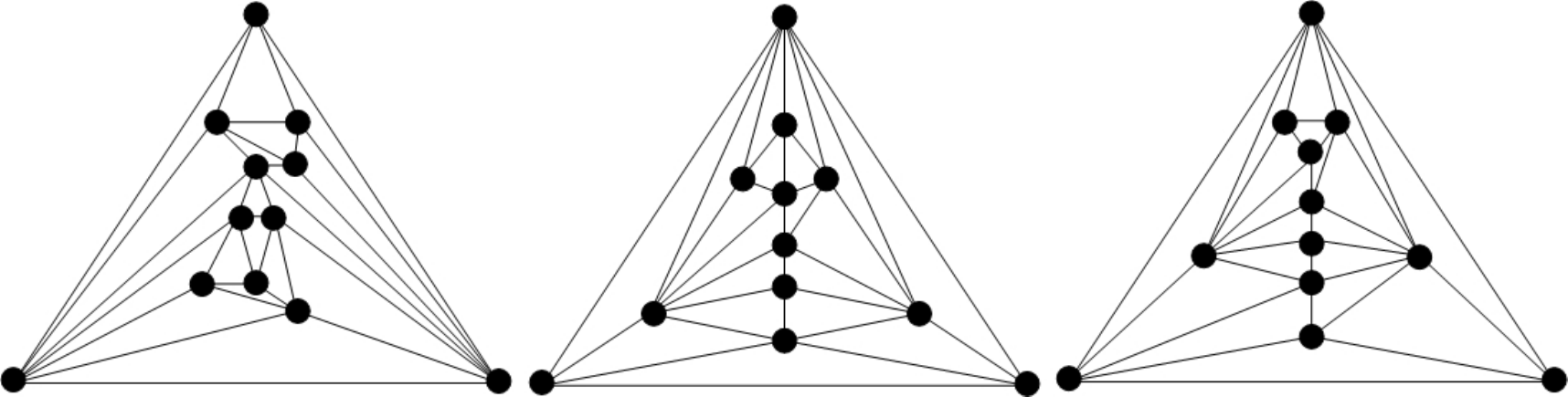}
    \includegraphics [width=340pt]{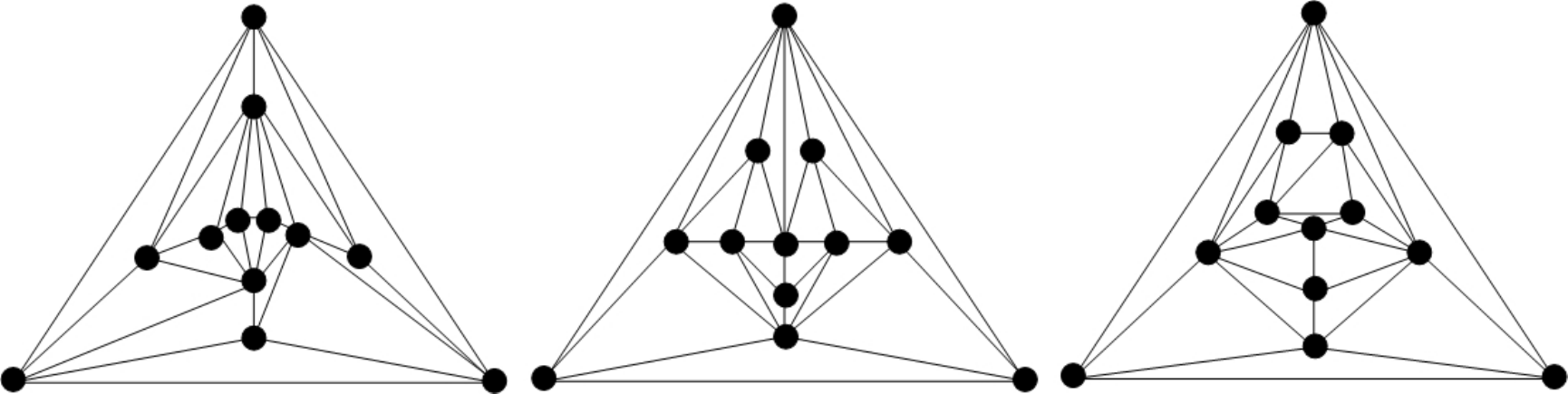}
    \includegraphics [width=340pt]{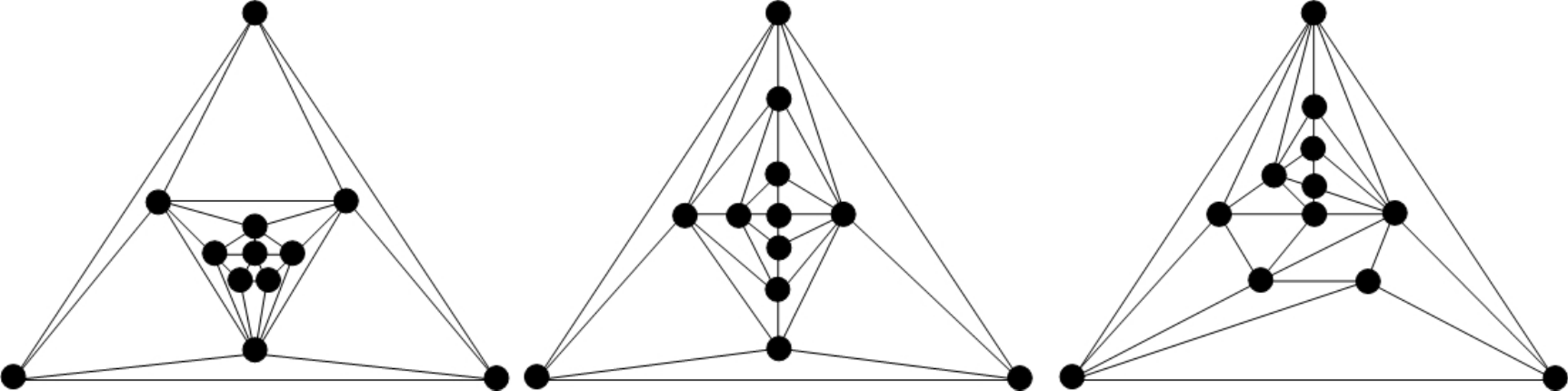}
    \includegraphics [width=340pt]{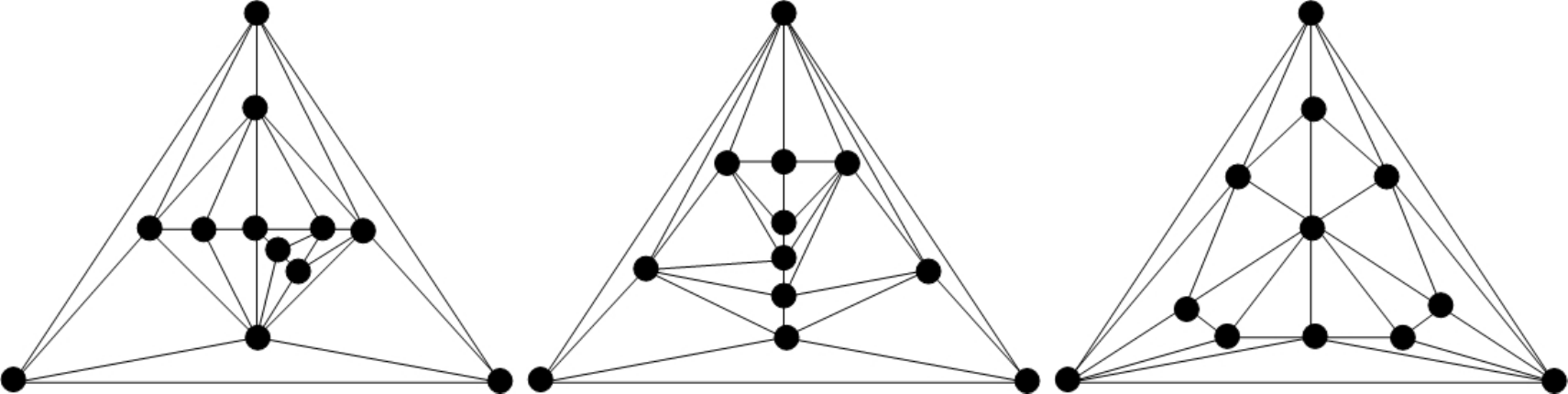}
    \includegraphics [width=340pt]{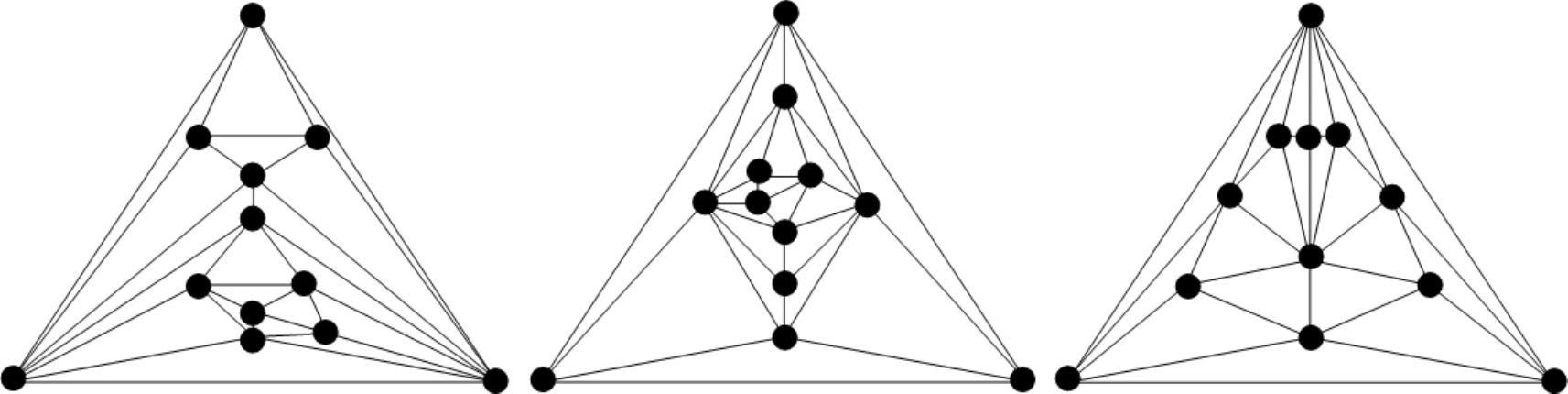}
    \includegraphics [width=340pt]{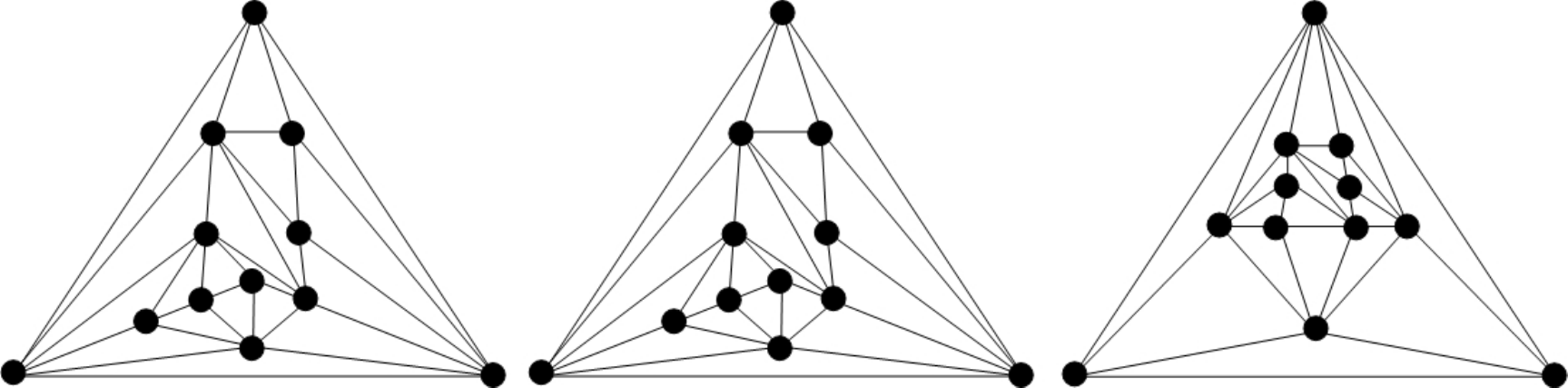}
    \includegraphics [width=340pt]{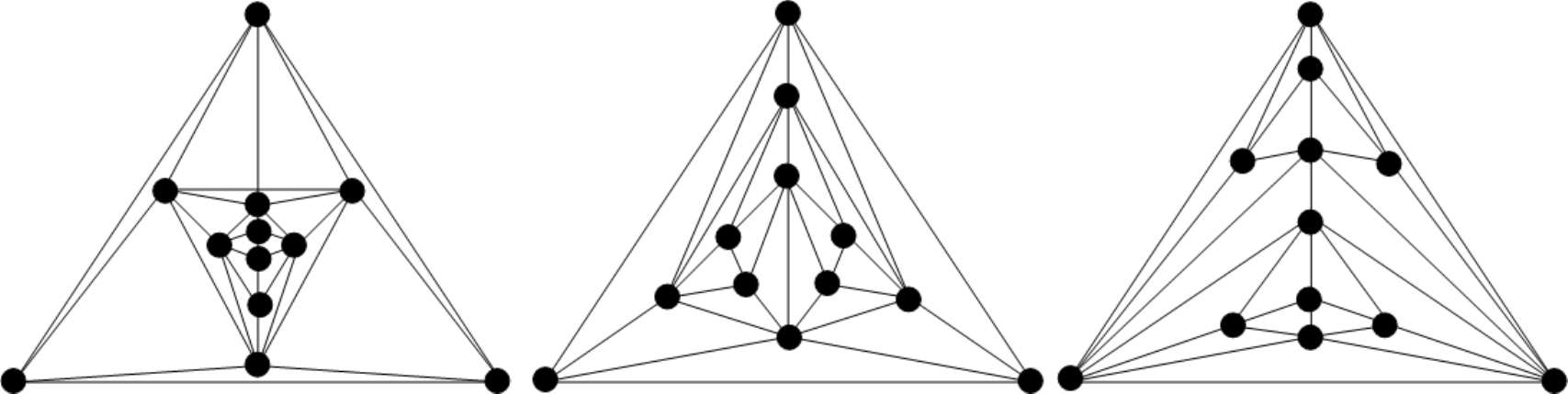}
    \includegraphics [width=340pt]{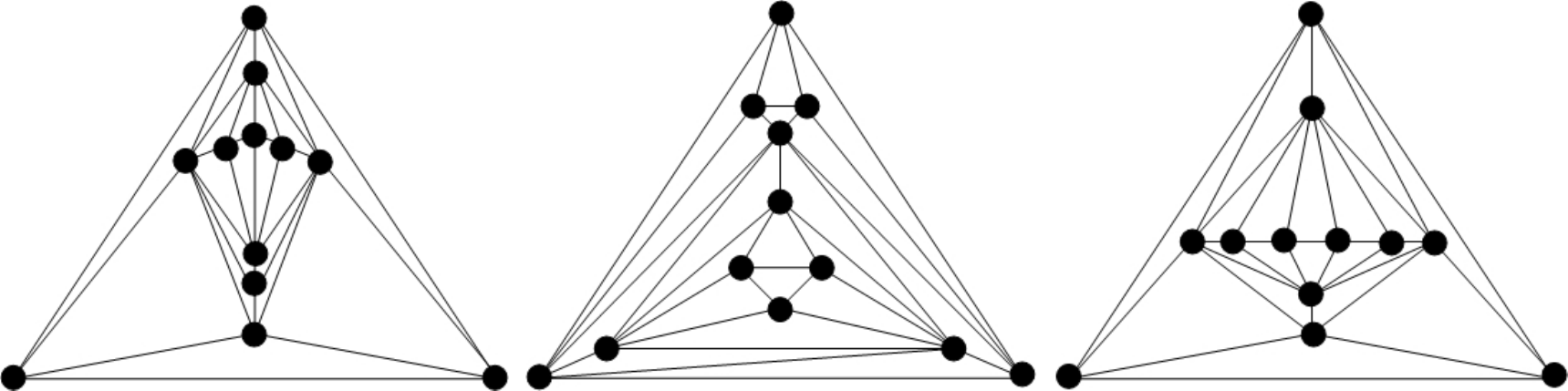}
    \includegraphics [width=340pt]{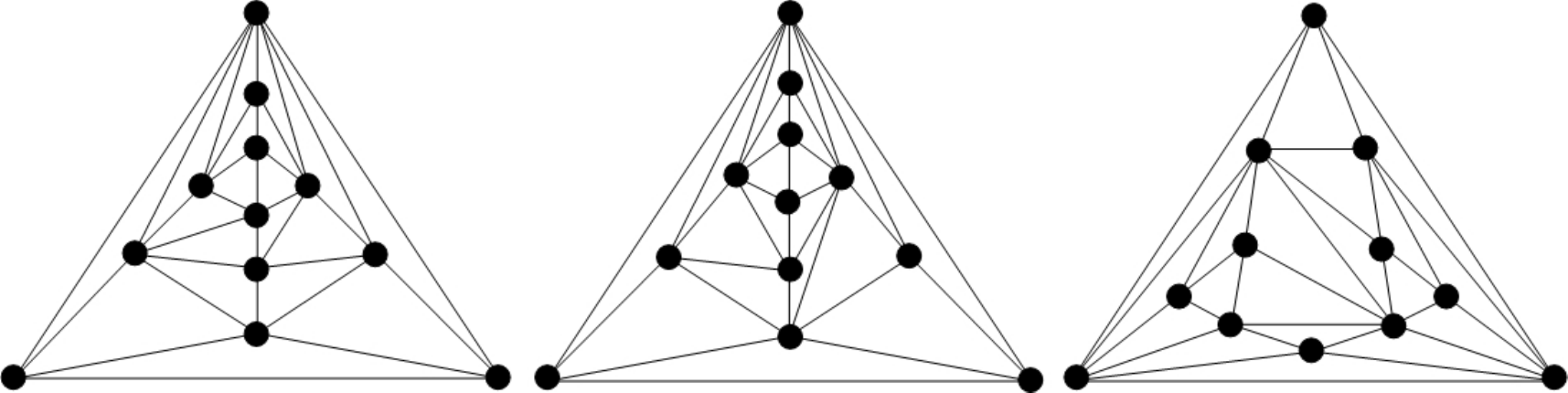}
    \includegraphics [width=340pt]{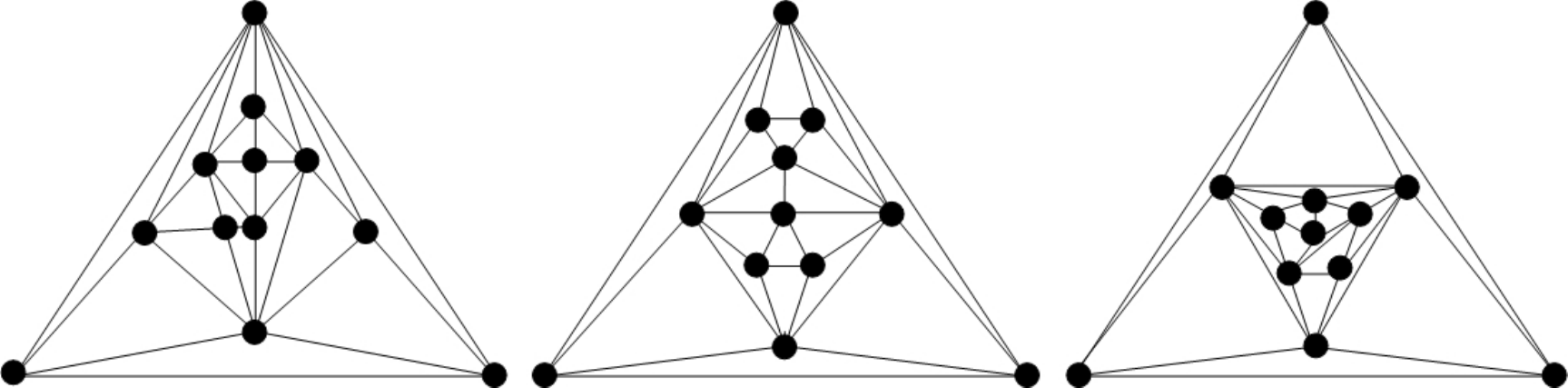}
    \includegraphics [width=340pt]{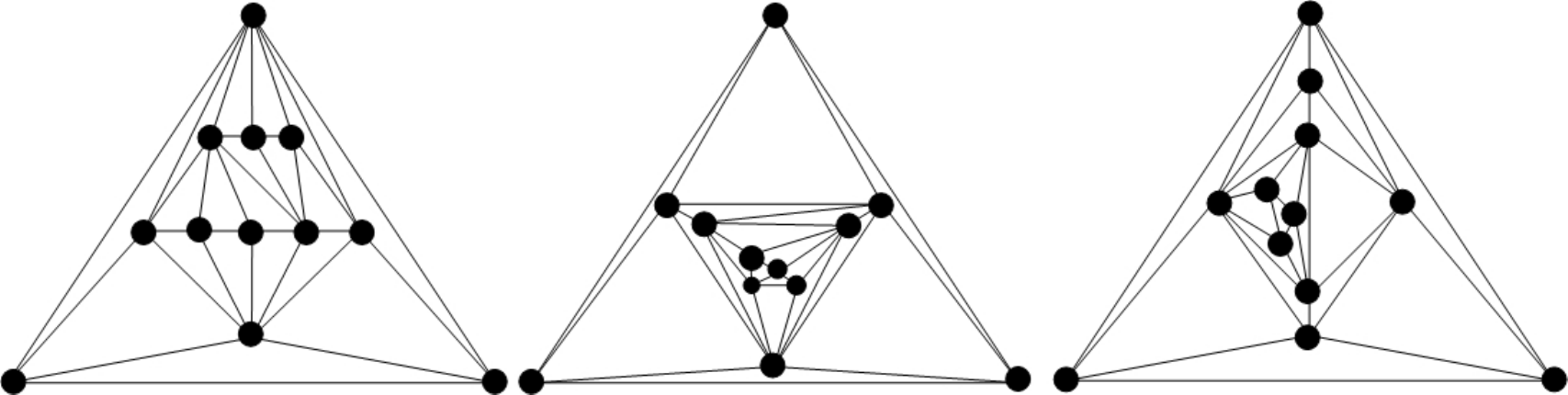}
    \includegraphics [width=340pt]{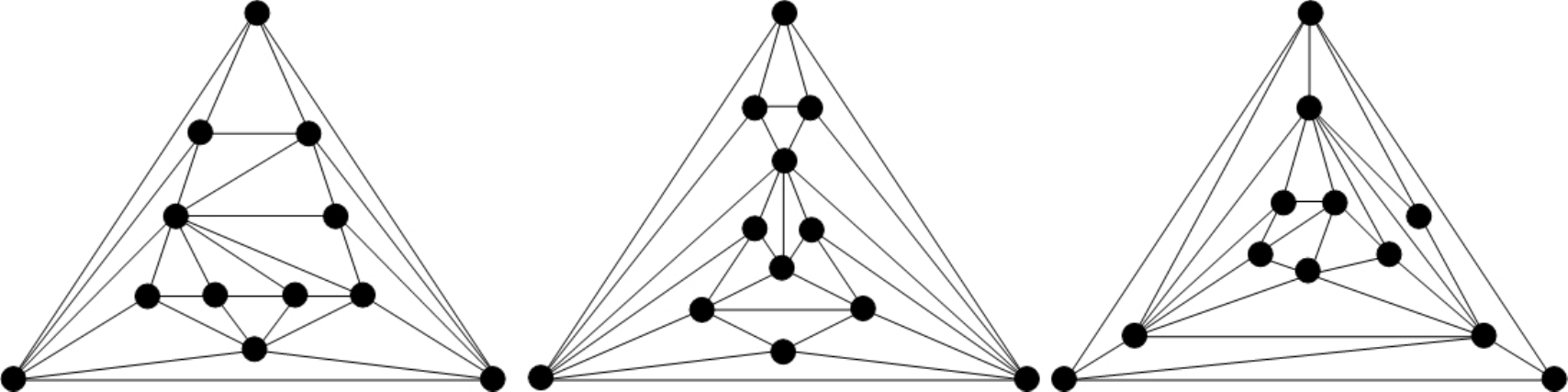}
    \includegraphics [width=340pt]{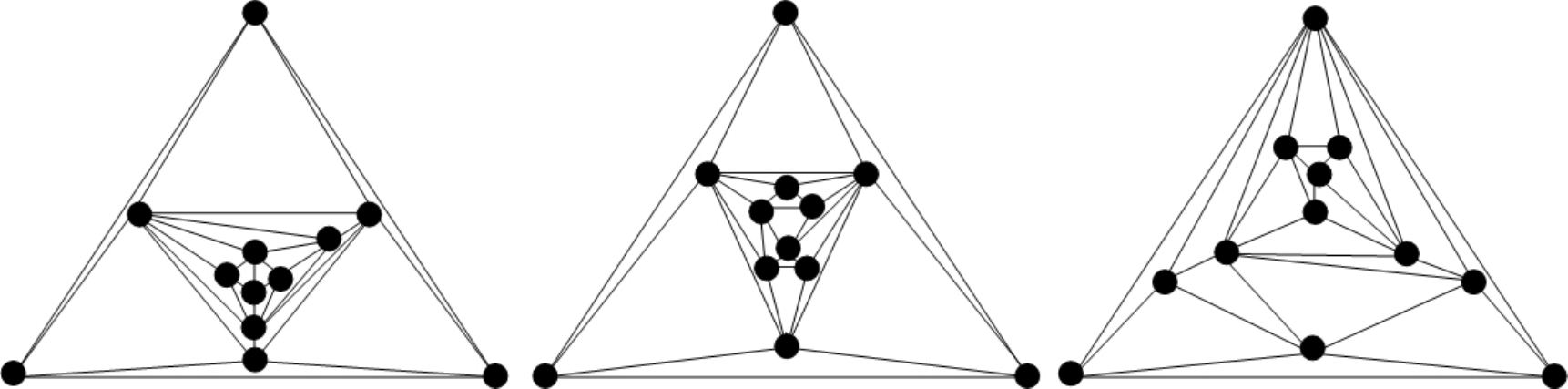}
    \includegraphics [width=340pt]{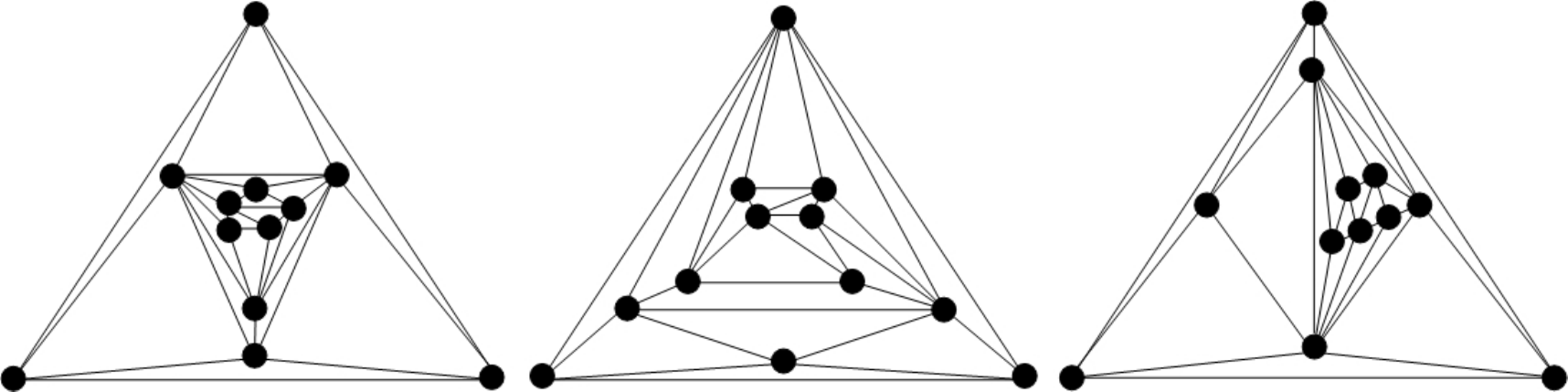}
    \includegraphics [width=340pt]{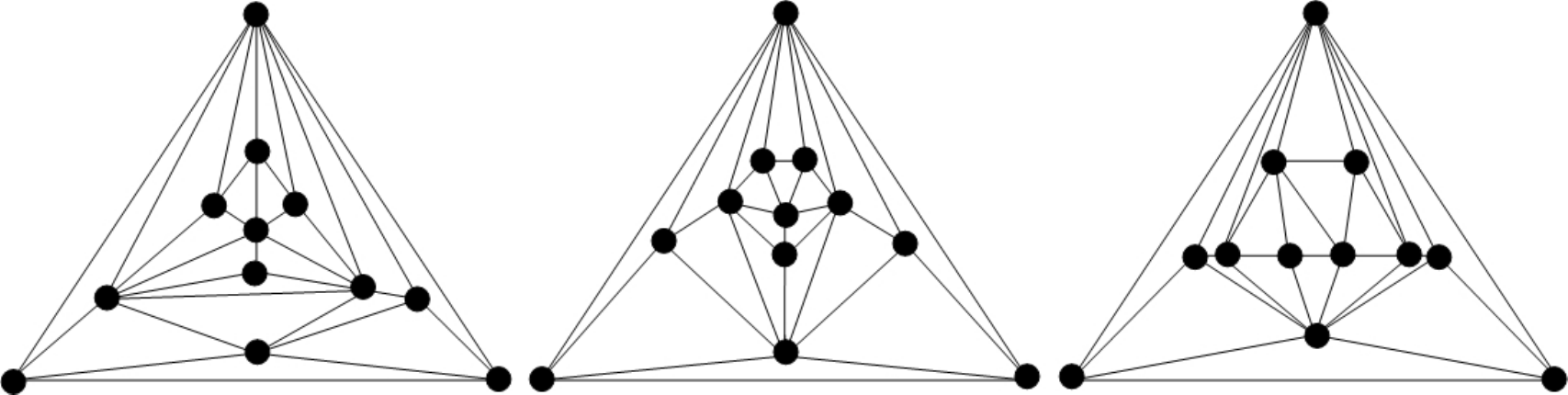}
    \includegraphics [width=340pt]{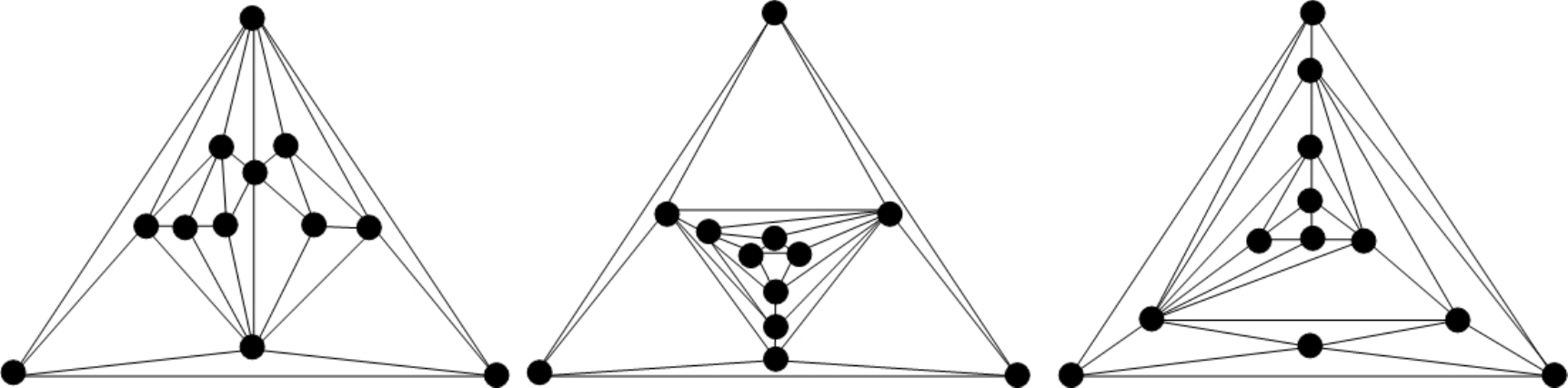}
    \includegraphics [width=340pt]{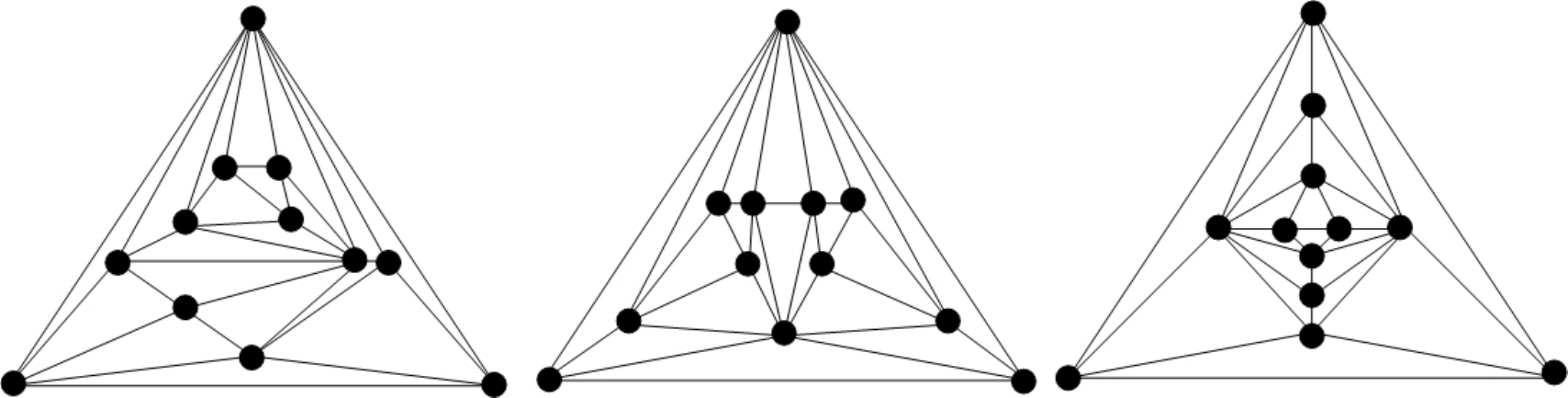}
    \includegraphics [width=340pt]{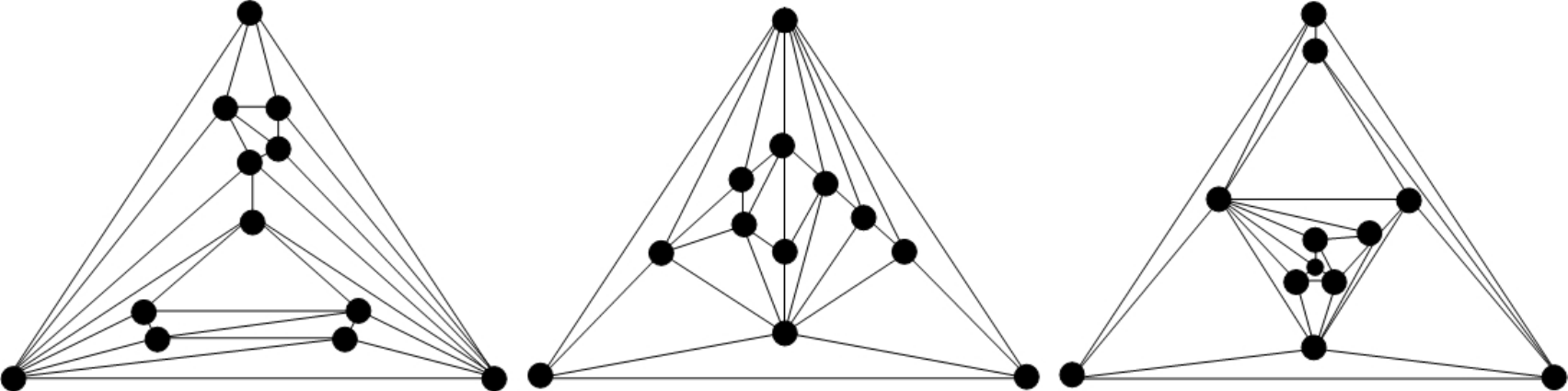}
    \includegraphics [width=340pt]{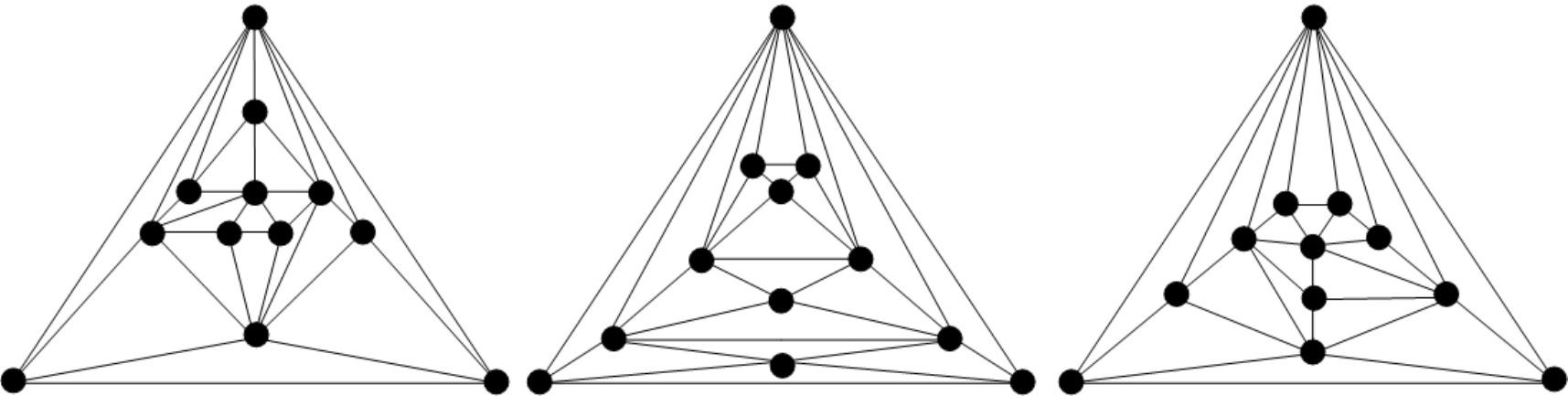}
    \includegraphics [width=340pt]{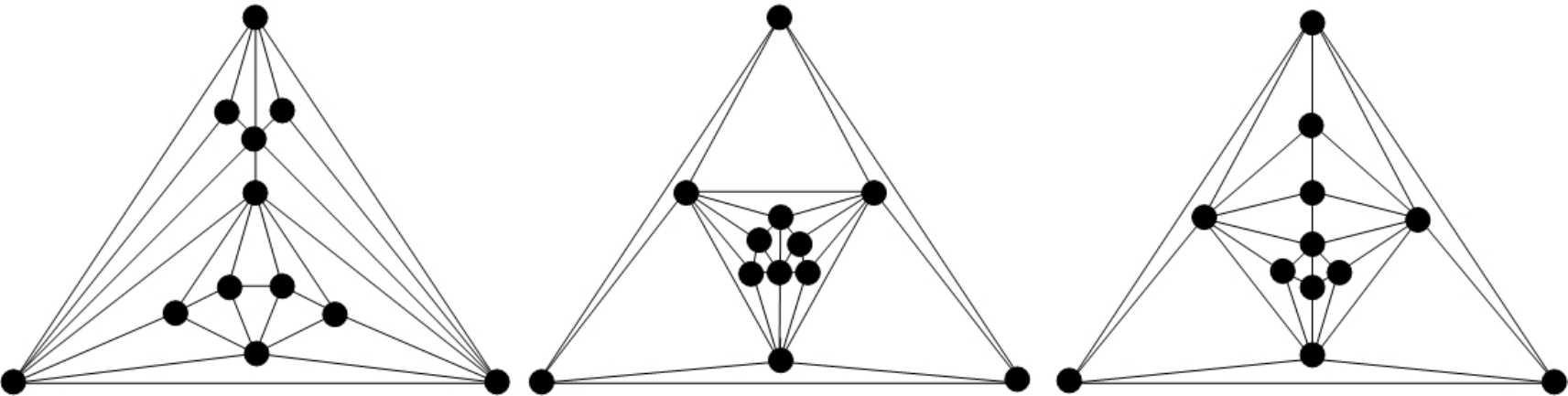}
    \includegraphics [width=340pt]{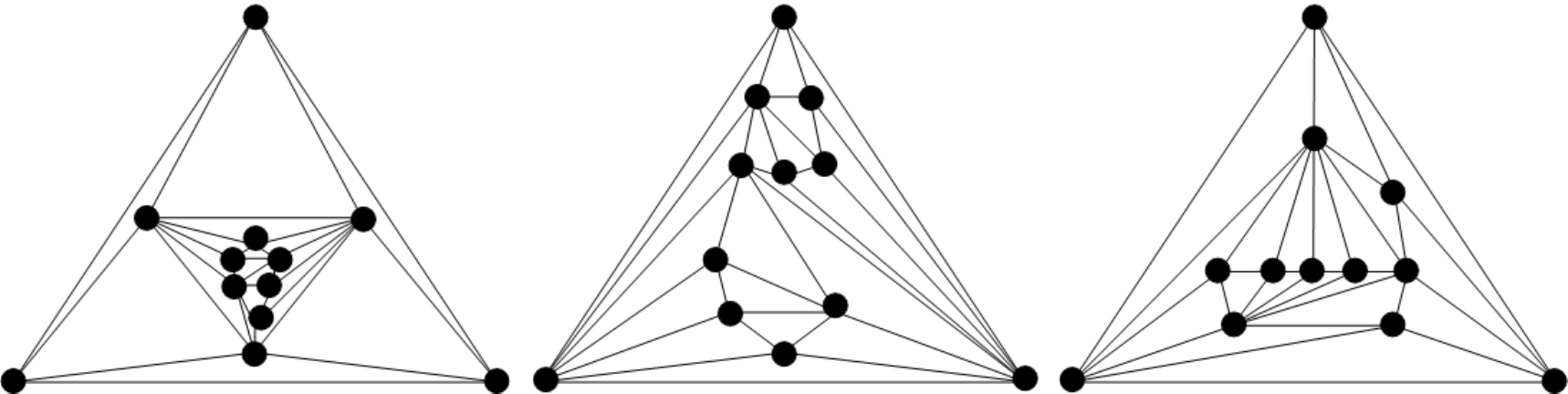}
    \includegraphics [width=340pt]{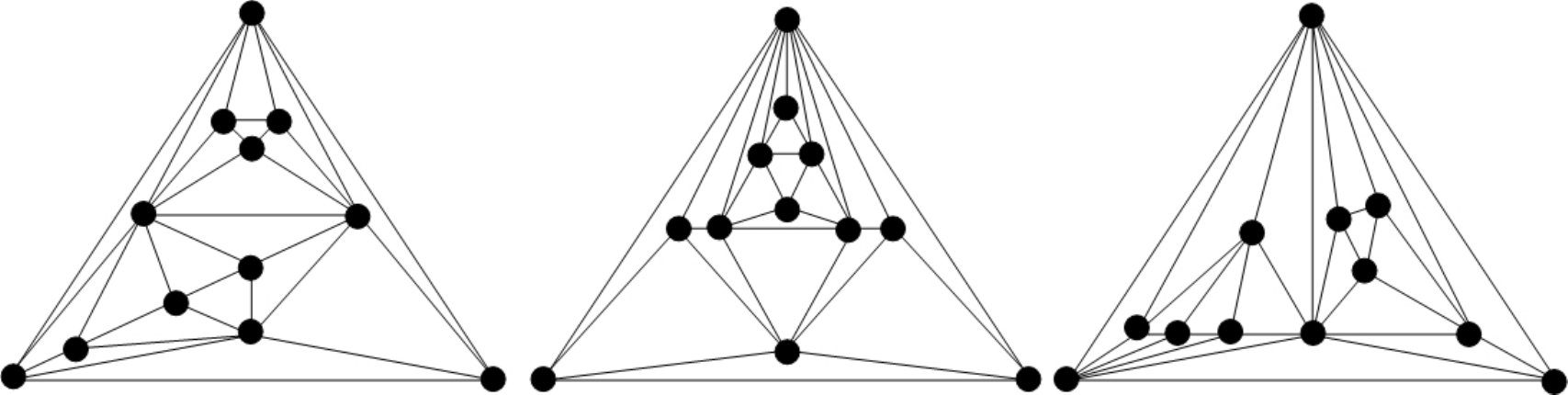}
    \includegraphics [width=340pt]{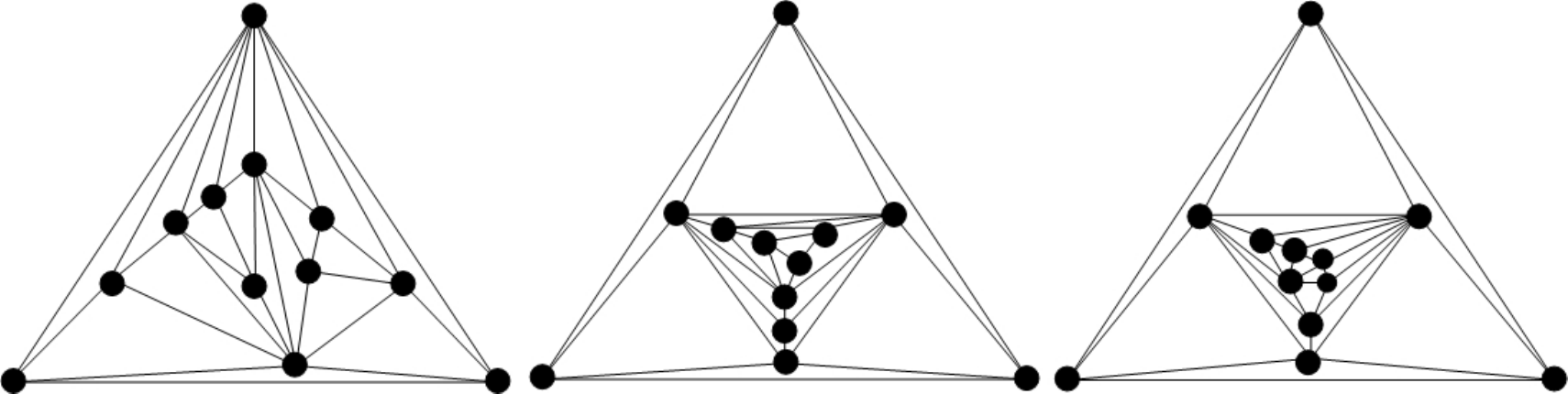}
    \includegraphics [width=340pt]{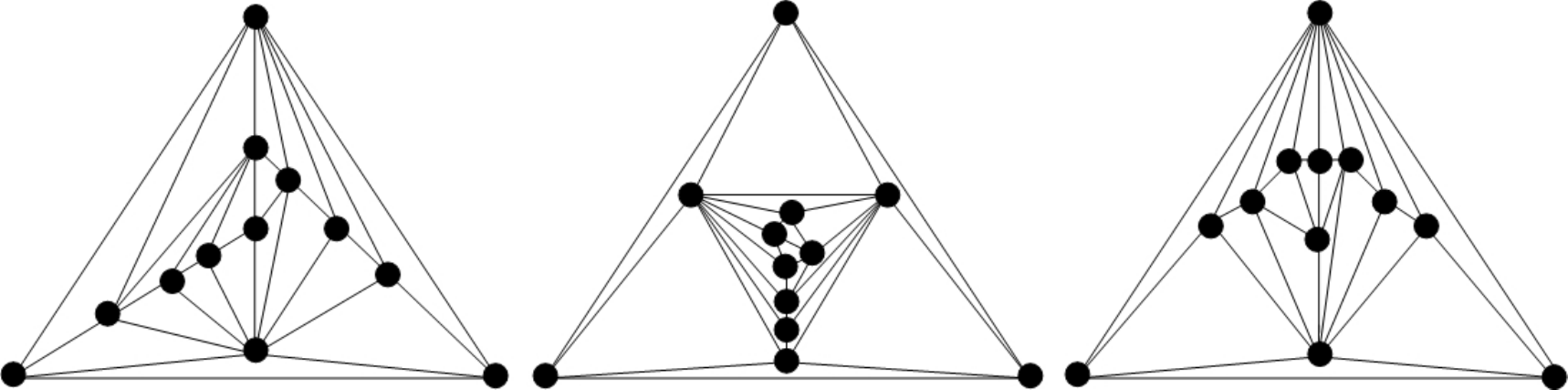}
    \includegraphics [width=340pt]{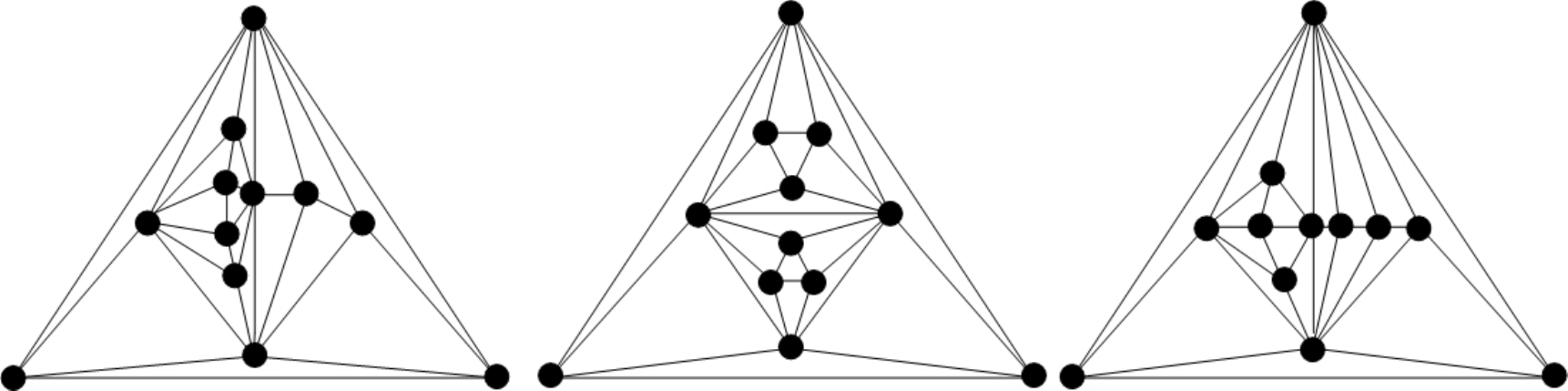}
    \includegraphics [width=340pt]{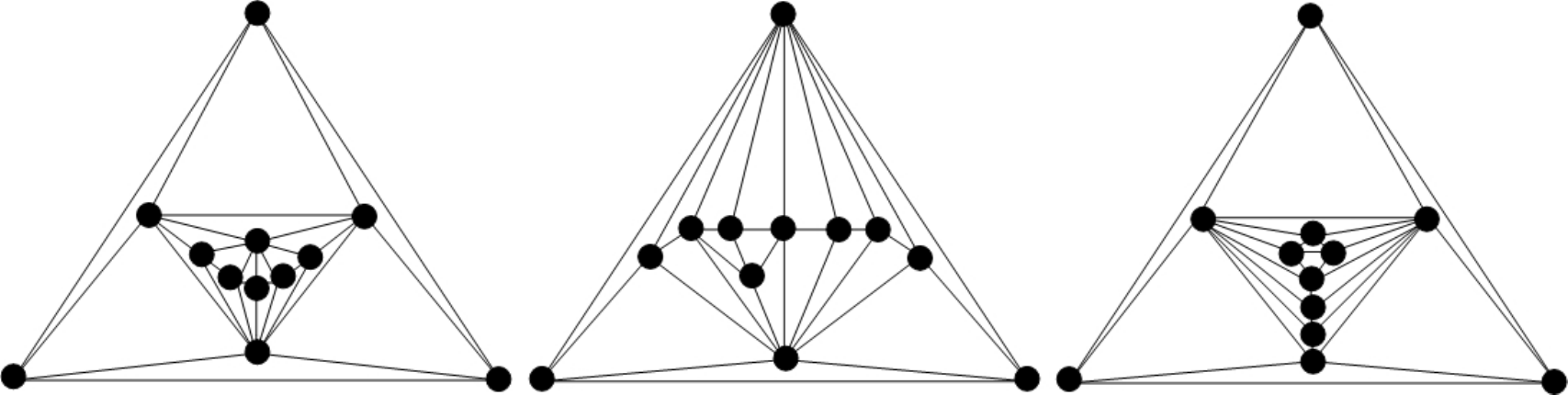}
    \includegraphics [width=110pt]{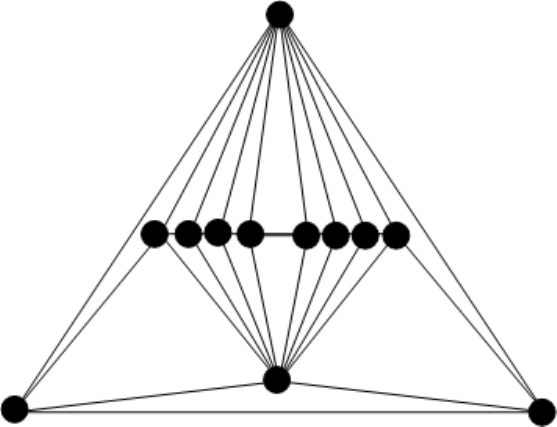}

\textbf{Figure 4.11.} All 130 maximal planar graphs with
$\delta(G)\geq 4$ whose order is 12
\end{center}

\subsection{Extending $k$-wheel operation and contracting $k$-wheel operation based on coloring}
Based on the last two subsections, in this subsection, we introduce
the \textbf{contracting $k$-wheel operation based on a coloring}
for a the maximal planar graph and the relevant inverse operations--
\textbf{the extending $k$-wheel operation based on a coloring},
where $k=2,3,4,5$, and give the related basic properties.

The contracting 2-wheel operation and the extending 2-wheel
operation based on a coloring are almost as the same as the
contracting 2-wheel operation and the extending 2-wheel operation
without coloring, only concerning how to assign a color to the
center of wheel. This process is very simple, so no more discussion
here. Please see section 4.2.

Let $G$ be a 4-colorable maximal planar graph. If $v$ is a 3-degree
vertex of $G$, and $\Gamma(v)=\{v_1,v_2,v_3\}$, then  $\forall f\in
C^{0}_{4}(G)$ ,\textbf{ the contracting 3-wheel operation based on
coloring $f$} about the 3-degree vertex $v$, means that delete
vertex $v$ from graph $G$. Naturally, the resulting graph is $G-v$,
and it is still a 4-colorable maximal planar graph. Meanwhile, the
extending 3-wheel operation on face $v_1-v_2-v_3$ based on coloring
$f\in C^{0}_{4}(G-v)$ is: add a new vertex  $v$ on that face, and
make $v$ adjacent to  $v_1,v_2,v_3$ respectively, and then assign to
$v$ a color different from $f(v_1),f(v_2),f(v_3)$.

Let $G$ be a 4-colorable maximal planar graph, $v$ be a 4-degree
vertex of $G$, and  $\Gamma(v)=\{v_1,v_2,v_3,v_4\}$. When $f\in
C^{0}_{4}(G)$, we have  $f(v_1)=f(v_3)$ or $f(v_2)=f(v_4)$. From now
on, we always assume that $f(v_1)=f(v_3)$ holds(see Figure 4.12(a)).
Then \textbf{ contracting 4-wheel operation based on coloring $f$}
about the 4-degree vertex $v$ means: deleting vertex $v$ from graph
$G$, and doing contracting operation on a pair of vertices
$\{v_1,v_3\}$, shown in Figure 4.12(b). \textbf{ Extending 4-wheel
operation based on coloring  $f$} means: for a 2-path $xuy$ of a
4-colorable maximal planar graph, first, doing the extending 4-wheel
operation without coloring(see Figure 4.2); second, assigning to the
new center $v$ of wheel a different color from $f(x),f(u),f(y)$.

\begin{center}
\includegraphics [width=280pt]{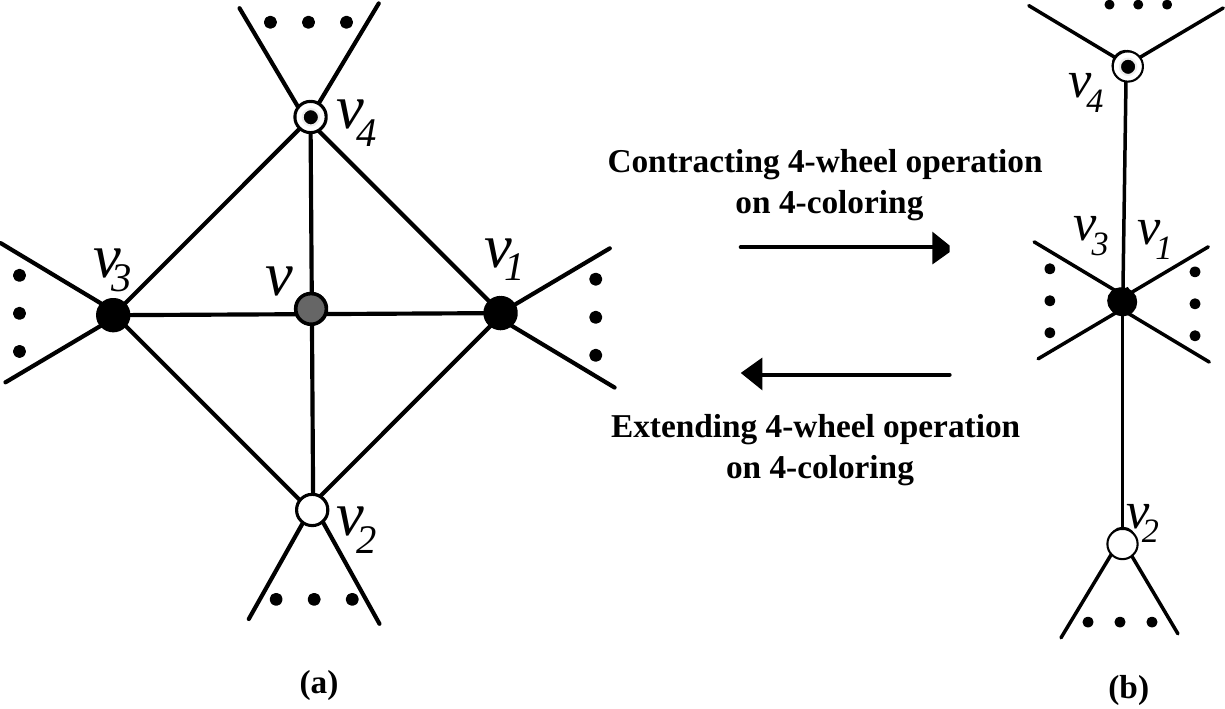}

\textbf{Figure 4.12.} Schematic diagram of contracting 4-wheel and
extending 4-wheel operations based on coloring
\end{center}

Let $v$ be a 5-degree vertex of 4-colorable maximal planar graph,
and $\Gamma(v)=\{v_1,v_2,v_3,v_4,v_5\}$. When $f\in C^{0}_{4}(G)$,
without lose of generality, let $f(v_1)=f(v_3)$, $f(v_2)=f(v_5)$,
shown in Figure 4.13(a). Then the \textbf{contracting 5-wheel
operation based on coloring $f$} means: deleting vertex $v$ from
graph $G$, and doing contracting operation on one pair of vertices
$\{v_2,v_5\}$ or $\{v_1,v_3\}$. Here we choose to $\{v_2,v_5\}$,
shown in Figure 4.13(b). For a 4-colorable maximal planar graph $G$,
the object of the \textbf{extending 5-wheel operation} is funnel
$L$, for which both the top and one of bottoms receive the same
color under a 4-coloring $f$ of $G$, shown in Figure 4.13(b),
$f(v_1)=f(v_3)$. The specific steps of the \textbf{extending 5-wheel
operation} is: first, doing the extending 5-wheel operation without
coloring (see Figure 4.4); second, assigning to the new center $v$
of wheel a different color from $f(v_1)=f(v_3),f(v_2),f(v_4)$(see
Figure 4.13).

\begin{center}
\includegraphics [width=320pt]{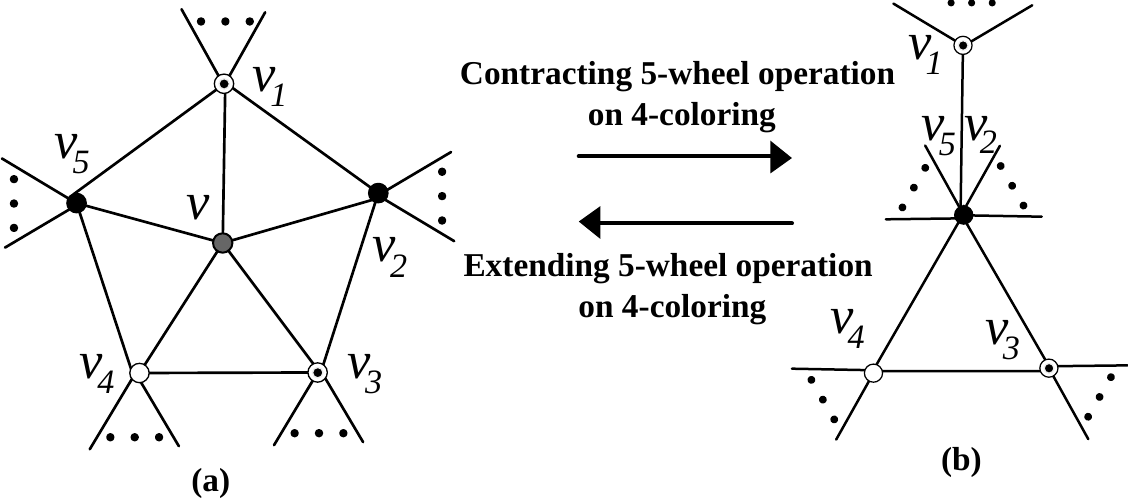}

\textbf{Figure 4.13.} Schematic diagram of contracting 5-wheel and
extending 5-wheel operations
\end{center}

So far, we have given the contracting $k$-wheel and the extending
$k$-wheel operations of a 4-colorable maximal planar graph
$G$($2\leq k\leq 5$) based on  a given coloring $f$ of $G$. It is
easy to see that, when $2\leq k \leq 5$, the extending $k$-wheel
operation and the contracting $k$-wheel operation are in one-one
correspondence. But when $k\geq 6$, the similar extending operation
and contracting operation are not in one to one correspondence, but
in one to many correspondence. In the following, we first give the
definition of  extending 6-wheel operation and contracting 6-wheel
operation of a 4-colorable maximal planar graph $G$ based on its
coloring $f$, then define extending $k$-wheel operation and
contracting $k$-wheel operation of $G$ under its coloring.

\begin{center}
  \includegraphics [width=340pt]{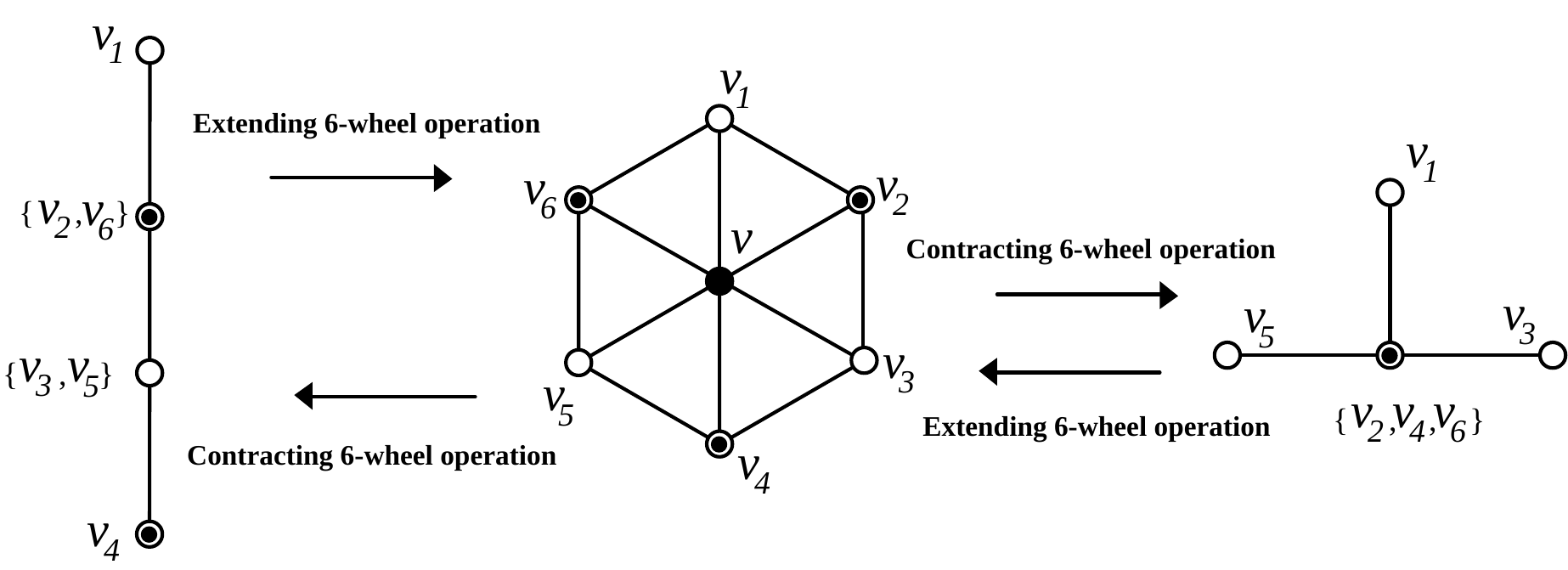}

  (a)

  \includegraphics [width=340pt]{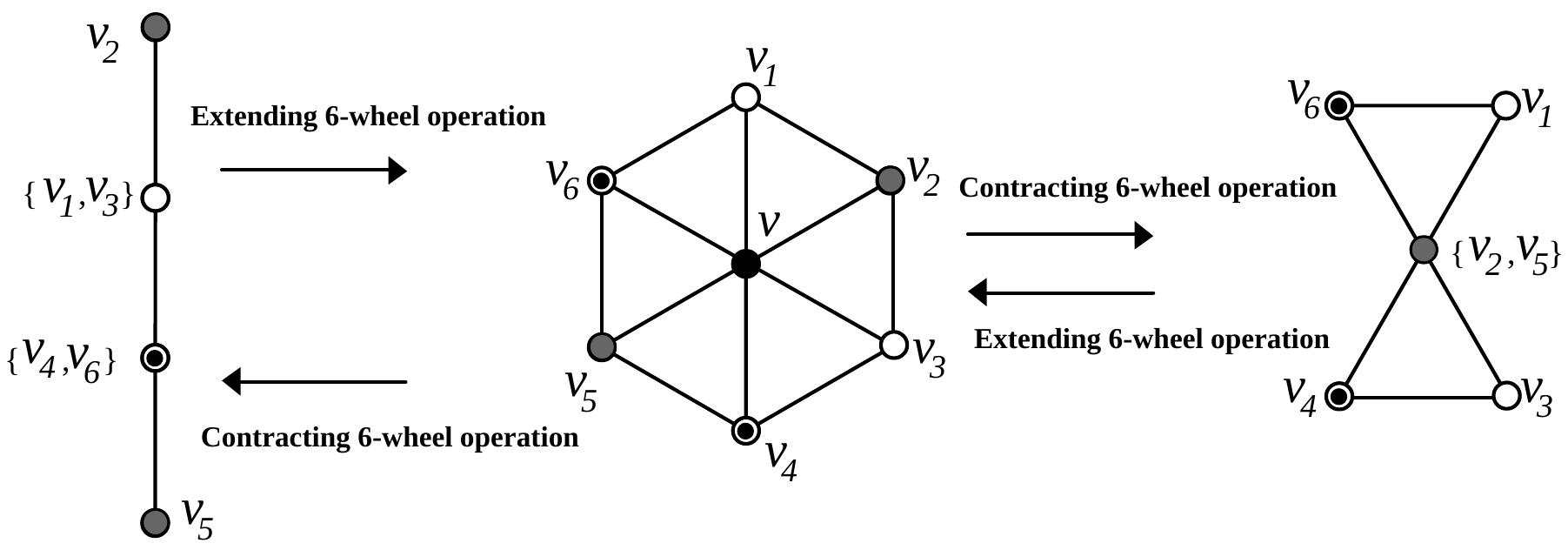}

  (b)

  \includegraphics [width=340pt]{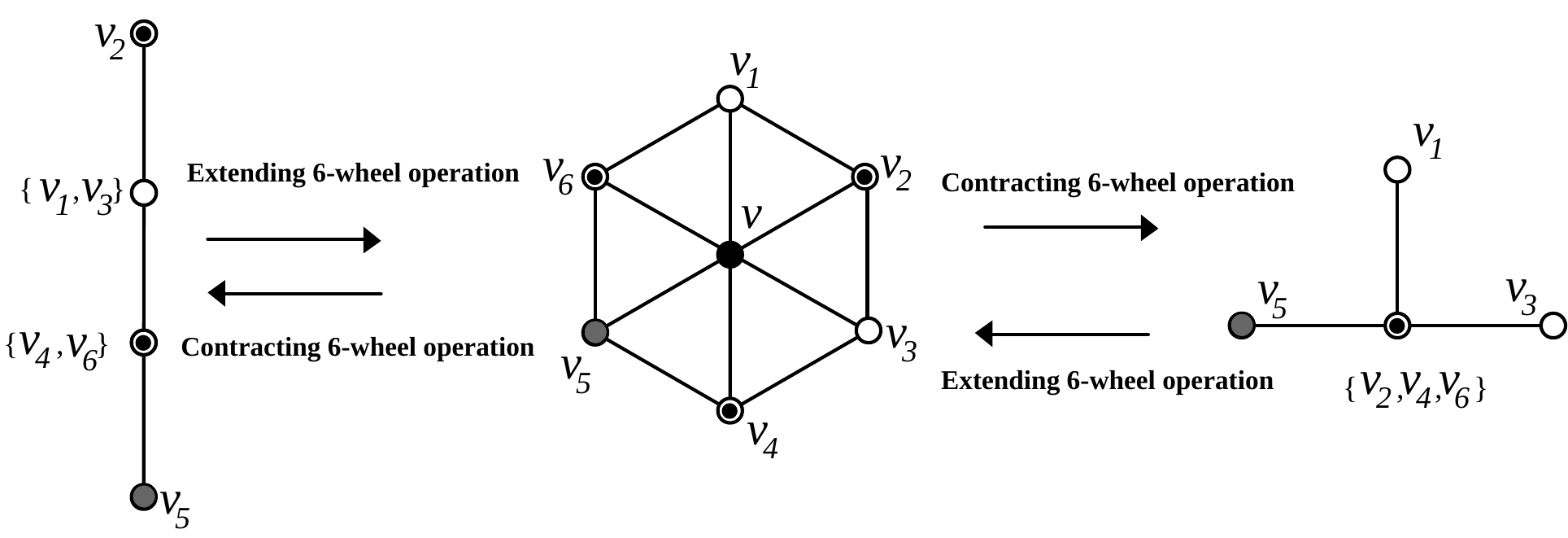}

  (c)

  \includegraphics [width=240pt]{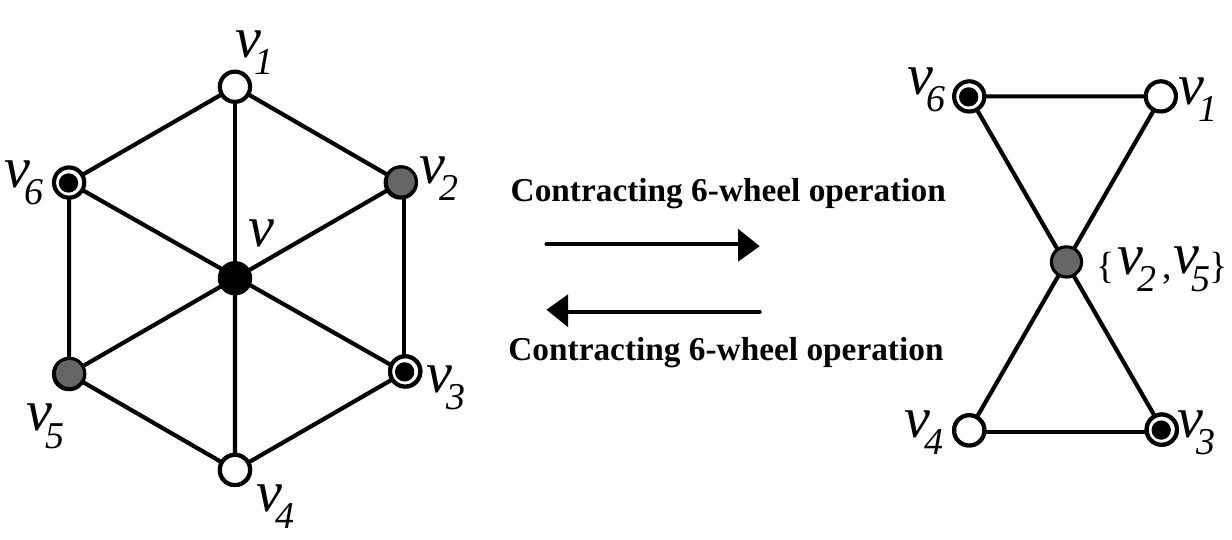}

  (d)

\textbf{Figure 4.14.} Schematic diagram of extending 6-wheel and
contracting 6-wheel operations
\end{center}

\begin{definition}
Let  $v$ be a 6-degree vertex of $G$, and
$\Gamma(v)=\{v_1,v_2,v_3,v_4,$ $v_5,v_6\}$. When $f\in
C^{0}_{4}(G)$, for $W_6=G[\Gamma(v)\cup \{v\}]$, there might be 4
kinds of colorings, shown in Figure 4.14.  For every kind of
colorings, the definition of the relevant extending 6-wheel
operation and contracting 6-wheel operation is clear from the shown
graphs.
\end{definition}

\begin{definition}
Let graph $G$ be 4-colorable, and  $f$ be one of its colorings. Let
$v$ be a vertex with degree $k\geq 3$ of $G$, and
$\Gamma(v)=\{v_1,v_2,\cdots,v_k\}$. The so-called \textbf{
contracting $k$-wheel operation based on $f$}  as to wheel
$W_k=G[\Gamma(v)\cup \{v\}]$ of graph $G$ is: deleting  vertex $v$,
then merging its neighbors received the same color into one vertex.
If the resulting graph is a maximal planar graph, then the
contracting $k$-wheel operation is completed; otherwise, the new
formed face with degree(the number of vertices incident to the face)
at least 4, for which there must have not less than two vertices
receiving the same color, then we merge these vertices again.
Repeat this process  until the resulting graph is a maximal planar
graph. According to the above definition of contracting $k$-wheel
operation, the process that does the extending operations step by
step reversely is the so-called\textbf{ extending $k$-wheel
operation based on coloring}.
\end{definition}

We sometimes use  $\zeta^{-}_{k}(G,v)$ to denote the graphs obtained
from contracting $k$-wheel operation on a degree-$k$ vertex $v$ of
the maximal planar graph $G$, in the case of no confusion, it can be
abbreviated as $\zeta^{-}_{v}(G)$. Use $\zeta^{+}_{k}(G,v)$ to
denote the graphs obtained from  extending $k$-wheel operation for a
maximal planar graph $G$, in the case of no  confusion, abbreviated
as $\zeta^{+}_{v}(G)$, where the added degree-$k$ vertex is denoted
by $v$.

Finally, a foundational result will be given in terms of the
operations of contracting $k$-wheel and extending $k$-wheel.

For a maximal planar graph $G$ of order $n(n\geq 5)$, then there
must exist at least one vertex with order three, four, or five. When
$G$ is a 4-colorable maximal planar graph, it must be obtained that
a maximal planar graph of order $n-1$, $n-2$ or $n-3$ through doing
the (combination) operations of contracting 3-wheel, 4-wheel,
5-wheel, 5-wheel and 3-wheel, 5-wheel and 2-wheel, or 4-wheel and
2-wheel from $G$. It means that any 4-colorable maximal planar graph
can be obtained by conducting the (combination) operations of
extending 3-wheel, 4-wheel, 5-wheel, 5-wheel and 3-wheel, 5-wheel
and 2-wheel or 4-wheel and 2-wheel.

\begin{theorem2}
For any given 4-colorable maximal planar graph $G$ of order $n(n\geq
8)$, it will be obtained by conducting the (combination) operations
of extending 3-wheel, 4-wheel, 5-wheel, 3-wheel and 5-wheel, 2-wheel
and 5-wheel or 2-wheel and 4-wheel from some graphs of order $n-1$,
$n-2$, or $n-3$.
\end{theorem2}

\section{Recursive maximal planar graphs}
The concept of recursive maximal planar graphs has been introduced in section 2: they can be obtained from $K_4$, embedding a 3-degree vertex in each triangular face continuously. The set of all the recursive maximal planar graphs is denoted by $\Lambda$, and the set of all those graphs with order $n$ is denoted by $\Lambda_n$, letting $\gamma_n=|\Lambda_n|$. Based on section 4, an exact definition of recursive maximal planar graphs is defined as follow: conducting extending 3-wheel operation continuously from $K_3$ or $K_4$, namely $\Lambda=(K_3,\zeta_3^+)$.

JT Conjecture states that a given 4-colorable maximal planar graph $G$ is uniquely 4-colorable if and only if it is a recursive maximal planar graph. So, the foundation to attack this conjecture is to further study recursive maximal planar graphs. This kind graphs are also
 called the FWF graphs. In the process of study, one class of graphs called the (2,2)-FWF
 graphs is actually the main graph class of recursive maximal planar graphs, which is indispensable in the proof of JT Conjecture.  Below, we give some related
 properties of the FWF graph, especially for the (2,2)-FWF graph.

 \subsection{Basic properties}
   \begin{theorem2}\label{th28}
    If $G$ is a FWF graph of order $n$, then it has at least two vertices of 3-degree.
 And when $n\geq 5$, any two vertices of degree 3 are not adjacent to each other.
   \end{theorem2}

   \begin{proof}
     Here we prove it by induction on the number of vertices, $n$. When $n=4,5,6$, $\gamma_{4}=\gamma_{5}=\gamma_{6}=1$;
 And the corresponding graphs are shown in the Figure 2.2.
 So the result is true obviously.

 Assume that the theorem holds when the number of vertices is $n$. That is, for any FWF graph $G$
 with $n$ vertices, it has at least two $3$-degree vertices, and all the vertices of $3$-degree
 are not adjacent to each other. Now we consider the case that the number of vertices is
 $n+1$.

     If a graph $G \in \Lambda_{n+1}$, then $G$ is constructed by adding a $3$-degree vertex $v$
 in any triangular face of a FWF graph with $n$ vertices, assuming $\Gamma_{G}(v)=\{v_{1},v_{2},v_{3}\}$. By induction, there are at least two $3$-degree vertices, and all the vertices of $3$-degree
 are not adjacent to each other in $G-v$. For $G-v$, if 3-degree vertices are included in $\{v_1,v_2,v_3\}$, then one exists at most, saying $v_1$. Obviously, there is at least another 3-degree vertex except $v_1$ in $G-v$, and all those $3$-degree vertices
 are not adjacent to each other. For $v$ is a 3-degree vertex of $G$ and it is not adjacent to any other vertices except $\{v_1,v_2,v_3\}$, so  there are also at least two 3-degree vertices in $G$, and all those $3$-degree vertices
 are not adjacent to each other. Thus, the conclusion holds. For $G-v$, if 3-degree vertices are not in $\{v_1,v_2,v_3\}$, the conclusion holds with the same discussion.
 \end{proof}

 \begin{theorem2}\label{th29}

 (1) There exists no maximal planar graph exactly having two adjacent vertices of degree 3.

 (2) There exists no maximal planar graph exactly having three vertices of degree 3,
  and each two of them are adjacent.
 \end{theorem2}
 \begin{proof}
   By contradiction.
   Assume that $G$ is a maximal planar graph with two adjacent vertices $u,v\in V(G)$ exactly, satisfying $d(u)=d(v)=3$.
 Since $u$ is also a 3-degree vertex, $\Gamma(u)=\{v,x,y\}$. Notice that $G$ is a maximal planar
 graph and  $u$ must be in the triangular face which consists of the vertices $v$, $x$ and $y$. In other words,
 $v$ is adjacent to vertices $x$ and $y$. These four vertices can form a subgraph $K_{4}$(shown in Figure 5.1).
 Since $G$ is a maximal planar graph, if there exist any other vertex, then it can
 form a triangle with $u$ or $v$. It contradicts  $d(u)=d(v)=3$. Otherwise,
 if there exists no other vertices, then $G$ is isomorphism to $K_{4}$ with
 four vertices of 3-degree. Therefore, there exists no maximal planar graph with two adjacent vertices
 of 3-degree exactly.
 \begin{center}
 \includegraphics[width=130pt]{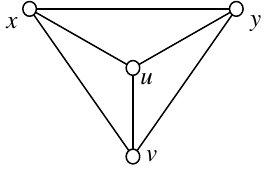}

 \textbf{Figure 5.1.} The schematic for the proof of Theorem 5.2
 \end{center}

    Similarly, we can conclude that there exists no maximal planar graph having exactly three vertices
 of 3-degree, and any two of them are adjacent. Let $u,v,x \in V(G)$, $d(u)=d(v)=d(x)=3$, and
 $uv, ux, vx \in E(G)$. There exist three vertices in $\Gamma(u)$, two of which are $v$ and $x$. Let $y$ denote the other vertices adjacent to $u$. So $\Gamma(u)=\{v,x,y\}$. Since a face can be constructed with three vertices $v$, $x$
 and $y$, any two of which are adjacent, so they can form a induced subgraph $K_{4}$, shown in Figure 5.1. Since $G$ is a maximal
 planar graph, if there exist any other vertices, it can form triangles with vertices $u$, $v$ or $x$. It contradicts
 the fact that $d(u)=d(v)=3$. Otherwise, $G\cong K_{4}$. Obviously, $K_{4}$ contains four vertices of 3-degree.
 Thus, there exists no maximal planar graph with three vertices of 3-degree, any two of which are adjacent.
 \end{proof}

 \begin{theorem2}\label{th30}
     If $G$ be a maximal planar graph having only one vertex of $3$-degree, then a subgraph without any
 $3$-degree vertex can be obtained by deleting $3$-degree vertices repeatedly.
 \end{theorem2}
 \begin{proof}
    Let $v$ be a unique vertex of 3-degree in graph $G$, and $\Gamma_{G}(v)=\{u_{1},u_{2},u_{3}\}$.
 Thus, these three vertices can form a triangle, any two of which are adjacent. Let $G_{1}$
 denote $G-v$, also a maximal planar graph. There may exist four cases as follows:

    (1)  $\delta(G_{1})\geq 4$;

    (2)  There exists only one 3-degree vertex;

    (3)  There exactly exist  two 3-degree vertices;

    (4)  There exactly exist three 3-degree vertices.

    For case (1), the theorem holds naturally. For case (3) and (4), we know they do not exist by Theorem 5.2.
 So we just need to consider about case (2). In this case, there exists a 3-degree vertex in subgraph $G_{1}$,
 denoted by $v_{1}$. Let $G_{2}=G_{1}-v_{1}$. Like the method mentioned above, if $\delta(G_{2})\geq 4$, then the
 theorem holds.
 Otherwise, the graph $G_{2}$ must contain a 3-degree vertex. In this way, we can get $\delta(G_{m})\geq 4$
 within finite $m$ steps; Otherwise, $G_{m} \cong K_{4}$ when $G_{m}$ contains only four vertices.
 It means the graph $G$ is a FWF graph. But there is only one 3-degree vertex in $G$.  It contradicts Theorem
 5.1.
 Therefore, this theorem holds.
 \end{proof}

 \subsection{(2,2)-FWF graphs}
    In this section, we introduce and study the (2,2)-FWF graph,  which is a special kind of the FWF graphs. A FWF graph is called the \textbf{(2,2)-FWF graph} if it contains only two vertices of 3-degree,
 and the distance between them is 2. It is easy to prove that there exist only one $(2,2)$-FWF graph with order 5 and 6 respectively, shown in Figure 5.3(a) and 5.3(b).

 To understand the structure of a (2,2)-FWF graph, the complete graph $K_{4}$
is divided into three regions, and then its vertices are labeled
correspondingly. Shown in Figure 5.2, the triangle
 is called \textbf{the outside triangle} when its vertices are labeled by 1, 2, 3, and the vertex $u$ (also labeled by 4) is called
 \textbf{the central vertex}. Here we define that the vertices 1--4 are colored by yellow, green, blue and red respectively. The four
 vertices and their corresponding colorings are called \textbf{the basic axes} in \textbf{the color-coordinate system} of a (2,2)-FWF graph.
 Four color axes are 1 (yellow), 2 (green), 3 (blue) and $u$ (red). Obviously, there exists no (2,2)-FWF graph of order 4; and
 there is only one (2,2)-FWF graph with 5 vertices under isomorphism, which can be obtained by embedding a 3-degree vertex in the region
 $I$, $II$ or $III$ of the graph $K_{4}$(shown in Figure 5.2).

 \begin{center}
    \includegraphics[width=240pt]{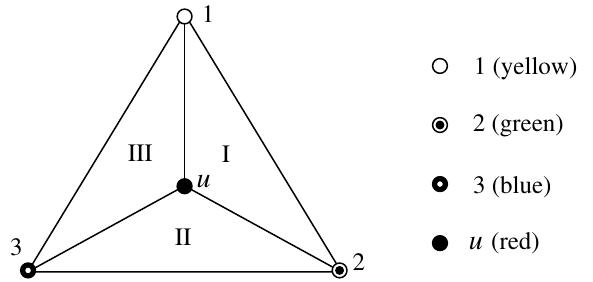}

    \textbf{Figure 5.2.}  The basic framework of the color-coordinate system
  \end{center}
    Without loss of the generality, we make the agreement that new vertices are only added in the region $II$. Thus, the vertex is colored by yellow (Figure 5.3 (a));
  the non-isomorphic (2,2)-FWF graphs of order $6$ can be obtained by embedding a 3-degree vertex in any region of the (2,2)-FWF graph of
 5-order. It is easy to prove that this kind of graphs with $6$ vertices obtained by embedding a new vertex in any face are isomorphic.
 Therefore, the (2,2)-FWF graphs of order $6$  are also unique. In general, we make an agreement that the 6th vertex is embedded in
 the face composed of the vertices 2, 4, 5 (i.e. the sub-region $I$ of the region $II$), which is colored by blue. (Figure 5.3 (b)).
 Further, for the (2,2)-FWF
 graphs with higher order, we restrict that new vertices are only added in the regions $I$ and $II$, but not in the region $III$.
 \begin{center}
   \includegraphics[width=320pt]{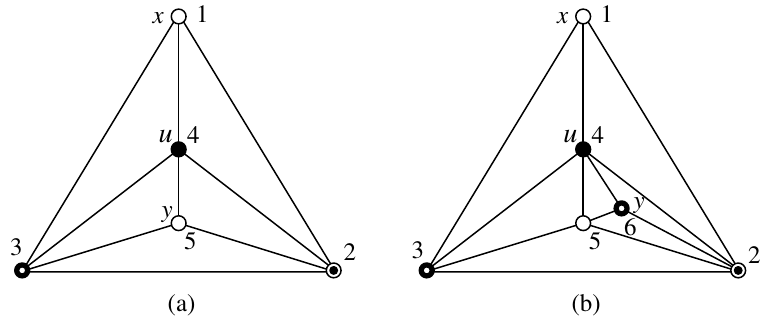}

    \textbf{Figure 5.3.}  Two (2,2)-FWF graphs \\ (a) a graph of order 5; (b) a graph of order 6
  \end{center}

    Based on the agreement above, we discuss about the classification of (2,2)-FWF graphs.
 Two methods are introduced as follows.

    The first is based on the region where the 3-degree vertices are embedded: (1) The (2,2)-FWF graphs are obtained by
    successively
 embedding the 3-degree vertices only in the region $II$. The graphs belong to this type are shown in Figure 5.4; (2) They are obtained by
 successively and randomly embedding the 3-degree vertices in the region $I$ and $II$, shown in Figure 5.5.  For planar graphs, there is a straightforward fact that
    \begin{Prop}$^{[54]}$
    Any face in the (maximal) planar graph can become the infinite outside face.
    \end{Prop}

    That is, the (2,2)-FWF graphs mentioned above are obtained by embedding 3-degree vertices randomly in the region $I$ and $II$.  We can transform
 any one 3-degree vertex in the region $I$ or $II$ to the outside triangular face by proposition 5.4, which is equivalent to the first classification.
 It means that this kind of (2,2)-FWF graphs are obtained by successively embedding 3-degree vertices only in the region $II$. Therefore, we only consider
 this kind of graphs in the following sections.

  \begin{center}
   \includegraphics[width=380pt]{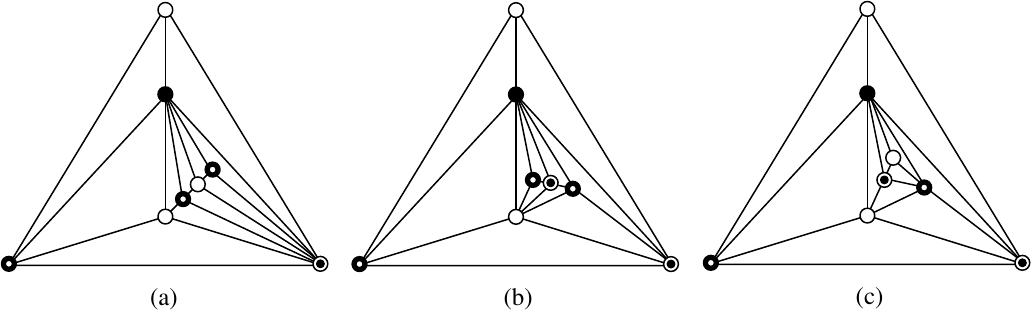}

    \textbf{Figure 5.4.}  The (2,2)-FWF graphs obtained by embedding 3-degree vertices only in the region $II$,
    \\(a) the adjacent type;(b) and (c) the non-adjacent type
   \end{center}
    \begin{center}
    \includegraphics[width=300pt]{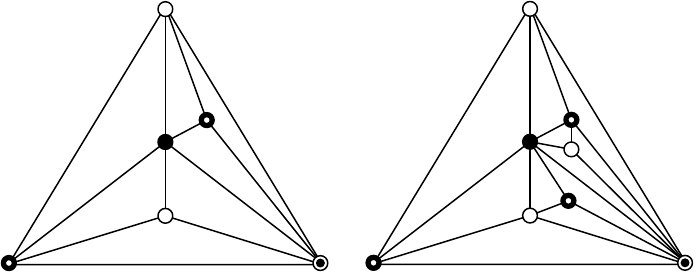}

    \textbf{Figure 5.5.} The (2,2)-FWF graphs obtained by embedding the 3-degree vertices in the region $I$ and $II$ randomly
    \end{center}

    The second classification is based on whether there exists a common edge between the two triangular surfaces of two $3$-degree vertices or not.
  It is called the adjacent type if there is a common edge; otherwise, the non-adjacent type. As shown in Figure 5.4, the first graph belongs to the adjacent type, while the last two graphs belong to the non-adjacent type.

  From the two classification methods above, all $(2,2)$-FWF graphs can be divided into the adjacent type of region $II$ and  the nonadjacent type of region $II$.

  From the Figure 5.3(a), the $(2,2)$-FWF graph of order 5 is a \textbf{double-center wheel}, and the degree of vertices in the neighbor of each center of the wheel is 4. When the order of $(2,2)$-FWF graph is not less than 6, we have the following result.

  \begin{theorem2}\label{th31}
    (1)Let $G$ be a $(2,2)$-FWF graph with order $n(n\geq 6)$, then for each 3-degree vertex $v$ in $G$, there only exists one vertex with order 4 in $\Gamma(v)$;  (2) Every (2,2)-FWF graph $G$ of nonadjacent type with order $n(n\geq 5)$ has one and only one $(n-1)$-degree vertex, and it is called \textbf{the central vertex} of the graph $G$,
   denoted by $u$.
    Further, in any partitions of color class in the (2,2)-FWF graph $G$, only the central vertex is colored by red; (3) For
   the (2,2)-FWF graphs of adjacent type obtained by embedding the 3-degree vertices only in the region $II$, only its central vertex is colored by red and also only
   its color axis $2$ is colored by green.
 \end{theorem2}
   \begin{proof}
    By induction. There is only one maximal planar graph of order 5(shown in Figure 5.3(a)), also a double-center wheel, so all triangular faces are equivalent.  Therefore, up to the isomorphism, there also only exist one FWF graph with order 6, also a $(2,2)$-FWF graph(shown in Figure 5.3(b)). Thus, the conclusion holds when $n=6$.

   Assume that the theorem holds when $n(n\geq 6)$, let us consider $(2,2)$-FWF graph $G$ with order $n+1$. Suppose $v$ is a vertex with 3-degree in $G$, two  cases as follows:

   Firstly, two or three vertices of 4-degree
are included in $\Gamma(v)$, then $G-v$ is also a FWF graph with order at least 5 which contains two or three vertices with 3-degree adjacent to each other, which contradicts theorem 5.1.

  Secondly,  the vertices of degree-4
is not included in $\Gamma(v)$, then $G-v$ is also a FWF graph with order at least 5 which contains only one vertex of 3-degree, which contradicts theorem 5.1.

In conclusion, we have proved that only one vertex of degree-4 is included in $\Gamma(v)$, the neighbor of the vertex $v$ of 3-degree.

Further, $G$ is a $(2,2)$-FWF graph of order $n$, then it exactly contains two vertices of 3-degree and the distance between them is two. Hence, there must be a vertex $u\in V(G)$ making any other vertices of $G$ adjacent to $u$,
namely $d(u)=n-1$, which can be proved by the gradual construction of $(2,2)$-FWF graphs.
   \end{proof}

On the basis of the theorem above, we now define some special triangular faces as follows: for a triangular face containing a vertex of 3-degree, if the degrees of three vertices in this triangular face are 3,4 and $n-1$ respectively, then this triangular face is called \textbf{$I$-type face}; if the degrees of three vertices in this triangular face are 3, $x$ and $n-1$ respectively, then this triangular face is called \textbf{$II$-type face}; if the degrees of three vertices in this triangular face are 3,4 and $x$ respectively, then this triangular face is called \textbf{$III$-type face}; where $5\leq x\leq n-1$.

\begin{theorem2}\label{thnew1}
$\gamma_5=\gamma_6=1,\gamma_7=2, \gamma_8=3, \gamma_9=6$
\end{theorem2}
The corresponding $(2,2)$-FWF graphs are shown in Figure 5.3, 5.6,5.7 and 5.8 respectively.

\begin{center}
    \includegraphics[width=320pt]{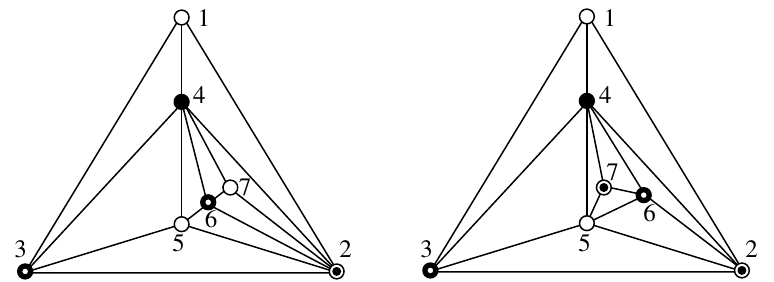}

   \textbf{Figure 5.6.} All of the two (2,2)-FWF graphs with
   order 7
\end{center}

\begin{center}
    \includegraphics[width=380pt]{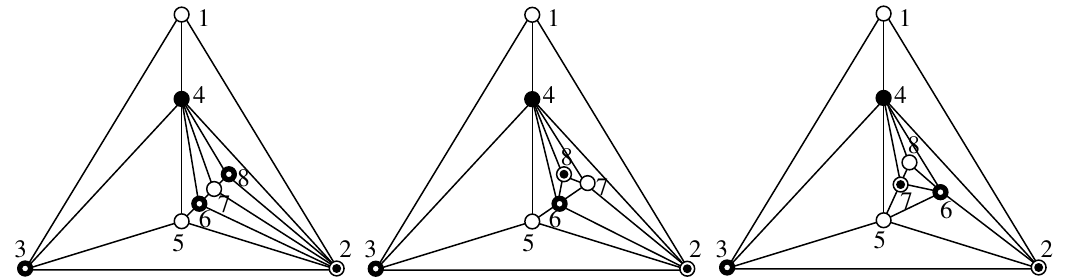}

   \textbf{Figure 5.7.} All of the three (2,2)-FWF graphs with
   order 8
\end{center}

\begin{center}
    \includegraphics[width=340pt]{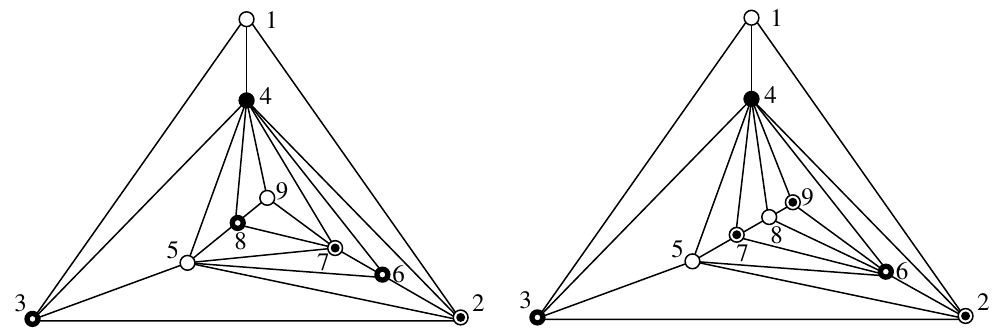}
     \includegraphics[width=340pt]{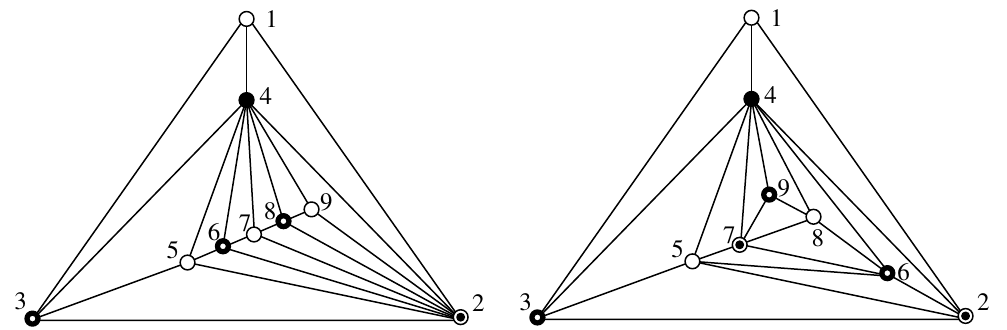}
     \includegraphics[width=340pt]{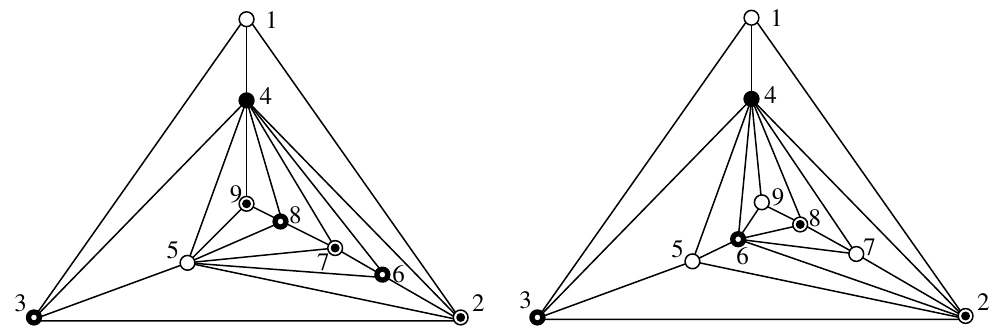}

   \textbf{Figure 5.8.} All of the six (2,2)-FWF graphs with
   order 9
\end{center}

\subsection{The color sequence of a (2,2)-FWF graph}
    Without loss of the generality, we can always assume that the (2,2)-FWF graph $G$ is obtained by embedding the $3$-degree vertices only in the region $II$
 in the following discussion. Thus, a (2,2)-FWF graph $G$ can be uniquely represented by its color sequence. The specific method is shown as follow:

    Let $V(G)=\{1,2,3,4=u,5,\ldots,n\}$, where vertex $1$($x$) indicates the first fixed  vertex of 3-degree, while the vertex $n(y)$ indicates the second
  vertex of 3-degree; the vertices $1(x)$, 2, 3, and $4(u)$ indicate the 1st, 2nd, 3rd, 4th color axis respectively; while the vertex $4(u)$ is the central
 vertex; the vertex $n-1$ signifies the  vertex of 3-degree of the subgraph $G_{n-1}=G-n$; the vertex $n-2$ denotes the  vertex of 3-degree of the subgraph
 $G_{n-1}-(n-1)$; the rest may be deduced by analogy. The sequence $c_{1}c_{2},\ldots,c_{n}$ is used to indicate the corresponding color sequence of
 the sequence $(1,2,3,4=u,5,\ldots,n)$, and the parameter $c_{i}$ is the color of the vertex $i$ in the (2,2)-FWF graph $G$. So we can obtain that

    $$c_{i}\in \{1=y(yellow),2=g(green),3=b(blue),4=r(red)\}$$

    According to the definition of the (2,2)-FWF graph, we can know that this representation also determine the structure of a graph uniquely.
 This structure starts from $K_{4}$ (shown in Figure 5.2), and selects a triangular face embedded the vertices according to the coloring of each vertex.

 \textbf{Example 5.1.} For the color sequence $c_{1}c_{2}c_{3}c_{4}c_{5}c_{6}c_{7}c_{8}c_{9}=ygbrybgyg$, its corresponding (2,2)-FWF graph can be analyzed
  easily, shown in Figure 5.9.

    \begin{center}
    \includegraphics[width=200pt]{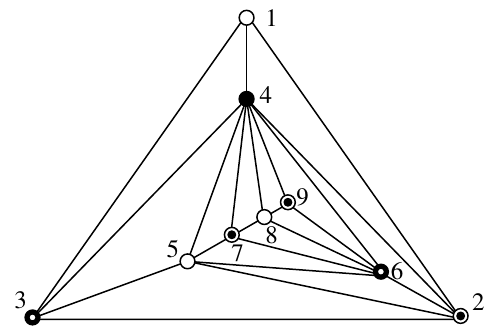}

    \textbf{Figure 5.9. } A color sequence $c_{1}c_{2}c_{3}c_{4}c_{5}c_{6}c_{7}c_{8}c_{9}=ygbrybgyg$ and its corresponding (2,2)-FWF graph
    \end{center}

    For the color sequence of a (2,2)-FWF graph, we can obtain the following theorem:
    \begin{theorem2}
        Let $c_{1}c_{2}\ldots c_{n}$ be the color sequence of a (2,2)-FWF graph. With the agreement in section 5.2,
    the colors of the first six vertices in this sequence is determined, namely $c_{1}=y, c_{2}=g, c_{3}=b, c_{4}=r, c_{5}=y, c_{6}=b$;
     if $G$ belongs to the adjacent type, then $c_{7}=y$; otherwise, $c_{7}=g$.
    \end{theorem2}

 \subsection{ Chromaticity of induced graphs by extending 4-wheel operation}

    In this section, we just discuss the vertex coloring problem of the induced graph from
 a (2,2)-FWF graph by extending $4$-wheel operation. We know that a given (2,2)-FWF graph is uniquely $4$-colorable, and according to the definition of the
 color-coordinate system in the section $5.3$, every vertex can be also colored determinately.
    \begin{definition}
 Let $G$ be a (2,2)-FWF graph, $f$ be the unique $4$-coloring of it,  and $xuy$ be a path of length 2 in $G$.  Obviously, there exist a coloring $f^*$ of graph $G*xuy$ that is induced from $G$ by extending 4-wheel operation on the path $xuy$, and
$$
f^*(x)=\left\{
\begin{array}{cc}
   f(u)& if~ x=u';  \\
   4&  if~ x=v;\\
   f(x)& otherwise.
\end{array}\right. \eqno{(5.1)}
$$
Namely, vertices $u$ and $u'$ are assigned the same color under $f^*$, and the new added vertex $v$ is assigned the different color 4 from vertices $x,y,u$, while the color of the rest vertices remain unchanged. We refer to $f^{*}$ as the  \textbf{natural 4-coloring} of graph  $G*xuy$.
\end{definition}

    Naturally, one question is proposed about whether the induced graph $G*xuy$ obtained by extending $4$-wheel operation graph
 is uniquely 4-colorable or not. This question is discussed as a key problem in this section. Definitely, the answer is negative, that is,
 $|C_{4}^0(G*xuy)|>1$.

 Here the definition of the color neighbor is introduced as follow:
 \begin{definition}
    Let $G$ be a $k$-chromatic graph, and $f \in C_{k}(G), u\in V(G)$. The \textbf{color neighbor} of vertex $u$ on
 coloring $f$ is the set which consists of all the colors assigned to $\Gamma(u)$ under $f$, denoted as $C(f, \Gamma(u))$.
 \end{definition}

 \begin{theorem2}
    Let $G$ be a (2,2)-FWF graph with order $n$ and $f$ be the unique 4-coloring of it. The vertices $x$, $y$ are two vertices of 3-degree and the vertex $u$ is
 the central vertex. Then, the induced graph $G*xuy$ is not uniquely $4$-colorable.
 \end{theorem2}
 \begin{proof}
    Obviously, if $f(x)=f(y)$,
 the vertex $v$ in the graph $G*xuy$ has two possible colors to choose when both vertices $u$ and $u'$ are colored by red. Hence, the graph $G*xuy$ is not uniquely 4-colorable.
 So we only need to consider the case that $f(x)\neq f(y)$.

    According to the classification in the section $5.2$, all the (2,2)-FWF graphs can be classified into two types:
     the adjacent type of region $II$ and  non-adjacent type of region $II$. The discussion as follows, respectively:

    \textbf{Case 1 :} The $(2,2)$-FWF graph $G$ of the adjacent type of region $II$.

     Based on the Theorem 5.5, we know that vertices $x$,$2$,$3$ are coordinate axes $1,2,3$, colored by yellow, green and blue respectively;
 and the central vertex $u$ is colored by red. Since all the 3-degree vertices can only be embedded in the subregion $I$ of
 the region $II$, the vertex $y$ can be colored by yellow or blue, illustrated by Figure 5.10(a). But when the vertex $y$ is colored
 by yellow, which is the same with vertex $x$, this case is not needed  considering. So we only discuss the case that the
 vertex $y$ is colored by blue. With the definition of extending 4-wheel
 operation,
an extending 4-wheel operation on the path $x-u-y$ can be done and the graph $G*xuy$ is obtained.

  \begin{center}
   \vspace{5mm}
     \includegraphics[width=360pt]{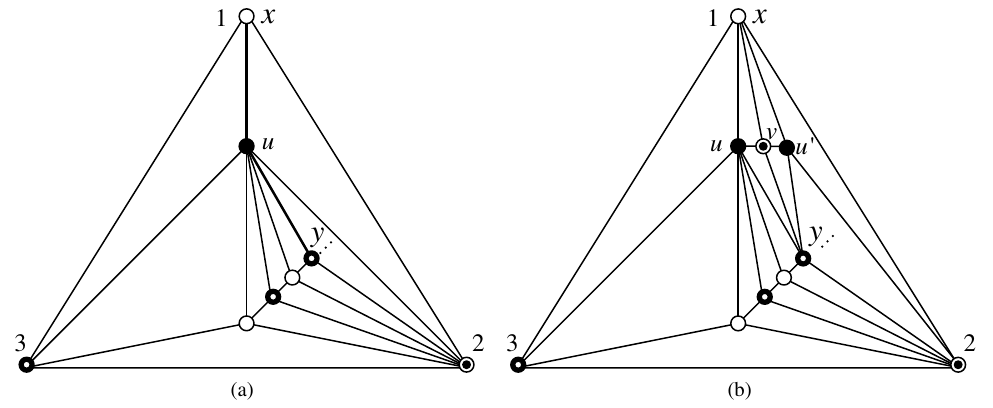}

     \includegraphics[width=180pt]{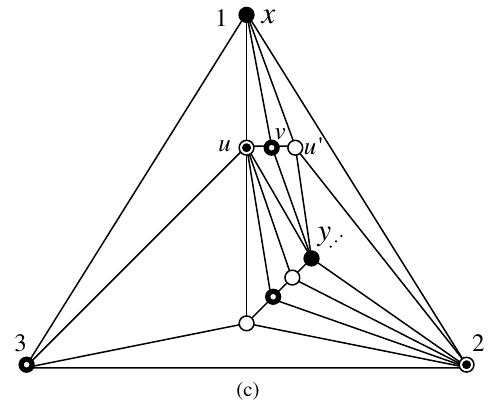}

  \textbf{Figure 5.10.} A graph of the adjacent type of region $II$ and two colorings of its induced graph by extending 4-wheel operation
  \end{center}

    It is easy to obtain two colorings of the graph $G*xuy$: one is the natural 4-coloring $f$ in which
 the vertex $u'$ is colored by red and the vertex $v$ embedded newly is colored by green.
 And the colorings of other vertices remain
 unchanged. Obviously, it is a coloring of the graph $G*xuy$, illustrated in figure 5.10(b). Besides,
from the discussion above, there is only one vertex $v$ colored by green in $\Gamma(u)$ and only one vertex $x$ colored by yellow in $\Gamma(u')$, under $f$.
Further, for graph $G*xuy$, it is  only two vertices $u$ and $u'$ that are colored by red under its natural coloring $f$, so we can obtain a new 4-coloring $f'$ of the graph $G*xuy$:
let the vertices
 $u$,$u'$,$x$,
 $y$ and $v$ be recolored  by green, yellow, red, red and blue respectively, other vertices remain unchanged on the basis of $f$.

 Since there is only one vertex $v$ colored by green in $\Gamma(u)$ under $f$,
 then change the color assigned to vertex $u$ from red to green, and only vertex $v$ receives the same green color.
 Similarly, since there is only
 one vertex $x$ colored by yellow $\Gamma(u')$ under $f$,
then change the color assigned to vertex $u'$ from red to yellow,
 and only vertex $x$ receives the same yellow color.
 After that, recolor nonadjacent vertices $x$ and $y$ by red, and change the color assigned to vertex $v$ from green to blue. Then we will obtain the new 4-coloring $f'$ of graph $G*xuy$ when remain the colors of any other vertices
 unchanged, shown in Figure 5.10(c).
 These two colorings $f$ and $f'$ are different apparently. Hence, the case 1 is proved.

    \textbf{Case 2:} The $(2,2)$-FWF graph $G$ of the nonadjacent type of region $II$.

    According to Theorem 5.7, the first six vertices of all the (2,2)-FWF graphs are colored in the same way, illustrated as follows:
   $$\left(
    \begin{array}{ccccccc}
        1 & 2 & 3 & 4 & 5 & 6 & \ldots \\
        y & g & b & r & y & b & \ldots \\
    \end{array}
   \right)\eqno{(5.2)}$$
 Namely, the color coordinate axes vertex 1(or vertex $x$), 2, 3 and 4(or vertex $u$) receive yellow, green, blue, and red colors respectively.
 Vertex 1 is a vertex of 3-degree and adjacent to the central vertex 4 colored by red, the  vertex 2 colored by green and the  vertex 3 colored by blue.
 Therefore, vertex 1 is a vertex with degree-5 in the graph $G*xuy$. And in the natural coloring of $G*xuy$, it is adjacent to the vertex 2 colored by green, the vertex 3 colored by blue, the vertex 4(or vertex $u$)colored by red, the vertex $u'$ colored by red and the vertex $v$ colored by blue respectively.

    Since the graph $G$ belongs to the nonadjacent type, so the vertices $7,8,\ldots,n$ must be added in the triangular face formed by vertices 4, 5 and 6, shown in the Figure $5.6(c)$ and Figure $5.11(a)$. According to Theorem 5.7, the 7th vertex  can only be colored by green, illustrated as
 follow:
   $$\left(
    \begin{array}{cccccccc}
        1 & 2 & 3 & 4 & 5 & 6 & 7 & \ldots \\
        y & g & b & r & y & b & g & \ldots \\
    \end{array}
   \right)\eqno{(5.3)}$$

    Then, it can be known easily that the vertex 2 colored by green, which is also a color coordinate axis and a vertex of 5-degree. The neighbors of vertex 2 are
 vertex 1(yellow), vertex 3(blue), vertex 5(yellow), vertex $u'$ (red) and vertex 6(blue). They are shown as follow.
    $$C(f,\Gamma(1))=\{g\{2\},b\{3,v\},r\{u,u'\}\} \eqno{(5.4)}$$
    $$C(f,\Gamma(2))=\{y\{1,5\},b\{3,6\},r\{u'\}\} \eqno{(5.5)}$$

    Now we take the representative graph in Figure 5.11(a) for example, the detailed steps of a new $4$-coloring induced by the natural 4-coloring  of graph $G*xuy$(Figure 5.11(b)) is given as follows:

    First, change the color assigned to vertex 1 from yellow to green. By formula (5.4), the two
 endpoints  of edge $\{1,2\}$ are both colored by green, and it is the unique pseudo color edge. Other vertices are colored properly, illustrated
 by Figure 5.11(c).

    Second, change the color assigned to vertex 2 from green to red. Thus,  the coloring of the two endpoints of edge $\{1,2\}$ becomes proper, while
  $\{u',2\}$ becomes a pseudo color edge, for its two endpoints are both colored by red.  The coloring of other vertices are all proper,
 illustrated by Figure 5.11(d).

    Third, the vertex $u'$ is colored yellow instead of red. Thus, the pseudo color edge $\{u',2\}$ becomes proper. There may be
 several vertices colored by yellow in the neighbor of vertex $u'$, which can form a set $C_{4}(u',yellow)$. Therefore,
 this step generates several pseudo color edges whose number is $|C_{4}(u',yellow)|$. Obviously, other edges are all proper,
 illustrated by Figure 5.11(e).

    Fourth, all the vertices in $C_{4}(u',yellow)$ are changed to red. Since in the $4$-coloring of the third step, only two
 vertices $u$ and 2 are colored by red. Obviously, $C_{4}(u', yellow)\subset \Gamma(u')$. So all the vertices in $C_{4}(u', yellow)$
 are not adjacent to the vertex $u$. In the vertex set adjacent to vertex 2 after the third step, the vertices 1, 3 and 6 are colored
 by green, blue and blue. Although the vertex 5 is colored by yellow, it is a vertex of 5-degree and not in $C_{4}(u', yellow)$. Therefore,
 the edges between vertex 2 and all the red vertices in $C_{4}(u', yellow)$ are proper. Moreover, the vertices in
 $C_{4}(u', yellow)$ form an independent set of the graph. So they can not generate pseudo color edges by themselves. This step is
 illustrated by Figure 5.11(f).

    Thus, based on the natural 4-coloring of the graph $G*xuy$, we can obtain a new coloring different from the natural $4$-coloring.
 It means that the induced graph constructed from nonadjacent (2,2)-FWF graph by extending 4-wheel operation is not uniquely $4$-colorable.

    To sum up the case 1 and 2, this theorem is proved.
 \end{proof}
 \begin{center}

  \includegraphics[width=240pt]{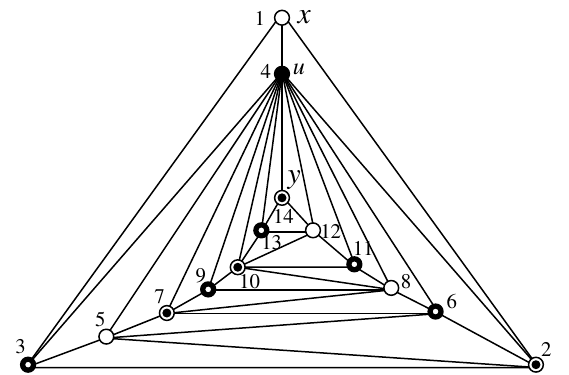}

  \textbf{(a)}  A representative (2,2)-FWF $G$

  \includegraphics[width=240pt]{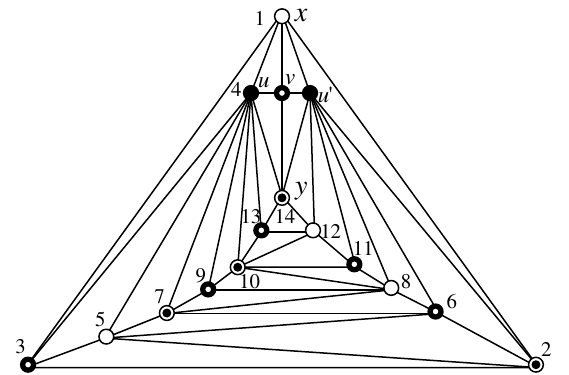}

  \textbf{(b)} The natural 4-coloring of the induced graph $G*xuy$ by extending 4-wheel operation of $G$

  \includegraphics[width=240pt]{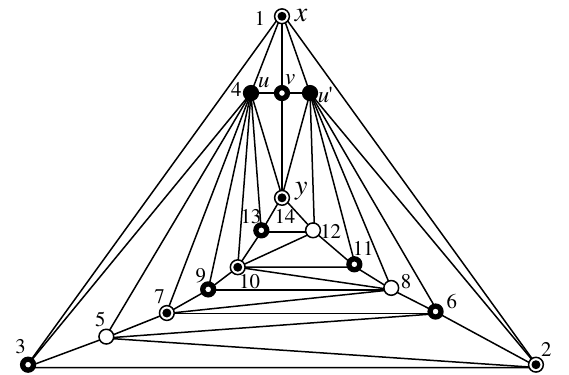}

  \textbf{(c)} The coloring of the induced graph $G*xuy$ when vertex 1 is colored by green instead of yellow,
  which generates a pseudo color edge \{1,2\}

  \includegraphics[width=240pt]{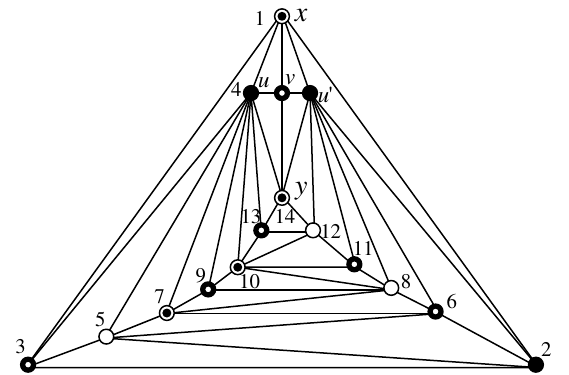}

  \textbf{(d) } The coloring of the induced graph $G*xuy$ based on the step(c) when vertex 2 is colored by red instead of green,
  which generates a pseudo color edge \{$u'$,2\}

  \includegraphics[width=240pt]{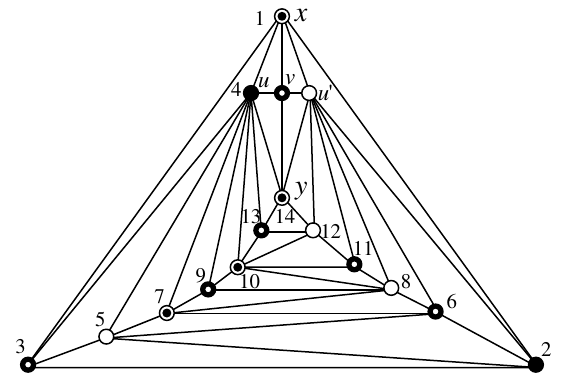}

  \textbf{(e)} The coloring of the induced graph $G*xuy$ based on the step (d) when vertex $u'$ is colored by yellow instead of red,
  which has several pseudo color edges $\{\{u',u^{\prime \prime}\},u^{\prime \prime}\in C_{4}(u',yellow)\}$

  \includegraphics[width=240pt]{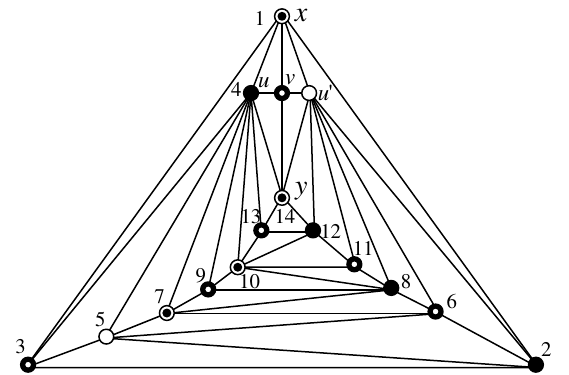}

  \textbf{(f) } A new coloring of the induced graph $G*xuy$ based on the step(e) when the vertices in the set $C_{4}$($u'$,yellow)
 are colored by red

 \textbf{Figure 5.11.} The illustration: the induced graph $G*xuy$ by extending 4-wheel operation of the $(2,2)$-FWF graph $G$ belonging to the nonadjacent type of region $II$, is not unique 4-colorable
  \end{center}

\section{Coloring-structure of maximal planar graph}
    Let $G$ be a 4-colorable maximal planar graph, $C(4)=\{1,2,3,4\}$ the color set of $G$ and $C_4^0(G)$ the set consisting of all the non-isomorphic 4-colorings of $G$. Then for $\forall f\in C_4^0(G)$, denote by $G[i,j]=G[V_i\cup V_j]$ the \textbf{bicolored induced subgraph} of $G$(or for short \textbf{bicolored subgraph}), where $V_i$ is the set that contains all the vertices of $G$ colored by  color $i$ based on $f$, $i,j=1,2,3,4$ and $i\neq j$. Obviously, the number of bicolored subgraphs of $G$ is altogether six. If each of these six bicolored subgraphs is a tree, then $f$ is called a \textbf{tree-coloring} of $G$; otherwise, there is at least one bicolored subgraph that contains a cycle, and then $f$ is called a \textbf{cycle-coloring} of $G$. In fact, the essential error, which appeared in the proof of four color conjecture by $Kempe$ in 1879, was that he did not make clear the basic structure of the 4-coloring of planar graphs. Although $Heawood$ found this error in 1890, he didn't gave a correction for this problem. However, he proved the five color theorem by means of $Kempe$'s idea. The results of this section seem say: the real reason that later many scholars could not correct $Kempe$'s error might be that the coloring-structure of 4-colorable maximal planar graphs didn't be studied in depth.

    The aim of this section is to make clear the structure of graph coloring corresponding to  a 4-coloring $f$ of a 4-colorable maximal planar graph $G$. The specific method is: $\textcircled{1}$ delete the vertices from $G$ that belong to one of the same independent set generated by a 4-coloring $f$ of $G$, thus the 4-coloring problem of a maximal planar graph transform a 3-coloring problem of a planar graph  correspondingly, and the structural problem of six bicolored subgraphs reduce to three bicolored subgraphs' structural problem correspondingly. So, not only do computation reduce largely, but also the structure becomes simple and easy to study. $\textcircled{2}$ Furthermore, in the process of researching the three bicolored subgraphs, we study the union structure of them and any two of them, respectively. The findings show that it is very important to study tree-colorings in $C_4^0(G)$ for JT-conjecture and other problems of graph coloring, so the tree-coloring and cycle-coloring are studied preliminary in this section.

\subsection{Cycle-coloring and tree-coloring}

\quad\quad A \textbf{cycle-coloring} $f$ of a 4-colorable maximal planar graph $G$ is a 4-vertex-coloring of $G$ such that there exists a cycle $C_{2m}={v_1v_2\cdots v_{2m}v_1}(m\geq 2)$ in $G$ satisfying $|\{f(v_1),f(v_2),\cdots, f(v_{2m})\}|=2$, where $V(C_{2m})=\{v_1,v_2,\cdots, v_{2m}\}$. It refer to $C_{2m}$ as a \textbf{bicolored cycle} of $f$, or say $f$ contains a \textbf{bicolored cycle}. Suppose $G$ is a 4-colorable maximal planar graph and $f\in C_4^0(G)$, if $f$ don't contain a bicolored cycle, then $f$ is called a \textbf{tree-coloring} of $G$. It is clear to know from the definition of cycle-coloring and tree-coloring that for any a 4-colorable planar graph $G$ and $f\in C_4^0(G)$, $f$ is either a cycle-coloring or a tree-coloring.

For example, there are a total of eight 4-colorings for the graph shown in Figure $6.1$, and all these eight 4-colorings are cycle-colorings; for the 4-colorings shown in Figure $6.2$, $f_1,f_2$ are cycle-colorings and $f_3,f_4$ are tree-colorings; Figure $6.3$ gives all of the ten 4-colorings of the icosahedron and they are all tree-colorings. Naturally, we can know a fact that all the 4-colorings of a maximal planar graph maybe contain only tree-colorings, or only cycle-colorings, or both tree-colorings and cycle-colorings. Then, which of graphs contain only tree-colorings? Which of graphs only cycle-colorings? Which of graphs both tree-colorings and cycle-colorings? Obviously, these problems are the basis of studying the coloring properties of 4-colorable maximal planar graphs.

Considering the above three examples, the maximal planar graphs can be divided into three categories according to cycle-coloring and tree-coloring: $\textcircled{1}$ \textbf{pure cycle-coloring} graphs, namely these graphs have only cycle-colorings;  $\textcircled{2}$ \textbf{pure tree-coloring} graphs, namely have only tree-colorings; $\textcircled{3}$ \textbf{impure coloring} graphs, namely have both cycle-colorings and tree-colorings.

From now on, we use four different icons(shown in Figure 6.1(a)) to denote color 1,2,3,4, respectively.
  \begin{center}
        \includegraphics [width=220pt]{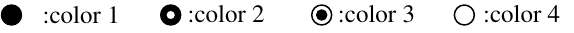}
        \vspace{1cm}

        \textbf{Figure 6.1(a).} The check figure between icons and colors.
  \end{center}
  \begin{center}
        \includegraphics [width=340pt]{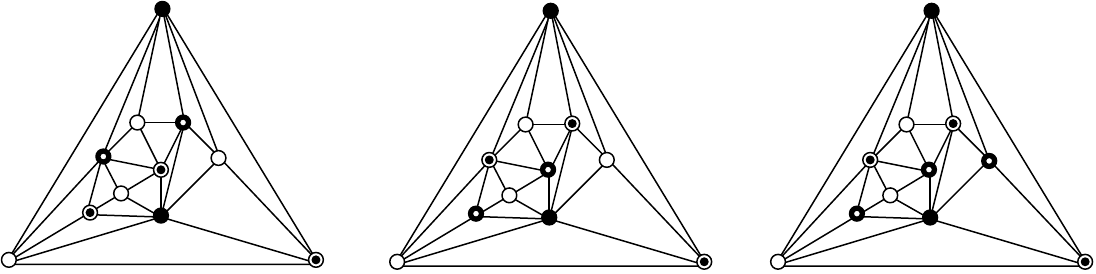}

       \small{ \hspace{0.8cm}$f_1$:1-4 cycle \hspace{1cm}$f_2$:1-3,1-4 and 3-4 cycles  \hspace{0.4cm}$f_3$:1-3 and 3-4 cycles}

        \includegraphics [width=340pt]{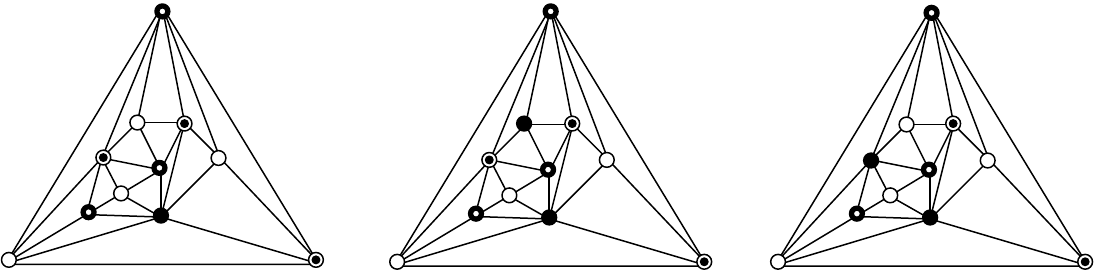}

        \small{\hspace{0cm}$f_4$:2-3,2-4 and 3-4 cycles \hspace{0.7cm}$f_5$:2-3 cycle \hspace{1cm}$f_6$:1-2,1-4 and 2-4 cycles}

        \includegraphics [width=260pt]{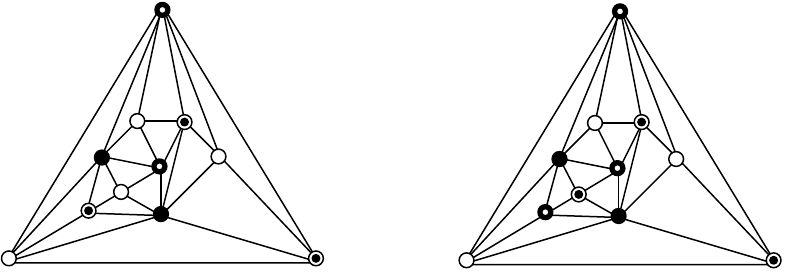}

         \small{\hspace{0cm}$f_7$:1-4 cycle \hspace{3cm}$f_8$:1-2 cycle}\\
        \textbf{Figure 6.1(b).} All the eight 4-colorings of a  maximal planar graph of order 11
  \end{center}

  \begin{center}

        \includegraphics [width=260pt]{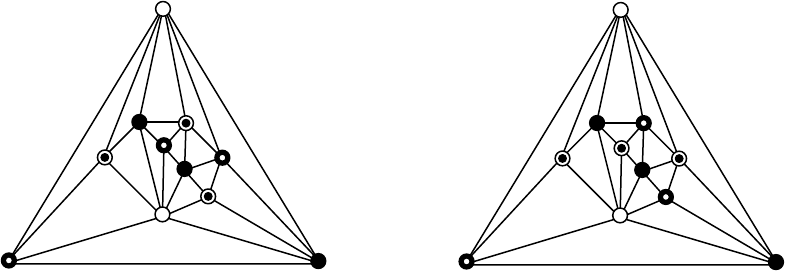}

       \small{\hspace{0cm}$f_1$:1-4 cycle \hspace{3cm}$f_2$:1-4 cycle}

        \includegraphics [width=260pt]{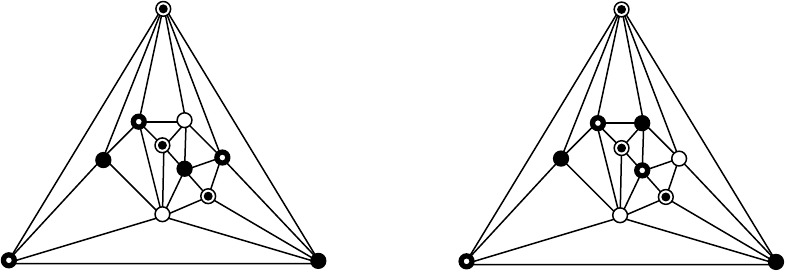}

         \small{\hspace{0cm}$f_3$:no bicolored cycle \hspace{2cm}$f_4$:no bicolored cycle }\\

        \textbf{Figure 6.2.} All the 4-colorings of a maximal planar graph with order 11
  \end{center}

  \begin{center}

        \includegraphics [width=320pt]{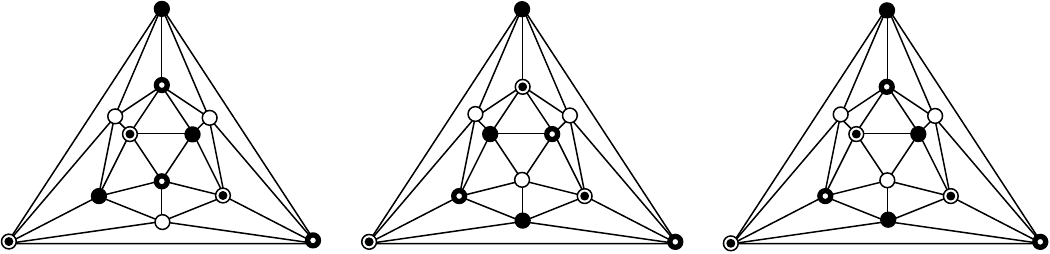}

        \includegraphics [width=320pt]{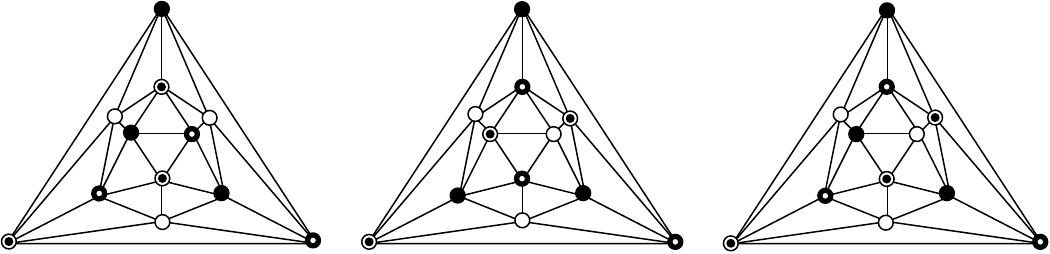}

        \includegraphics [width=320pt]{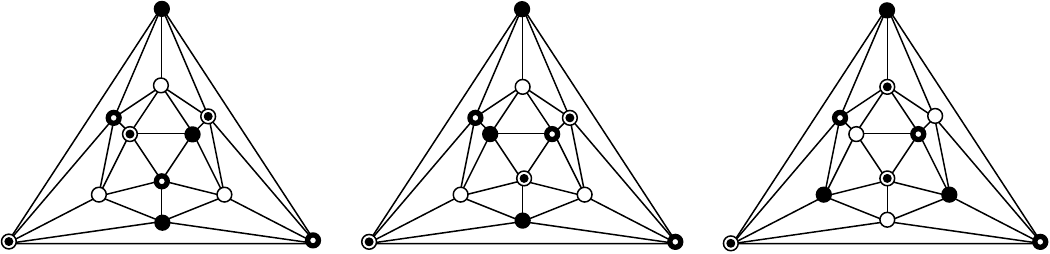}

        \includegraphics [width=105pt]{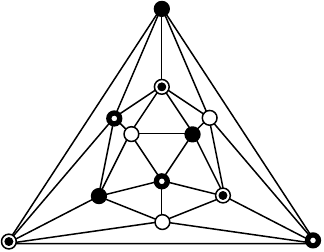}

        \textbf{Figure 6.3.} All the ten 4-colorings of icosahedron
  \end{center}

In terms of the relationship between 4-colorings and the structure of a maximal planar graph, the following results is obvious.

  \begin{theorem2}\label{th6.1} Let $G$ be a double-center wheel graph with $\delta(G)\geq 4$. Then $G$ is a pure cycle-coloring graph.
  \end{theorem2}
  \begin{proof} Let $u,v$ be the wheel-center vertices of $G$ and $f$ a 4-coloring of $G$. If $u,v$ are assigned the same color based on $f$, it is easy to infer that $f$ contains at least a bicolored cycle of length 4. Otherwise if $u,v$ are assigned different color, because $G\setminus \{u,v\}$ is a cycle $C$, then the length of $C$ must be even and the vertices of $C$ can be dyed only by two colors, so $f$ also contain a bicolored cycle.
  \end{proof}

  \begin{theorem2}\label{th6.2}  For the maximal planar graph $G_1$ and $G_2$  shown in Figure 6.4(a) and 6.4(b) respectively: $\textcircled{1}$ when $l$ is even,  $G_1$ has only one tree-coloring; $\textcircled{2}$ when $l$ is odd, $G_2$ has only one tree-coloring.
  \end{theorem2}

\begin{proof}
The proof is easy, so omit here. Now we give some examples of this theorem: the first and fourth graphs shown in Figure 6.5 illustrate the first case of this theorem and the sixth graph shown in Figure 6.5 illustrate its second case.
\end{proof}

\begin{center}

        \includegraphics [width=300pt]{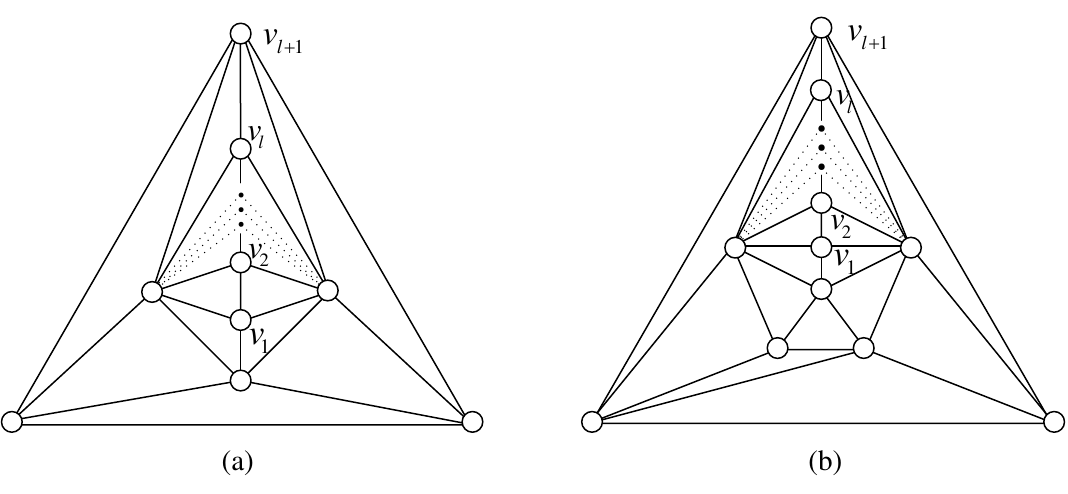}

        \textbf{Figure 6.4.} Two types of maximal planar graph with impure coloring
\end{center}

When a maximal planar graph $G$ contains 3-degree vertices, there has not any changed about the
coloring properties of the graph obtained by deleting these 3-degree vertices from $G$
comparing with the original graph $G$. So, we  need only consider the maximal planar graphs with minimum degree 4 or 5 when we study the coloring properties of them.

Considering the colorings of maximal planar graphs whose orders are from 7 to 11 and whose minimum degrees are not less than 4, we can obtain such a fact that tree-colorings are very few comparing with cycle-colorings. We can see these graphs and their 4-colorings in Appendix $I$, where there is only one graph with order 7 and it is  a double-center wheel graph, so it has not tree-colorings; two graphs with order 8: one is a double-center wheel graph and 3-colorable, the other has two cycle-colorings and one tree-coloring; five graphs with order 9: the first is 3-colorable and divisible, the second has only six cycle-colorings, the third which we refer to as \textbf{9-mirror graph} has only two tree-colorings, the fourth is double-center graph, the fifth has only seven cycle-colorings; thirteen graphs with order 10: just the second and twelfth have tree-colorings, and only one respectively, the others do not have tree-colorings; thirty-four graphs with order 11: only the tenth, twelfth, twenty-sixth, thirty-second and thirty-fourth have  tree-colorings and each of them at most contain two tree-colorings, the others do not contain tree-colorings.

\begin{table}[tbp]\label{table6.1}
\caption{\small{The cycle-coloring number and tree-coloring number of the maximal
planar graphs of order 7 to 11, the minimum degrees of which are 4 or 5.}}
\begin{tabular}{|c|c|c|c|c|c|c|c|c|c|c|c|}  
\hline
\small{GL} &7&$8_1$&$8_2$&$9_1$&$9_2$&$9_3$&$9_4$&$9_5$&$10_1$&$10_2$&$10_3$\\ \hline 
\small{CN} &5&$\ast$&2&$\ast$&6&0&17&7&7&10&$\ast$\\ \hline          
\small{TN} &0&$\ast$&1&$\ast$&0&2&0&0&0&1&$\ast$\\ \hline   
\small{GL} &$10_4$&$10_5$&$10_6$&$10_7$&$10_8$&$10_9$&$10_{10}$&$10_{11}$&$10_{12}$&$10_{13}$&$11_1$\\ \hline
\small{CN} &14&6&$\ast$&6&8&$\ast$&14&13&4&$\ast$&8\\ \hline          
\small{TN} &0&0&$\ast$&0&0&$\ast$&0&0&1&$\ast$&0\\ \hline         
\small{GL} &$11_2$&$11_3$&$11_4$&$11_5$&$11_6$&$11_7$&$11_8$&$11_9$&$11_{10}$&$11_{11}$&$11_{12}$\\ \hline
\small{CN} &$\ast$&25&29&41&85&14&10&10&5&8&2\\ \hline
\small{TN} &$\ast$&0&0&0&0&0&0&0&1&0&2\\ \hline
\small{GL} &$11_{13}$&$11_{14}$&$11_{15}$&$11_{16}$&$11_{17}$&$11_{18}$&$11_{19}$&$11_{20}$&$11_{21}$&$11_{22}$&$11_{23}$\\ \hline
\small{CN} &$\ast$&21&13&10&$\ast$&$\ast$&$\ast$&$\ast$&$\ast$&$\ast$&10 \\ \hline
\small{TN} &$\ast$&0&0&0&$\ast$&$\ast$&$\ast$&$\ast$&$\ast$&$\ast$&0 \\ \hline
\small{GL} &$11_{24}$&$11_{25}$&$11_{26}$&$11_{27}$&$11_{28}$&$11_{29}$&$11_{30}$&$11_{31}$&$11_{32}$&$11_{33}$&$11_{34}$\\ \hline
\small{CN} &$\ast$&8&5&$\ast$&11&13&17&9&12&16&6\\ \hline
\small{TN} &$\ast$&0&1&$\ast$&0&0&0&0&1&0&1\\ \hline
\end{tabular}
\vspace{0.2cm}

 \small{Where $GL$ denotes graph label; CN the number of cycle-coloring; TN the number of tree-coloring; $i_j$ the $j$th graph with order $i$ in appendix $I$; $\ast$ denotes the corresponding graph is 3-colorable or divisible.}
\end{table}

\newpage
From table $6.1$ we can see that there has much more cycle-colorings than tree-colorings. It is in total of fifty-five maximal planar graphs whose orders are from 7 to 11 and whose minimum degrees are 4 or 5, however only one is pure tree-coloring graph and we refer to this graph as \textbf{9-mirror graph}(the fourth graph in Figure 6.5); nine graphs contain at least one tree-coloring(see Figure 6.5). In addition, there are thirty pure cycle-coloring graphs, eight impure coloring graphs and sixteen divisible(or 3-colorable) graphs.  Apart from 3-colorable and divisible graphs in these fifty-five graphs, the number of 4-colorings of the remaining graphs is 518, but the number of tree-colorings is just 11 and share the proportion of 2.12\%, rarely!

\begin{center}

        \includegraphics [width=320pt]{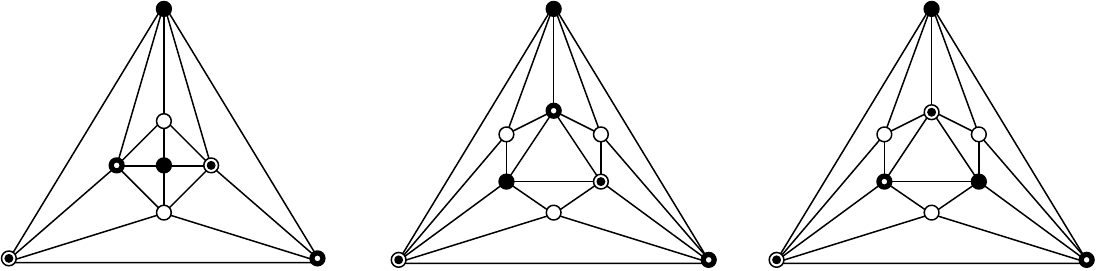}

        \includegraphics [width=320pt]{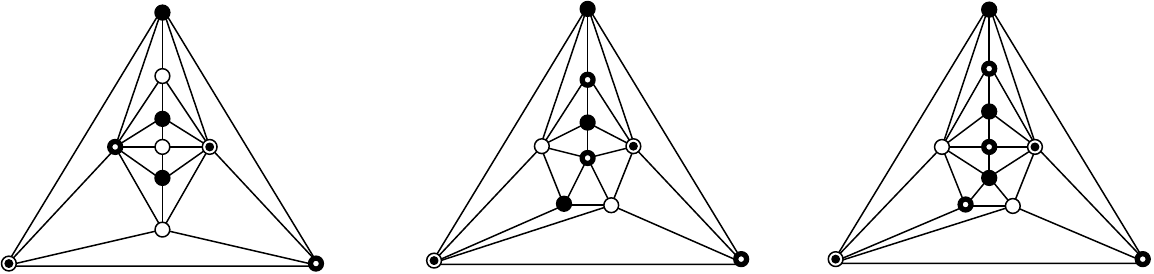}

        \includegraphics [width=320pt]{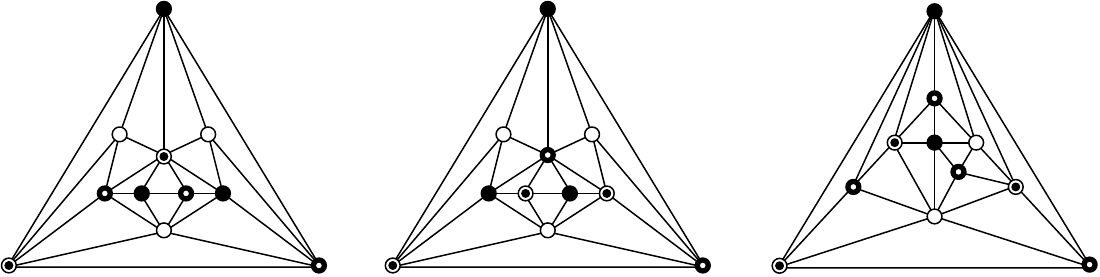}

        \includegraphics [width=210pt]{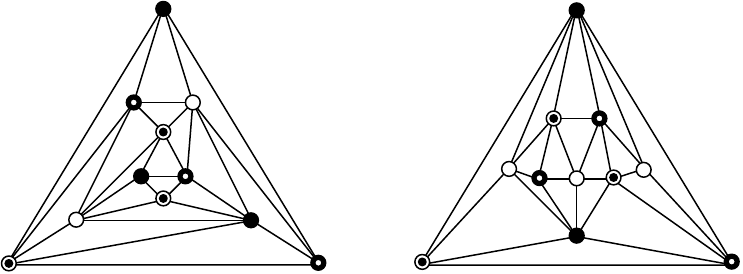}

        \textbf{Figure 6.5.} All the tree-colorings of the maximal planar graphs with orders from 7 to 11 and $\delta\geq 4$
\end{center}

For the pure tree-coloring graphs, there is an evident fact as follow:

\begin{Prop} If $G$ is a uniquely 4-colorable maximal planar graph, then $G$ is a pure tree-coloring graph.

\end{Prop}

So, all of the recursive maximal planar graphs are pure tree-coloring graphs. Up to now, we have found three pure tree-coloring graphs with minimum degree 4 or 5: one is the 9-mirror graph and the other two are the icosahedron(see Figure 6.3) and 13-mirror graph(see Figure 6.6), respectively.

\begin{center}

        \includegraphics [width=260pt]{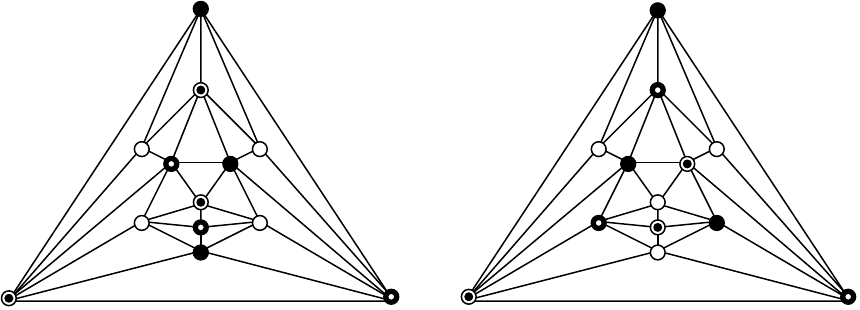}

        \includegraphics [width=260pt]{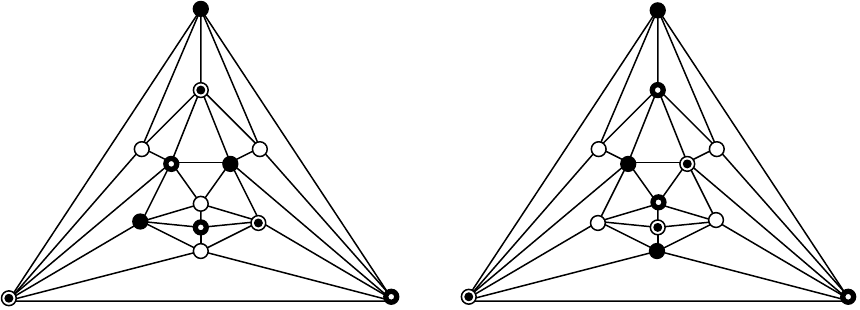}

        \textbf{Figure 6.6.} The third pure tree-coloring graph(13-mirror graph) and its 4-colorings
\end{center}

Naturally, an important problem will be proposed as follow:

\begin{problem} What is the characteristic of the pure tree-coloring graph whose minimum degree is not less than 4? And how many such graphs are there?
\end{problem}

Surprisingly! Only three maximal planar graphs, whose minimum degrees are not less than 4, are the pure thee-coloring graphs: the 9-mirror graph, 13-mirror graph and icosahedron. The detailed proof on this problem will be given in the following sections.

For the pure cycle-coloring graphs and the impure coloring graphs, we also propose two problems as follows:

\begin{problem}
What are the necessary and sufficient conditions that a maximal planar graph is a pure cycle-coloring graphs?
\end{problem}
\begin{problem}
What are the necessary and sufficient conditions that a maximal planar graph is a  impure coloring graph?
\end{problem}

In the same way, the detailed discussions about these two problems will also be given subsequently in the following sections.

\subsection{Equivalency of coloring  between tricolored induced subgraphs and maximal planar graphs}

\quad\quad For a given 4-colorable maximal planar graph $G$, let $f$ be a 4-coloring of $G$ and the color classes partition of $f$ is $\{V_1,V_2,V_3,V_4\}$, where $V_i$ denotes the set of vertices assigned color $i$. Obviously, when three color classes of them are determined, the last one is also determined uniquely. So, we only need make clear the 3-coloring structure of the tricolored induced subgraph that are induced by any three classes partition, such as $V_1,V_2,V_3$.

For the sake of convenient expression, here we introduce a definition of big-cycle. Let $C$ be a subgraph of a  planar graph $G$. If $C$ is a facial cycle and the length of it is not less than 4, then we call $C$ a \textbf{big-cycle}.

The following gives some examples that illustrate the equivalency of colorings between a maximal planar graph and its tricolored induced graph keeping that each of the big-cycles is colored by three colors. For the first graph shown in Figure 6.7, it has three different 4-colorings totally(see Figure 6.7(a),(b),(c)). If denote by $V_4$ the set consisting of vertices received by color 1, then the three 3-colorings of $G-V_4$ corresponding to Figure 6.7(a),(b) and (c) are shown as Figure 6.7(a'),(b') and (c'). For the first graph shown in Figure 6.8(a), when delete the vertices received by color 1, we can obtain its tricolored induced subgraph  and the corresponding 3-coloring(see Figure 6.8(b)). Figure 6.8(c) and (d) exhibit two 4-colorings of this subgraph satisfying that each of its  big-cycles  is colored by at most three colors(here just three).

\begin{center}

        \includegraphics [width=320pt]{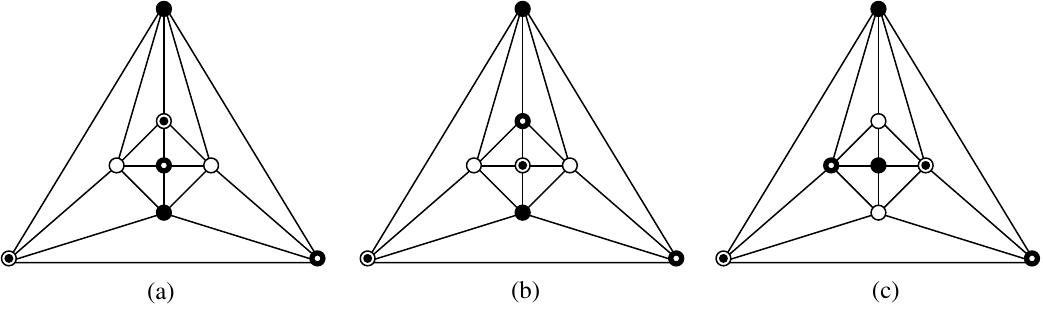}
        \includegraphics [width=320pt]{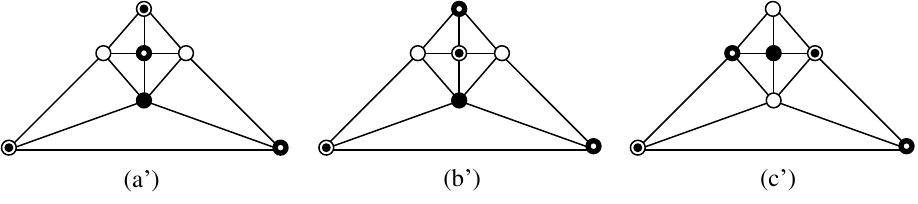}

        \textbf{Figure 6.7.} All the 4-colorings of a maximal planar graph of order 8 and  the corresponding 3-colorings of their tricolored induced subgraphs.
\end{center}

\begin{center}

         \includegraphics [width=340pt]{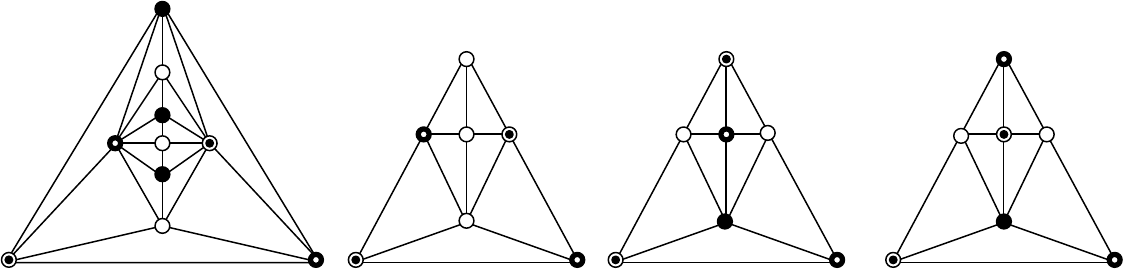}

        \textbf{Figure 6.8.} A 4-coloring of a maximal planar graph of order 10 and  the corresponding colorings of their tricolored induced subgraphs.
\end{center}

It is clear from the above two examples that there has the intuitive and understandable advantage when we study 4-coloring problem from $G-V_4$ comparing to from $G$ directly. So, when we study the 4-coloring problem of a maximal planar graph $G$, we need only research the 4-coloring problem of $G-V_4$. However, in the process of studying this problem, we need to pay attention to three points as follows:

The first, in terms of the choice of $V_4$, we should try to choose such a color class that contain the most vertices  as $V_4$, so that $G-V_4$ can become simple as much as possible.

The second, because the number of triangles in a maximal planar graph of order $n$ is $2n-4$, the number of triangles in $G-V_4$ is
$$
2n-4-\sum\limits_{v\in V_4}d_G(v)
 \eqno{(6.1)}
$$
For example, Figure 6.7 exhibit a maximal planar graph with order $n=8$ and the sum of the degrees of the two vertices in $V_4$ is 10, so the number of triangles in $G-V_4$ is $2\times 8-4-10=2$;  for the maximal planar graph of order 10 shown in Figure 6.8, similarly, we can calculate the number of triangles in $G-V_4$ is 3 by Formula 6.1.

The third, following equality holds:
$$
G[V_1\cup V_2\cup V_3]=G[1,2]\cup G[1,3]\cup G[2,3]
 \eqno{(6.2)}
$$
So, we should study the structure and property of tricolored induced subgraph $G[V_1\cup V_2\cup V_3]$ from any two bicolored induced subgraph, such as $G[1,2]\cup G[1,3]$. Studying the structure and property of $G[V_1\cup V_2\cup V_3]$ gradually is the basic idea in this section.

\begin{theorem2}Suppose $G$ is a 4-colorable maximal planar graph and $f$ is a 4-coloring of $G$, if the color classes partition of $f$ is $\{V_1,V_2,V_3,V_4\}$, then

\textcircled{1} $f$ is a tree-coloring of $G$ if and only if $f$ is a tree-coloring  limited to $G[V_1\cup V_2\cup V_3]$;

\textcircled{2} $f$ is a cycle-coloring of $G$ if and only if $f$ is a cycle-coloring( or a disconnected coloring) limited to $G[V_1\cup V_2\cup V_3]$;\\
where``limited to $G[V_1\cup V_2\cup V_3]$" refer to the color classes partition of $f$ only presenting on $G[V_1\cup V_2\cup V_3]$; the disconnected coloring refers there are disconnected bicolored induced subgraph in $G[V_1\cup V_2\cup V_3]$ under the coloring $f$.
\end{theorem2}

\begin{proof}
These two results is easy to prove, so omit here. Obviously, this theorem can bring us some convenience when we judge whether a coloring is a tree-coloring. So, in the following study, we will mainly consider the coloring structure of tricolored induced subgraph.
\end{proof}

\subsection{Union structure of two bicolored induced subgraphs}
\quad\quad This subsection introduce the concept of fence, and prove that the union of any two
bicolored induced subgraphs of 4-colorable maximal planar graph is a fence.
Further, some special fences and the characteristic of maximal planar graphs
corresponding to them are discussed.

\subsubsection{General theory}

\quad\quad The degree of a face is the number of
edges in its boundary, cut edges being counted twice. Let $G$ be a planar graph. If
the degree of every face of $G$ is even and not less then 4, then $G$ is called the \textbf{fence}.
The graphs shown in Figure 6.9(a), (c), (d) are fences, but the graph shown in Figure 6.9(b)
is not a fence, because there exists an odd cycle in this graph. For a fence $G$,  it may not have any suspending vertices, of course, it  may also contain several suspending vertices. Where the suspending vertices refer to the vertices with degree 0 or 1.
If there exist a suspending vertex $v$ in $G$, then the subtree containing $v$  maybe connect with a cycle  by a common vertex, say $u$ and called \textbf{weld-vertex}. Denote by $t$
the distance between $u$ and $v$, then there exists a path of length $t$ between them.
Choose a maximum $t$ and refer to $G$ as a \textbf{$t$-fence}.  If $G$ has no suspending vertex, it is called a $0$-fence.
The graph shown in
Figure 6.9(a) is a 1-fence; the graph shown in Figure 6.9(c) is a 0-fence; the graph shown in
Figure 6.9(d) is a 2-fence; the graph shown in Figure 6.9(e) is a 3-fence. If there exists no
path between a suspending vertex and any cycles of $G$, that is to say, the graph $G$ is disconnected,
then $G$ is called a $\infty$-fence. The graph shown in Figure 6.9(f) is a $\infty$-fence.

\begin{center}
        \includegraphics [width=260pt]{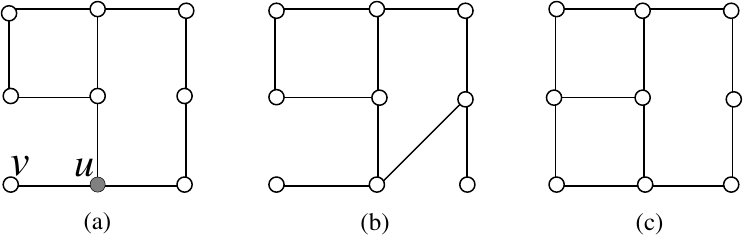}

        \includegraphics [width=260pt]{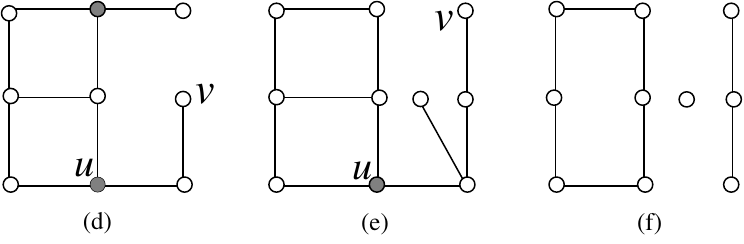}

        \textbf{Figure 6.9.} The illustration of the concept of fence.
\end{center}

\begin{theorem2}\label{th}
Let $G$ be a 4-colorable maximal planar graph, $C(4)=\{1,2,3,4\}$ the color set. For any a 4-coloring $f$ and the union of any two bicolored induced subgraphs with a common color, say $G[1,2]\cup G[1,3]$, there are following results.

 (1) $G[1,2]\cup G[1,3]$ has no odd cycle;

 (2) If the order of both $G[1,2]$ and $G[1,3]$ are at most 3, then for graph $G[1,2]\cup G[1,3]$, it is a cycle under only one case and a tree in other cases;

(3) If $f$ is a tree-coloring and the order of  $G[1,2]$ or $G[1,3]$ is at least 4, then $G[1,2]\cup G[1,3]$
is a 1-fence or 0-fence;

(4)  If $f$ is a tree-coloring, then every suspending vertex of $G[1,2]\cup G[1,3]$ must adjacent to the vertices
colored by the common color 1.
\end{theorem2}
 \begin{proof}
 (1) Assume that $G[1,2]\cup G[1,3]$ has a odd cycle $C$. Since the set of colors appeared on the cycle $C$ must contain color 1,
2 and 3, so $C$ contains not only the 1-2 edges and 1-3 edges, but also the 2-3 edges. But $G[1,2]\cup G[1,3]$ has no edges in $G[2,3]$, it is a contradiction. Where $i-j$ edge denotes the edge whose two end vertices are colored by color $i$ and $j$ respectively.

 (2) If the order of both $G[1,2]$ and $G[1,3]$ are not more than 3, then all 4 cases are shown in Figure 6.10. It is easy
 to see that only in one case is $G[1,2]\cup G[1,3]$ a cycle and a tree in other cases.

  \begin{center}

        \includegraphics [width=260pt]{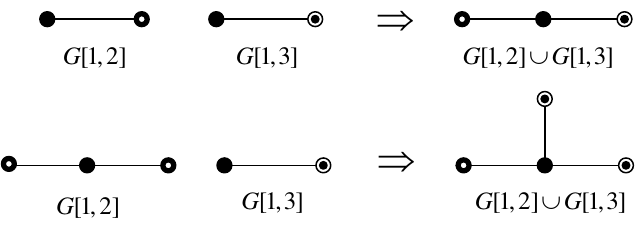}

        \includegraphics [width=260pt]{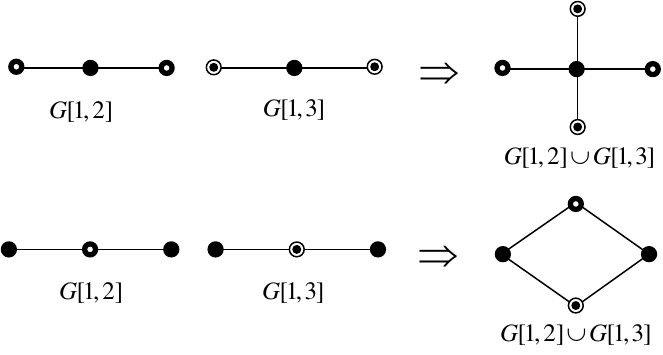}

        \textbf{Figure 6.10.} Four cases of $G[1,2]\cup G[1,3]$ in Theorem 6.5(2).
  \end{center}

(3) Based on the case (1) and (2), we study the case of the order of  $G[1,2]$ or $G[1,3]$ is at least 4.
Suppose $v$ is a suspending vertex of $G[1,2]\cup G[1,3]$ and the distance of $v$ to the nearest cycle in $G[1,2]\cup G[1,3]$ is 2. The unique vertex adjacent to $v$ in $G[1,2]\cup G[1,3]$ is denoted by $v^{\prime}$.
 If $f(v)=1$, since $f$ is a tree-coloring, both $G[1,2]$ and $G[1,3]$ are connected, so $v$ is adjacent to some vertices
 colored by color 2 and 3 respectively. Thus, $v$ is adjacent to at least two
 vertices in $G[1,2]\cup G[1,3]$, it is a contradiction to the fact that $v$ is a suspending vertex.   If $f(v)=2$(or $f(v)=3$),
 then $f(v^{\prime})=1$ and the vertex $w$($\neq v$) adjacent to $v^{\prime}$ may be assigned by two colors. If $f(w)=3$,
 then $G[1,2]$ is disconnected; if $f(w)=2$, then $G[1,3]$ is disconnected. There is a contradiction for both of these two cases.  So, $G[1,2]\cup G[1,3]$
is a either 1-fence or 0-fence.

(4) This case is obvious, so omitted here.
\end{proof}

The above results show that: for a 4-coloring $f$ of a 4-colorable maximal planar graph $G$,  $G[1,2]\cup G[1,3]$ don't contain odd cycles, that is to say, it is a tree or contain only even cycles; furthermore, for the structure including the even cycles, we prove that $G[1,2]\cup G[1,3]$ is a 1-fence or 0-fence when $f$ is a tree-coloring. In addition, for a fence, if it do  contain suspending vertices, they must be adjacent to the vertices colored by the common color 1.

In Figure 6.11(b), it is easy to see that the resulted graph by adding a new vertex $v$ to every face of degree at least 4 and connecting
$v$ to all vertices on the cycle of the face including $v$ is a maximal planar graph $G$ with $\delta(G)\geq 4$.  The graph shown in Figure 6.11(a) exhibit the case that the union of two bicolored induced subgraph of a
maximal planar graph $G$ with $\delta(G)\geq 4$ has no suspending vertex.

\begin{center}

        \includegraphics [width=280pt]{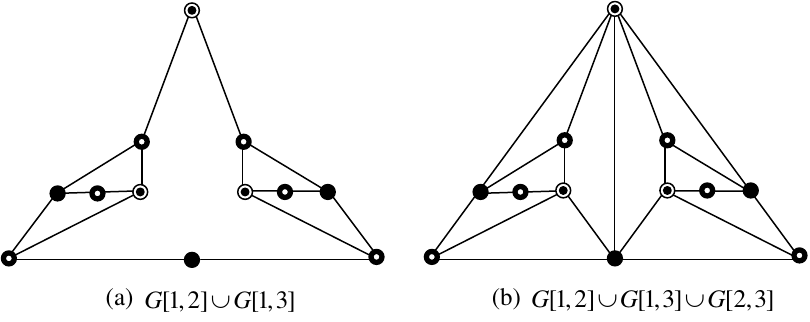}

        \textbf{Figure 6.11.} The illustration for the union of two bicolored induced subgraph has no suspending vertex.
  \end{center}

\subsubsection{Structure of the bicolored induced subgraphs and their union of the pure tree-coloring graph}

The necessary and sufficient condition, which a maximal planar graph $G$ is a pure
tree-coloring graph is that $G$ is the icosahedron, 9-mirror graph or 13-mirror graph. This statement will be discussed in the following sections.
Now, we analyze the structure of the union of bicolored induced subgraphs of these three maximal planar graphs.

For the 9-mirror graph(see the 2th, 3th graph in Figure 6.5), it is easy to prove that its all bicolored induced subgraphs are paths,
and the length of them has only two kinds: one is 4 and another is 3. Concerning the structure of the union of two bicolored induced subgraphs with a common color, there are three cases: $\textcircled{1}$ the union consist of two paths of length 4(see Figure 6.12(a)): it is a fence without suspending vertices and includes two cycles of length 6 and one cycle of
length 4. $\textcircled{2}$ The union consist of two paths with length 3 and 4 respectively: it is a 1-fence with 3 suspending vertices and
a cycle of length 4(see Figure 6.12(b)). $\textcircled{3}$ The union consist of two paths of length 3: this union is a 1-fence with 2 suspending vertices and a cycle of length 4(see Figure 6.12(c)).

For the icosahedron(shown in Figure 6.3),
all of its bicolored induced subgraphs are isomorphic and each of them is a path of length 5(see Figure 6.12(d)). The union of two bicolored induced subgraphs, which consist of two paths of length 5 with a common color(see Figure 6.12(d)), is formed by adding two suspending vertices on the basis of the graph shown in Figure 6.12(a).

\begin{center}

         \includegraphics [width=260pt]{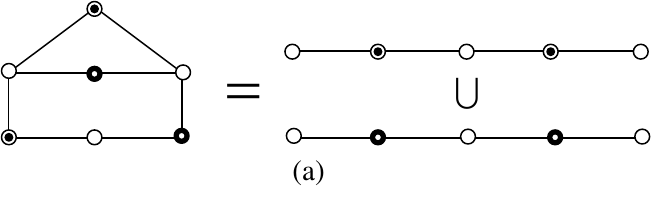}

         \includegraphics [width=260pt]{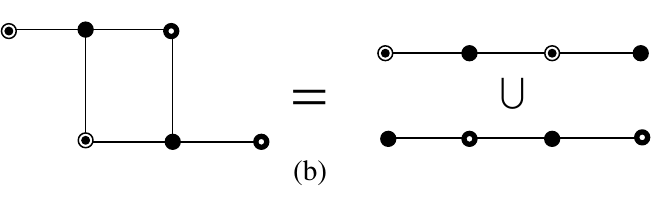}

         \includegraphics [width=260pt]{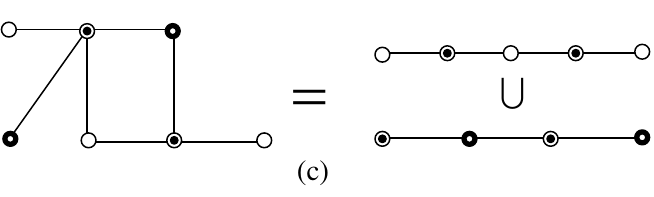}

         \includegraphics [width=260pt]{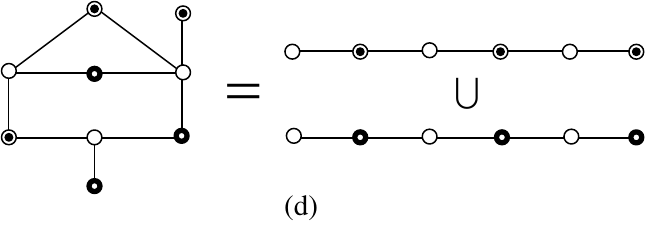}

         \includegraphics [width=260pt]{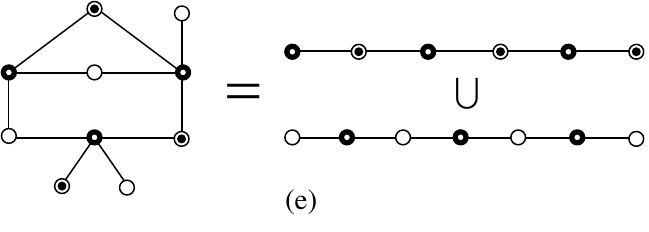}

         \includegraphics [width=280pt]{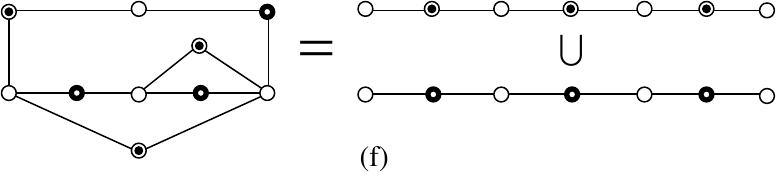}

         \includegraphics [width=280pt]{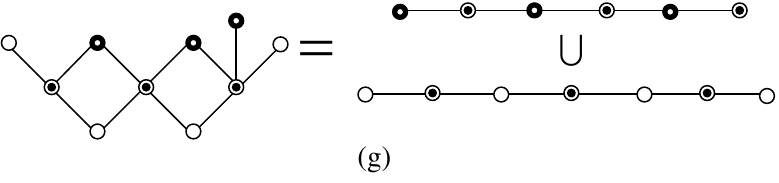}

        \textbf{Figure 6.12.} Structure of all bicolored induced subgraph of the icosahedron, the 9-mirror graph  and 13-mirror graph.
  \end{center}

 For the 13-mirror graph(shown in Figure 6.6),  all of its bicolored induced subgraphs  are paths and the length of them has only two kinds: one is 5 and another is 6. Concerning the structure of the union of two bicolored induced subgraphs with a common color, there are four cases shown  in Figure 6.12(d), 6.12(e), 6.12(f), 6.12(g) respectively.

Further, the icosahedron, the 9-mirror graph and 13-mirror graph have several properties as follows:

 $\textcircled{1}$ Each of the bicolored induced subgraph of the icosahedron is a path of length 5; each of the bicolored
 induced subgraph of the 9-mirror graph  is also a path, the length of which is 5, 4 or 3; each of the bicolored
 induced subgraph of the 13-mirror graph is also a path, the length of which is 6 or 5, and there are two paths of length 6.
Of course, all of the bicolored induced subgraph of these three graphs have no vertex with degree at least 3.

$\textcircled{2}$ For these three graphs, there are two cycles of length 4, two cycles of length 6 and only one cycle of length 8 in the union of  the two bicolored
induced subgraphs with a common color.

 $\textcircled{3}$ These three graphs are symmetrical strongly.

The foregoing discussions inspire us to excogitate an idea that  prove
the necessary and sufficient condition: for a maximal planar graph $G$, it is a pure
tree-coloring graph if and only if $G$ is the icosahedron, the 9-mirror graph or 13-mirror graph. The details as follows:

The first, if some  bicolored induced subgraph of $G$ has a vertex of degree at least 3, then $G$ isn't a pure tree-coloring graph;

The second, if there is a bicolored path with length at least 7 in $G$,
  then $G$ isn't a pure tree-coloring graph.

  These two aspects tell us that: if $G$ is a pure tree-coloring maximal planar graph, then
  the order of $G$ is at most $7+3+3=13$, and each of the bicolored induced subgraph of $G$ is a path. In fact, for a tree-coloring $f$ of $G$ and a given bicolored path $P$ of length 6, suppose the vertices of $P$ are assigned color 1 and 2, where the number of vertices colored by 1 is four, then the length of the path  whose vertices are assigned color 3 and 4 is at most 5. Otherwise, there must is a path whose vertices are assigned color 1 and 3(or 1 and 4) with length at least 7, contradiction!  So the order of $G$ is at most 13.

  The third is to seek the pure tree-coloring graphs among the maximal planar graphs with order at most 13 and
  this task is easy to complete.

 \subsubsection{Structure on the union of two bicolored induced subgraphs whose non-suspending vertices are totally in a cycle}

This subsection study the structure of a class of specific graphs that the union of two bicolored induced subgraphs
whose  non-suspending vertices are totally in a cycle.

Let $G$ be a 4-coloring maximal planar graph, $f$ a tree-coloring of $G$, $G[1,2]$ and $G[1,3]$  any two bicolored induced subgraphs having a common color,
then $G[1,2]\cup G[1,3]$ is a 1-fence or a 0-fence.
For the graph $G[1,2]\cup G[1,3]$, suppose the number of vertices assigned  color 1, 2 and 3 are $a, b$ and $c$, respectively.
If all of the non-suspending vertices are in a cycle, then $G[1,2]\cup G[1,3]$ has several properties as follows:

\textbf{Property 1.} Besides the edge incident with the suspending vertex, each edge is in either one or two cycles of length 4;

\textbf{Property 2.} There are $a-i$ cycles of length $2i+2$, $i=1,2,\cdots,a-1$, where the cycles of length 4 are adjacent in turn and any
adjacent two have just a common edge;

\textbf{Property 3.}  The number of vertices with degree 1 is $b+c-a$.

It is easy to prove that:
 \begin{theorem2}\label{th6.6}
Let $G$ be a 4-colorable maximal planar graph, $C(4)=\{1,2,3,4\}$ the color set and $f$ a 4-coloring of $G$.
If the order of
$G[1,2]$ or $G[1,3]$ is 3, say $G[1,2]$, then  $\delta(G)=3$ in condition of odd suspending vertices in $G[1,2]\cup G[1,3]$, or two vertices receiving color 1 in $G[1,2]$.
 \end{theorem2}

\subsection{Construction of tricolored induced subgraphs}

The Theorem 6.4 tells us that for a 4-coloring $f$ of a maximal planar graph, its properties can be  characterized by its tricolored
induced subgraph completely.
In fact, this subgraph is the union of three bicolored induced subgraphs. It have made clear the structure on the union of two bicolored induced subgraphs, which is the basis of studying the structure of tricolored induced graphs.

Without loss of the generality, we still denote by $G[1,2]\cup G[1,3]$ the union of any two bicolored induced subgraphs here and we will construct the tricolored induced subgraphs based on it by connecting the edges between the vertices receiving color 2 and color 3.
We should remark that if $G[1,2]\cup G[1,3]$ is connected, there is always at least a path $uu_1^{1}u_2^{2-3}u_3^{1}u_4^{2-3}\cdots u_k^{1}u^{\prime}$ between any pair of vertices $u$ and $u'$ that are assigned color 2 and 3 respectively.
Where $u_i^{1} (1\leq i \leq k)$  denotes the number of vertices colored 1, $u_i^{2-3} (2\leq i \leq k)$ denotes
the number of vertices colored by 2 or 3. Obviously, there are odd vertices on the path between the vertices $u$ and $u^\prime$. So, we have

\begin{theorem2}\label{th}
Let $G$ be a 4-colorable maximal planar graph, $C(4)=\{1,2,3,4\}$  the color set, $f\in C_4^0(G)$ and $\{u,u^\prime\}$ a pair of vertices colored by 2 and 3 respectively, then any path between
$u$ and $u^\prime$ has odd vertices in the connected subgraph $G[1,2]\cup G[1,3]$ .
\end{theorem2}

This theorem shows that: in order to construct $G[V_1\cup V_2\cup V_3]=G[1,2]\cup G[1,3]\cup G[2,3]$ based on
$G[1,2]\cup G[1,3]$, when $G[2,3]$ is a tree with $q$ edges, we need only to connect $q$ edges which can form a tree between the vertices colored by 2 and 3 in $G[1,2]\cup G[1,3]$.

Since $G[1,2]\cup G[1,3]$ has only even cycles and the length of any pair of vertices colored by 2 and 3 is an odd number, so every edge of $G[2,3]$ contributes to $G[V_1\cup V_2\cup V_3]$ at least one odd cycle. Now it is discussed in detail
by three cases.

 \textbf{Case 1.} Cycle-cycle edge: the edge in $G[2,3]$ and the two ends $u$ and $u'$ of this edge are on some cycles in $G[1,2]\cup G[1,3]$. For this case, by Theorem 6.7
 we know that there are at least two different paths of even length. So $G[V_1\cup V_2\cup V_3]$
 is contributed at least two odd cycles when  $u$ and $u'$ are joined in $G[1,2]\cup G[1,3]$.

\textbf{Case 2.} Cycle-suspending edge: the edge in $G[2,3]$ and is formed by joining a vertex on a cycle and a suspending vertex in $G[1,2]\cup G[1,3]$.

\textbf{Case 3.} Suspending-suspending edge: the edge in $G[2,3]$ and is formed by joining two suspending vertices in $G[1,2]\cup G[1,3]$.

\textbf{Example 6.1.} For the 4-colorable maximal planar graph $G$(shown in Figure 6.13(a)), let $f$ be a 4-coloring of $G$. The subgraph shown in Figure 6.13(b) is obtained by deleting the subset consisting of all vertices assigned color 4; the subgraph shown in Figure 6.13(c) is $G[1,2]\cup G[1,3]$, in which there are only two suspending vertices $v_3$ and $v_8$. Thus, the edges $v_4v_8$ and $v_2v_3$ are so-called
cycle-suspending edges, and the edge $v_3v_8$ is so-called suspending-suspending edge(see Figure 6.13(b)).

Since there exist two paths $v_4v_7v_8$ and $v_4v_1v_2v_7v_8$ from $v_4$ to $v_8$ in graph $G[1,2]\cup G[1,3]$, so the edge $v_4v_8$
contributes  to $G[V_1\cup V_2\cup V_3]$ two cycles $v_4v_8v_7v_4$ and $v_4v_1v_2v_7v_8v_4$. This also illustrate Theorem 6.7.
Similarly, the edge $v_2v_3$ contributes to $G[V_1\cup V_2\cup V_3]$  two odd cycles: $v_2v_1v_3v_2$ and  $v_2v_7v_4v_1v_3v_2$. In addition, because there exist two paths $v_3v_1v_4v_7v_8$ and $v_3v_1v_2v_7v_8$ from $v_3$ to $v_8$ in graph $G[1,2]\cup G[1,3]$,
the edge $v_3v_8$ contributes to $G[V_1\cup V_2\cup V_3]$ two cycles $v_3v_1v_4v_7v_8v_3$ and $v_3v_1v_2v_7v_8v_3$.

In graph $G[1,2]\cup G[1,3]$, since there exist two paths $v_3v_1v_4$ and $v_3v_1v_2v_7v_4$ from $v_3$ to $v_4$, so the two edges in $G[2,3]$
$v_4v_8$ and $v_8v_3$  contribute to $G[V_1\cup V_2\cup V_3]$ two even cycles $v_3v_1v_4v_8v_3$ and $v_3v_1v_2v_7v_4v_8v_3$. Similarly, because there exist two paths $v_2v_7v_8$ and $v_2v_1v_4v_7v_8$ from $v_2$ to $v_8$,  the two edges in $G[2,3]$
$v_2v_3$ and $v_3v_8$  contribute to $G[V_1\cup V_2\cup V_3]$ two even cycles $v_2v_7v_8v_3v_2$ and $v_2v_1v_4v_7v_3v_8v_2$.

\begin{center}

         \includegraphics [width=380pt]{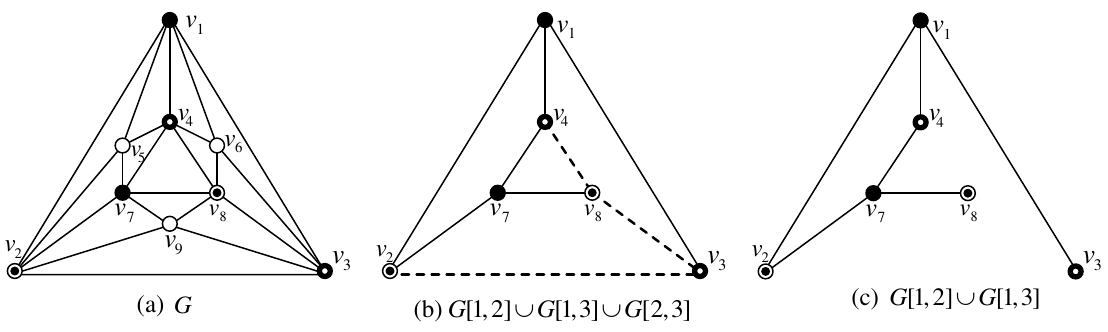}

        \textbf{Figure 6.13.}  The illustration of the concepts of cycle-cycle edge, cycle-suspending edge and suspending-suspending edge
  \end{center}

Furthermore, we consider the three edges $v_4v_8$, $v_8v_3$ and $v_3v_2$ in $G[2,3]$. Since there exist two paths $v_2v_7v_4$ and $v_2v_1v_4$
from $v_2$ to $v_4$ in graph $G[1,2]\cup G[1,3]$, the three edges in $G[2,3]$
$v_4v_8$, $v_8v_3$ and $v_3v_2$ contribute to $G[V_1\cup V_2\cup V_3]$ two odd cycles $v_2v_7v_4v_8v_3v_2$ and $v_2v_1v_4v_8v_3v_2$.

From the discussion of Example 6.1, we can obtain the following theorem.

\begin{theorem2}\label{th}
Let $G$ be a 4-colorable maximal planar graph, $C(4)=\{1,2,3,4\}$ the color set and $f\in C_4^0(G)$. Suppose the number of the paths
from the vertex $u$ to $u^\prime$ is $q$ in $G[1,2]\cup G[1,3]$, for any a path $P$ of length $p$ in $G[2,3]$, then

$\textcircled{1}$ $q$ is a even number;

$\textcircled{2}$ $G[1,2]\cup G[1,3]\cup P$  has $q$ odd(even) cycles including
$P$ more than $G[1,2]\cup G[1,3]$ when $p$ is a odd(even) number.
 \end{theorem2}

This theorem is easy to prove, so omitted here.

 It is clear from the Theorem 6.8 that the coloring structure of a 4-colorable maximal planar graph $G$ under a 4-coloring $f$.
Let $V_1,V_2,V_3,V_4$ be four independent sets of $G$ based on $f$, then

$\textcircled{1}$ the coloring structure of $G$ corresponding to $f$ is equivalent to $G[V_1\cup V_2\cup V_3]$;

$\textcircled{2}$ $G[1,2]\cup G[1,3]$ is a fence;

$\textcircled{3}$ $G[V_1\cup V_2\cup V_3]$ is obtained by adding continuously the edges of $G[2,3]$ in $G[1,2]\cup G[1,3]$.
$G[V_1\cup V_2\cup V_3]$ is contributed several odd cycles when a edge of $G[2,3]$ is added to $G[1,2]\cup G[1,3]$, furthermore, $G[V_1\cup V_2\cup V_3]$ is contributed
 several odd(even) cycles including $P$ when a odd(even) path $P$ of $G[2,3]$ is added to $G[1,2]\cup G[1,3]$.

\section{Black-White coloring, and necessary and sufficient conditions for 2-colorable cycle}

In this section, we will propose a new coloring method, saying Black-White coloring.
For a maximal planar graph $G$, the Black-White coloring of $G$ is to assign only two colors, black and white, to the vertices of $G$. At the beginning of an even cycle, either black or white is colored to the  vertices  according to a definite rule. Based on this method, if all the vertices of $G$ are colored by black or white, then we can deduce that $G$ is 4-colorable if no odd-cycles are included in the induced subgraphs by the vertices that are colored the same color(black or white). In fact, the Black-White coloring provides a subset of 4-colorings for maximal planar graphs.

Depending on the Black-White coloring, we plan to yield a necessary and sufficient condition that an even cycle in a maximal planar graph is 2-colorable. For this we put forward some new concepts, like closed-maximal planar graphs, opened-maximal planar graphs, semi-maximal planar graphs and 2-colorable cycles. Further, we conduct research deeply in the characteristics of even-cycles in a maximal planar graph and discuss the enumeration of even-cycles.

Let $G$ be a maximal planar graph and $C$ a cycle in $G$ with length not less than 4. We refer to the subgraph of $G$ that is induced by the vertex set consisting of the vertices of $C$ and the inside component of $C$, as the \textbf{semi-maximal planar graph}. In other words, the so-called semi-maximal planar graphs are just a kind of special planar graphs, in which the boundary of their infinite faces are cycles with length not less than 4 and other faces are triangles. We also call this sort of graphs \textbf{the semi-maximal planar graphs on $C$ }, which is written simply as $G^C$. Obviously, for $G$ and $C$, there are just two semi-maximal planar graphs on $C$, and we call them the semi-maximal planar graphs on $C$ of $G$.

\subsection{Characteristics and distribution of the even-cycles}

Even-cycle is the basic element of cycle colorings, so considering the characteristics and distribution of the even-cycles in a maximal planar graph is very important to research whether it contains cycle colorings. Suppose that $C=v_1v_2\cdots v_mv_1$ is a cycle with length $m$. A $\textbf{chord} $ on a cycle $C$ is an edge that link two nonadjacent vertices $v_i,v_j$ of $C$. We refer to the cycle containing chords as $\textbf{chord-cycle}$. Both the longest cycles in Figure 7.1(b) and 7.1(c) are chord-cycles.

Thus it is clear that for a maximal planar graph and its any connected subgraph $H$, the subgraph induced by the neighbor set of $H$ contains either a cycle or  a chord-cycle \emph{except for a tree}. For example,
in Figure 7.1(d), the cycle $C=v_1v_2v_3v_5v_{10}v_1$ is a chord-cycle(see Figure 7.1(e)). In addition, there are a special kind of chord-cycles, called \textbf{one-side cycle}, which contain no vertices inside or outside(in Figure 7.1(f), the cycle $C=v_1v_2v_3v_4v_1$ in Figure 7.1(d) is a one-side cycle).

 \begin{center}
      \includegraphics [width=180pt]{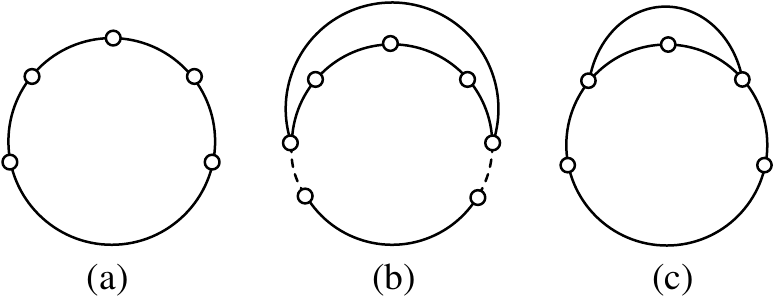}

       \includegraphics [width=320pt]{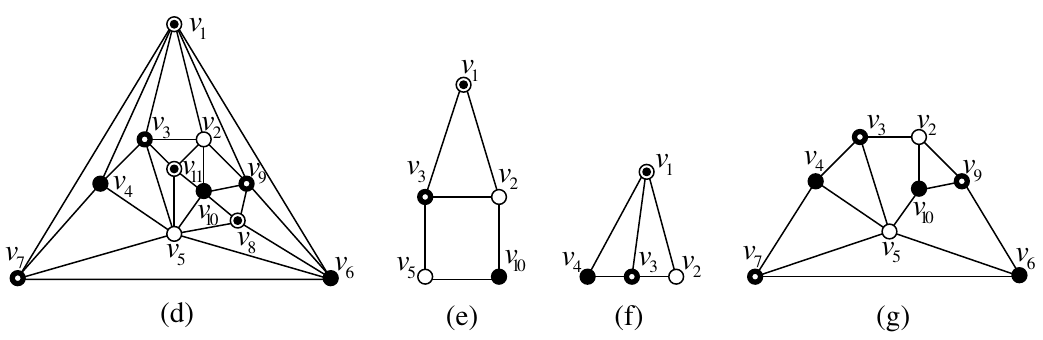}

\textbf{Figure 7.1.} Graphic expression of the definition of cycle, chord-cycle and one-side cycle.
  \end{center}

Here, we make an agreement that $P=v_1v_2\cdots v_m$ is called a \textbf{basic path} of a maximal planar graph $G$ only if the subgraph induced by $\{v_1,v_2,\cdots, v_m\}$ is a path, and  $C=v_1v_2\cdots v_mv_1$ is called a \textbf{basic cycle} of a maximal planar graph $G$ only if the subgraph induced by  $\{v_1,v_2,\cdots, v_m\}$ is a cycle with length $m\geq 4$. Obviously, for a cycle $C$ in a maximal planar graph $G$, $G[V(C)]$ is either a basic cycle or a chord-cycle.

In a maximal planar graph, we are more concerned about the structure, distribution and enumeration of cycles. So, on the basis of the definition defined above, now we are discussing these problems in depth.

\begin{theorem2}\label{th7.1}
Suppose that $G$ is a maximal planar graph with $\delta(G)\geq 4$ and $G$ is not a divisible graph, then the subgraph induced by the set of neighbors of each vertex $v\in V(G)$ is just a cycle with length $d(v)$.
\end{theorem2}

\begin{proof}
Let $\Gamma(v)$ be the neighbor set of $v$. Then there are three possible cases for the induced subgraph $G[\Gamma(v)]$ as follows:

Case 1. $G[\Gamma(v)]$ is a cycle;

Case 2. $G[\Gamma(v)]$ contains triangle;

Case 3. $G[\Gamma(v)]$ contains at least three cycles with length not less than 4.

Case 1 is just the result of this theorem; for the Case 2 and Case 3, it is easy to prove $G$ contains vertices with degree 3 and $G$ is a divisible graph, respectively.
\end{proof}

Denote by $\varsigma(G)$ the set containing all of the cycles with the lengths at least 4 of $G$; $\varsigma^{1}(G)$ the set containing all of the odd-cycles with the lengths at least 5 of $G$; $\varsigma(G)^{2}$ the set containing all of the even-cycles. Where, odd-cycle denotes such a cycle that the number of its length is odd, so does the even-cycle. Obviously,
$$
\varsigma(G)=\varsigma^{1}(G) \cup \varsigma(G)^{2}.
$$

Actually, Theorem 7.1 says in the condition that  $G$ is a maximal planar graph with $\delta(G)\geq 4$ and a nondivisible graph, each vertex $v$ just contributes a cycle with length $d(v)$ for $\varsigma(G)$. Concretely, if $d(v)$ is an odd number, then $v$ just contributes an odd-cycle for $\varsigma^{1}(G)$; if $d(v)$ is an even number, then $v$ contributes an even-cycle for $\varsigma^{2}(G)$.

Naturally, for a maximal planar graph $G$, a correlative problem will be asked about the structure of the subgraph induced by the neighbor set of two adjacent vertices $u,v$ or a connected subgraph with more vertices $v_1,v_2,\cdots, v_m$. What are the structures of $G[\Gamma(u,v)]$  and $G[\Gamma(v_1,v_2,\cdots,v_m)]$? Trees? Cycles (odd-cycle or even-cycle)? Chord-cycles or one-side cycles?

Anyway, the length of cycles should be considered mainly.

\begin{theorem2}\label{th7.2}
Suppose that $G$ is a maximal planar graph and $P=v_1v_2$ $\cdots v_m$ is a basic path of $G$. If the subgraph induced by the neighbor set of $P$ is a cycle, denoted $C=G[\Gamma(v_1,v_2,\cdots, v_m)]$, then the length of $C$ is equal to
$$
|C|=\sum\limits_{i=1}^{m}d(v_i)-4(m-1) \eqno{(7.1)}
$$
\end{theorem2}

\begin{proof}
By induction on $m$. When $m=1$, if the subgraph, denoted $H$, which is induced by the neighbors of some vertex $v$ of $G$ is a cycle, then the length of $H$ is equal to $d(v)$ by Theorem 7.1 and the assertion holds. Suppose that it is true for all paths of $G$ with fewer than $m$ vertices, where $m\geq 2$, and let $P=v_1v_2\cdots v_{m-1}v_m$ be a path of $G$ with order $m$ and the subgraph induced by the neighbor set of $P$ be a cycle, denoted $C$. Choose edge $e=v_{m-1}v_m$ on $P$ and contract $e$ in $G$. Here, we denote the two ends $v_{m-1},v_m$ of $e$ by a new vertex $v$. Then $G/e$ is a maximal planar graph with a path $P^{\prime}=v_1v_2\cdots v_{m-2}v$ and the subgraph in $G/e$, induced by the neighbor set of $P'$ is also the cycle $C$ because the neighbors of $P$ in $G$ and the neighbors of $P^{\prime}$ in $G/e$ are  identical. By the induction hypothesis,
$$
|C|=\sum\limits_{i=1}^{m-2}d_{G/ e}(v_i)+d_{G/ e}(v)-4(m-2) \eqno{(7.2)}
$$
Using the relations
$$
d_{G/e}(v_i)=d_{G}(v_i),i=1,2,\cdots,m-2, d_{G/e}(v)=d_G(v_{m-1})+d_G(v_{m})-4
$$
we obtain
$$
|C|=\sum\limits_{i=1}^{m}d_G(v_i)-4(m-1)
$$
The theorem follows by induction.
\end{proof}
From Formula 7.1, we can see that the parity of cycle $C$'s length only depends on the number of vertices with odd-degree in path $P$.

\begin{corollary}\label{coro7.3}
Suppose that $G$ is a maximal planar graph and $P=v_1v_2\cdots v_m$ is a basic path of $G$. If the subgraph induced by the neighbor set of $P$ is a cycle, denoted $C$, then $C$ is an even-cycle if and only if there are even number of vertices with odd-degree among $v_1,v_2,\cdots, v_m$.
\end{corollary}

More generally, we are seeing about such cycles in $G$ that are induced by the neighbor set of a connected subgraph of $G$.

\begin{theorem2}\label{th7.4}
Suppose that $G$ is a maximal planar graph and $H=G[\{v_1,$ $v_2,\cdots, v_m\}]$ is a connected subgraph of $G$. If the subgraph induced by the neighbor set of $H$ is a cycle, denoted $C=G[\Gamma(v_1,v_2,\cdots, v_m)]$, then
$$
|C|=\sum\limits_{i=1}^{m}d_G(v_i)-\sum\limits_{i=1}^{m}d_H(v_i)-b(H) \eqno{(7.3)}
$$
where  $b(H)$ denotes the number of edges on the boundary of $H$, in which the cut edges being counted twice.
\end{theorem2}

It follows by the similar proof as Theorem 7.2, so omit the detailed process.

\begin{corollary}\label{coro7.5}
Suppose that $G$ is a maximal planar graph and $H=G[\{v_1,v_2,\cdots, v_m\}]$ is a connected subgraph of $G$. If the subgraph induced by the neighbor set of $H$ is a cycle, denoted $C=G[\Gamma(v_1,v_2,\cdots, v_m)]$, then $C$ is an even-cycle if and only if $b(H)$ has the same parity with the number of vertices with odd-degree among $v_1,v_2,\cdots, v_m$.

\end{corollary}

Let $G$ be a maximal planar graph and $H=G[\{v_1,v_2,\cdots, v_m\}]$ be a connected subgraph of $G$. If the subgraph induced by the neighbor set of $H$ is a cycle or a chord-cycle, then the resulting graph $G-C-H$, denoted by $\bar{H}_{C}$, is called the \textbf{complement of $H$ on $C$} in $G$. Namely
$$
G-C-H{\triangleq} \bar{H}_{C} \eqno{(7.4)}
$$

\begin{Prop}\label{pro7.6}
Suppose that $G$ is a maximal planar graph and $H=G[\{v_1,v_2,\cdots, v_m\}]$ is a connected subgraph of $G$. If the subgraph induced by the neighbor set of $H$ is a cycle or a chord-cycle, then when $\bar{H}_{C}$ is connected, the subgraph induced by the neighbors of $\bar{H}_{C}$ is a cycle.

\end{Prop}

On the basis of Proposition 7.6, if the subgraph induced by the neighbors of $\bar{H}_{C}$ is written simply as $C^{\prime}$, then the following theorem is true.

\begin{theorem2}
Let $G$ be a maximal planar graph and $H=G[v_1,v_2,\cdots, v_m]$ be a connected subgraph of $G$. If the subgraph induced by the neighbor set of $H$ is not a tree, denoted $C$, then

I.  If $\bar{H}_{C}$ is unconnected, then $C$ is a chord-cycle;

II. If  $\bar{H}_{C}$ is connected, then either $C^{\prime}=C$ that shows $C$ is a basic cycle of $G$, or $|V(C^{\prime})|<|V(C)|$ and $C$ is a chord-cycle that each chord is contained in a triangle.
\end{theorem2}

\subsection{Enumeration of even-cycles}

In this subsection, we intend to have a try on  finding a necessary and sufficient condition for a 2-colorable cycle. Obviously, it is very significant to research how many even-cycles are contained in a maximal planar graph.

Suppose that $G$ is a maximal planar graph with order $n$ and $\delta(G)\geq 4$. Denote by $\pi(G)=(d_1,d_2,\cdots,d_n)$ the degree sequence of $G$.  Now, considering the number of cycles in $G$. If $G$ is indivisible, then according to Theorem 7.1, $G$ contains $n$ basic cycles that are induced only by the neighbors of a vertex of $G$, and we denote the set that contains  all of such cycles by $C_1$, obviously $|C_1|=n$. Let $m$ be the number of odd-degree vertices of $G$, then $C_1$ contains $m$ odd-cycles; naturally, the number of even-cycles in $C_1$ is $n-m$. Especially, when $m=n$, there are no even-cycles in $C_1$.  However, the longest cycle in the subgraph of $G$ induced by any two adjacent vertices is either a basic cycle or a chord-cycle,
so there are $3n-6$ such basic cycles and chord-cycles, in which the length of both basic cycles and chord-cycles are even by Theorem 7.2, denoted by $C_2$ and $|C_2|=3n-6$.

In a maximal planar graph $G$, the induced subgraph by the vertices on the boundary of each semi-maximal planar graph, $H$, of $G$ is either a basic cycle or a chord-cycle. We say $H$ corresponds to a basic cycle or chord-cycle. Conversely, we also say a basic cycle or a chord-cycle corresponds to a semi-maximal planar graph.
In addition, it is easy to see that each basic cycle or chord-cycle in $G$ just corresponds to two semi-maximal planar graph of $G$. Let $Ha(G)$ be the set of all the semi-maximal planar graphs of $G$, $Ch(G)$ the set of all the chord-cycles of $G$, $Cy(G)$ the set of all the basic cycles with length not less than 4 of $G$. Then it follows

\begin{theorem2} Suppose that $G$ is a maximal planar graph, then
$$
|Cy(G)|\leq\frac{1}{2}|Ha(G)|-|Ch(G)| \eqno{(7.5)}
$$
\end{theorem2}

\textbf{Remark:} Conducting research in the number of cycles in a maximal planar graph is not the focus of this article, so we would like to discuss it at length by another paper.

\subsection{Black-White coloring operation}

Suppose that $G^{C}$ is a semi-maximal planar graph on even-cycle $C=v_1v_2\cdots v_{2m}$\\$v_1(m\geq 2)$. Denote by $\Gamma^{*}(C)$ the vertex-set consisting of such vertices of $V(G^{C})-V(C)$ that adjacent to both the vertices of $C$ with odd-subscript and even-subscript. The so-called Black-White coloring for $G^{C}$, denoted $f_{bw}:V(G)\rightarrow \{b,w\}$, is to divide $V(G^{C})$ into two subsets, $B$ and $W$ that are called black vertex-set and white vertex-set, respectively, $V(G^{C})=B \cup W, B,W\neq \emptyset$. Where all of the vertices in $B$ and $W$ are colored by black and white, respectively. Further, for a Black-White coloring $f_{bw}=(B,W)$ of $G^{C}$, if both $G^{C}[B]$ and $G^{C}[W]$ contains no odd-cycles, then we refer to $f_{bw}=(B,W)$ as a \textbf{proper Black-White coloring}; otherwise, an \textbf{improper Black-White coloring}.

Now, for a semi-maximal planar graph $G^{C}$, we introduce a kind of operation on $C$, called Black-White coloring operation, which is closely related to 4-colorings. Following gives the detailed steps.

Step 1. Color vertices of $C$ by black;

Step 2. Color vertices of $\Gamma^{*}(C)$ by white;

Step 3. Let $\Gamma^{*}(\Gamma^{*}(C))\triangleq \Gamma^{2*}(C)$,  then color  vertices of $\Gamma^{2*}(C)$ by black;

$\cdots\cdots\cdots\cdots\cdots\cdots\cdots\cdots\cdots\cdots\cdots\cdots\cdots\cdots\cdots\cdots\cdots\cdots\cdots\cdots$

Step $2i$. Color vertices of $\Gamma^{(2i-1)*}(C)$ by white;

Step $2i+1$. Color vertices of $\Gamma^{(2i)*}(C)$ by black;

$\cdots\cdots\cdots\cdots\cdots\cdots\cdots\cdots\cdots\cdots\cdots\cdots\cdots\cdots\cdots\cdots\cdots\cdots\cdots\cdots$

Until

Step $t+1.$  $\Gamma^{(t)*}(C)=\emptyset$;

Step $t+2.$  Color all the vertices in $G^{C}$ that are not colored by black or white by grey, and $A$ denotes the set containing all of the vertices colored by grey;

Step $t+3.$ For $\forall v\in V(A)$, when $v$ is colored only by black(or white), there will be odd-cycles in $G[B]$(or $G[W]$), then we call $v$ the \textbf{fixed-vertex} and color it by white(or black);

Step $t+4.$ For $\forall u\in V(A)$, no matter which color(black or white) is colored to it, there always odd-cycles in $G[B]$(or $G[W]$), then we call vertex $u$ the \textbf{petal-vertex} and color it by black or white optionally. Then, if there are fixed-vertex in $A$, go back to Step $t+3$; otherwise, stop.

Here, if $A\neq \emptyset$ when the operation stops, then we also use $f_{bw}=(B,W,A)$ to denote the Black-White coloring operation. According to this operation, an obvious result can be obtained as follow.
\begin{theorem2}\label{th7.9}
Let $G^{C}$ be a semi-maximal planar graph on even-cycle $C$, and $f_{bw}=(B,W,A)$ be the Black-White coloring operation of $G^{C}$ on $C$, then $G^{C}[B]$ or $G^{C}[W]$ contains odd-cycles if and only if there appears petal-vertices at Step $t+4$ in the process of $f_{bw}=(B,W,A)$.

\end{theorem2}
Obviously, the petal-vertex is adjacent to at least two vertices of $B$ and $W$, respectively. Figure 7.2(a) and 7.2(b) show a structural characteristic partially of the petal-vertex. And the vertex $u$, shown in Figure 7.2(c), is a petal-vertex.

 \begin{center}
       \includegraphics [width=240pt]{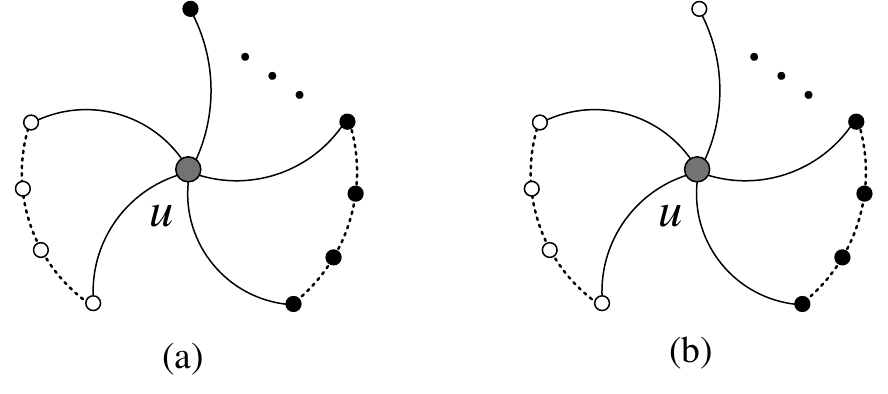}

        \includegraphics [width=320pt]{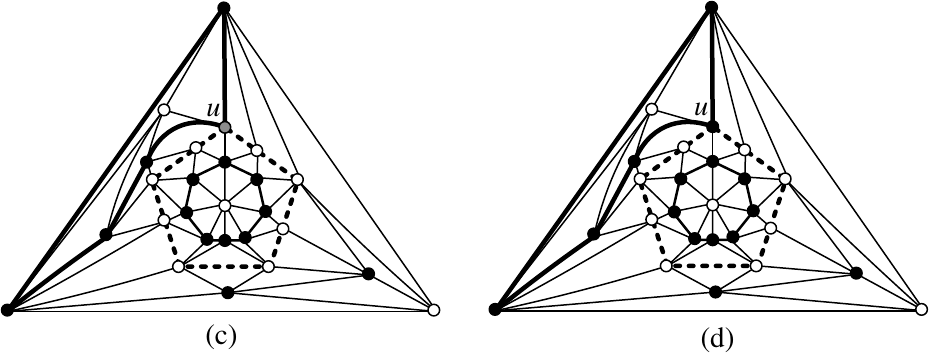}

  \textbf{Figure 7.2.} Schematic diagram of petal-vertices
  \end{center}

 For a maximal planar graph $G$ and a cycle $C$ of $G$, the Black-White coloring operation on $C$ for $G$ is to conduct Black-White coloring operation on $C$ for the two semi-maximal planar graphs on $C$ of $G$. For example, Figure 7.3(a) is the resulted graph after conducting Black-White coloring operation on cycle $v_1v_3v_5v_6v_1$ for the maximal planar graph shown in Figure 7.1; Figure 7.3(b) is the resulted graph after conducting Black-White coloring operation on cycle $v_1v_2v_3v_4v_1$ for the maximal planar graph shown in Figure 6.1; Figure 7.3(c) is the resulted graph after conducting Black-White coloring operation on cycle $v_1v_2v_3v_4v_5v_6v_1$ for icosahedron. In addition, for the semi-maximal planar graph shown in Figure 7.3(d), when we conduct Black-White coloring operation on the inside cycle $C$, $\Gamma^{3*}(C)=\emptyset$.

 \begin{center}
       \includegraphics [width=340pt]{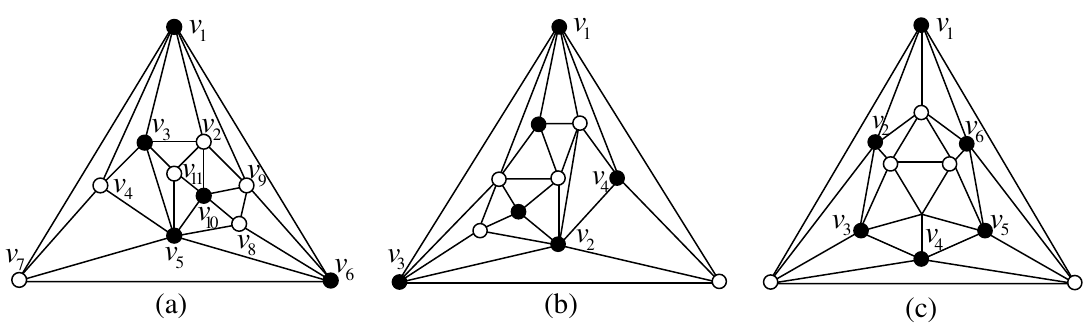}

        \includegraphics [width=200pt]{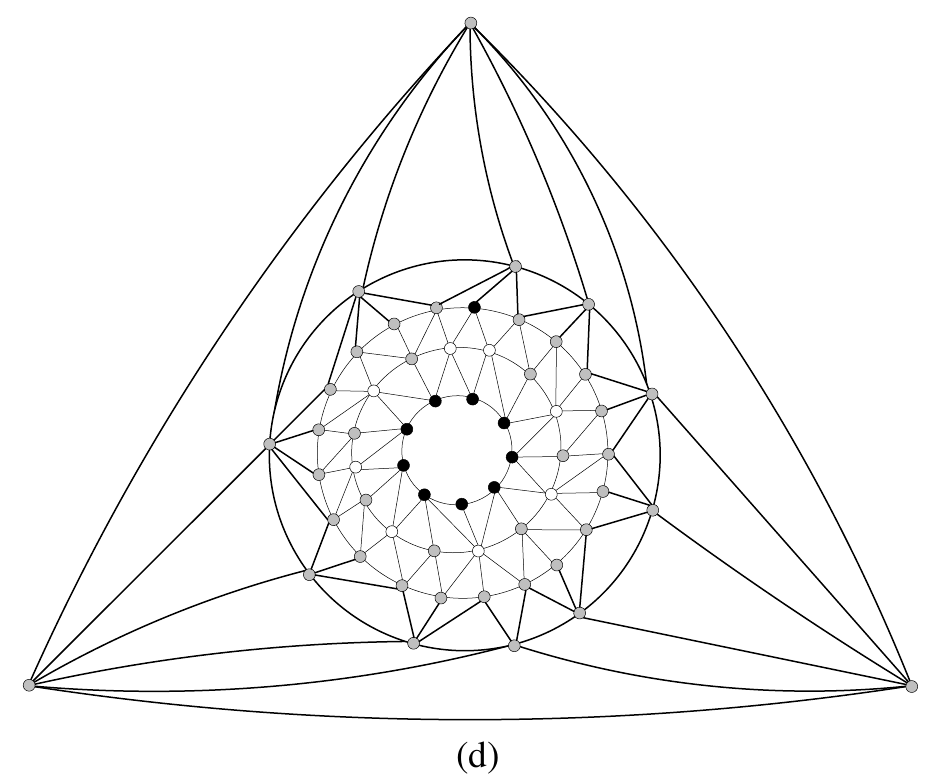}

       \textbf{Figure 7.3.} An example of Black-White coloring operation
  \end{center}

 For some maximal planar graphs, we can know from the above examples that after we conduct Black-White coloring operations on cycles for them, some may have vertices colored by grey and others may not have such vertices. We refer to Black-White colorings  for the former as \textbf{unique Black-White colorings} on cycles, and for the latter as \textbf{non-unique Black-White colorings}. For example, the coloring in Figure 7.3(a),7.3(b) and 7.3(c) are unique Black-White colorings, but in Figure 7.3(d) is a non-unique Black-White coloring.

 \subsection{The necessary and sufficient condition of the 2-colorable cycle based on petal-syndrome}

\quad\quad It is easy to prove the following result:
\begin{theorem2}\label{th7.10} Let $G$ be a maximal planar graph, $C$ be an even-cycle of $G$. Suppose $f_{bw}$ is a Black-White coloring on $C$ of $G$. If $f_{bw}$ is unique, then $C$ is 2-colorable if and only if $f_{bw}$ is proper.

  \end{theorem2}

According to the Theorem 7.10, we can clearly judge that for the maximal planar graph $G$ shown in Figure 7.1(d), cycle $C=v_1v_3v_5v_6v_1$ is 2-colorable because the Black-White coloring on $C$ of $G$ is unique, and $G[B]$ and $G[W]$ are forest and 1-fence, respectively(see Figure 7.3(a)); similarly, for the maximal planar graph $G$ shown in Figure 7.1(b), cycle $C=v_1v_2v_3v_4v_1$ is also a 2-colorable cycle ; however, for icosahedron, the 6-cycles induced by the neighbors of any two adjacent vertices
is not 2-colorable because $G[W]$ contains two triangles(see Figure 7.3(c)).

Let $G^{C}$ be a semi-maximal planar graph on even-cycle $C$. We conduct the Black-White coloring operation $f_{bw}$ on $C$, which partitions $V(G^{C})$ into three subsets: black vertex-set $B$, white vertex-set $W$ and prey vertex-set $A$. Denoted by $$f_{bw}=(B,W,A) \eqno{(7.6)}$$

Let $f_{bw}=(B,W,A)$ be a Black-White coloring on $C$ of $G^C$, in which $A \neq \emptyset$. If $C$ is an even-cycle and both $G^C[B]$ and $G^C[W]$ contain no odd-cycles, then we recolor any vertex $v \in A$ by black or white, and remain the colors of other vertices in $A$ unchanged.  Denote by $f^\prime_{bw}=(B^\prime,W^\prime,A^\prime)$ the new Black-White coloring. Obviously, both $G^C[B^\prime]$ and $G^C[W^\prime]$ still contain no odd-cycles. And then, we call the vertices of $A$ the \textbf{free vertices}.

For a Black-White coloring $f_{bw}=(B,W,A)$ on $C$ of $G^C$. Suppose that $|A|\geq 2$ and $u,v \in A$, if $G^C[B\cup \{u,v\}]$ or $G^C[W\cup \{u,v\}]$ contains odd-cycles including $u$ and $v$ when they are recolored by black(white), then $\{u,v\}$ are called the \textbf{petal-pair}.  In addition, let $S\subseteq A$, if any pair of vertices $u$ and $v$ in $S$ is a petal-pair, then $S$ is called the \textbf{petal-set}. For example, both the set $S=\{u_1,u_2,u_3,u_4\}$ in Figure 7.4(a) and the set $S=\{u_1,u_2,u_3\}$ in Figure 7.4(b) are petal-sets.

We can easily obtain the following theorem by the fact that any planar graph contains no $K_5$ and its subdivision.

\begin{theorem2}\label{th7.11} Let $G^C$ be a semi-maximal planar graph on the even-cycle $C$, and $f_{bw}=(B,W,A)$ be a Black-White coloring on $C$ of $G^C$. If $A\neq \emptyset$ and $S\subseteq A$ is a petal-set of $G^C$, then
$$|S|\leq 4 \eqno{(7.7)}$$
\end{theorem2}

An edge $uv$ is called \textbf{petal-edge} if $\{u,v\}$ is a petal-pair. In Figure 7.4(b), edges $u_1u_2$, $u_2u_3$, and $u_1u_3$ are petal-edges. A path $P$ in $G^C[A]$ is called the \textbf{special petal-path} if each edge of $P$ is a petal-edge; analogously, a cycle $C^*$ in $G^C[A]$ is called the \textbf{special petal-cycle} if each edge of $C^*$ is a petal-edge. The subgraph induced by a sequence of $p$ vertices $x_1,x_2,\cdots,x_p$ in $G^C[A]$ is called the \textbf{general petal-path} if $\{x_i,x_{i+1}\}$ is a petal-vertex pair, $i=1,2,\cdots,p-1$; further, if $\{x_1,x_{p}\}$ is a petal-pair, then the subgraph is called the \textbf{general petal-cycle}. Actually, a special petal-path is a special case of the general petal-path; similarly, a special petal-cycle is a special case of the general petal-cycle. So the general petal-path and general petal-cycle are called straightly the \textbf{petal-path} and \textbf{petal-cycle}, respectively. For example, the cycle $u_1,u_2,u_3,u_1$ in Figure 7.4(b) is a petal-cycle(special).

\begin{center}
         \includegraphics [width=360pt]{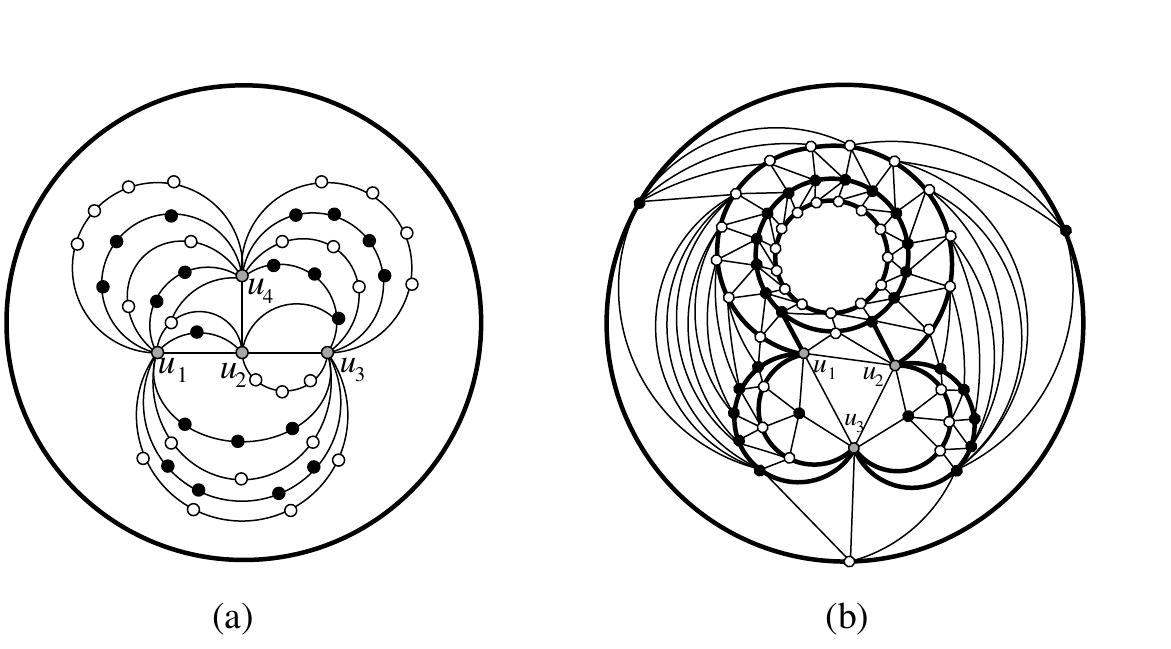}

        \textbf{Figure 7.4.} The illustrations of the petal-vertex set, petal-edge and petal-cycle
\end{center}

On the basis of the above arguments, we define a new graph---petal-graph. For a semi-maximal planar graph $G^C$, we conduct the Black-White coloring operation on $C$, denoted $f_{bw}=(B,W,A)$, $A\neq \emptyset$. The so-called \textbf{petal-graph $G_S$ on coloring $f_{bw}$} is that: its vertex set $\{x_1,x_2,\cdots,x_n\}\subseteq A$ and each vertex at least belongs to one petal-pair, in which the vertices $x_i$ and $x_j$ are adjacent if and only if $\{x_i,x_j\}$ is a petal-pair, where $i,j=1,2,\cdots,n$, $i\neq j$. Obviously, a petal-graph is planar and contains no isolated vertices. For example,if $G^{C}[V(G_S)]$ is  a petal-path, then $G_S$ is a path; if $G^{C}[V(G_S)]$ is a petal-cycle, then $G_S$ is a cycle; and if $V(G_S)$ is a petal-set, then $G_S$ is a complete graph.

Let $G^C$ be a semi-maximal planar graph, $f_{bw}=(B,W,A)$ a Black-White coloring on $C$, $A\neq \emptyset$ and $G_S$ the petal-graph on $f_{bw}$. If $G_S$ contains no odd-cycles, namely a bipartite graph with two independent sets $X$ and $Y$, then $G_S$ is called the \textbf{exclusive petal-graph on coloring $f_{bw}$} if the following two conditions are satisfied:

(1) Both $G^C[B\cup X]$ and $G^C[B\cup Y]$, or both $G^C[W\cup X]$ and $G^C[W\cup X]$ contain odd-cycles;

(2) Both $G^C[B\cup X]$ and $G^C[W\cup X]$, or both $G^C[B\cup Y]$ and $G^C[W\cup X]$ contain odd-cycles.\\
As shown in Figure 7.5(a),(b), if there exist cycles $C_1$ and $C_2$, or $C_3$ and $C_4$ in $G^C$, then $G^C[X\cup Y]$ is a  exclusive petal-graph.

When there exists a path $P=u_1u_2\cdots u_l$ in $G^C[A]$ such that both $G^C[B\cup \{u_{i}u_{i+1}\}]$ and $G^C[W\cup \{u_{i+1}u_{i+2}\}]$ contain odd-cycles, or both $G^C[W\cup \{u_{i}u_{i+1}\}]$ and $G^C[B\cup \{u_{i+1}u_{i+2}\}]$ contain odd-cycles, $i=1,2,\cdots,l-2$, then $P$ is called the \textbf{Black-White path}. For the sake of convenience, we always assume that $G^C[B\cup \{u_{i}u_{i+1}\}]$ and $G^C[W\cup \{u_{i+1}u_{i+2}\}]$ contain odd-cycles. Obviously, if a vertex $u$(not the ends) of path $P$ is recolored by black or white, then there may appear several fixed-vertices in $A$. Let $D^w_b$ and  $D^w_w$(contains $u$) be the sets of vertices in $A$ recolored by black and white respectively when $u$ is recolored by white, and let $D^b_b$(contains $u$) and  $D^b_w$ be the sets of vertices in $A$ recolored by black and white respectively when $u$ is recolored by black. A Black-White path $P$ is called the \textbf{exclusive Black-White path} if one of the following conditions is satisfied:

 (1) $G^C[ D^b_b ]$ or $G^C[ D^b_w]$, and $G^C[ D^w_b]$ or $G^C[ D^w_w]$ contain simultaneously odd-cycles;

 (2) $G^C[ D^b_b]$ or $G^C[ D^b_w]$ contains two odd edge-disjoint cycles;

 (3) $G^C[ D^w_b]$ or $G^C[D^w_w]$ contains two odd edge-disjoint cycles.\\
 It is easy to see that if $G^C[A]$ contains a exclusive Black-White path $P$, then there must exist a unicolor odd-cycle  whatever colors(black or white) are assigned to the vertices of $P$. In Figure 7.5(c), when the vertex $u$ is recolored by black,  $D^b_b=\{u,u_2,u_4\}$, $D^b_w=\{u_1,u_3,u_5\}$ and $G^C[D^b_w]$ contains an odd-cycle; when the vertex $u$ is recolored by white, $D^w_b=\{u_6,u_8,u_{10}\}$, $D^w_w=\{u,u_7,u_9,u_{11}\}$, and $G^C[D^w_w]$ contains an odd-cycles, shown in Figure 7.5(d).

\begin{center}
             \includegraphics [width=260pt]{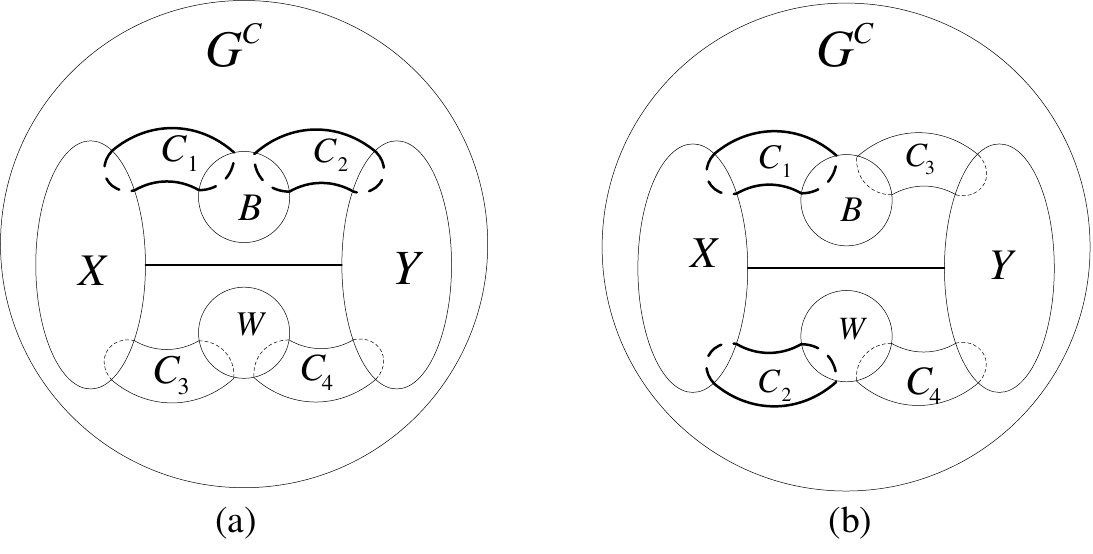}

                            \hspace{1cm}
             \includegraphics [width=320pt]{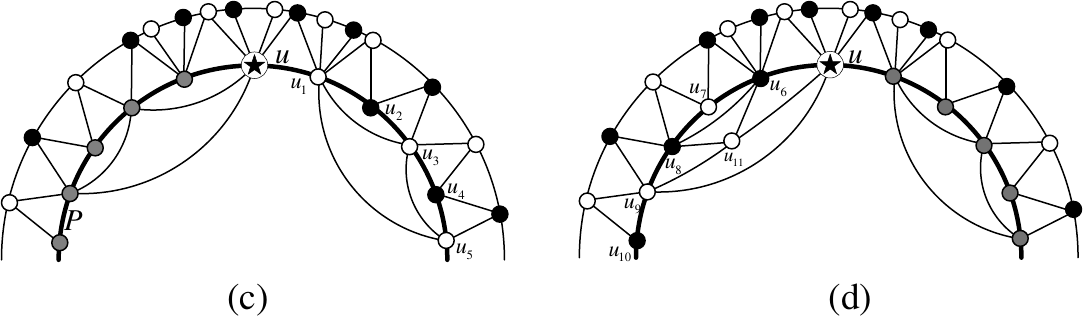}

        \textbf{Figure 7.5.} The illustrations of the exclusive petal-graph and the exclusive Black-White path
\end{center}

For a vertex $u$ of $G^C[A]$, when we recolored $u$ by black or white, some fixed-vertices will appear in $\Gamma(u)$, denoted by $D_1$. Then color them and another set of fixed-vertices in $\Gamma(D_1)$, $D_2$, will appear correspondingly. Further, when we recolored $D_2$, $D_3$ will appear correspondingly. In this way, continue this action until some $D_k$, for which no fixed-vertices are appeared in $\Gamma(D_k)$. Then we call $D=D_1\cup D_2 \cup \cdots \cup D_k$ the \textbf{fixed-set} of $u$.

Supposed that $u$ is a vertex of $G^C[A]$. When $u$ is recolored by white, we denote by $D^w_b$ and $D^w_w$(contains $u$) the vertex-sets fixed by $u$, in which vertices are recolored by black and white, respectively.
 When $u$ is recolored by black,  we denote by  $D^b_b$(contains $u$) and  $D^b_w$ the vertex-sets fixed by $u$, in which vertices are recolored by black and white, respectively.
 When $u$ is recolored by white,
 let $f^\prime_{bw}=(B^\prime,W^\prime,A^\prime)$ be the extended coloring of $f_{bw}=(B,W,A)$ in the condition that $u$ is recolored by white, in which $B^\prime=B\cup D^w_b$, $W^\prime=W\cup D^w_w$ and $A^\prime=A-(D^w_b\cup D^w_w)$.  When $u$ is recolored by black, Let $f^{\prime\prime}_{bw}=(B^{\prime\prime},W^{\prime\prime},A^{\prime\prime})$ be the extended coloring of $f_{bw}=(B,W,A)$ in the condition that $u$ is recolored by black, in which $B^{\prime\prime}=B\cup D^b_b$, $W^{\prime\prime}=W\cup D^b_w$ and $A^{\prime\prime}=A-(D^b_b\cup D^b_w)$. Then, $u$ is called the \textbf{general petal-vertex} if one of the following four conditions is satisfied:

(1) $G^C[B^\prime]$ or $G^C[W^\prime]$ contains odd-cycles, and  $G^C[B^{\prime\prime}]$ or $G^C[W^{\prime\prime}]$ contains odd-cycles;

(2) Both the petal-graphs on $f^\prime_{bw}$ and $f^{\prime\prime}_{bw}$ of $G^C$ contain odd-cycles;

(3) At least one of $G^C[B^\prime]$ and $G^C[W^\prime]$ contains odd-cycles, and the petal-graph on $f^{\prime\prime}_{bw}$ of $G^C$ contains cycles;

(4) At least one of $G^C[B^{\prime\prime}]$ and $G^C[W^{\prime\prime}]$ contains odd-cycles, and the petal-graph on $f^{\prime}_{bw}$ of $G^C$ contains cycles.\\
For example, each vertex of a exclusive petal-graph is a general petal-vertex; in Figure 7.5(a) and (b), the vertex $u$ on the exclusive Black-White path is a general petal-vertex; if a petal-graph contains an odd-cycle, then all of the vertices on the odd-cycle are general petal-vertices. Based on the above arguments, we refer to the phenomenon that $G^C[A]$ contains general petal-vertices as the \textbf{petal-syndrome}.

\begin{theorem2}\label{th7.12} Let $G^C$ be a semi-maximal planar graph on the even-cycle $C$, and $f_{bw}=(B,W,A)$ be a Black-White coloring on $C$ of $G^C$. If $A\neq \emptyset$ and $G^C[A]$ contains petal-syndrome, then
$C$ is not a 2-colorable cycle.
  \end{theorem2}

 \begin{lemma}\label{th1}
 Let $G^C$ be a semi-maximal planar graph on the even-cycle $C$, and $f_{bw}=(B,W,A)$ be a Black-White coloring on $C$ of $G^C$. If $A\neq \emptyset$ and on the promise that no unicolor odd-cycles appear, then for any two adjacent vertices $u,v\in A$:

 (1) $u,v$ are recolored only by two different colors if and only if $uv$ is a petal-edge;

 (2) $u,v$ are recolored only by the same color if and only if $u,v$ are on a petal-path and the distance between them is even.
  \end{lemma}
\begin{proof}
(1) According the definition of the petal-edge, the sufficient condition holds. Conversely, suppose that $u,v$ are recolored only by two different colors, it shows that both $G^C[B\cup \{u,v\}]$ and $G^C[W\cup \{u,v\}]$ contain  odd-cycles. Therefore, $uv$ is a petal-edge.

(2) Suppose that $u,v$ are on a petal-path and the distance of them on the path is even, then the colors of $u$ and $v$ must be recolored by the some color. Conversely, suppose that $u,v$ are recolored only by the some color, namely when one of them, says $u$, is recolored by black(white), then $v$ must be in the vertex-set $D$ fixed by $u$, in which vertices are recolored by black(white). So there exists a petal-path between $u$ and $v$ and the distance of them on the path is even.
\end{proof}

For a semi-maximal planar graph $G^C$ on the even-cycle $C$, in order to deal with the possible problem of the petal-syndrome after conducting  Black-White coloring operation on $C$ for $G^C$, now we give an \textbf{improved operation of Black-White coloring on $C$} of $G^C$ as follow.

Step 1. Color vertices of $C$ by black;

Step 2. Color vertices of $\Gamma^{*}(C)$ by white;

Step 3. Let $\Gamma^{*}(\Gamma^{*}(C))\triangleq \Gamma^{2*}(C)$,  then color  vertices of $\Gamma^{2*}(C)$ by black;

$\cdots\cdots\cdots\cdots\cdots\cdots\cdots\cdots\cdots\cdots\cdots\cdots\cdots\cdots\cdots\cdots\cdots\cdots\cdots\cdots$

Step $2i$. Color vertices of $\Gamma^{(2i-1)*}(C)$ by white;

Step $2i+1$. Color vertices of $\Gamma^{(2i)*}(C)$ by black;

$\cdots\cdots\cdots\cdots\cdots\cdots\cdots\cdots\cdots\cdots\cdots\cdots\cdots\cdots\cdots\cdots\cdots\cdots\cdots\cdots$

Until

Step $t+1$.  $\Gamma^{(t)*}(C)=\emptyset$;

Step $t+2$.  If the subgraph, induced by the set of vertices colored by black(or white), contains odd-cycles, then recolor any vertex of the odd-cycles grey and stop. Otherwise, if all vertices of $G^C$ are colored by black or white, stop; else if there are vertices in $G^{C}$ that are not colored by black or white, then color these vertices grey;

For any prey vertex $u$ in $G^{C}$ that is not a general petal-vertex,
if only one color, say white, is assigned to $u$,
the subgraph, induced by the set of vertices colored by black or white, contains odd-cycles, after recoloring the vertices in fixed-set of $u$ by black or white properly, then $u$ is called the \textbf{restricted-vertex}. Of course, when we recolor $u$ black, no unicolor odd-cycles appear.

Step $t+3$. If there are restricted-vertices in the set of prey vertices, then properly color them and the vertices of their fixed-sets correspondingly;

Step $t+4$. If there exist petal-vertices or petal-syndrome, or the subgraph, induced by the set of all black vertices or all white vertices, contains odd-cycles, then stop;

Step $t+5$. If there are restricted-vertices in the set of prey vertices, then go back to Step $t+3$; otherwise, go to next step;

Step $t+6$. Choose a prey vertex $v$, called the \textbf{sign-vertex}, which has the most neighbors colored by black or white, and color $v$ black \emph{and remark it in sequence}. If there is no prey vertex, stop; otherwise, go to the next step;

Step $t+7$. If there are restricted-vertices in the set of prey vertices, then properly color them and the vertices of their fixed-sets correspondingly;

Step $t+8$. If there is no prey vertex, stop. Otherwise, if there exist petal-vertices or petal-syndrome, or the subgraph, induced by the set of all black vertices or all white vertices, contains odd-cycles, when there are black sign-vertices, choosing the latest black sign-vertex, denoted $w$, then we assign prey to
the vertices, which were colored by black or white after $w$ was colored by black. At the same time, recolor $w$ by white and go back to Step $t+7$; when there are no black sign-vertices, stop. Else if there are not petal-vertices or petal-syndrome, or the subgraph, induced by the set of all black vertices or all white vertices, contains no odd-cycles, go to the next step;

Step $t+9$. If there are restricted-vertices in the set of prey vertices, then go back to Step $t+7$; otherwise, go back to Step $t+6$.

Next, we will give an example to illustrate the process of an improved operation of Black-White coloring. Let $G$ be a maximal planar graph and $G^C_1, G^C_2$ are two semi-maximal planar graphs on the cycle $C$, shown in Figure 7.6(a). Now, we are conducting the improved operation of Black-White coloring on $C$ of $G^C_1$.

First, color all vertices of $C$ by black, and conduct the operation until Step 3, then $\Gamma^{2*}(C)=\emptyset$. Because not all vertices in $G^C_1$ are colored by black or white and the subgraph, induced by the set of the vertices colored by black(or white), contains no odd-cycles, we color other vertices by prey, shown in Figure 7.6(b).

Second, conduct Step $t+3$($t=2$) of the operation. Since neither prey vertex is a restricted-vertex, so conduct Step $t+4$. Because there are no petal-vertices, petal-syndrome, sign-vertices, and the subgraph, induced by the set of all black vertices or all white vertices, contains no odd-cycle, conduct Step $t+6$ directly. Choose a sign-vertex $v_1$ and color it by black; Then conduct Step $t+7$ because $G^C_1$ still contains prey vertices. Here, $v_{11}$ is a restricted-vertex and $\{v_{12}\}$ is the fixed-set of $v_{11}$. Color vertices $v_{11}$ and  $v_{12}$ by the corresponding colors, shown in 7.6(b').

Third, because there are no petal-vertices, petal-syndrome, sign-vertices, and the subgraph, induced by the set of all black vertices or all white vertices, contains no odd-cycles, go back to Step $t+6$. Choose a sign-vertex $v_2$, and color it by black. Here, $v_{21}$ is a restricted-vertex, and $v_{22}$ is the fixed-set of  $v_{21}$. Color $v_{21}$ and  $v_{22}$ by corresponding colors. According to the operation and the above arguments, we choose a sign-vertex $v_3$, then $v_{31}, v_{33}$ are two restricted-vertices and $\{v_{32},v_{34},v_{35},v_{36},v_{37},v_{38},v_{39},v_{310}\}$ is the fixed-set of them. We color these vertices properly, then the operation stops for there is no prey vertex in $G^C_1$(see Figure 7.6(b")). So $C$ is 2-colorable in $G^C_1$ for $G^C_1$ contains no unicolor odd-cycle colored.

Analogously, we conduct the improved operation of Black-White coloring on $C$ of $G^C_2$, and the resulting coloring is shown in Figure 7.6(c) and (c'). On the process of the coloring, $v_1,v_2,v_3,v_4$ and $v_5$ are five sign-vertices, and  $\{v_{11},v_{12}\}$, $\{v_{21},v_{22},\cdots,v_{210}\}$, $\{v_{31},v_{32},v_{33}\}$, $\{v_{41},v_{42}\}$, and $\{v_{51},v_{52}\}$ are the fixed-set of them. We can see that $G^C_2$ contains no prey vertex, no unicolor odd-cycle colored  when the operation stops, so $C$ is 2-colorable in $G^C_2$.

Hence, $C$ is a 2-colorable cycle of $G$ and the proper Black-White coloring on $C$ of $G$ is shown in Figure 7.6(b") and (c').

\begin{center}
            \includegraphics [width=240pt]{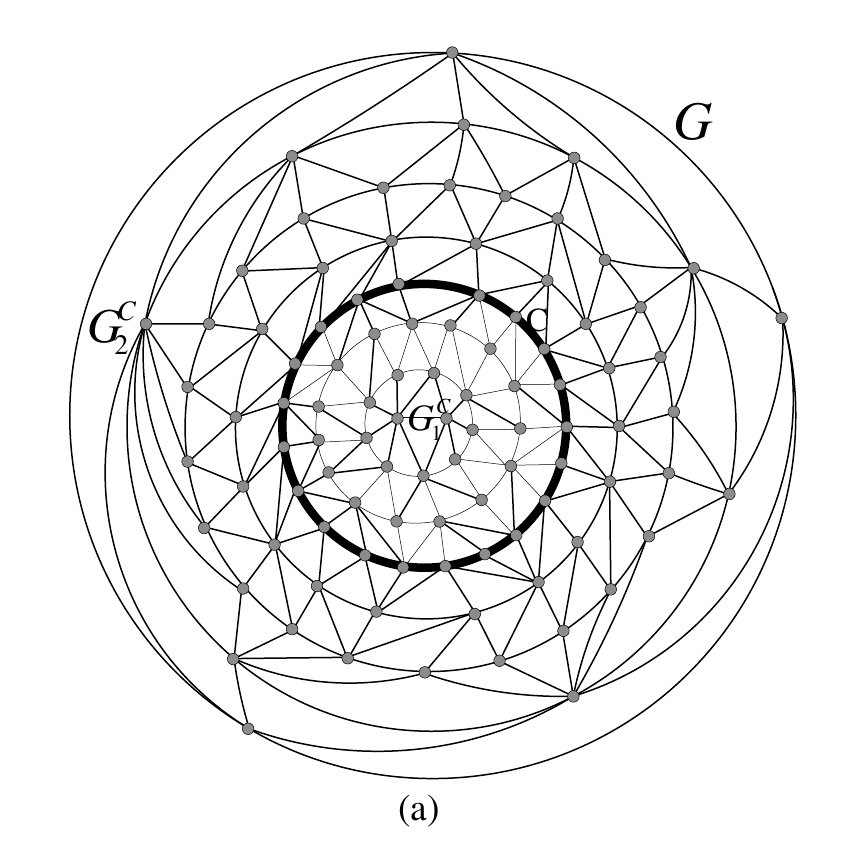}

             \includegraphics [width=380pt]{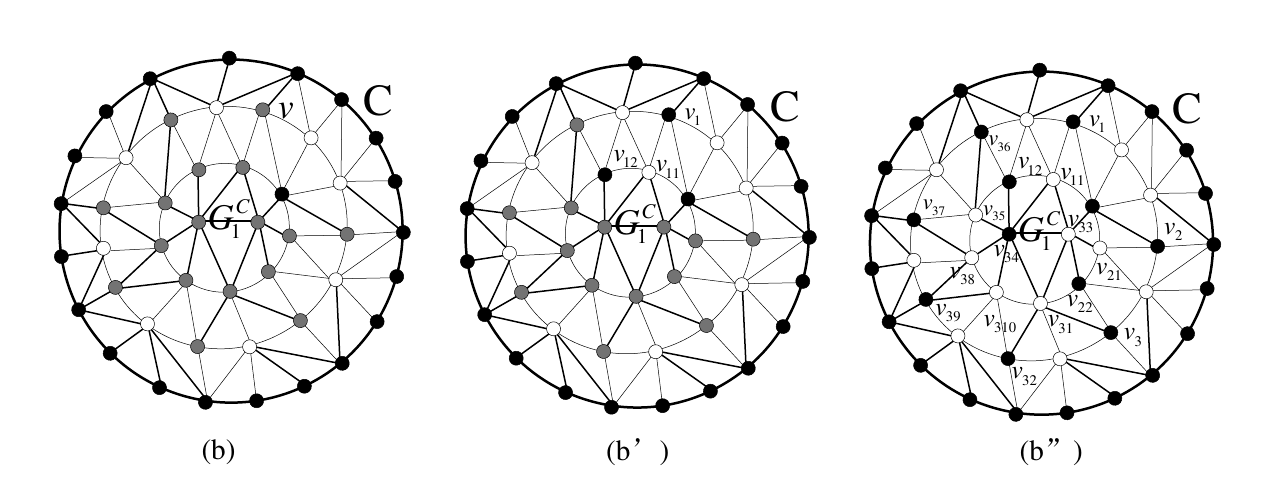}

             \includegraphics [width=360pt]{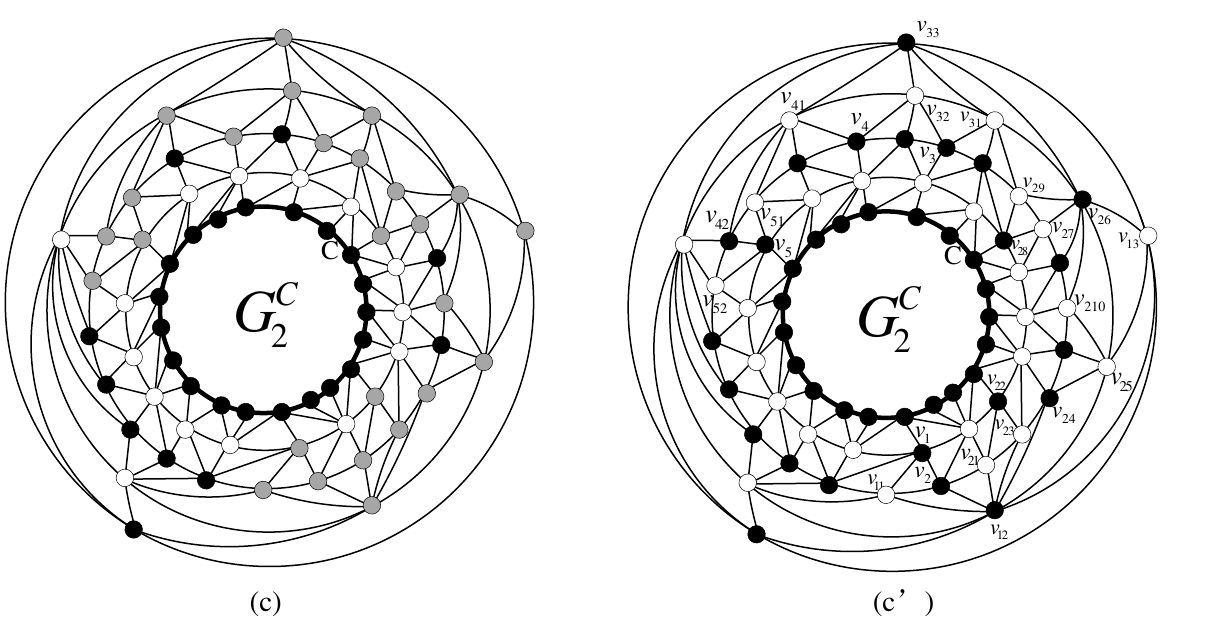}

        \textbf{Figure 7.6.} A maximal planar graph $G$ and two semi-maximal planar graphs $G^C_1, G^C_2$ on $C$
\end{center}

From the improved operation of Black-White coloring, we can obtain the following result.

\begin{theorem2}\label{th7.14} Let $G$ be a maximal planar graph with $\delta(G)\geq 4$, $C$ be an even-cycle of $G$, and $C$ splits $G$ into two semi-maximal planar graphs $G^C_1$, $G^C_2$.  Then $C$ is 2-colorable if and only if neither $G^C_1$ nor $G^C_2$ contains prey vertices after conducting the improved operation of Black-White coloring on $C$ for $G^C_1$ and $G^C_2$, respectively.
  \end{theorem2}
  \begin{proof}
  Suppose that $C$ is 2-colorable in $G$. If we conduct the improved operation of Black-White coloring on $C$ for $G^C_1$ and $G^C_2$ respectively, at least one of $G^C_1$ and $G^C_2$, say $G^C_1$, contains prey vertices when the operations stop, then there appears petal-syndrome in the process of conducting the improved operation of Black-White coloring on $C$ for $G^C_1$. That is to say, no matter how to color the vertices of $G^C_1$, there always exist odd-cycles colored by the same color. Contradict the assumption $C$ is 2-colorable.

  Conversely, suppose that both $G^C_1$ and $G^C_2$ contain no prey vertices after conducting the improved operation of Black-White coloring on $C$ for $G^C_1$ and $G^C_2$, respectively, then there appears no odd-cycle in the process of conducting the improved operation of Black-White coloring on $C$ for $G^C_1$ or $G^C_2$. Therefore, $C$ is 2-colorable both in $G^C_1$ and $G^C_2$. So $C$ is 2-colorable.
  \end{proof}

\subsection{Necessary and sufficient conditions of 2-colorable cycles based on structure}

On the basis of petal-syndrome, subsection 7.4 has given a necessary and sufficient condition of 2-colorable cycles. Sometimes, finding or judging petal-syndrome is a very tough task, so this subsection will study the characteristics of 2-colorable cycles on the structure of a maximal planar graph.

Suppose that $G^{C}$ is a semi-maximal planar graph on even-cycle $C$, and $V'=V(G^{C}-V(C)$. Define $G^{C}[V']\triangleq G'$ as the subgraphs of $G^{C}$, induced by the inner vertices of $C$.

\begin{Prop} Suppose that $G^{C}$ is a semi-maximal planar graph on even-cycle $C$, then $C$ is 2-colorable if and only if $G^{C}$ can be partitioned into two bipartite subgraphs, and $C$ is included in one of them.
\end{Prop}

Let $G^{C}$ be a semi-maximal planar graph on even-cycle $C$. Suppose that $C$ is 2-colorable, and $G_1,G_2$ are the two bipartite subgraphs of $G^{C}$. Without loss of generality, we can assume $C\in G_1$. Thus, we can discuss the structure of $G^C$ through dividing even-cycle $C$ into two categories according to the relationship between $\Gamma(C)$ and $\Gamma^{*}(C)$.

\textbf{Type 1. Closed-cycles}

If  $\Gamma(C)=\Gamma^{*}(C)$, then we call cycle $C$ the closed-cycle of $G^{C}$. Further, we can divide it into three subcases to consider in detail.

\textbf{Case 1.1. Closed cycle-cycle type}, namely $G^C[\Gamma^{*}(C)]$ is a cycle, denoted $C^*=G^C[\Gamma^{*}(C)]$. We say cycle $C^*$ enclose $C$. Correspondingly, we say the semi-maximal planar graph $G^{C}$ with closed-cycle $C$ is a semi-maximal planar graph of closed cycle-cycle type. If the Black-White coloring $f_{bw}$ on $C$ for $G^{C}$ is unique, and the subgraphs, induced by the defined-vertices in each step of $f_{bw}$, is a closed-cycle, then we refer to $G^{C}$ as a \textbf{closed type semi-maximal planar graph on $C$}.

For example, for the semi-maximal planar graph $G^{C_6}$ on 6-cycle $C_6$, shown in Figure 7.7, $C_6=v_1v_2v_3v_4v_5v_6v_1$ is a closed-cycle of $G^{C_6}$. Because $\Gamma(C_6)=\Gamma^*(C_6)=\{v_1',v_2',v_3',v_4',v_5',v_6'\}$ and $G^{C_6}[\Gamma^*(C_6)]=C_6'$ is also a closed-cycle for $\Gamma(C_6')=\Gamma^*(C_6')=\{x,y,z\}$, then $G^{C_6}$ is a closed type semi-maximal planar graph on 6-cycle $C_6$.

\begin{center}
         \includegraphics [width=280pt]{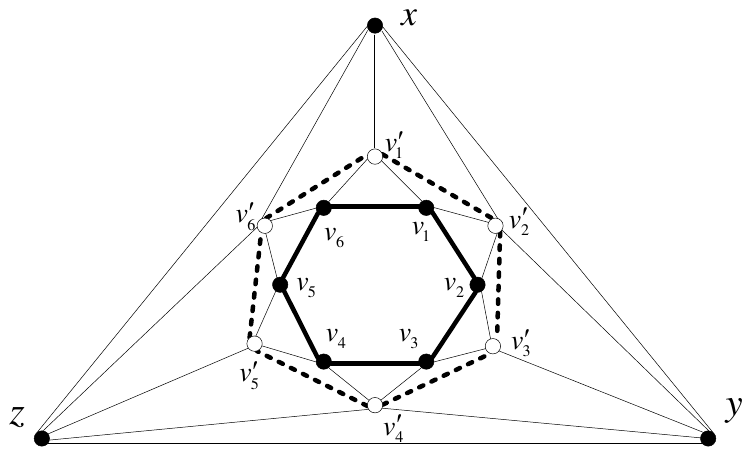}

        \textbf{Figure 7.7.} Schematic diagram for showing closed-cycle and closed type semi-maximal planar graph
\end{center}

For closed type semi-maximal planar graphs, it follows an obvious result as follow.

\begin{lemma}
Suppose that $G^{C}$ is a semi-maximal planar graph on even-cycle $C$, and $\Gamma(C)=\Gamma^{*}(C)$. If $G^C[\Gamma^{*}(C)]$ is a cycle, denoted $C^{*}$, then $C$ is 2-colorable in $G^{C}$ if and only if $C^{*}$ is 2-colorable in $G^{C}-V(C)$.
\end{lemma}

\begin{proof} If $C^{*}$ is 2-colorable in $G^{C}-C$, then $C$ is also 2-colorable in $G^{C}$, obviously. Conversely, if $C$ is a 2-colorable cycle in $G^{C}$, then there exists a coloring $f$ of $G^{C}$ satisfying $|f(C)|=2$. Without loss of generality, we may assume that $f(C)=\{1,2\}$. By the definition of $C^{*}$, we have $f(v)\neq 1,2$ for any $v\in V(C^*)$. So $f(C^*)=\{3,4\}$, namely $C^{*}$ is 2-colorable in $G^{C}$. Thus, $C^{*}$ is also 2-colorable in $G^{C}-V(C)$.
\end{proof}

For example, in Figure 7.7, because cycle $xyzx$ is not a 2-colorable cycle,  then $C_6$ is not a 2-colorable cycle.

\textbf{Case 1.2. Closed cycle-tree type}, namely $G^C[\Gamma^{*}(C)]=G'$ and $G'$ a tree. In this case, $G_1=C$, $G_2$ is a tree, and we say  $G^C$ is a semi-maximal planar graph of closed cycle-tree type. Obviously, $C$ is a 2-colorable cycle.

\textbf{Case 1.3. Closed cycle-fence type}, namely $G^C[\Gamma^{*}(C)]\neq G'$ and $\Gamma^{*}(C)$ is a fence. Then, we say that $G^C$ is a semi-maximal planar graph of closed cycle-fence type. In fact, if  $\Gamma^{*}(C)$ isn't a fence, namely $\Gamma^{*}(C)$ contains odd-cycles, then $C$ isn't a 2-colorable cycle obviously according to \emph{Lemma 7.2}. For a semi-maximal planar graph of closed cycle-fence type $G^C$, $\Gamma^{*}(C)$ is a fence that contains one or more even-cycles. Figure 7.8 gives two examples, in which $\Gamma^{*}(C)$ contains one even-cycle. For this case, obviously, $C$ is 2-colorable if and only if all of the even-cycles in $\Gamma^{*}(C)$ are 2-colorable.

\begin{center}
          \includegraphics [width=260pt]{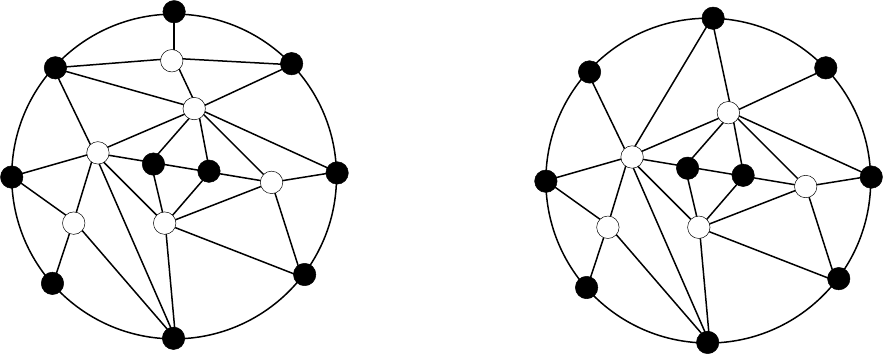}

        \textbf{Figure 7.8.} Two semi-maximal planar graphs of closed cycle-fence type
\end{center}

\textbf{Type 2. Opened-cycles}

If  $\Gamma(C) \neq \Gamma^{*}(C)$, namely $\Gamma^{*}(C)\subset \Gamma(C)$,  then we call cycle $C$ the \textbf{opened-cycle} of $G^{C}$, and $O(C)=\Gamma(C)-\Gamma^{*}(C)$ the \textbf{opened-vertex set} of $C$,  in which the vertices are called \textbf{opened-vertices}. In this case, we refer to $G^{C}$ as \textbf{opened type semi-maximal planar graph}. For example, in Figure 7.1(d), consider the 4-cycle $C=v_1v_3v_5v_6v_1$ and the its inner components, which is a semi-maximal planar graph on $C$. Obviously, $C$ is a opened-cycle, and the opened-vertex set $O(C)=\{v_{10}\}$; in Figure 7.3(b), for the semi-maximal planar graph including the 4-cycle $C=v_1v_3v_2v_4v_1$ and its inner components, $C$ is also a opened-cycle.

 If $C$ is a opened-cycle, then $G_1$ can only contain the opened-vertices of $\Gamma(C)$. Similar as the condition of closed-cycles, we can also subdivide the opened type semi-maximal planar graphs into six subtypes as follows:

\textbf{Case 2.1. Fence-tree type.} $G_1$ is a fence that contains only a even-cycle $C$, and $G_2=G^{C}-G_1$ is a tree. In this case, we call $G^{C}$ the semi-maximal planar graph of fence-tree type. Consider a special situation that $G_1$ is a 0-fence and $G_2$ is a tree, which is similar as case 1.2.

\textbf{Case 2.2. Fence-cycle type.} The connected components including $C$ of $G_1$ is a fence $G_F$ that contains only a even-cycle $C$, and $G_2=G^{C}-G_1$ contains even-cycles, all of which are 2-colorable. Then,  we call $G^{C}$ the semi-maximal planar graph of fence-cycle type(see Figure 7.9(a)).

\textbf{Case 2.3. Scycles-forest type.} Besides $C$, $G_1$ also contains other even-cycles, which either connect mutually through a path, or connect to $C$ through a path containing opened-vertex, and $G_2=G^{C}-G_1$ is a forest. In this case, we call $G^{C}$ the semi-maximal planar graph of scycles-forest type(see Figure 7.9(b)).

\begin{center}
         \includegraphics [width=260pt]{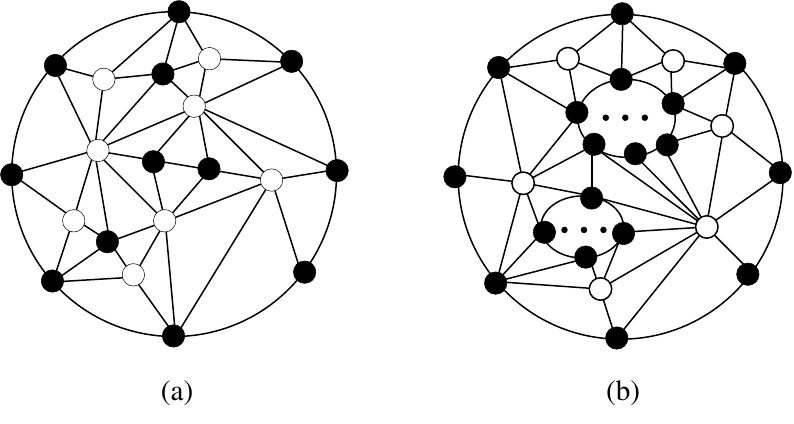}

        \textbf{Figure 7.9.} Semi-maximal planar graphs of fence-cycle, cycles-forest, scycles-scycles type
\end{center}

\textbf{Case 2.4. Scycles-Scycles type.} Besides $C$, $G_1$ also contains other even-cycles. Let $G_
s$ be the connected components including $C_1$ in $G_1$, then $G_2=G^{C}-G_1$ also contains some even-cycles, which are 2-colorable in $G^{C}-G_S$. In this case, we call $G^{C}$ the semi-maximal planar graph of scycles-scycles type.

\textbf{Case 2.5. Intersected cycles-forest type.} $G_1$ consists of some even-cycles $C,C_1,\cdots, C_m$, each of which at least has two vertices of $C$. Namely, the vertices of $C$ are partitioned into $m$ subsets so that each subset are included in at least one even-cycle. At same time, $G_2=G^{C}-G_1$ is a forest. In this case, we call $G^{C}$ the semi-maximal planar graph of intersected cycles-forest type(see Figure 7.10).

\textbf{Case 2.6. Intersected cycles-cycles type.} $G_1$ contains $q$ cycles $C,C_1,\cdots, C_q$. Among them there are $m$ even-cycles $C_1,\cdots, C_m$, each of which at least has two vertices of $C$.
Namely, the vertices of $C$ are partitioned into $m$ subsets so that each subset are included in at least one even-cycle. At same time, $G_2=G^{C}-G_1$ also contains even-cycles that are 2-colorable. In this case, we call $G^{C}$ the semi-maximal planar graph of intersected cycles-cycles type.

\begin{center}
             \includegraphics [width=260pt]{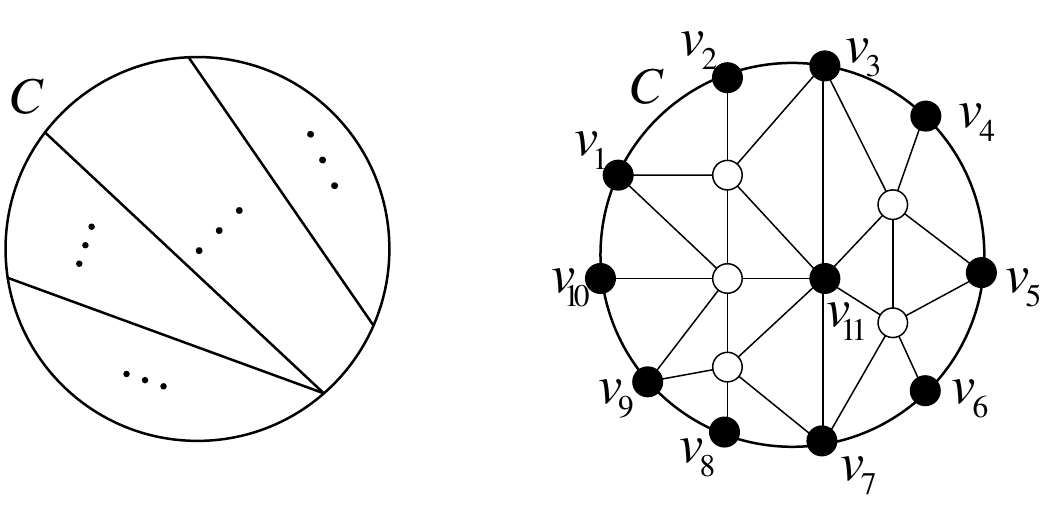}

        \textbf{Figure 7.10.} Semi-maximal planar graphs of intersected cycles-cycles type
\end{center}

The above arguments of nine cases on closed-cycles and opened-cycles, in fact, give a necessary and sufficient condition that $C$ is 2-colorable. In addition, these nine cases are all of the categories when we try to give a classification for semi-maximal planar graphs according to closed-cycle and opened-cycle. So, we can obtain the following result.

\begin{theorem2}
Suppose that $G^{C}$ is a semi-maximal planar graph on even-cycle $C$. Then $C$ is 2-colorable if and only if $G^{C}$ belongs to one of the following items:

\textbf{\textcircled{1}} Closed cycle-cycle type, in which $G^{C}[\Gamma^{*}(C)]$ is 2-colorable in $G^{C}-C$;

\textbf{\textcircled{2}} Closed cycle-tree type;

\textbf{\textcircled{3}} Closed cycle-fence type, in which the even-cycles in the fence are 2-colorable;

\textbf{\textcircled{4}} Fence-tree type;

\textbf{\textcircled{5}} Fence-cycle type;

\textbf{\textcircled{6}} Scycles-forest type;

\textbf{\textcircled{7}} Scycles-Scycles type;

\textbf{\textcircled{8}} Intersected cycles-forest type;

\textbf{\textcircled{9}} Intersected cycles-cycles type.
\end{theorem2}

\subsection{Construction of semi-maximal planar graphs with 2-colorable cycles}

In a semi-maximal planar graph $G^{C}$ on $C$, when we conduct the Black-White coloring on $C$, if $C$ is a bicolored cycle, and there is no other bicolored cycles in the component inside $C$, then we call $C$ the \textbf{basic bicolored cycle} of $G^{C}$ and say $G^{C}$ is a \textbf{basic type}; otherwise, the \textbf{compound bicolored cycle} and  say $G^{C}$ is a \textbf{compound type}. In Theorem 7.16, only in the second and the fourth case, \emph{closed cycle-tree and fence-tree type},  $C$ is a basic bicolored cycle. Obviously, compound bicolored cycles can be gained from basic bicolored cycles through some given operations. Then, what are the operations? This subsection will  reply this question. In fact, just three operations involved: one is cycle-spliced operation, and another two are the bicolored path-split operation and its inverse operation, bicolored cycle-contracted operation.  At the end of this subsection, we study the characteristics and structure of a semi-maximal planar graph $G^{C}$ with a basic bicolored cycle $C$. For this, we introduce two new operations: \textbf{folded operation on even-cycles} and its inverse operation, \textbf{unfolded operation on even-cycles}.  Further, we give the structural characteristics of a semi-maximal planar graph that belongs to  fence-tree type.

\subsubsection{Generating operation system of semi-maximal planar graphs with compound bicolored cycles}
Denote by $G_B^C$ and $G_C^C$ the semi-maximal planar graphs with basic bicolored cycles and compound bicolored cycles, respectively. This subsubsection will give a generating operation system of semi-maximal planar graphs with compound bicolored cycles, denoted $\zeta{(G_C^C)}$, as follow:
\begin{center}
    \textbf{$\zeta{(G_C^C)}$=($G_B^C$, $S$)}
\end{center}
in which $S$ contains three operations: cycle-spliced operation, bicolored path-split operation and bicolored cycle-contracted operation. Now, we will describe these three operations at length.

\textbf{I. Cycle-spliced operation}

Suppose that $G^{C_1}$ and $G^{C_2}$ are two semi-maximal planar graphs with 2-colorable cycle $C_1$ and $C_2$, respectively. Choose two paths, $P_1$ and $P_2$, with length $x$ on $C_1$ and $C_2$, satisfying
$$
1\leq x\leq \left\{\begin{array}{c}
              \min\{|C_1|,|C_2|\}-1, |C_1|\neq |C_2|; \\
             \hspace{0.8cm} |C_1|-2, \hspace{1cm} |C_1|=|C_2|.
            \end{array}\right.
$$
Then, when we merge $P_1$ and $P_2$ into a path $P$ with length $x$, $G^{C_1}$ and $G^{C_2}$ will form a new semi-maximal planar graph $G^{C}$ on $C$, whose length is
$$
|C|=|C_1|+|C_2|-2x.
$$
This process is shown in Figure 7.11. Easily, it follows that $|C|\geq 4$, and $C$ is a 2-colorable cycle of  $G^{C}$.

\begin{center}
          \includegraphics [width=160pt]{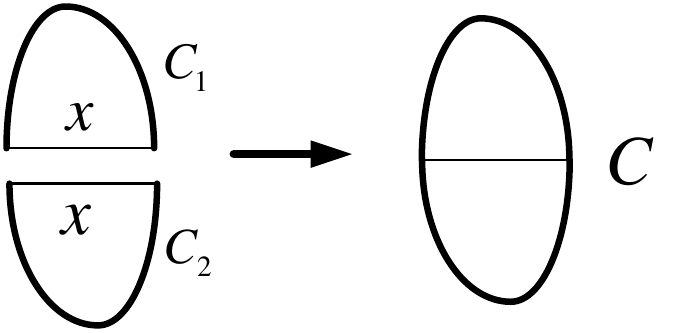}

        \textbf{Figure 7.11.} Schematic diagram of showing cycle-spliced operation
\end{center}

\textbf{II. Bicolored path-split operation}

Let that $G^C$ be a 4-colorable semi-maximal planar graph on even-cycle $C$, which is a 2-colorable cycle.  Suppose that $f\in C_4^0(G^C)$, and the vertices of $C$ are colored by 1,2 alternately. Under $f$, let $P$ be a bicolored path with length $l(\geq 2)$, the vertices of which are colored by 1,2(or 3,4) alternately(see Figure 7.12(a)). Conducting bicolored path-split operation on $P$ in $G^C$, we can gain a new semi-maximal planar graph $G_1^C$ on $C$, and $C$ is also a 2-colorable cycle in $G_1^C$. Following, we will give the detailed description of bicolored path-split operation.

\textbf{(1)} Similar as the extending $4$-wheel operation, along the direction from one end to another end of $P$, cut a crack inside the vertices and edges in accordance with edge-vertex-$\cdots$-vertex-edge order. That is to say except the two ends of $P$, other vertices and edges are cut a crack from their inner side. In this way, each vertex $v$ of $P$ (except the ends) reproduces a new vertex inheriting its color, namely a copy of $v$; and each edge of $P$ reproduces a new edge correspondingly(see Figure 7.12(b)).

\textbf{(2)} Extend path $P$ from the crack and a bicolored cycle $C^{\prime}$ with length $2(l-1)$ will be obtained(see Figure 7.12(c)).

\textbf{(3)} Finally, add a tree $T$ to the inside of $C^{\prime}$ properly, and connect the vertices  between $T$ and $C^{\prime}$ properly so that all the faces inside $C^{\prime}$ are triangle, and no vertices have degrees less than 4.

\textbf{III. Bicolored cycle-contracted operation}

Let that $G^C$ be a 4-colorable semi-maximal planar graph on even-cycle $C$, which is a 2-colorable cycle.  Suppose that $f\in C_4^0(G^C)$, and the vertices of $C$ are colored by 1,2 alternately.

\textbf{(1)} Under $f$, suppose that $C_p$ is a bicolored cycle that is different from $C$, and the vertices of $C$ are colored by 1,2(or 3,4) alternately. Choose a vertex on $C_p$ arbitrarily, say $v_1$, and then choose another vertex on $C_p$, say $v_2$, which has the longest distance with $v_1$(see Figure 7.12(c)).

\textbf{(2)} Delete all of the vertices inside $C_p$(see Figure 7.12(c)).

\textbf{(3)} Starting from $v_1$, we identify any pair of vertices that have the same distance with $v_1$(see Figure 7.12(a),(b)). Thus, a new 4-colorable semi-maximal planar graph, $G^{C}_{C}$, is to obtained, and $C$ is also a 2-colorable in $G^{C}_{C}$.

\begin{center}
        \includegraphics [width=260pt]{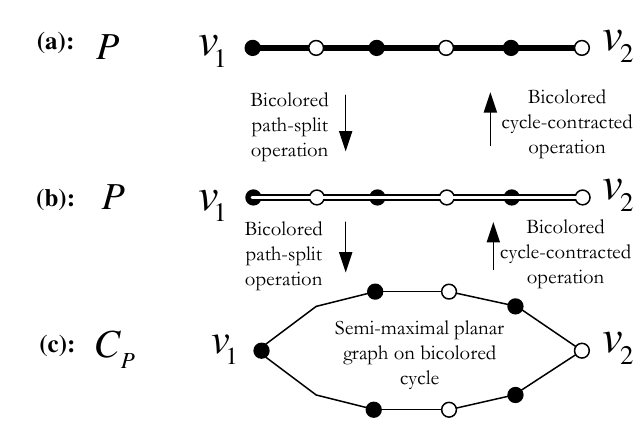}

        \textbf{Figure 7.12.} Bicolored path-split operation and bicolored cycle-contracted operation
\end{center}

On the basis of the above three operations, now we give a method to construct a semi-maximal planar graph with 2-colorable cycles as follow.

Step 1. Given some semi-maximal planar graphs, $G^{C_1},G^{C_2},\cdots,G^{C_m}$, on basic bicolored cycles $C_1,C_2,\cdots,C_m(m\geq 2)$. Namely, $G^{C_1},G^{C_2},\cdots,G^{C_m}$ belong to closed-tree type or fence-tree type;

Step 2. Conducting cycle-spliced operation among  $G^{C_1},G^{C_2},\cdots,G^{C_m}$, we can produce some new semi-maximal planar graphs with 2-colorable cycles, which belong to intersected-cycle type.

Step 3. The semi-maximal planar graphs with 2-colorable cycles, which belong to closed-cycle type and closed-fence type, can be constructed by conducting bicolored path-split operation on the semi-maximal planar graphs with 2-colorable cycles, which belong to closed-path type and closed-tree type, respectively;

Step 4. The semi-maximal planar graphs with 2-colorable cycles, which belong to fence-cycle type, scycle-forest type and scycle-scycle type, can be constructed by conducting bicolored path-split operation on the semi-maximal planar graphs with 2-colorable cycles, which belong to fence-path type;

Step 5. For any semi-maximal planar graphs with 2-colorable cycles belonging to one of the nine types of Theorem 7.16, we can always produce various semi-maximal planar graph with 2-colorable cycles by conducting bicolored path-split operation and cycle-spliced operation, simultaneously.

Step 6. For any semi-maximal planar graph $G_C^{C}$ on compound bicolored cycles, when we conduct  bicolored cycle-contracted operation on $G_C^{C}$, a semi-maximal planar graph $G_C^{C}$ on a basic bicolored cycle will be produced.

\subsubsection{Folded operation, unfolded operation and the characteristics of 2-colorable cycles}

We have known that semi-maximal planar graphs $G^{C}$ on the basic bicolored cycle $C$ belongs to either closed-tree type or fence-tree type. For the former, its structure is very clear that the subgraph $G^{\prime}$ induced by the vertices inside cycle $C$ is a tree; but for the latter, $G^{C}$ can be viewed as the union of a fence $G_F$ on cycle $C$ and a tree $T$. Namely, $G^{C}=G_F\cup T$ and $G_F\cap T=\emptyset$. So, $G^{\prime}$ contains not only the vertices in $T$, but also the vertices in $G_F-C$. Then, what are the characteristics of $G^{\prime}$? The following will discuss this problem, for which we need a pair of operation on even-cycles: \textbf{folded operation} and \textbf{unfolded operation}.

\textbf{I. Folded operation}

Actually, the so-called folded operation is to conduct the action of identifying some succussive pairs of vertices. Suppose that $G^{C}$ is a semi-maximal planar graph on cycle $C$, which belongs to closed-tree type, and $f\in C_4^0(G^C)$. When $|C|\geq 6$, choose a vertex on $C$ arbitrarily, say $u$, and identify the two vertices adjacent to $u$ on $C$, say $v_1,v_1'$. Define that the new 4-colorable semi-maximal planar graph$ G^{C}\circ\{v_1,v_1'\}\triangleq G_1^{C_1}$,  $\{v_1,v_1'\}\triangleq v_1$ and $C_1$ is the new outer cycle(see Figure 7.13).
If $|C_1|\geq 6$, then identify the two vertices adjacent to $v_1$ on $C_1$, say $v_2,v_2'$, and we will produce another new semi-maximal planar graph $G_2^{C_2}$ on cycle $C_2$. This procedure can be continued until some new outer cycle $C_i$ satisfying \emph{$|C_i|=4$}.

\textbf{II. Unfolded operation}

In fact, the so-called unfolded operation is the inverse operation of folded operation. Different from the bicolored path-split operation, for a bicolored path $P=v_1v_2\cdots v_l$, we can split one of its ends $v_l$, the terminal end,  but it is not permitted in bicolored path-split operation. The detailed description as follows:

Step 1. For a semi-maximal planar graph $G^{C_1}$ on cycle $C_1$, which belongs to fence-tree type, and $f\in C_4^0(G^{C_1})$.  Suppose $G_F$ and $T$ are the fence and tree, respectively. Namely, $G_F\cup T=G^{C_1}$, and $C\subseteq G_F$. We choose a bicolored $t$-path in $G_{F}$, denoted $P=u-v_1$, which only contains one vertex $v_1$ of $C_1$, namely $v_1$ is one end of $P$;

Step 2.  Starting from vertex $u$, along the direction from $u$ to $v_1$ on path $P$, cut a crack inside the vertices and edges in accordance with edge-vertex-$\cdots$-edge -vertex order. That is to say except the end $u$, other vertices and edges are cut a crack from their inner side. In this way, each vertex $v$ of $P$ (except $u$) reproduces a new vertex inheriting its color, namely a copy of $v$; and each edge of $P$ reproduces a new edge correspondingly;

Step 3. Extend path $P$ from the crack and a new semi-maximal planar graph will be obtained.

Obviously, when we conduct unfolded operation for all the paths in $G_F$, which only contain one vertex of $C_1$, the resulting graph is going to be a semi-maximal planar graph belonging to closed-tree type.
For example, Figure 7.13(b) is a semi-maximal planar graph belonging to fence-tree type. Of course, it is opened and $u$ is a opened-vertex. Conducting the unfolded operation on path $P=uv_1$, we can obtain the graph shown in Figure 7.13(a).

\begin{center}
         \includegraphics [width=280pt]{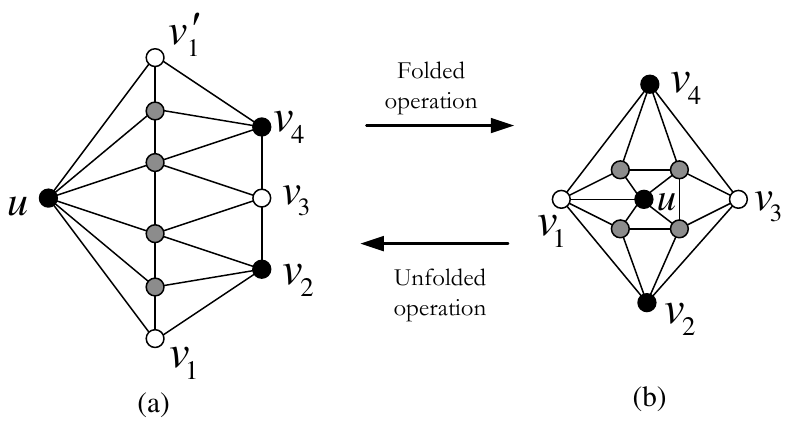}

        \textbf{Figure 7.13.} Schematic diagram of showing folded operation and unfolded operation
\end{center}

Considering the case that in $G_F$, some suspending vertices connect with the vertex $v_1$ of $C_1$ through a tree. That is to say, in $G_F\backslash (V(C_1)-\{v_1\})$, the connected branch containing $u$ is a tree, denoted $T'$, not a path. At this case, we need process it by stages. First, we should choose a path $P$ that starts from a vertex $v_1$ of $C_1$, which is adjacent to opened-vertices, to a vertex $u'$ of $T'$, which has degree not less than 3. Conducting unfolded operation on $P$, we can obtain a new semi-maximal planar graph $G^{C_2}$, which  belongs to fence-tree type. Obviously, if we remain colors appearing in the vertices in $G^{C_1}$ unchanged, including the copies of the vertices of $P$, then the 4-coloring of $G^{C_2}$ is also proper.
Next, similar as the above process, conducting the  unfolded operation in $G^{C_2}$, we can obtain another 4-colorable semi-maximal planar graph $G^{C_3}$. Continue this procedure until some semi-maximal planar graph $G^{C_i}$, which belongs to closed-tree type. Figure 7.14 gives the illustration of this procedure, in which Figure 7.14(a) gives a semi-maximal planar graph belonging to fence-tree type, Figure 7.14(b) gives the new semi-maximal planar graph after conducting unfolded operation on path $P=uv$, Figure 7.14(c) gives another semi-maximal planar graph after conducting unfolded operation on path $P=uu_1$, and Figure 7.14(d) gives another semi-maximal planar graph after conducting unfolded operation on path $P=uu_2$.

\begin{center}
        \includegraphics [width=280pt]{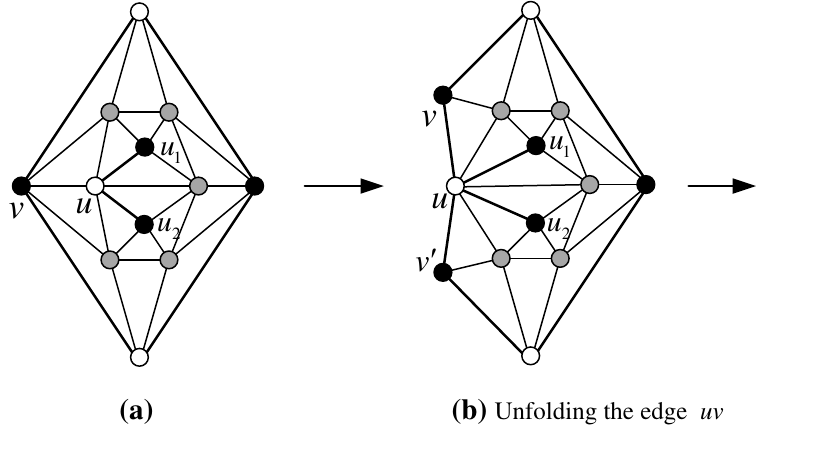}

         \includegraphics [width=300pt]{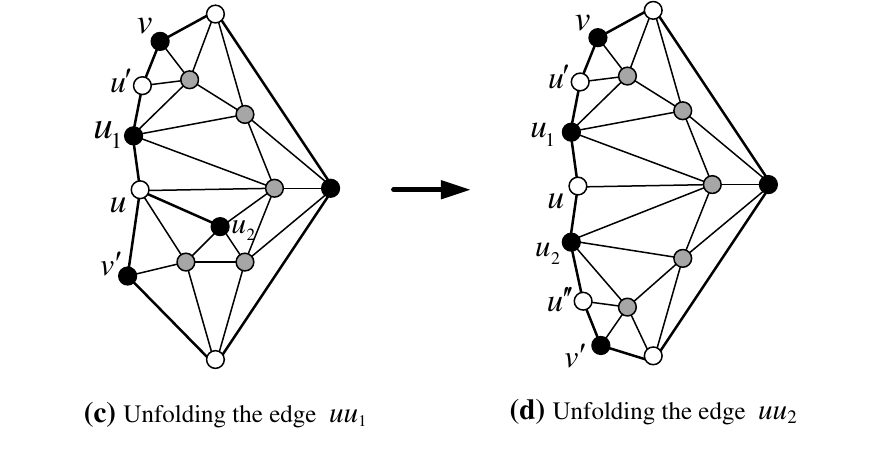}

        \textbf{Figure 7.14.} Schematic diagram of unfolded operation in special case
\end{center}

It follows the following fact from the above two operations.

\begin{theorem2}
Suppose that $G^{C}$ is a semi-maximal planar graph belonging to closed-tree type with $|C|\geq 6$, $V'=V(G^{C})-V(C)$, and $G[V']\triangleq G'$ is a tree. When conduct folded operation on $C$, we can obtain a new semi-maximal planar graph $G^{C_1}$, which belongs to fence-tree type. Denote by $G_F$ and $T$ the fence and tree of $G^{C_1}$ respectively, then $T=G'$. That is to say, when we conduct the folded operation in a semi-maximal planar graph belonging to closed-tree type, the structure of its fence gets changed, but the tree remains unchanged. Conversely, for a semi-maximal planar graph belonging to fence-tree type $G^{C_1}$, when we conduct unfolded operation in $G^{C_1}$ repeatedly, at the end we can obtain a semi-maximal planar graph belonging to closed-tree type, and the structure of tree will not be affected in the process of unfolded operation.
\end{theorem2}

This theorem actually tells us that fence-tree structures can be seen obtained from closed-tree structures by conducting folded operation in them, repeatedly. Because the tree structure remained unchanged in the process of folded operation, we need only study the change of edges on the cycle after conducting this operation.

Here, we give an example of constructing a semi-maximal planar graph on two paths. Suppose $P,P'$ are two paths with length not less than 3. Now, we will construct a semi-maximal planar graph $G^{C}$ from $P$ and $P'$(see Figure 7.15), and the resulting graph is called the semi-maximal planar graph of \textbf{cycle-path type}.

Step 1. Connect one end $u$ of $P'$ to the two ends of $P$, then connect another end of $P'$ to at least three vertices of $P$;

Step 2. For the vertices of $P'$, connect them to the vertices of $P$ such that each connected edge is in a triangle, and the degrees of vertices on $P$ except $u$ have to increase to 4 at least.

\begin{center}
         \includegraphics [width=360pt]{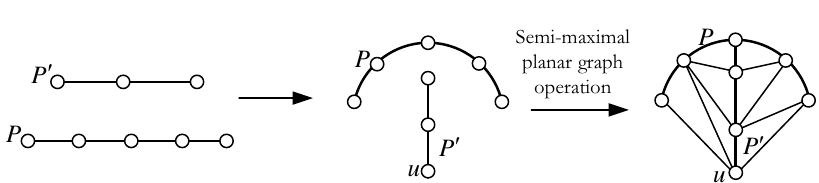}

        \textbf{Figure 7.15.} Schematic diagram of constructing a semi-maximal planar graph belonging to closed-tree type on two paths
\end{center}

In the process of constructing semi-maximal planar graphs on two paths, we can deem that the path $P'$ is obtained by identifying pairs of vertices with the same color in a semi-maximal planar graph $G^{C}$ belonging to closed-tree type repeatedly. That is to say, $P'$, generated in the process of conducting folded operation for $G^{C}$, is a path in a semi-maximal planar graph $G^{C_1}$ belonging to fence-tree type. If $G^{C_1}$ can be obtained only by conducting folded operation once in $G^{C}$, then $G^{C_1}[V_1' \cup\{u\}]$ is the join of a semi-maximal planar graph belonging to cycle-path type
and some possible trees, where $V'=V(G^{C_1})-V(C_1)$ and $u$ is the weld-vertex in the fence; if $G^{C_1}$ is obtained by conducting folded operation many times in $G^{C}$, then this process can be seen turned back as conducting cycle-spliced operation among semi-maximal planar graphs belonging to closed-path type, and then connecting some possible trees on the outer cycle. We call these graphs \textbf{barrette-structure} graphs. Thus, we have prove that if $G^{C_1}$ is a semi-maximal planar graph belonging to fence-tree type, then the subgraphs induced by the inside vertices of $C_1$ and the vertices on $C_1$ joined the trees in the fence is a barrette-structure graph.

\begin{theorem2}
Suppose that $G^{C_1}$ is a semi-maximal planar graph belonging to fence-tree type, $u_1,u_2,\cdots,u_k$ are all weld-vertices in the fence, and $V'=V(G^{C_1})-V(C_1)$, then $G^{C_1}[V_1' \cup \{u_1,u_2,\cdots,u_k\}]$ is a barrette-structure graph.
\end{theorem2}

\subsection{Summary}

The contents of this section mainly cover the following six aspects.

\textbf{First}, we point out that there have two categories of cycles in maximal planar graphs, basic cycles and chord cycles, and study their distribution and enumeration, which paves the way for the latter research.

\textbf{Second}, propose a new method, Black-White coloring, to study maximal planar graphs. The advantages of this method is that its process can be realized simply, and for a maximal planar graph, each of its Black-White coloring consists of a subset of its 4-colorings set. Especially, it is a power technique to study the 2-colorable cycles.

\textbf{Third}, set up the petal-syndrome, on which we find a necessary and sufficient condition that a even-cycle is a 2-colorable cycle in a maximal planar graph.

\textbf{Fourth}, make clear a necessary and sufficient condition of a 2-colorable cycle on structure. Namely, independently satisfy two basic types and seven compound types(see Theorem 7.15).

\textbf{Fifth}, prove that each compound type can be obtained from basic type by conducting cycle-spliced operation, bicolored path-split operation and bicolored cycle-contracted operation.

\textbf{Sixth}, make clear the relationship of two basic types through introducing the folded operation and unfolded operation, and depict the inner structure of the semi-maximal planar graphs belonging to fence-tree type: a barrette-structure graph. For this, the structure of the semi-maximal planar graphs with 2-colorable cycles are described deeply.

However, for the compound types, it is still a tough problem to judge which type they belong to in the seven cases of Theorem 7.15. In order to tackle this problem, we need argue it combined with the open-vertices properly.

The more in-depth study on this problem will be given in later sections. 
\section{Coloring operation system for maximal planar graphs(I)--protected-cycle operation}

In the fourth section, we proposed and deeply discussed the generating operation system of maximal planar graphs, from which we have clearly known the structure of maximal planar graphs. On the basis of the fourth section, the following chapters are to advance the \textbf{generating operation system of coloring} for maximal planar graphs, the aim of which is to research the relationship between any two 4-colorings of a maximal planar graph $G$. It intends to deduce other 4-colorings based on one 4-coloring of $G$. Namely make fundamentally clear how the 4-colorings are generated and the mutual relationship between two non-isomorphic 4-colorings. Specially, on the basis of the classification of cycle-colorings and tree-colorings of 4-colorable maximal planar graphs, we put forward the generating operation system of coloring for a maximal planar graph $G$, denoted $< f_0, \Phi>$, in which $f_0$ denotes a 4-coloring of $G$ and $\Phi$ the set of coloring operations. $\Phi$ contains three basic operators: \textbf{protected-cycle operation}, denoted $\sigma_{p-C}$; \textbf{broken-cycle operation}, denoted $\sigma_{b-C}$; \textbf{broken-tree operation}, denoted $\sigma_{b-T}$. The main task of this section is to study the first operator---protected-cycle operation. Here, we, especially, announce that all of the graphs involved in this section are the maximal planar graphs with $\delta \geq 4$.

\subsection{Definitions and properties}

Let $G$ be a 4-colorable maximal planar graph, $C(4)=\{1,2,3,4\}$ the color set, and $C=v_1v_2\cdots v_{2m}v_1$ a even-cycle of $G$. If $f\in C_4^0(G)$ is a cycle-coloring on cycle $C$, under which $f(v_i)=1$ or $f(v_i)=2$, $i=1,2,\cdots, 2m$, then the so-called \textbf{protected-cycle operation of $f$} on $C$ , denoted $\sigma_{p-C}$, is to exchange the colors (3 and 4) assigned to the vertices inside(or outside) cycle $C$, and remain the colors assigned to other vertices unchanged. Obviously, the protected-cycle operation of $f$ is to change it into another new cycle-coloring of $G$, which is denoted by $f'$, and $$
\sigma_{p-C}(f)=f'\eqno{(8.1)}
$$
We refer to $f'$ as the \textbf{complementary cycle-coloring} of $f$ on $C$ . Of course, $f$ is also a complementary cycle-coloring of $f'$ on $C$ . Namely, $f$ and $f'$ are a pair of complementary cycle-coloring on $C$.

For example, for the graph in Figure 6.1(b), there are 8 4-colorings $f_1,f_2,f_3,f_4$,
$f_5,f_6,f_7,f_8$. It is easy to see that on $1-4$ cycle(red-yellow cycle), $f_1'=f_2$; on $1-3$ cycle(red-green cycle), $f_3'=f_2$;
on $3-4$ cycle(green-yellow cycle), $f_4'=f_2$; on $2-3$ cycle(blue-green cycle), $f_5'=f_4$; on $2-4$ cycle(blue-yellow cycle), $f_7'=f_6$; on $1-2$ cycle(red-blue cycle), $f_8'=f_6$. Here, we should notice that when there are many bicolored cycles under a 4-coloring $f$, the complementary cycle-colorings of $f$ on different bicolored cycles are different. And the following result is obvious.

\begin{theorem2}
Suppose that $f$ is a cycle-coloring of a 4-colorable maximal planar graph $G$, and $C$ is a bicolored cycle of $f$, then
$$
\sigma_{p-C}(\sigma_{p-C}(f)) \triangleq \sigma_{p-C}^{2}(f)=f \eqno{(8.2)}
$$
$$\eqno{\Box}$$
\end{theorem2}

Based on protected-cycle operation of coloring, we will put forward a new concept, eigen graph of protected-cycle coloring, which helps to research the relationships among the 4-colorings in $C_4^{0}(G)$, clearly and intuitively.

\subsection{Eigen graph of protected-cycle coloring}

Let $G$ be a 4-colorable maximal planar graph, and $C_4^{0}(G)=\{f_1,f_2,\cdots,f_n\}$. Now, we construct a new graph $P(G)$, called \textbf{Eigen graph of protected-cycle coloring of $G$}, in which $V(P(G))=C_4^{0}(G)=\{f_1,f_2,\cdots,f_n\}$, and $f_i$ is adjacent to $f_j$ if and only if $f_i$ and $f_j$ are a pair of complementary cycle-coloring on some bicolored cycle, $i,j=1,2,\cdots, n$, $i\neq j$.

Obviously, for a maximal planar graph $G$, if one of its 4-coloring $f$ contains bicolored cycles, then $\sigma_{p-C}(f)$ also contains bicolored cycles. We call the connected branches that contain cycle-colorings in $P(G)$ \textbf{cycle branches}. Obviously, each of cycle branches contains at least two vertices. For any $f_1,f_2\in C_4^{0}(G)$, if they are  in the same connected branch in $P(G)$, then we say $f_1,f_2$ are \textbf{connected}, otherwise, say they are \textbf{disconnected}. In addition, we similarly call the connected branches that contain tree-colorings in $P(G)$ \textbf{tree branches}. Obviously, in $P(G)$, each of the tree branches just contains a isolated vertex.

\textbf{Example 8.1} For the graph $G$ shown in Figure 6.1(b), there are in total eight 4-colorings for $G$. Correspondingly, Figure 8.1(a) gives its eigen graph of protected-cycle coloring, and it is a tree, which indicates that any two of 4-colorings of $G$ can be deduced mutually by protected-cycle coloring operation.
Figure 6.2 has totally four 4-colorings: a pair of complementary cycle-coloring and two tree-colorings. Its eigen graph of protected-cycle coloring, shown in Figure 8.1(b), is a disconnected graph. In addition, for the icosahedron, there are in total ten 4-colorings(see Figure 6.3), each of which is a tree-coloring, so its  eigen graph of protected-cycle coloring is a complete null graph with ten isolated vertices.

\begin{center}
   \includegraphics [width=260pt]{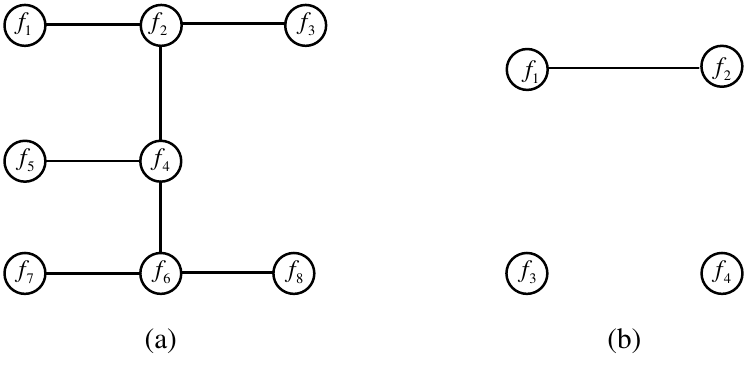}

  \textbf{Figure 8.1} Eigen graphs of protected-cycle coloring of graphs shown in Figure 6.1(b) and 6.2
\end{center}

However, inspired by the above example, we will naturally propose a problem as follow.

\begin{problem}
\textcircled{1} What kinds of the 4-colorable maximal planar graphs are these that the eigen graph of protected-cycle coloring of them are connected?

\textcircled{2} What is the necessary and sufficient condition that $P(G)$ is a $l$-path for a maximal planar graph $G$?
\end{problem}

It easily follows
\begin{theorem2}
Let $G$ be a 4-colorable maximal planar graph and $P(G)$ the eigen graph of protected-cycle coloring of $G$. It has

\textcircled{1} $G$ is a unique 4-colorable maximal planar graph if and only if $P(G)=K_1$;

\textcircled{2} for any vertex $f$ of $P(G)$, $d_{P(G)}(f)=k$ if and only if the coloring $f$ of $G$ contains $k$ bicolored cycles.
$$\eqno{\Box}$$
\end{theorem2}

\begin{corollary}
Suppose that $G$ is a 4-colorable maximal planar graph. If there is a 4-coloring $f\in C_4^0(G)$ that contains $k(\geq 1)$ bicolored cycles, namely $d_{P(G)}(f)=k$, then
$$
|C_4^0(G)|\geq k+1 \eqno{(8.3)}
$$
$$\eqno{\Box}$$
\end{corollary}

For a 4-colorable maximal planar graph $G$, the following point out the key issues that need studying from three aspects.

First, two cycle-colorings $f_1,f_2 \in C_4^0(G)$ may be disconnected, that is to say $f_1,f_2$ are not necessarily to be in the same connected branch in $P(G)$. For example, the cycle-colorings shown in Figure 6.1(b) and Figure 6.2 are connected, but the two cycle-colorings shown in Figure 8.2 must be disconnected because the first (denoted by $f$) only contains one bicolored cycle, 1-4 cycle(red-yellow cycle), and the new cycle-coloring, obtained by exchanging color 2 and 3(green and blue) on the vertices inside the cycle, also contains one bicolored cycle, 1-4 cycle(red-yellow cycle). We call such cycle-coloring $f$ the \textbf{bicolored cycle-unchanged coloring}, and its specific definition will be given in the next section. Obviously, in Figure 8.2, the second cycle-coloring can not be obtained by conducting the protected-cycle operation of $f$ on its unique bicolored cycle  1-4 cycle(red-yellow cycle), so these two cycle-colorings must be disconnected. In addition, if $G$ contains bicolored cycle-unchanged colorings, then such coloring must appear in pairs, and each pair constructs a connected branch $K_2$ in $P(G)$. Naturally, a question is likely to be asked whether there is a cycle with length not less than 3 in $P(G)$.

\begin{center}
\includegraphics [width=340pt]{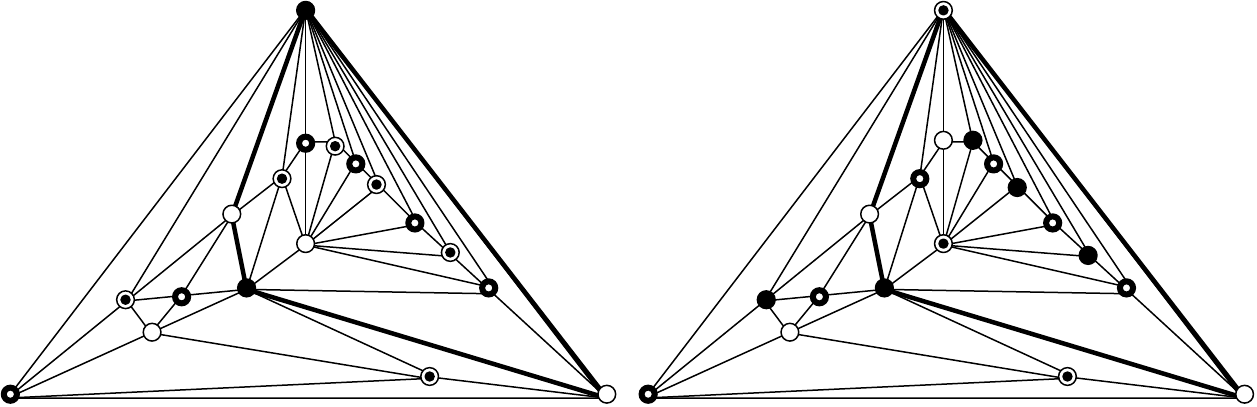}

\textbf{Figure 8.2} Two disconnected colorings of a maximal planar graph
\end{center}

\begin{problem}
Let $G$ be a 4-colorable maximal planar graph and $P(G)$ the eigen graph of protected-cycle coloring of $G$. What is the necessary and sufficient condition that there is a cycle with length not less than 3 in $P(G)$?
\end{problem}

Second, how to construct one component through another component in $P(G)$ so as to gain the whole $P(G)$? There are three aspects that deserve to be considered seriously. \textcircled{1} What methods or techniques are to be used to interchange between two cycle branches?   \textcircled{2} How to interchange between cycle branches and tree branches in $P(G)$? \textcircled{3} How to interchange among tree branches in $P(G)$?  For the first two cases, we depend on the so-called broken-cycle operation to study them deeply; for the last case, we introduce the broken-tree operation to study it.

Third, on the basis of the above two problems, how to give the estimated value or precise value of $|C_4^0(G)|$?

In the following chapters, we will discuss mainly around the above three issues. In the process of the discussion, it is necessary to research the bicolored cycle-unchanged colorings and their corresponding graphs.

\subsection{Properties of protected-cycle operation}

In the sixth section, there was a probable result that tree-colorings are very few compared to cycle-colorings for maximal planar graphs. Among all maximal planar graphs with order at most $11$ and $\delta\geq 4$, tree-colorings accounts for only about 2.12\%. So, for a given cycle-coloring $f$ of a maximal planar graph $G$, it is largely possible to find other 4-colorings of $G$ on the basis of $f$. For example, the eigen graph of protected-cycle coloring of the graph shown in Figure 6.1(b) is a connected graph, so if only we know one coloring of $C_4^0(G)$, all of others can be deduced by conducting protected-cycle operation. For the general case, it has

\begin{theorem2}
Suppose that $G$ is a 4-colorable maximal planar graph. If $P(G)$ is connected, then the difference is $|V(P(G))|-1$ between complexities of the two algorithms of finding out all of the 4-colorings and finding out a 4-coloring of $G$. That is to say these two algorithms are equivalent.
\hspace{1cm} $\Box$
\end{theorem2}

In addition, according to the definition of protected-cycle operation of a cycle-coloring $f$, it is in fact a destruction to other bicolored cycles of $f$ in some sense! That is to say based on $f$, some bicolored cycles of $f$ may disappear in $\sigma_{p-C}(f)$; likewise, some new bicolored cycles may also appear in $\sigma_{p-C}(f)$.

\begin{theorem2}
Suppose that $G$ is a 4-colorable maximal planar graph, and $f$ is a 4-coloring of $G$. $C$ is a bicolored cycle on $f$, and $f$ is a suspending vertex in $P(G)$. If the degree of $f'=\sigma_{p-C}(f)$ is not less than 2 in $P(G)$, then there must exist a coloring $f''\in C_4^0(G)$, in which $C$ is not a bicolored cycle of $f''$.
\end{theorem2}

\begin{proof}
Suppose that the color set $C(4)=\{1,2,3,4\}$. For cycle $C$, a bicolored cycle(2-3 cycle) of $f$(see Figure 8.3(a)), we conduct the protected-cycle operation of $f$ on $C$. Without loss of generality, we assume that the colors on the vertices inside $C$ are exchanged. According to the assumption of the theorem, $f'=\sigma_{p-C}(f)$ contains at least two bicolored cycles. So, under coloring $f$, there are $k$ vertices( say $v_1,v_2,\cdots,v_k$ ) that are colored by the same color(say 2) on $C$, which satisfies $v_{i},v_{i+1}$ connected by a bicolored path, 1-2 path(or 4-2 path), inside $C$ or a bicolored path, 4-2 path(or 1-2 path) outside $C$(see Figure 8.3(a), the case of $k=2$), $i=1,2,\cdots,k$. Here, the subscript of each vertex is taken modulo $k$ and in $\{1,2,\cdots,k\}$. Obviously, $C'$, $1-2$ cycle, and $C$, $2-3$ cycle, are two different bicolored cycles of $f'$, in which there are vertices of $C$ inside $C'$(see Figure 8.3(b)). Then, we conduct protected-cycle operation of $f'$ on $C'$. Likewise, we exchange the color 3 and 4 on the vertices inside $C'$ in the process of the operation, and obtain a new 4-coloring $f''$ of $G$, namely $\sigma_{b-C}(f')=f''$. So, $C$ is not a bicolored cycle of $f''$(see Figure 8.3(c)).

\end{proof}

\begin{center}
  \includegraphics [width=300pt]{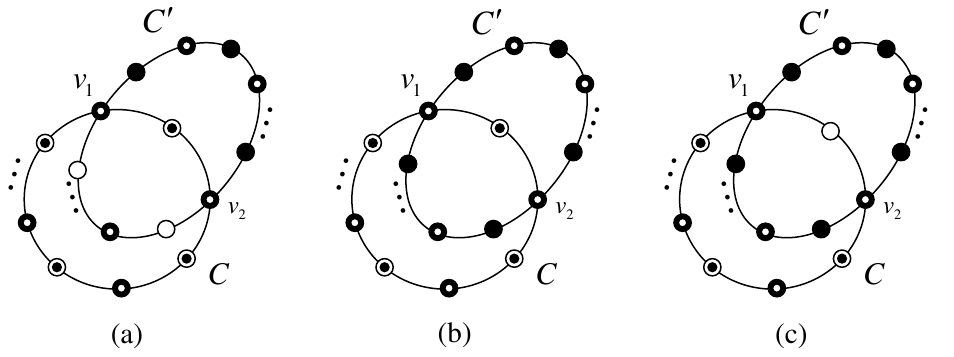}

\textbf{Figure 8.3} Schematic diagram of illustrating Theorem 8.5 when $k=2$
\end{center}

This theorem in fact shows that if $f$ is not a bicolored cycle-unchanged coloring on cycle $C$, then there must exist another 4-coloring $f''$, of which $C$ is not a bicolored cycle. In other words, we can also say that $f''$ is a \textbf{broken-cycle operation} of $f$ on $C$. Namely finding a 4-coloring $f''\in C_4^0(G)$ such that the number of colors assigned to the vertices of $C$ is not less than 3 under $f''$.

\begin{corollary}
Suppose that $G$ is a 4-colorable maximal planar graph, $f$ is a 4-coloring of $G$, and $C$ is a bicolored cycle on $f$. If $f$ is not a bicolored cycle-unchanged coloring on $C$, then there must exist a 4-coloring $f''\in C_4^0(G)$ satisfying that $C$ is not a bicolored cycle of $f''$.
$$
\eqno{\Box}
$$
\end{corollary}

From the above discussions, we can know that if $G$ is a 4-colorable maximal planar graph, $f$ is a 4-coloring of $G$, and $C$ is a bicolored cycle on $f$, then there must exist a  4-coloring $f''\in C_4^0(G)$ satisfying that $C$ is not a bicolored cycle of $f''$ when the number of vertices of the connected branch containing $f$ in $P(G)$ is not less than 3. We refer to the process of conducting the action that make $C$ be not a bicolored cycle of $f''$ from that $C$ be a bicolored cycle of $f$, as \textbf{broken-cycle operation of $f$ on $C$}, denoted $\sigma_{b-C}$. Namely
$$
\sigma_{b-C}(f)=f''\eqno{(8.4)}
$$
Further, we say \textbf{$C$ is destructible}. However, when $f$ is a bicolored cycle-unchanged coloring of $G$ on even-cycle $C$, is there a 4-coloring $f''\in C_4^0(G)$ satisfying that $C$ is not a bicolored cycle of $f''$? Namely, are this kind of cycles also destructible? So, we propose the following problem.

For any bicolored cycle $C$ of a 4-coloring of a maximal planar graph, is $C$ destructible? Of course, the answer is positive, we will prove this fact completely in next section. 
    \vspace{5mm}
    \begin{center}
    \textbf{Acknowledgements}
    \end{center}
    \vspace{5mm}

    The first draft of this paper was written in 1991, when I worked in the Mathematics Department of Shaanxi Normal
 University, my Alma Mater. I reported my work to several professors there, including Prof.Guojun Wang, Prof. Xiansun
 Wei, Prof. Zhongqiang Yang, Prof. Taihe Fan, Prof. Baolin Guo and Prof. Wanmin Zhang, etc. All of them gave me some
 useful advices, especially Prof. Zhongqiang Yang pointed out a fatal error, which was a necessary and
 sufficient condition for uniquely 4-colorable planar graphs. In other words, there were some errors in the proof of the FWF
 conjecture. To overcome it, I have since spent more than eighteen years on this problem until August 6th, 2009.
 I would like to thank my teachers Prof. Xinmin Wang, Prof. Hongke Du for their
 encouragement, caring and support on my work in 1991, especially Prof. Hongke Du, who has been
encouraging and directing me on it since then. Here I would like to
thank them all deeply.

    After the draft was completed in Nov. 9th, 2009, Prof. Jianfang Wang, Prof. Liang
Sun discussed with me directly. They pointed out that there was one
obvious leak on the proof of \"If a maximal planar graph has a
minimum degree of 5, it is not uniquely 4-colorable. \" Here I also
would like to thank them deeply. At the same time, many experts on
this field have looked over my paper, such as Prof. Zhongfu Zhang,
Prof. Bin Yao, Prof. Xiangen Chen, Prof. Hui Chen, Prof. Muchun Li,
Prof. Huiying Qiang, student Zepeng Li and many international
experts such as Prof.Jensen. All of them have pointed out the same
leak as Prof. Liang Sun and Prof. Jianfang Wang's. Here I would like
to thank them all deeply.

    During the drafting of first version and second version of this paper,
 my students Fang Xi, Mei Chen, Enqiang Zhu, Jingming Liu, Zhen Chen,
 Yufang Huang, Ziqi Wei and Dongming Zhao spent a lot of time on this paper,
 from translation, drawing, typewriting to proofreading, especially for
 Enqiang Zhu and Fang Xi. During the accomplishment of second version of my
 paper, Enqiang Zhu paid lots of attention on it. He discussed with me about
 some problems in my proof, completed the drawing of numbers of
 figures and English translation; Fang Xi was in charge of the
 overall work of final draft, who had maximal workload and often
 worked very late into the night. Moreover, Prof.chunling Quan validated all the
 graphs in Figure 5.7, 5.8 and 5.9 by electronic computer. Here I would like
 to thank them all deeply.

I would like to thank my colleges Prof. Daoheng Yu for his continuous encouragement and  
Prof. Tian Liu for his helps to improve our English presentation.

    Here, I also would like to thank my tutors Prof. Ziguo Wang, Prof. Yingluo Wang (Fellow of Chinese Academy of Engineering) and Prof. Zheng Bao
 (Fellow of Chinese Academy of Sciences) for their advice and help through years.


\section{Appendix}
    This appendix gives all 4-colorings of the maximal
 planar graphs whose orders are from 6 to 10 and
 $\delta(G)=4$.
\begin{table}[!htb]
\begin{center}
\caption{\label{tab:test}Lower Bound of the number of partitions of color
groups about all  maximal planar graphs whose orders are from 6 to
10 and contain no three adjacent vertices of degree 4.}

\vspace{5mm}
\begin{tabular}{ccccccccccccc}
  \hline
  Order && 6 && 7 && 8 && 9 && 10 \\
  \hline
  Lower Bound && 4 && 5 && 3 && 6 && 5 \\
  \hline
\end{tabular}
\end{center}
\end{table}

 1. There is only one maximal planar graph of order 6 whose minimal
 degree is 4. It has 4 different colorings, and its partitions of
 color group are shown as follow:

 $$ \{\{v_1\}\{v_2,v_6\}\{v_3,v_4\}\{v_5\}\}, \{\{v_1,v_6\}\{v_2\}\{v_3,v_4\}\{v_5\}\}$$
 $$ \{\{v_1,v_5\}\{v_2,v_6\}\{v_3\}\{v_4\}\}, \{\{v_1,v_5\}\{v_2,v_6\}\{v_3,v_4\}\}$$

 The following figures show the drawing of this graph and its four different
 4-colorings:

  \begin{center}
        \includegraphics [width=160pt]{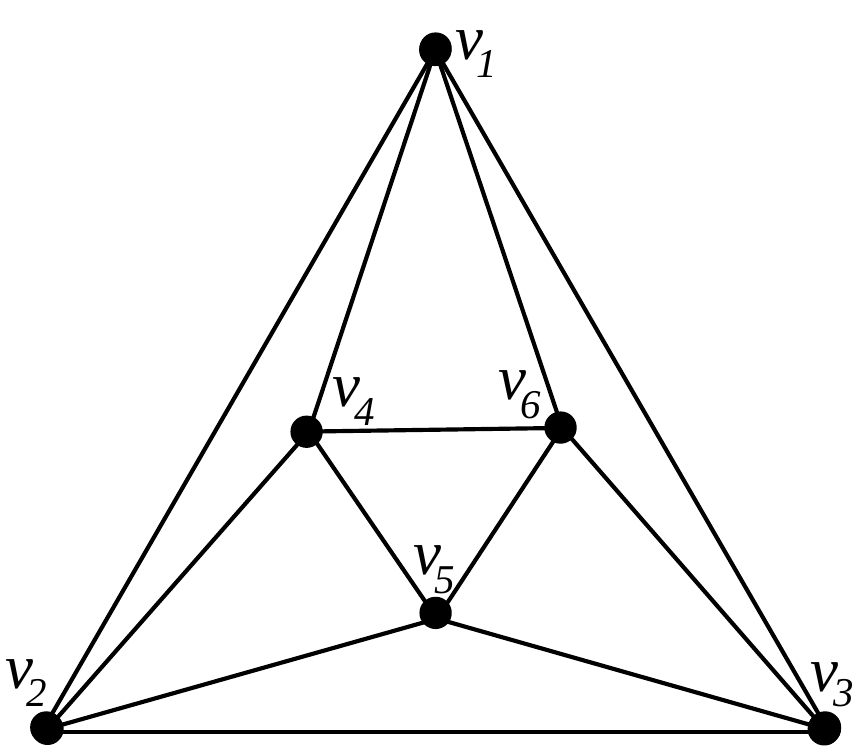}
         \vspace{5mm}

         \includegraphics [width=380pt]{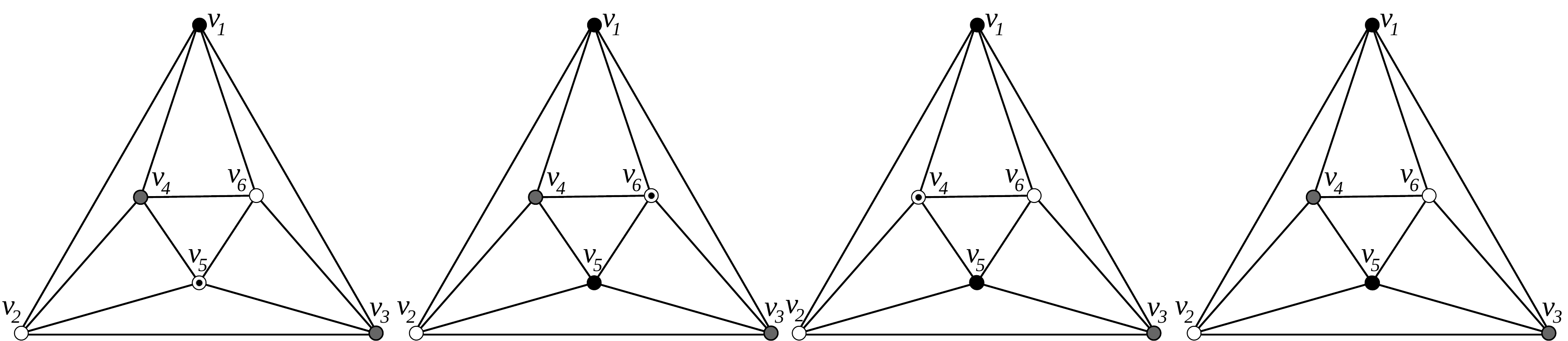}
  \end{center}

 2. There is only one maximal planar graph of order 7 whose minimal
 degree is 4. It has 5 different colorings.

 $$ \{\{v_1,v_7\}\{v_2,v_5\}\{v_3,v_4\}\{v_6\}\}, \{\{v_1,v_7\}\{v_2,v_6\}\{v_3,v_4\}\{v_5\}\}$$
 $$ \{\{v_1,v_7\}\{v_2,v_5\}\{v_3\}\{v_4,v_6\}\}, \{\{v_1,v_7\}\{v_2,v_6\}\{v_3,v_5\}\{v_4\}\}$$
 $$ \{\{v_1,v_7\}\{v_2\}\{v_3,v_5\}\{v_4,v_6\}\}$$

   \begin{center}
        \includegraphics [width=160pt]{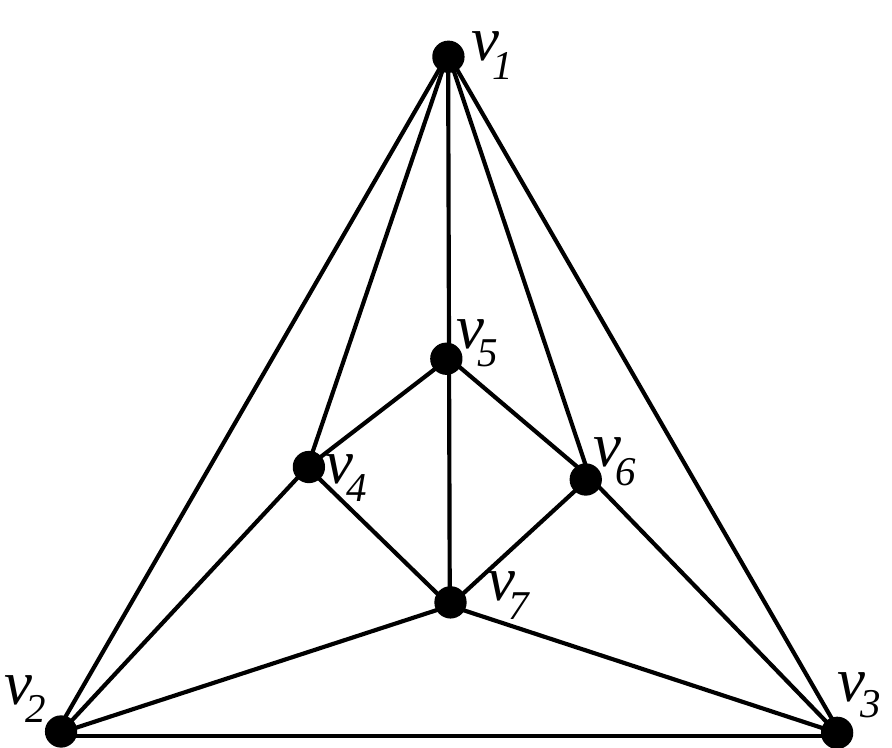}
         \vspace{5mm}

         \includegraphics [width=380pt]{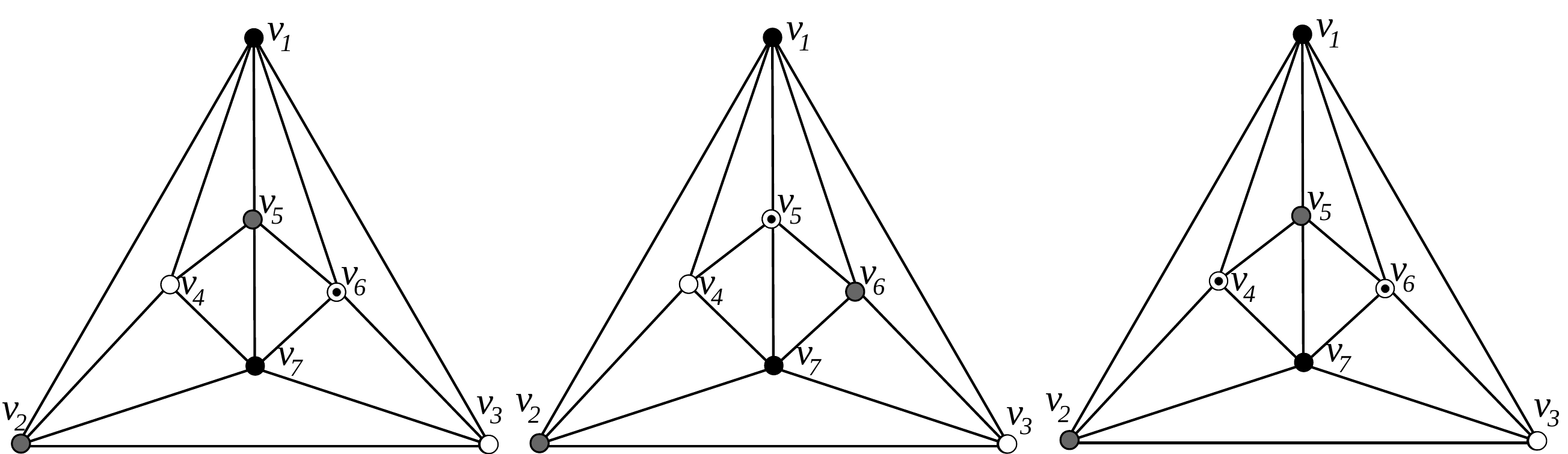}

         \vspace{8mm}
         \includegraphics [width=300pt]{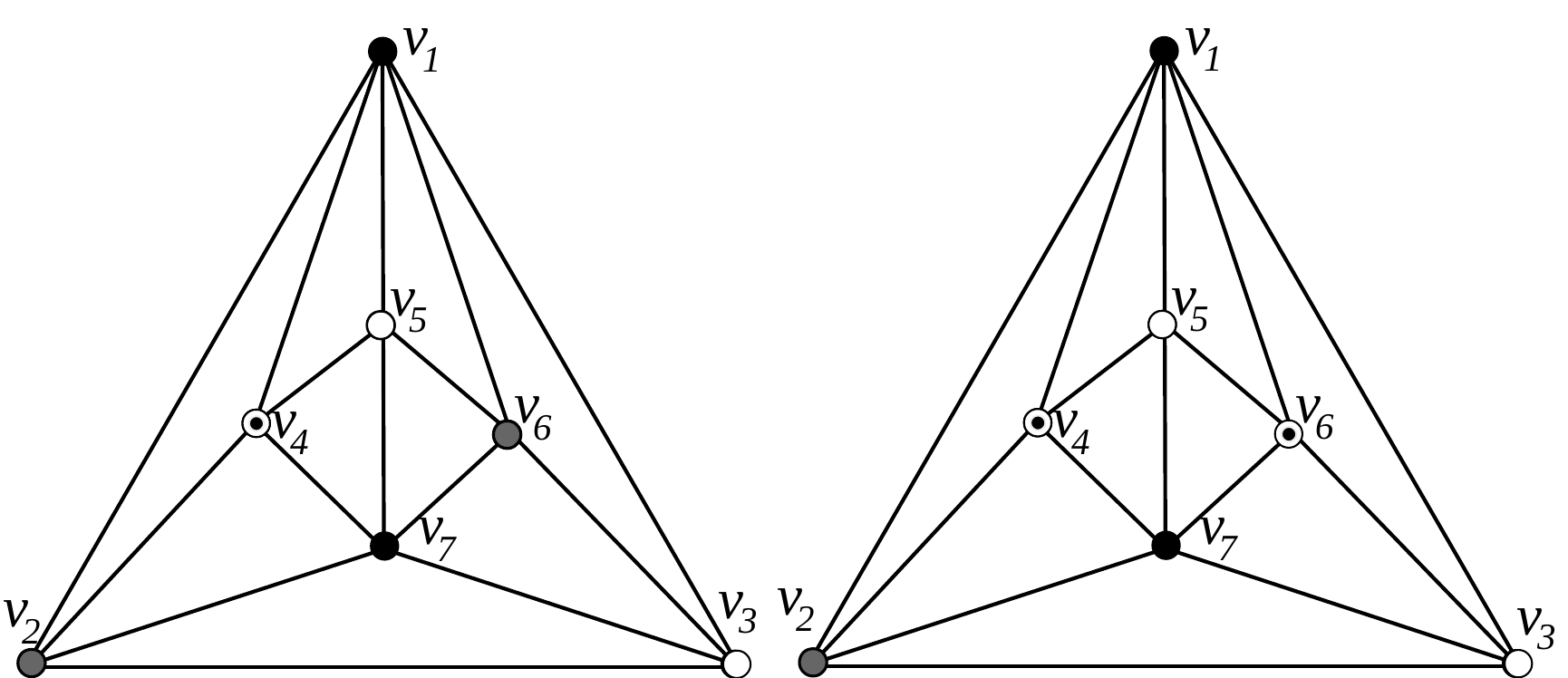}
  \end{center}

 3. There are two  maximal planar graphs of order 8 whose minimal
 degree are 4.
 \vspace{3mm}

 3.1 Degree sequence is 44444466, and it has 12 kinds of different colorings.
 \vspace{3mm}
 $$ \{\{v_1,v_8\}\{v_2,v_5,v_7\}\{v_3,v_4,v_6\}\}, \{\{v_1,v_8\} \{v_2,v_5\}\{v_3,v_4,v_6\}\{v_7\}\}$$
 $$ \{\{v_1,v_8\}\{v_2,v_5,v_7\}\{v_3,v_4\}\{v_6\}\}, \{\{v_1,v_8\} \{v_2,v_6\}\{v_3,v_4\}\{v_5,v_7\}\}$$
 $$ \{\{v_1,v_8\}\{v_2,v_7\}\{v_3,v_4,v_6\}\{v_5\}\}, \{\{v_1,v_8\} \{v_2\}\{v_3,v_4,v_6\}\{v_5,v_7\}\}$$
 $$ \{\{v_1,v_8\}\{v_2,v_5,v_7\}\{v_3,v_6\}\{v_4\}\}, \{\{v_1,v_8\} \{v_2,v_5\}\{v_3,v_6\}\{v_4,v_7\}\}$$
 $$ \{\{v_1,v_8\}\{v_2,v_5,v_7\}\{v_3,v_6\}\{v_4\}\}, \{\{v_1,v_8\} \{v_2,v_5\}\{v_3,v_6\}\{v_4,v_7\}\}$$
 $$ \{\{v_1,v_8\}\{v_2,v_5,v_7\}\{v_3\}\{v_4,v_6\}\}, \{\{v_1,v_8\} \{v_2,v_6\}\{v_3,v_5\}\{v_4,v_7\}\}$$
 $$ \{\{v_1,v_8\}\{v_2,v_7\}\{v_3,v_5\}\{v_4,v_6\}\}, \{\{v_1\} \{v_2,v_5,v_7\}\{v_3,v_4,v_6\}\{v_8\}\}$$
   \begin{center}
        \includegraphics [width=160pt]{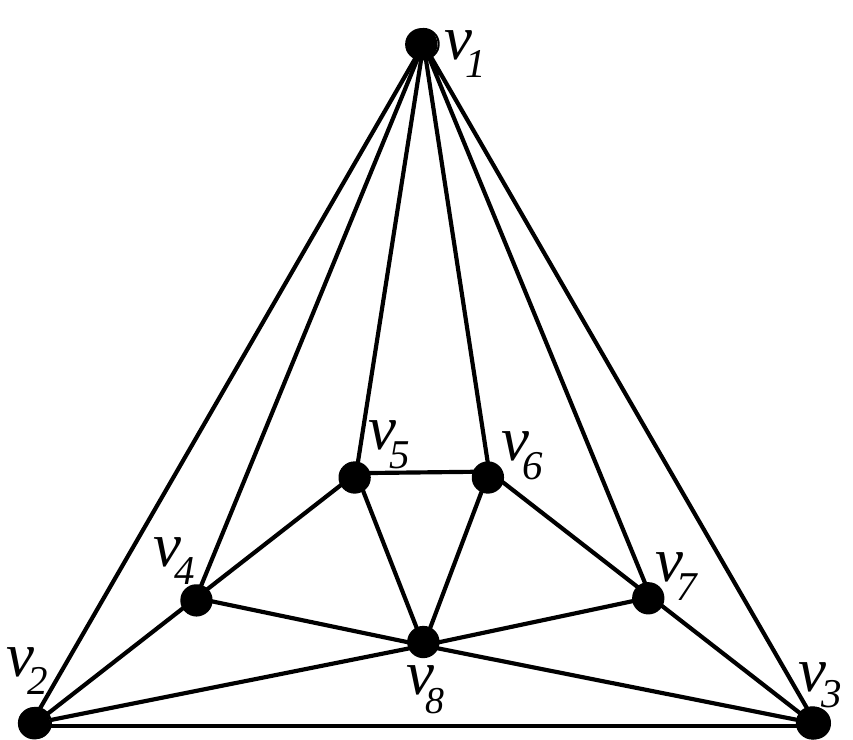}
         \vspace{5mm}

         \includegraphics [width=380pt]{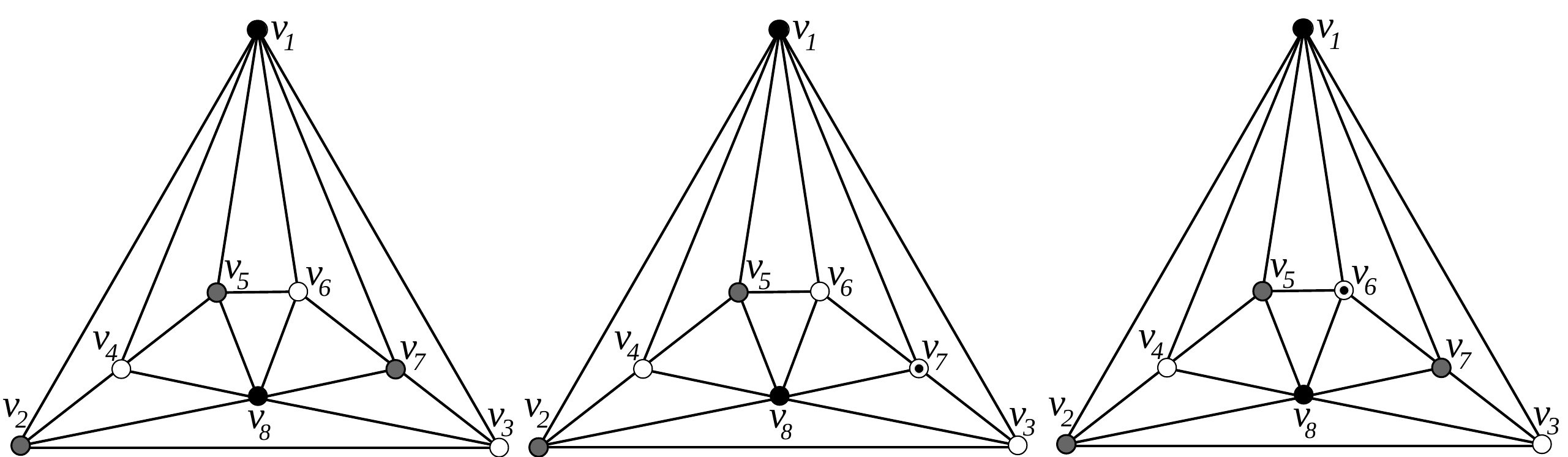}

         \vspace{5mm}
         \includegraphics [width=380pt]{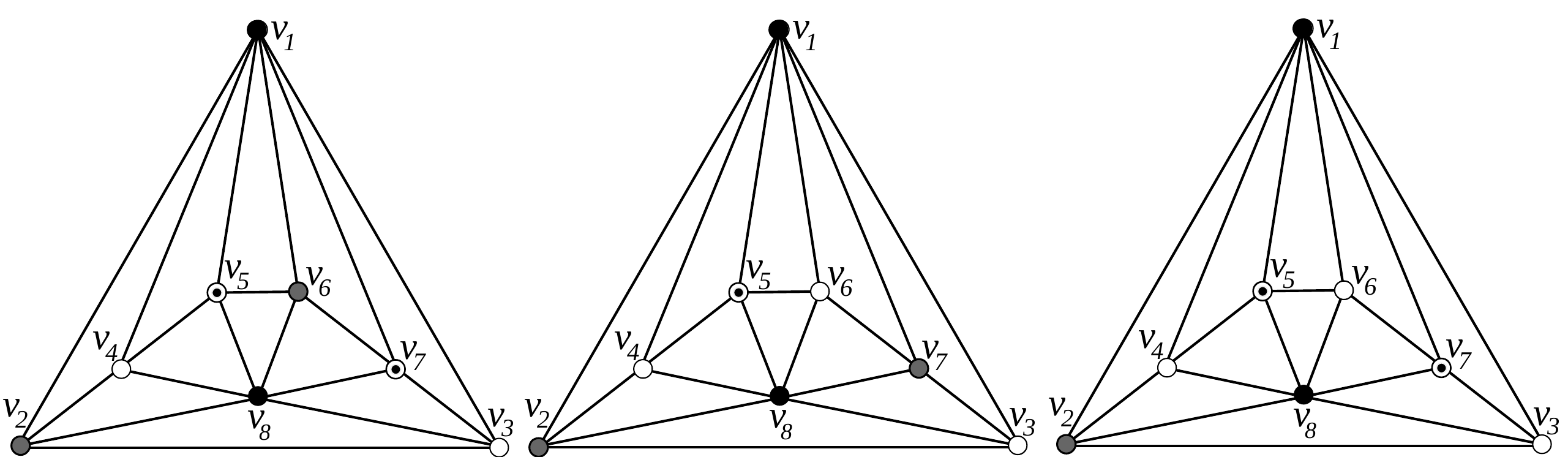}

          \vspace{5mm}
         \includegraphics [width=380pt]{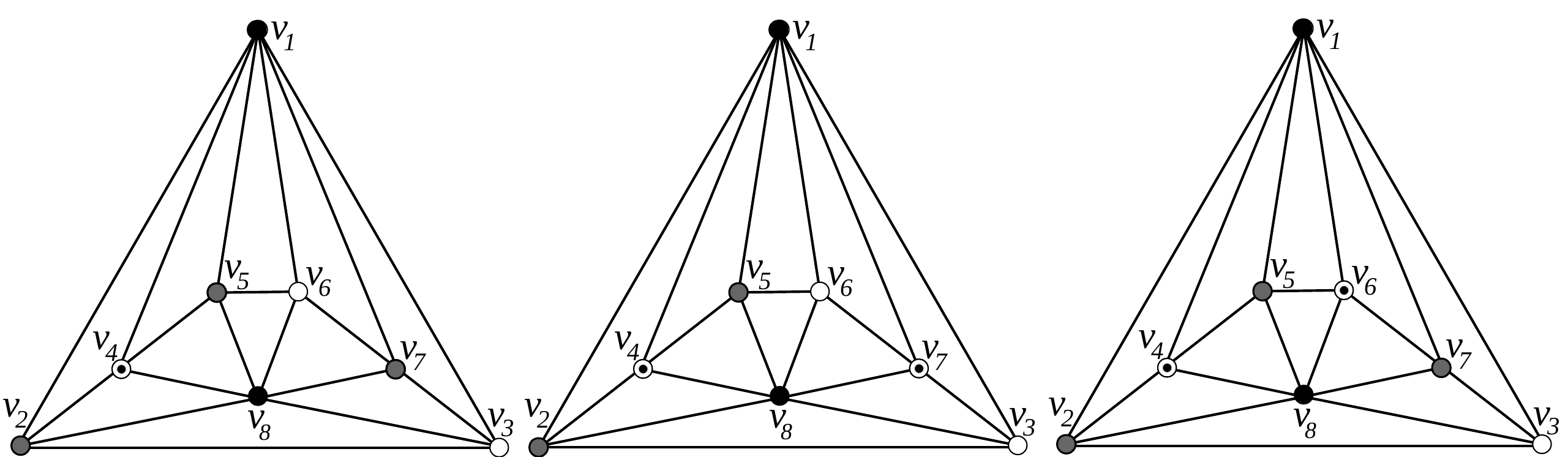}

         \vspace{5mm}
         \includegraphics [width=380pt]{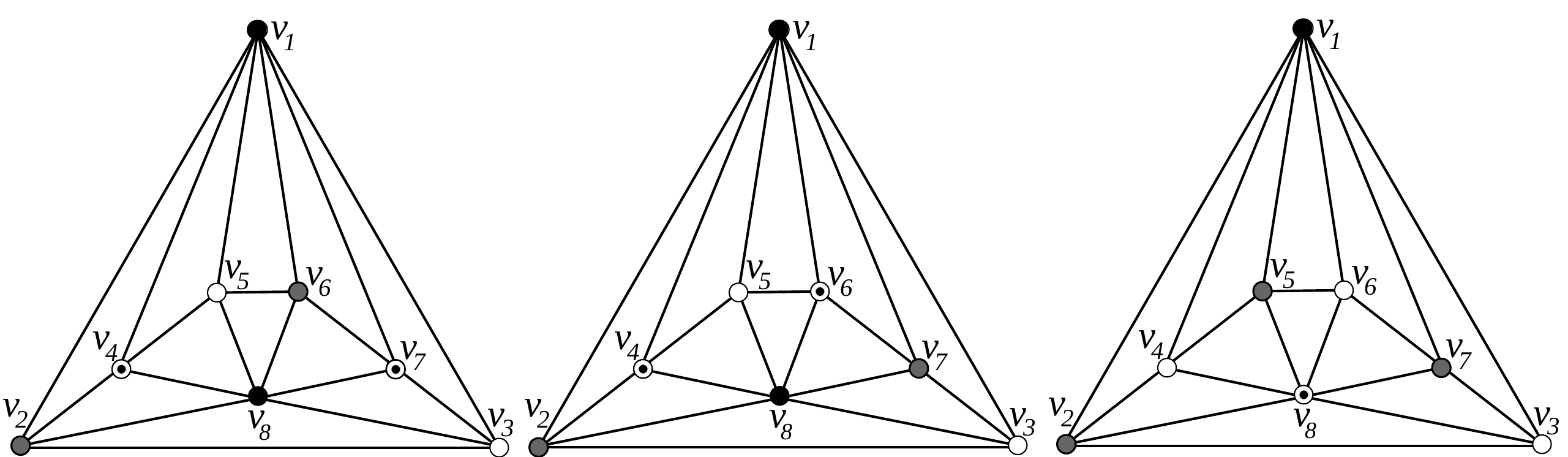}

         \vspace{8mm}
  \end{center}
 3.2 Degree sequence is 44445555, and it has 3 kinds of different colorings.

 $$ \{\{v_1,v_8\}\{v_2,v_4\}\{v_3,v_6\}\{v_5,v_7\}\}, \{\{v_1,v_8\} \{v_2,v_6\}\{v_3,v_4\}\{v_5,v_7\}\}$$
 $$ \{\{v_1,v_6\}\{v_2,v_7\}\{v_3,v_5\}\{v_6,v_8\}\}$$
 \begin{center}
        \includegraphics [width=160pt]{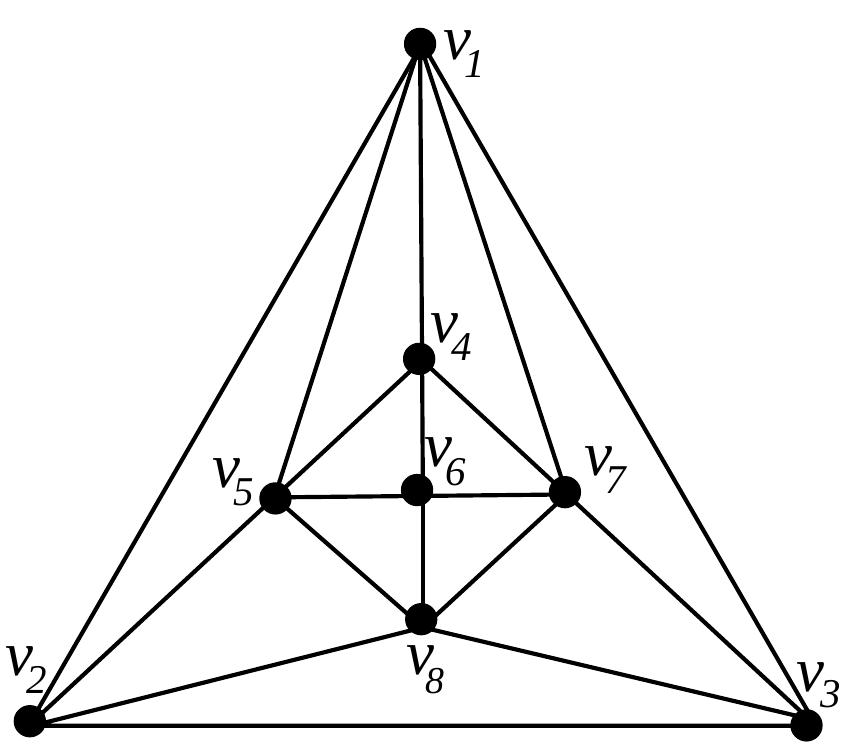}
         \vspace{5mm}

         \includegraphics [width=380pt]{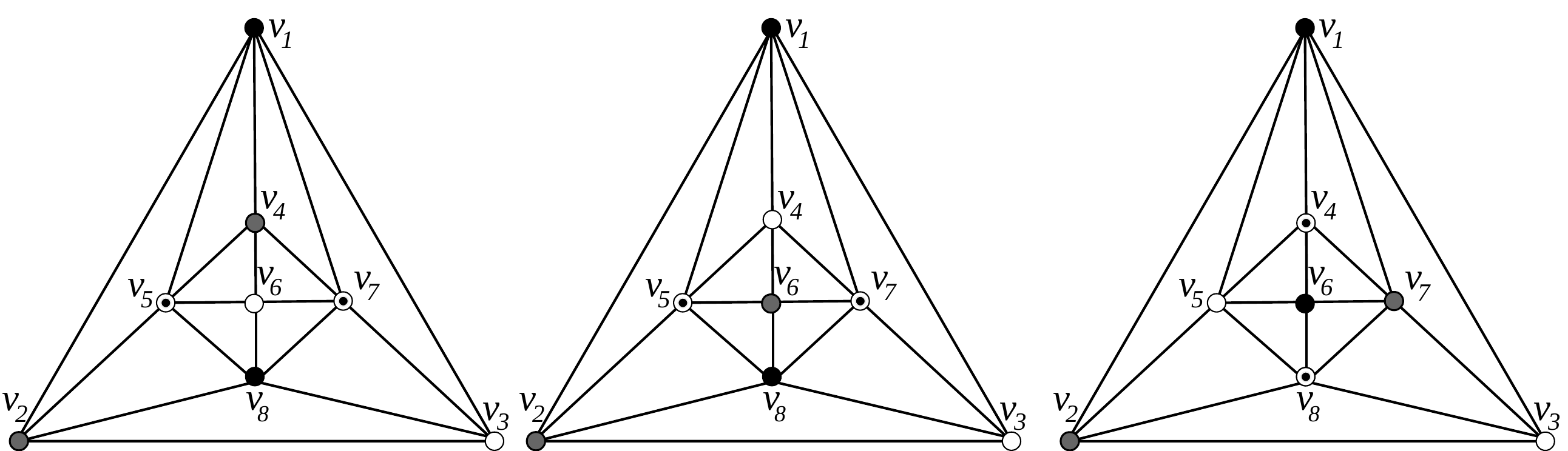}

         \vspace{8mm}
  \end{center}
 4. There are five  maximal planar graphs of order 9 whose minimal
 degree is 4.

 4.1 Degree sequence is 444444666, and it is 3-colorable.  The unique partitions of
 color group are shown as follow:
 \begin{center}
        \includegraphics [width=320pt]{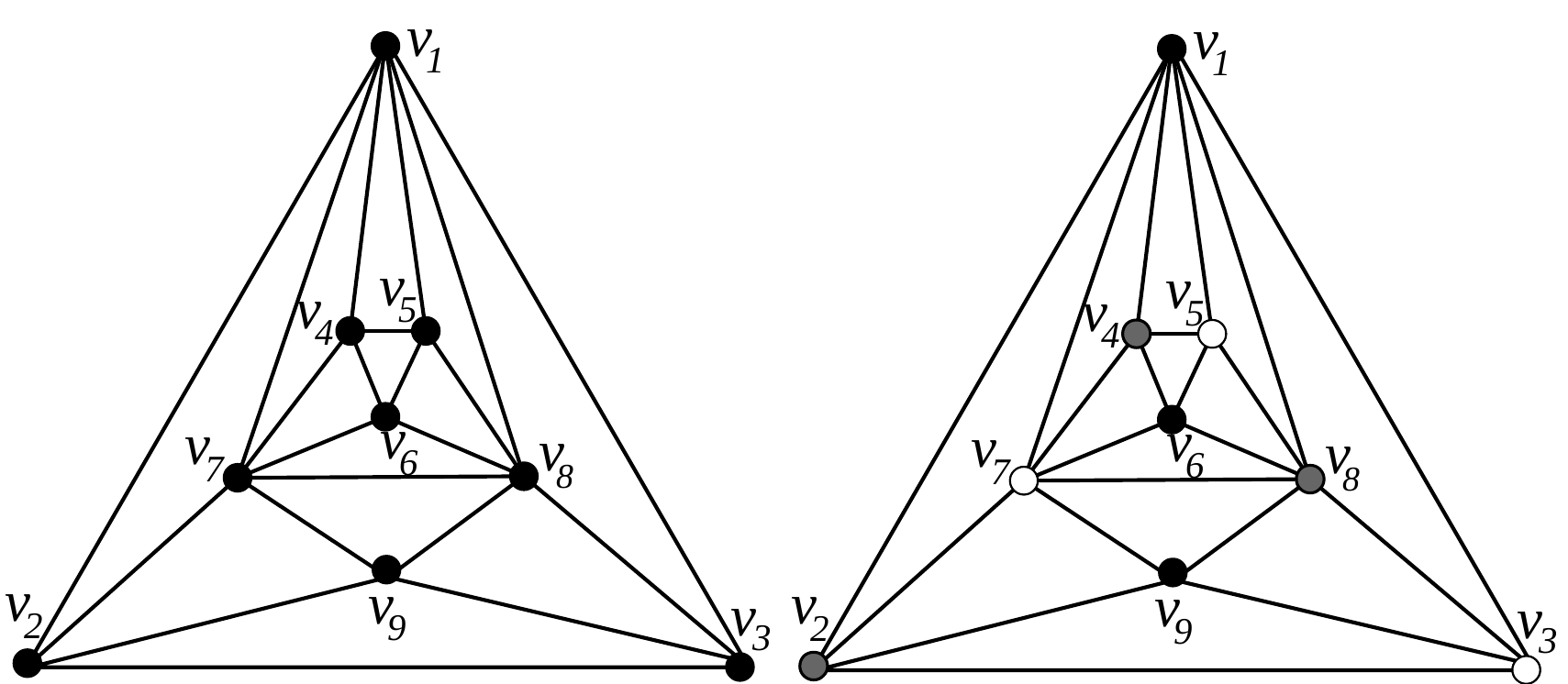}

         \vspace{8mm}
  \end{center}
 4.2 Degree sequence is 444455556, and it has 6 kinds of different colorings.

  $$ \{\{v_1,v_7,v_8\}\{v_2,v_6\}\{v_4,v_5,v_9\}\{v_3\}\}, \{\{v_1,v_7,v_8\} \{v_2\}\{v_4,v_5,v_9\}\{v_3,v_6\}\}$$
  $$ \{\{v_1,v_8\}\{v_2,v_6\}\{v_4,v_5,v_9\}\{v_3,v_7\}\}, \{\{v_1,v_7\} \{v_2,v_8\}\{v_3,v_6\}\{v_4,v_5,v_9\}\}$$
  $$ \{\{v_1,v_7,v_8\}\{v_2,v_6\}\{v_3,v_4\}\{v_5,v_9\}\}, \{\{v_1,v_7,v_8\} \{v_2,v_5\}\{v_3,v_6\}\{v_4,v_9\}\}$$
 \begin{center}
        \includegraphics [width=160pt]{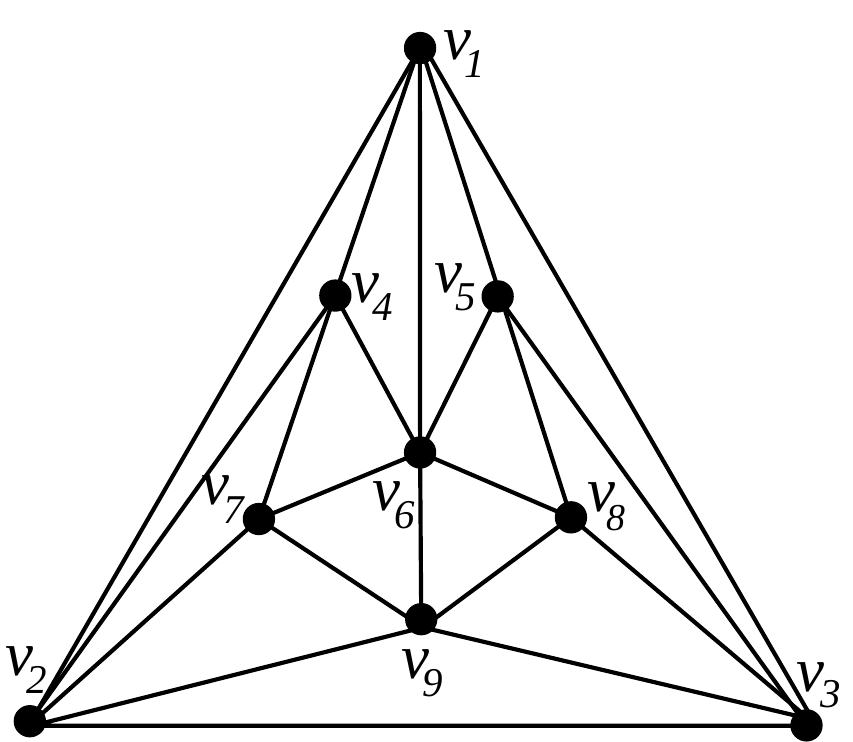}

         \vspace{5mm}
         \includegraphics [width=380pt]{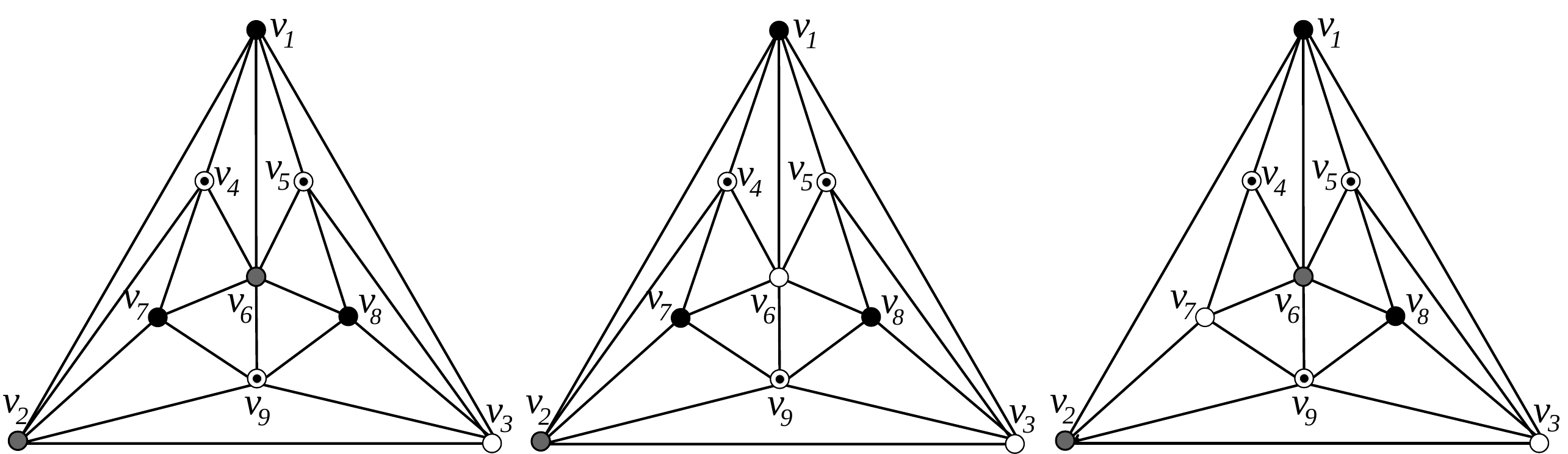}

         \vspace{5mm}
         \includegraphics [width=380pt]{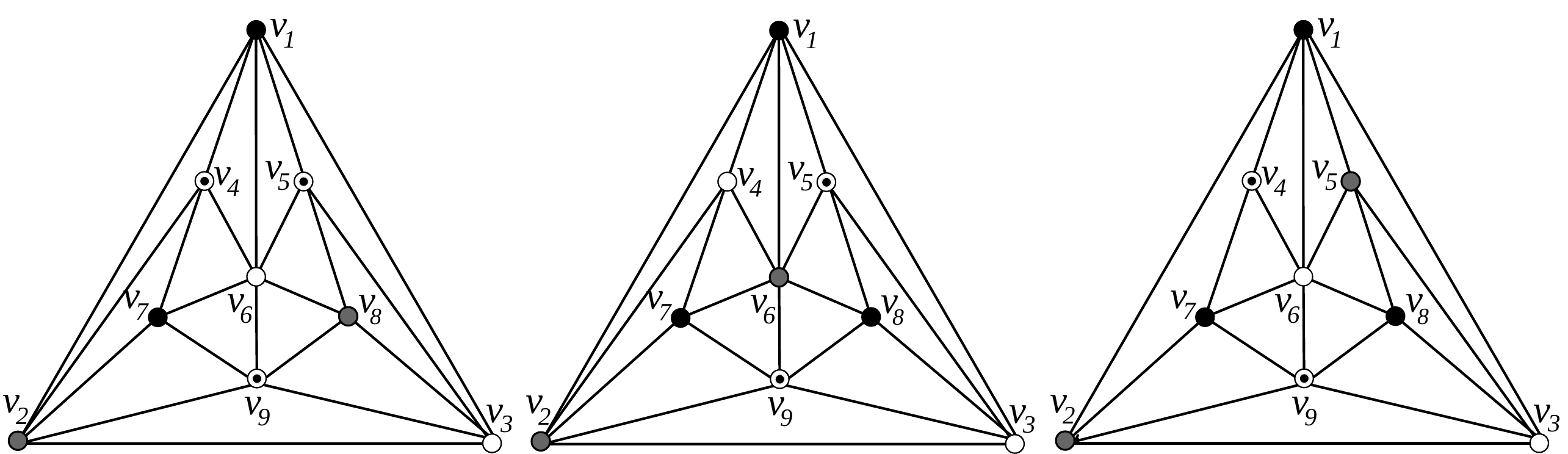}

         \vspace{8mm}
  \end{center}
 4.3 Degree sequence is 444555555, and it has 2 kinds of different colorings.

 $$ \{\{v_1,v_7\}\{v_2,v_8\}\{v_3,v_4\}\{v_5,v_6,v_9\}\}, \{\{v_1,v_8\} \{v_2,v_4\}\{v_3,v_7\}\{v_5,v_6,v_9\}\}$$
 \begin{center}
        \includegraphics [width=160pt]{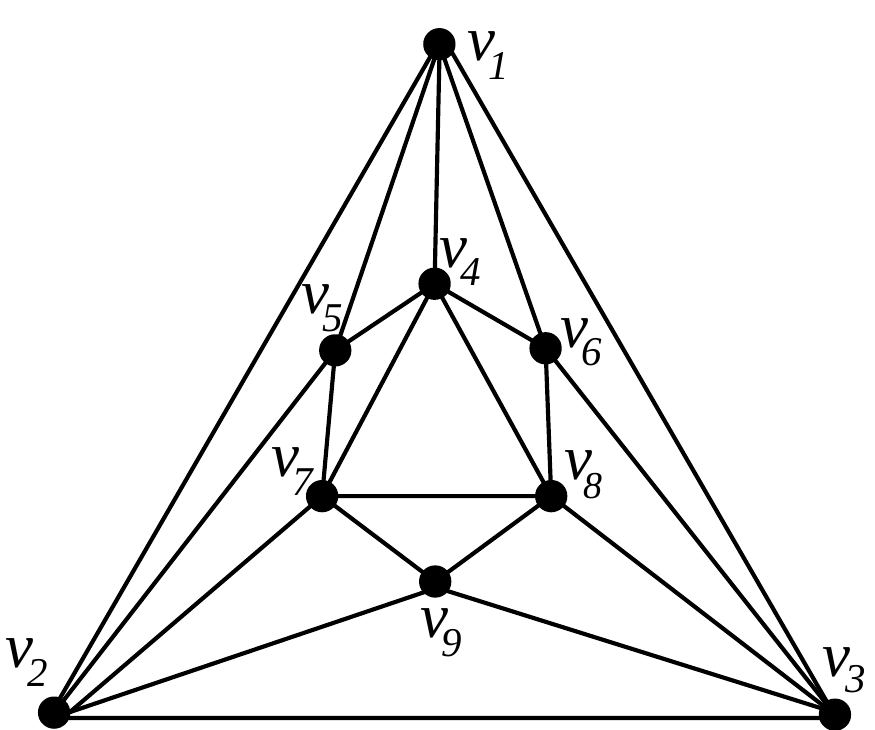}

         \vspace{5mm}
         \includegraphics [width=320pt]{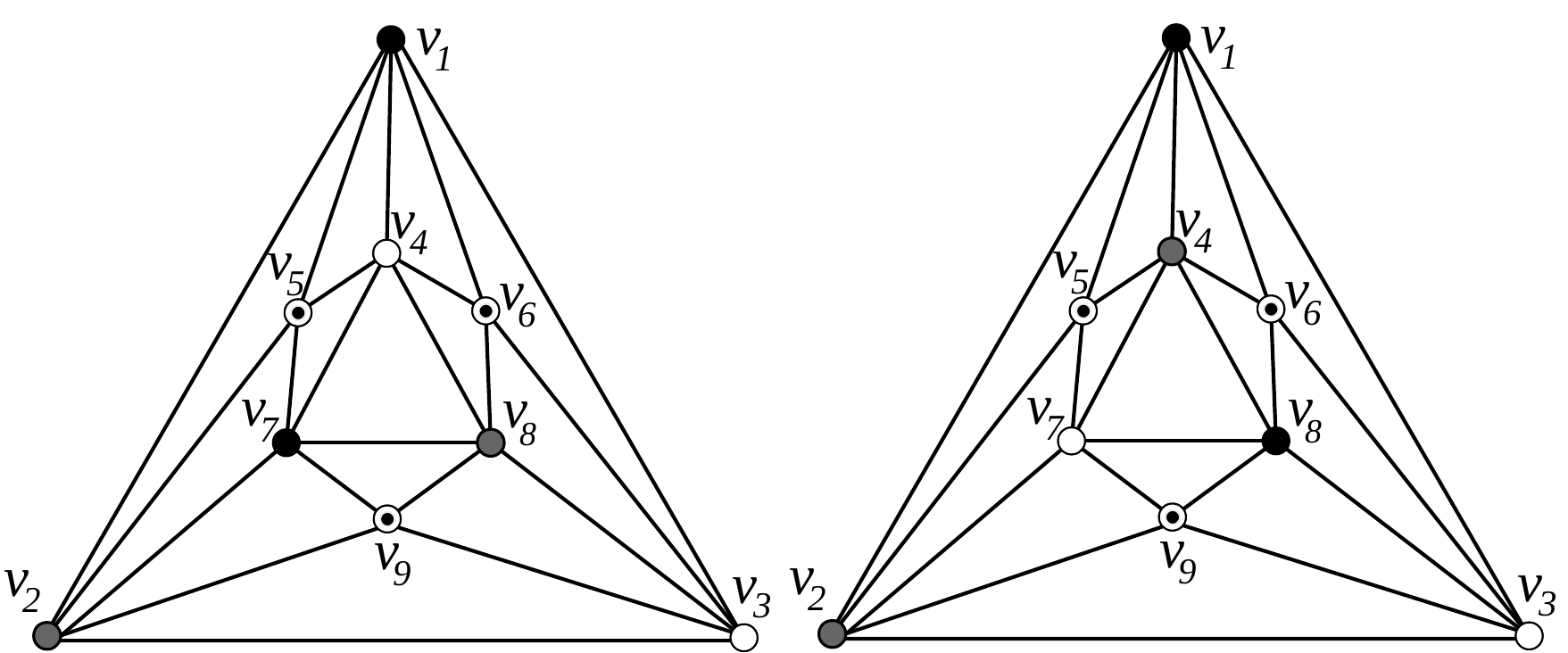}

         \vspace{8mm}
  \end{center}
 4.4 Degree sequence is 4444444477, and it has 17 kinds of different colorings.
 \vspace{3mm}

  $$ \{\{v_1,v_9\}\{v_2,v_5,v_8\}\{v_3,v_6,v_7\}\{v_4\}\}, \{\{v_1,v_9\} \{v_2,v_4,v_8\}\{v_3,v_6,v_7\}\{v_5\}\}$$
  $$ \{\{v_1,v_9\}\{v_2,v_5,v_8\}\{v_3,v_4,v_7\}\{v_6\}\}, \{\{v_1,v_9\} \{v_2,v_5,v_6\}\{v_3,v_4\}\{v_7,v_8\}\}$$
  $$ \{\{v_1,v_9\}\{v_2,v_4\}\{v_3,v_5,v_6\}\{v_7,v_8\}\}, \{\{v_1,v_9\} \{v_2,v_8\}\{v_3,v_4,v_7\}\{v_5,v_6\}\}$$
  $$ \{\{v_1,v_9\}\{v_2,v_4,v_8\}\{v_3,v_7\}\{v_5,v_6\}\}, \{\{v_1,v_9\} \{v_2\}\{v_3,v_5,v_6\}\{v_4,v_7,v_8\}\}$$
  $$ \{\{v_1,v_9\}\{v_2,v_5,v_6\}\{v_3\}\{v_4,v_7,v_8\}\}, \{\{v_1,v_9\} \{v_2,v_6\}\{v_3,v_7\}\{v_4,v_7,v_8\}\}$$
  $$ \{\{v_1,v_9\}\{v_2,v_5\}\{v_3,v_6\}\{v_4,v_7,v_8\}\}, \{\{v_1,v_9\} \{v_2,v_5\}\{v_3,v_6,v_7\}\{v_4,v_8\}\}$$
  $$ \{\{v_1,v_9\}\{v_2,v_5,v_6\}\{v_3,v_7\}\{v_4,v_8\}\}, \{\{v_1,v_9\} \{v_2,v_5,v_8\}\{v_3,v_6\}\{v_4,v_7\}\}$$
  $$ \{\{v_1,v_9\}\{v_2,v_8\}\{v_3,v_5,v_6\}\{v_4,v_7\}\}, \{\{v_1,v_9\} \{v_2,v_5,v_6\}\{v_3,v_4,v_7\}\{v_8\}\}$$
  $$ \{\{v_1,v_9\}\{v_2,v_4,v_8\}\{v_3,v_5,v_6\}\{v_7\}\}$$
 \begin{center}
        \includegraphics [width=160pt]{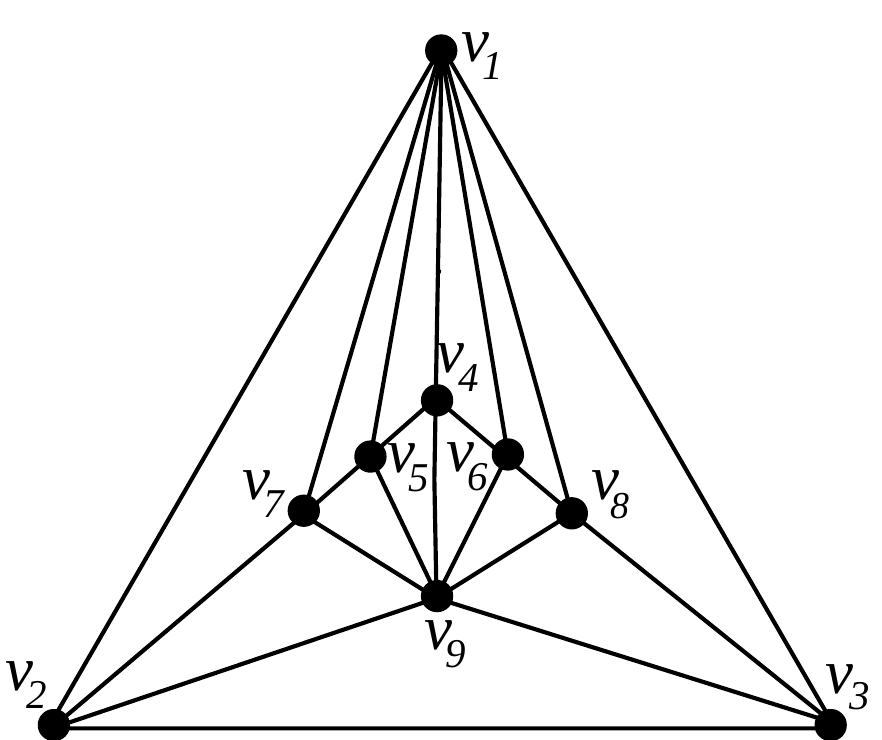}

         \vspace{5mm}
         \includegraphics [width=380pt]{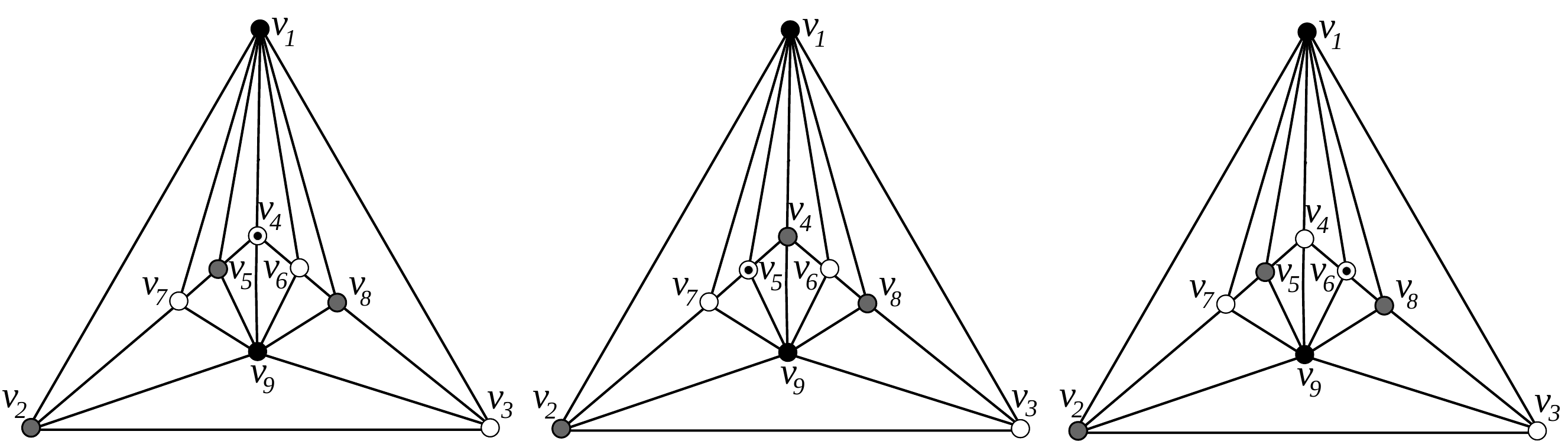}

         \vspace{5mm}
         \includegraphics [width=380pt]{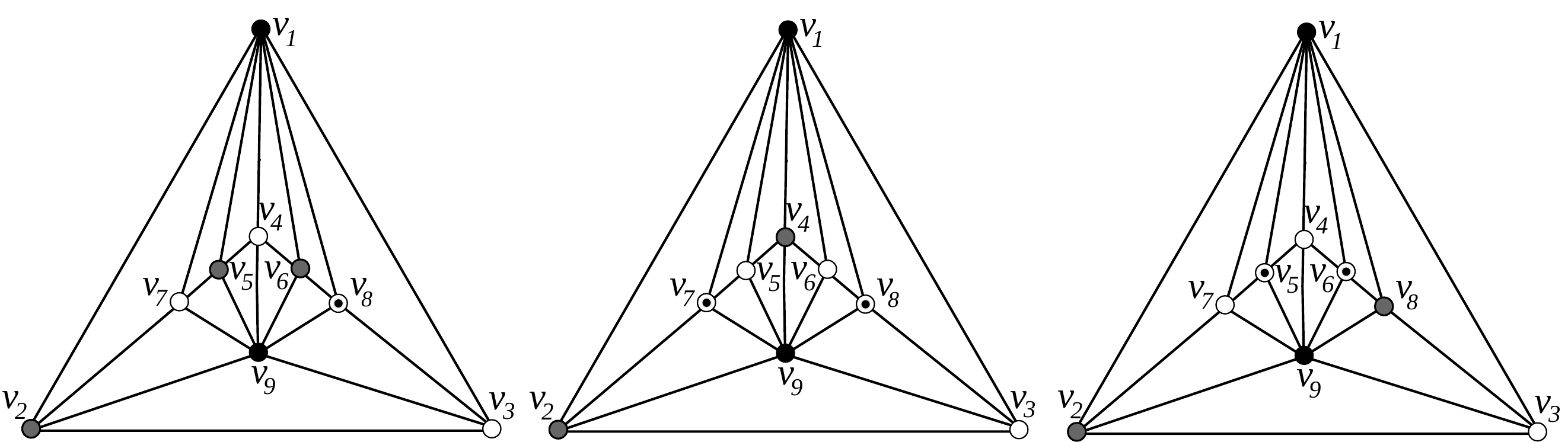}

         \vspace{5mm}
         \includegraphics [width=380pt]{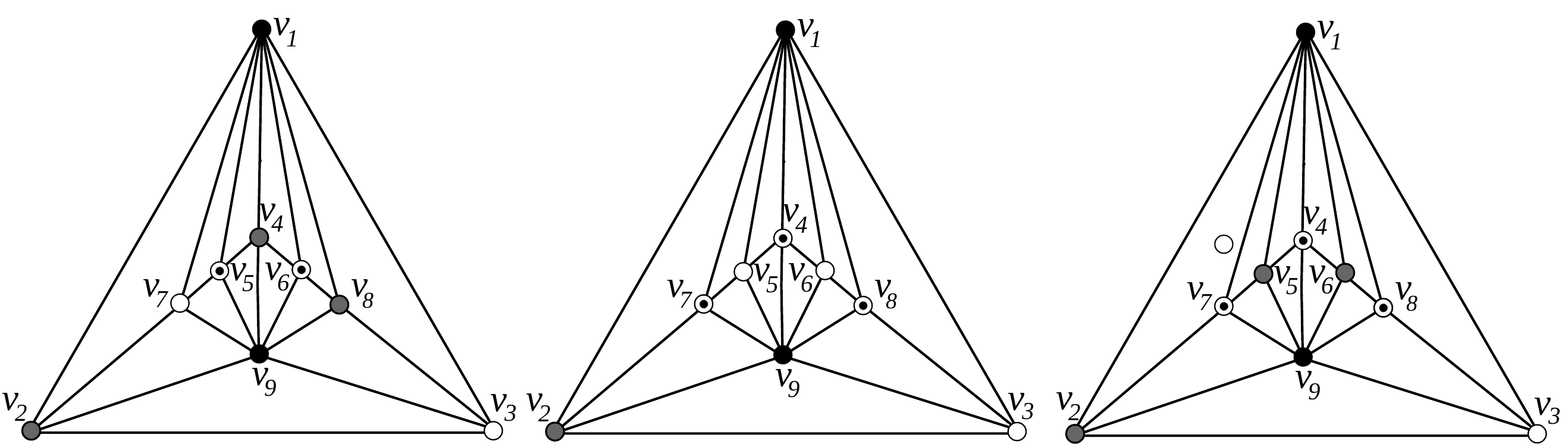}

         \vspace{5mm}
         \includegraphics [width=380pt]{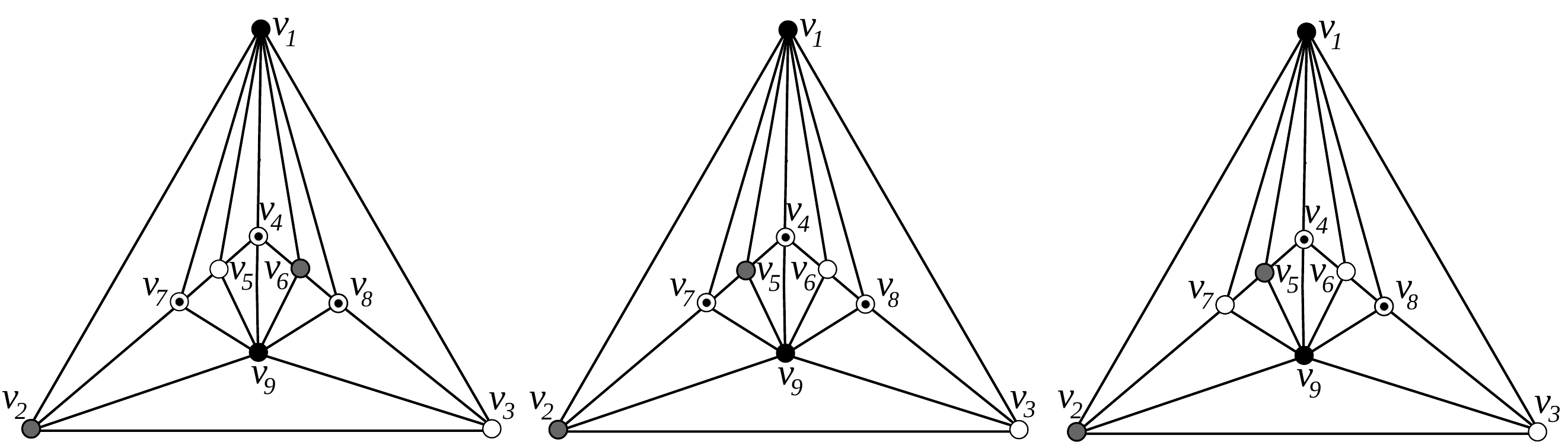}

         \vspace{5mm}
         \includegraphics [width=380pt]{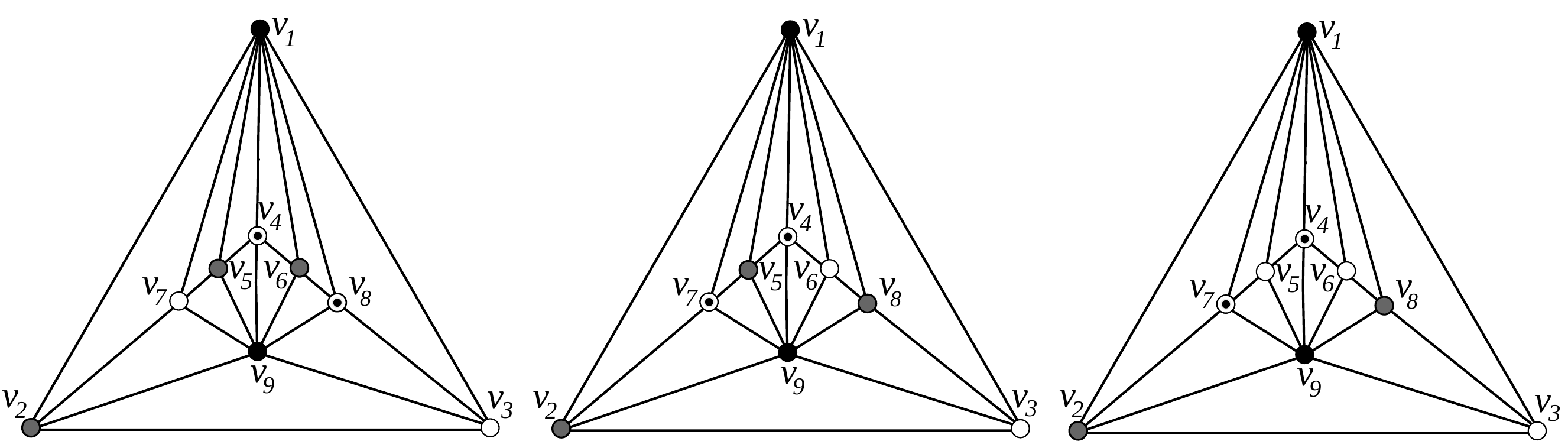}

         \vspace{5mm}
         \includegraphics [width=320pt]{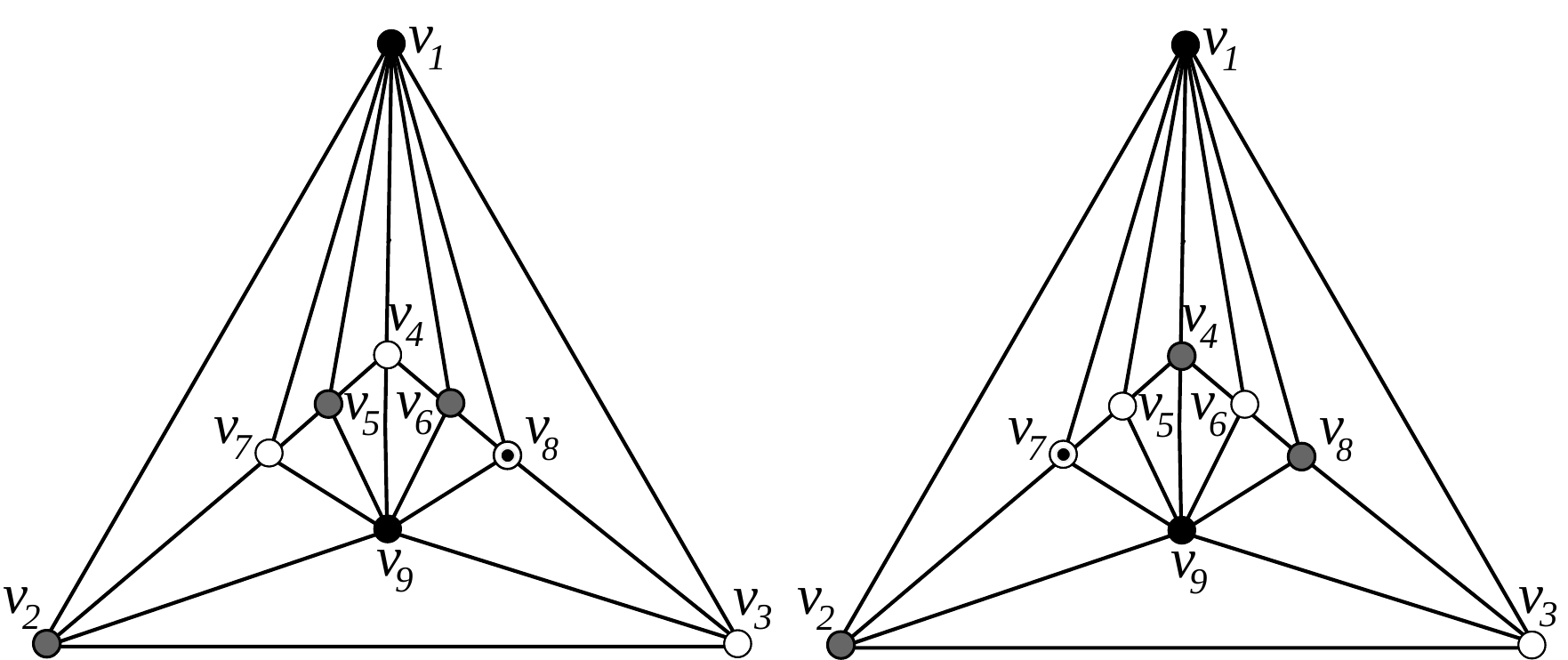}

         \vspace{8mm}
  \end{center}
 4.5 Degree sequence is 4444445566, and it has 7 kinds of different colorings.
 \vspace{3mm}

  $$ \{\{v_1,v_5,v_9\}\{v_2,v_8\}\{v_3,v_6\}\{v_4,v_7\}\}, \{\{v_1,v_5,v_9\} \{v_2,v_7\}\{v_3,v_4\}\{v_6,v_8\}\}$$
  $$ \{\{v_1,v_5,v_9\}\{v_2,v_4\}\{v_3,v_7\}\{v_6,v_8\}\}, \{\{v_1,v_5,v_9\} \{v_2\}\{v_3,v_4,v_7\}\{v_6,v_8\}\}$$
  $$ \{\{v_1,v_5,v_9\}\{v_2,v_4,v_7\}\{v_3\}\{v_6,v_8\}\}, \{\{v_1,v_5,v_9\} \{v_2,v_8\}\{v_6\}\{v_3,v_4,v_7\}\}$$
  $$ \{\{v_1,v_5,v_9\}\{v_2,v_4,v_7\}\{v_8\}\{v_3,v_6\}\}$$
 \begin{center}
        \includegraphics [width=160pt]{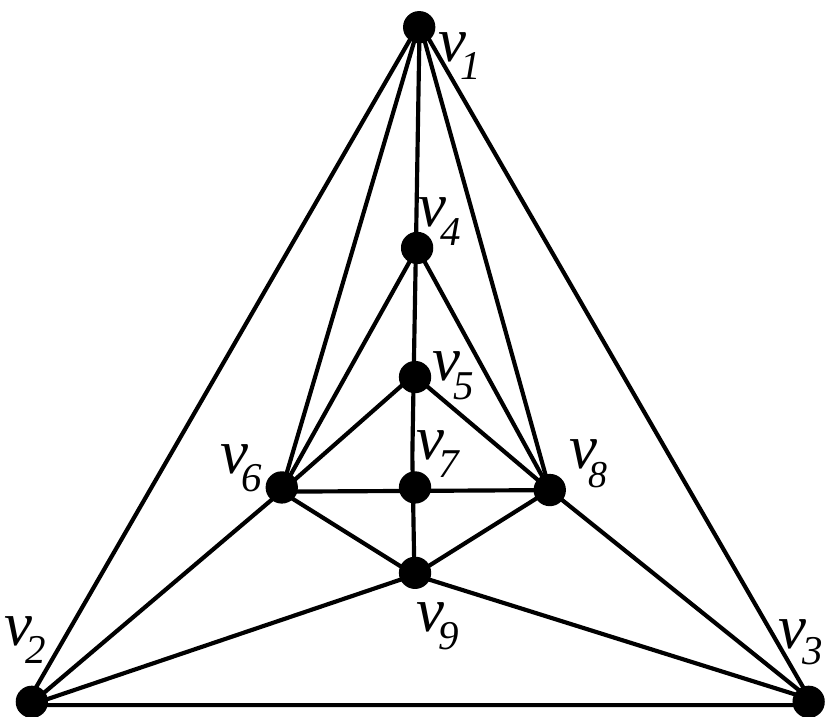}

         \vspace{5mm}
         \includegraphics [width=380pt]{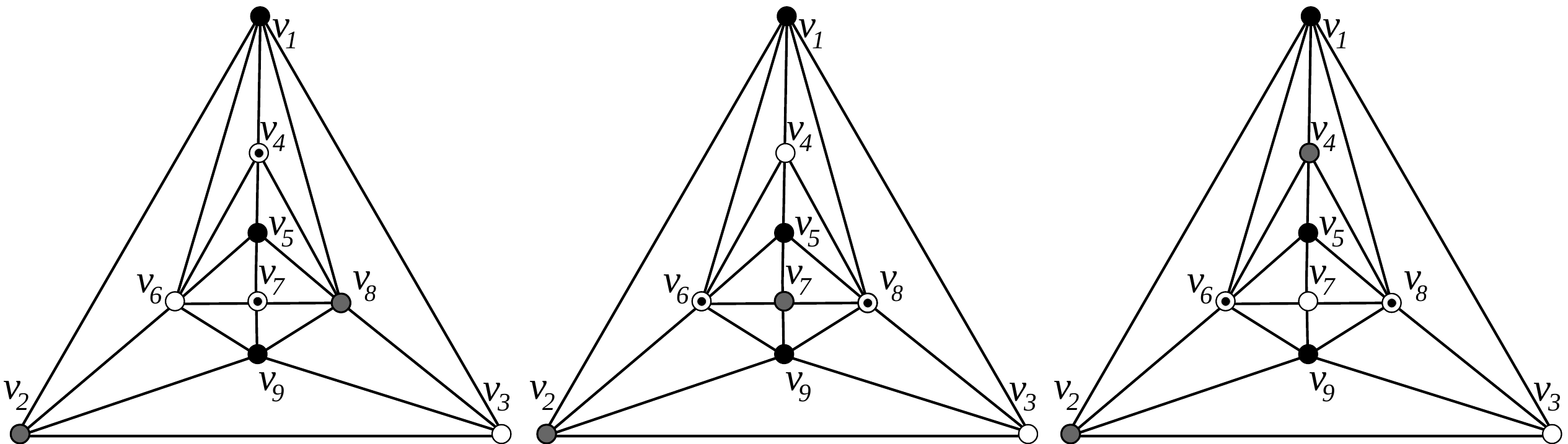}

         \vspace{5mm}
         \includegraphics [width=380pt]{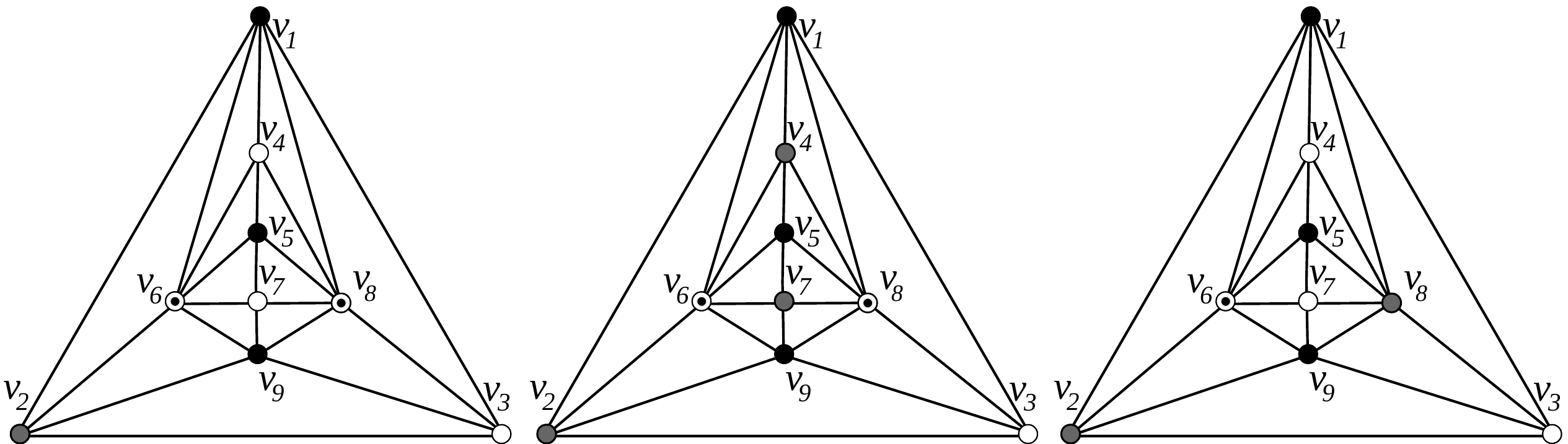}

         \vspace{5mm}
         \includegraphics [width=160pt]{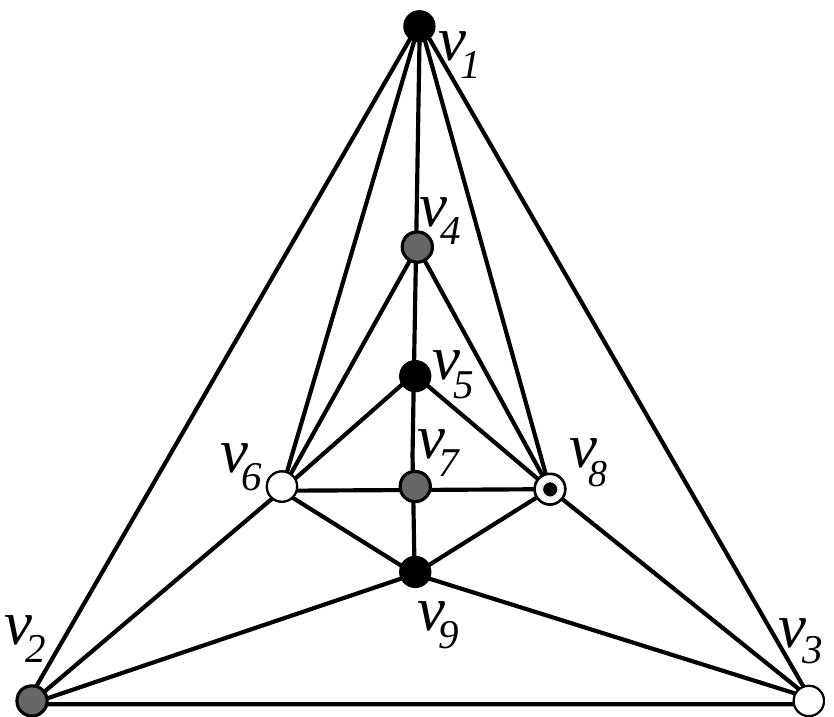}

         \vspace{8mm}
  \end{center}
 5. There are 13 maximal planar graphs of order 10 whose minimal
 degree is 4.

 5.1 Degree sequence is 4444455567, and it has 7 kinds of different colorings.
    $$\{\{v_1,v_9\}\{v_2,v_4\}\{v_5,v_8,v_{10}\}\{v_3,v_6,v_7\}\} \{\{v_1,v_9\}\{v_2,v_6,v_7\}\{v_5,v_8,v_{10}\}\{v_3,v_4\}\} $$ $$\{\{v_1,v_6,v_7\}\{v_2,v_5,v_{10}\}\{v_3,v_8\}\{v_4,v_9\}\} \{\{v_1,v_6,v_7\}\{v_2,v_{10}\}\{v_3,v_5,v_8\}\{v_4,v_9\}\} $$ $$\{\{v_1,v_6\}\{v_2,v_5,v_{10}\}\{v_3,v_7,v_8\}\{v_4,v_9\}\} \{\{v_1,v_7\}\{v_2,v_6,v_{10}\}\{v_3,v_5,v_8\}\{v_4,v_9\}\}$$ $$\{\{v_1,v_5\}\{v_2,v_6,v_{10}\}\{v_3,v_7,v_8\}\{v_4,v_9\}\}$$

     \begin{center}
        \includegraphics [width=160pt]{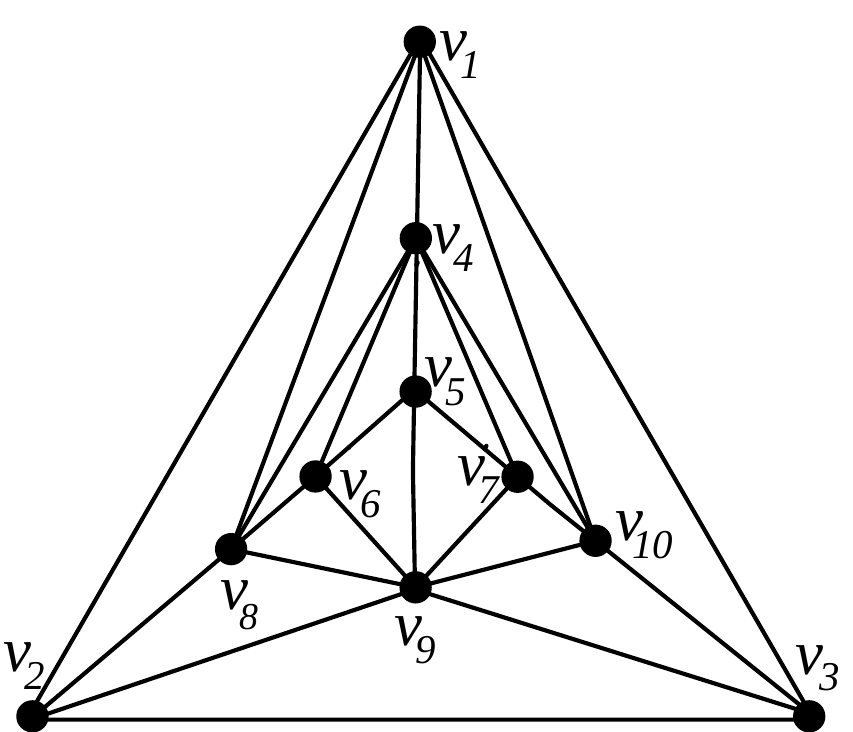}

         \vspace{5mm}
         \includegraphics [width=380pt]{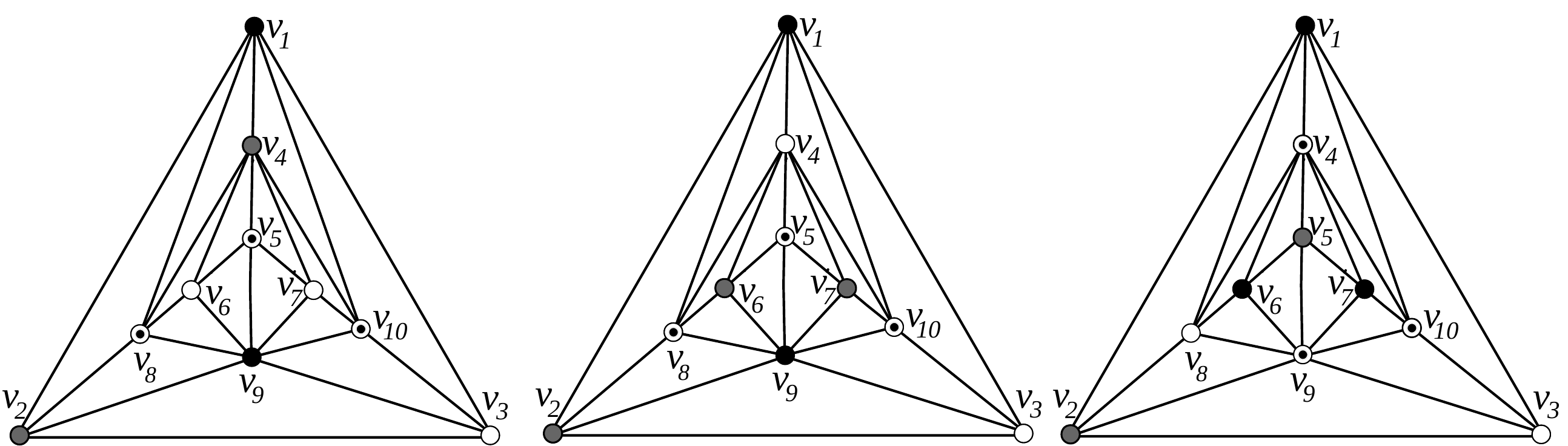}

         \vspace{5mm}
         \includegraphics [width=380pt]{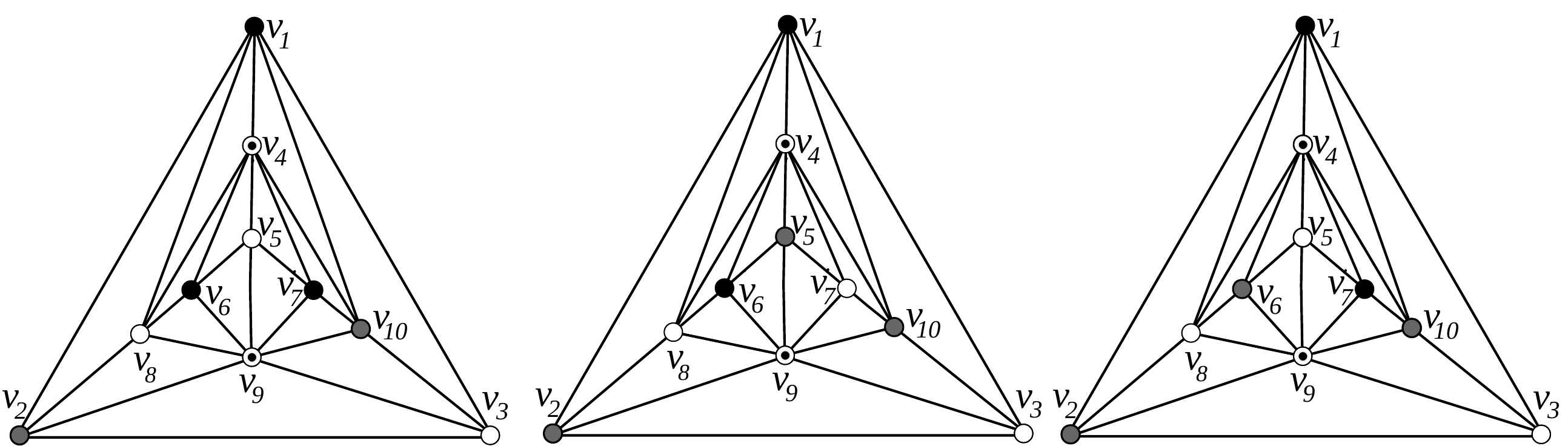}

         \vspace{5mm}
         \includegraphics [width=160pt]{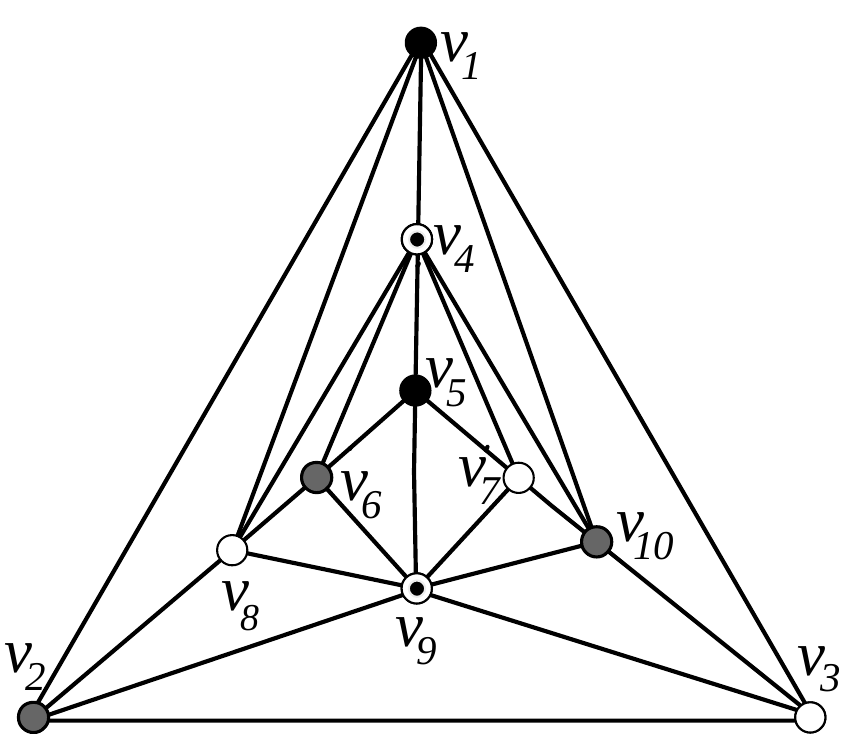}

         \vspace{8mm}
  \end{center}

 5.2 Degree sequence is 4444445577, and it has 10 kinds of different colorings.
    $$\{\{v_1,v_7,v_{10}\}\{v_2,v_4\}\{v_3,v_5,v_9\}\{v_6,v_8\}\} \{\{v_1,v_7,v_{10}\}\{v_2,v_4,v_9\}\{v_3,v_5\}\{v_6,v_8\}\}$$ $$\{\{v_1,v_{10}\}\{v_2,v_4,v_7\}\{v_3,v_5,v_9\}\{v_6,v_8\}\} \{\{v_1,v_5,v_{10}\}\{v_2,v_4,v_7\}\{v_3,v_9\}\{v_6,v_8\}\}$$ $$\{\{v_1,v_5,v_{10}\}\{v_2,v_7\}\{v_3,v_4,v_9\}\{v_6,v_8\}\} \{\{v_1,v_5,v_{10}\}\{v_2,v_9\}\{v_3,v_4,v_7\}\{v_6,v_8\}\}$$ $$\{\{v_1,v_7,v_{10}\}\{v_2,v_5\}\{v_3,v_4,v_9\}\{v_6,v_8\}\} \{\{v_1,v_7,v_{10}\}\{v_2,v_5,v_9\}\{v_3,v_4\}\{v_6,v_8\}\}$$ $$\{\{v_1,v_{10}\}\{v_2,v_5,v_9\}\{v_3,v_4,v_7\}\{v_6,v_8\}\} \{\{v_1,v_5,v_9\}\{v_2,v_8\}\{v_3,v_6\}\{v_4,v_5,v_{10}\}\}$$
     \begin{center}
        \includegraphics [width=160pt]{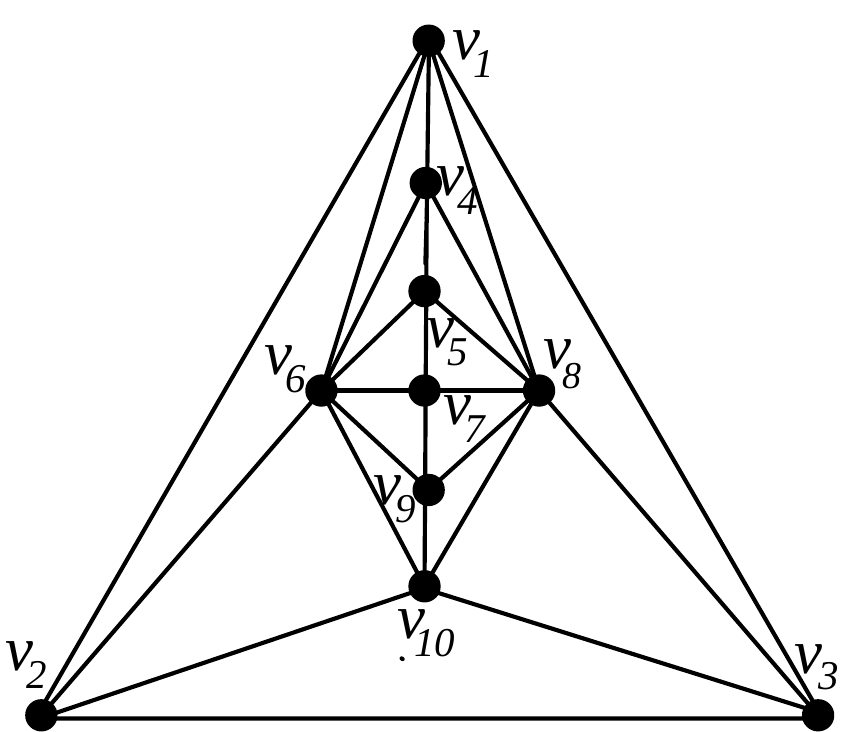}

         \vspace{5mm}
         \includegraphics [width=380pt]{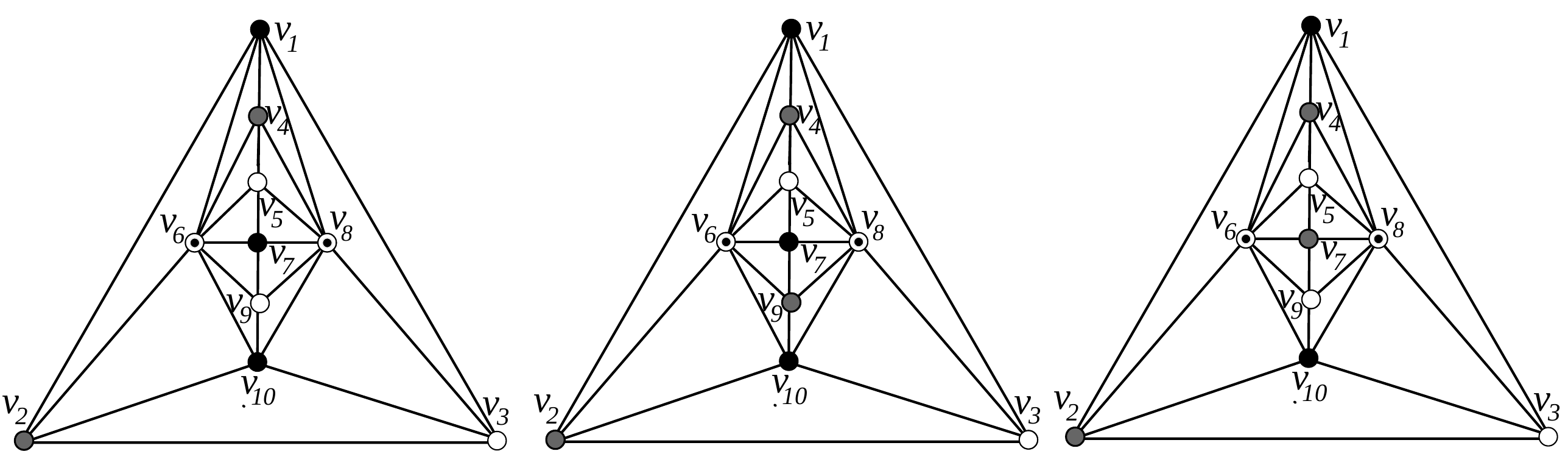}

         \vspace{5mm}
         \includegraphics [width=380pt]{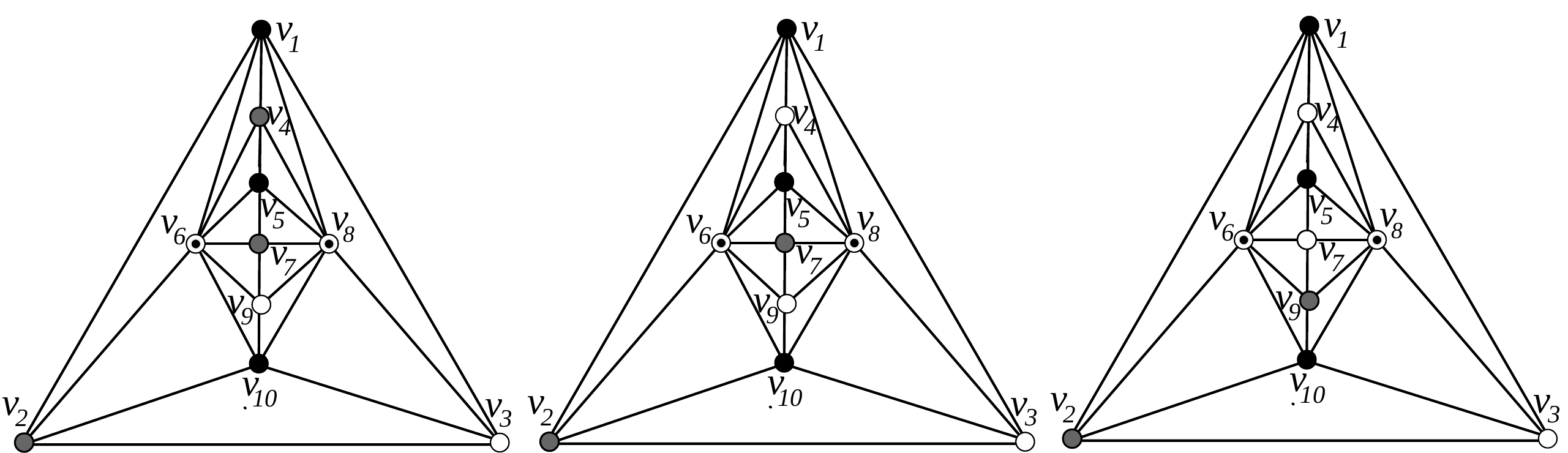}

         \vspace{5mm}
         \includegraphics [width=380pt]{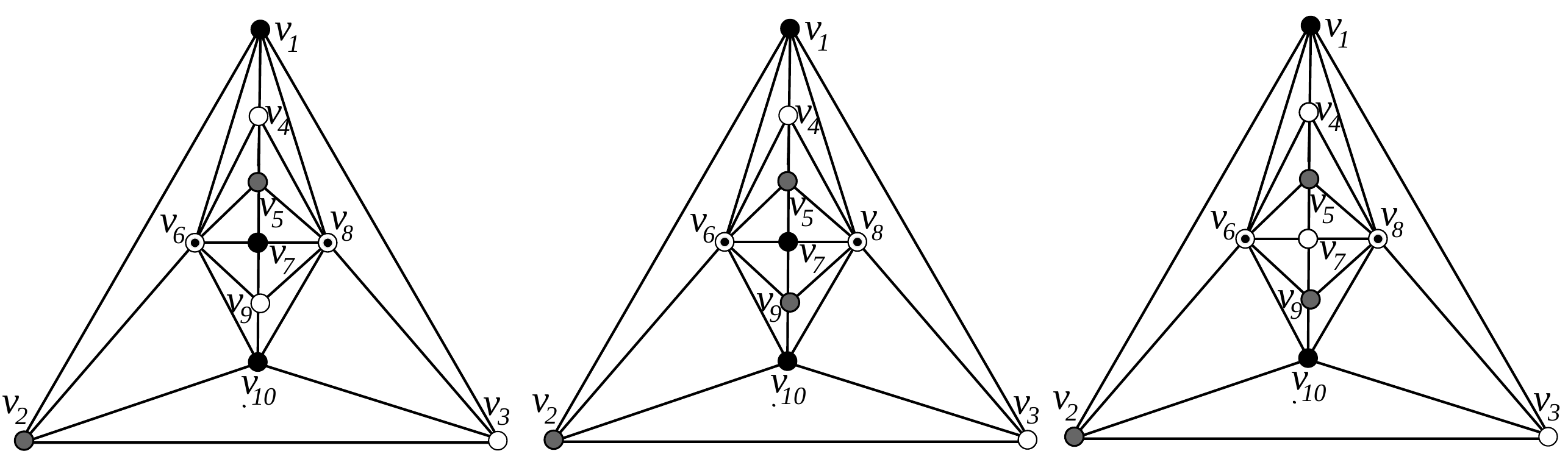}

         \vspace{5mm}
         \includegraphics [width=120pt]{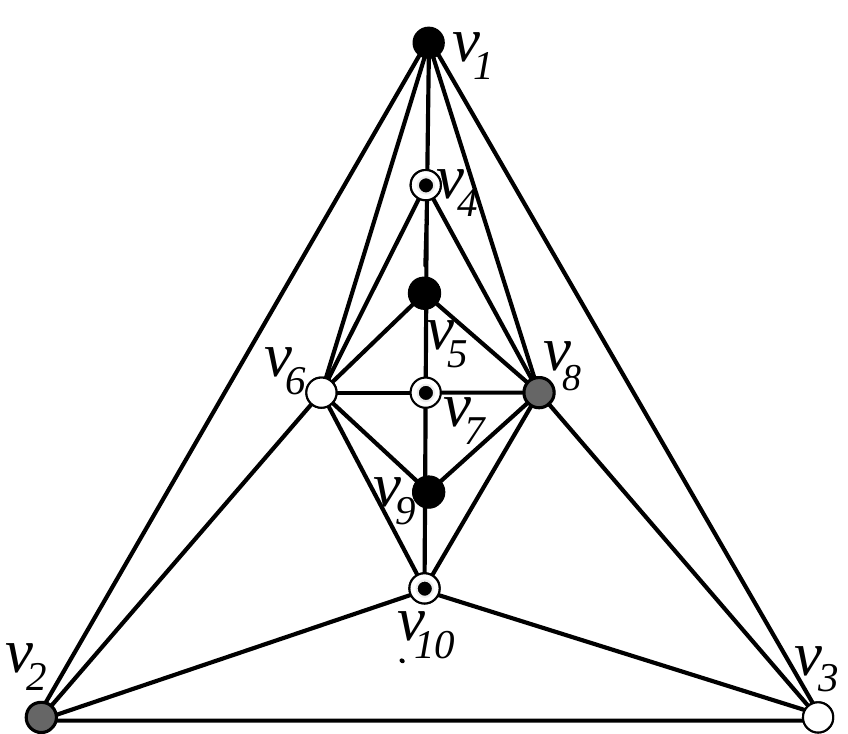}

         \vspace{8mm}
  \end{center}
 5.3 Degree sequence is 4444444488, and the unique partitions of color group are shown as follow:
   \begin{center}
        \includegraphics [width=320pt]{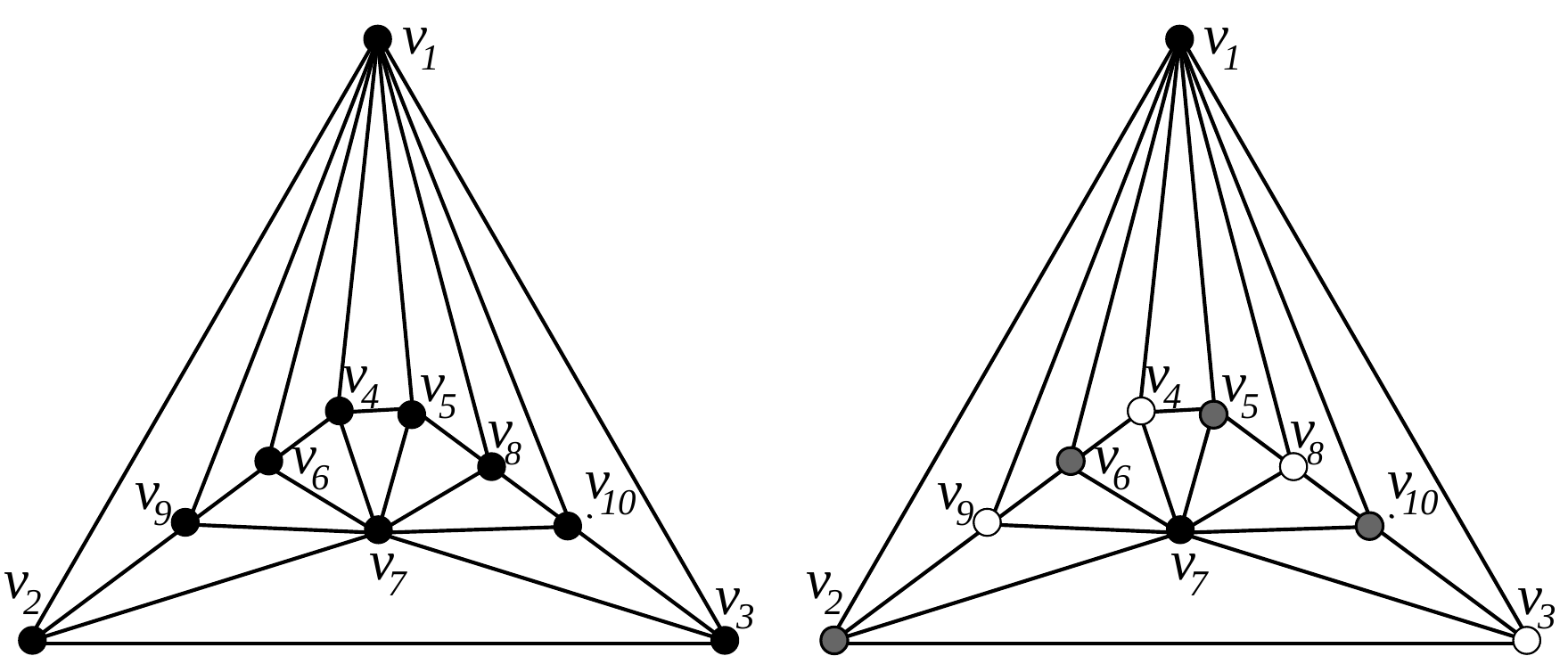}

         \vspace{8mm}
  \end{center}
 5.4 Degree sequence is 4444455666, and it has 14 kinds of different colorings.
    $$\{\{v_1,v_5,v_9\}\{v_2,v_4,v_6\}\{v_3,v_7,v_8\}\{v_{10}\}\} \{\{v_1,v_5,v_9\}\{v_2,v_4\}\{v_3,v_7,v_8\}\{v_6,v_{10}\}\}$$ $$\{\{v_1,v_5,v_9\}\{v_2,v_4,v_7\}\{v_3,v_8\}\{v_6,v_{10}\}\} \{\{v_1,v_5,v_9\}\{v_2,v_6,v_{10}\}\{v_3,v_8\}\{v_4,v_7\}\}$$ $$\{\{v_1,v_5,v_9\}\{v_2,v_6,v_{10}\}\{v_3,v_7,v_8\}\{v_4\}\} \{\{v_1,v_5,v_9\}\{v_2,v_{10}\}\{v_3,v_7,v_8\}\{v_4,v_6\}\}$$
     $$\{\{v_1,v_5,v_9\}\{v_2,v_4,v_7\}\{v_3,v_6\}\{v_8,v_{10}\}\} \{\{v_1,v_5,v_9\}\{v_2,v_4,v_6\}\{v_3,v_7\}\{v_8,v_{10}\}\}$$ $$\{\{v_1,v_5,v_9\}\{v_2,v_{10}\}\{v_3,v_4,v_6\}\{v_7,v_8\}\} \{\{v_1,v_5,v_9\}\{v_2,v_6,v_{10}\}\{v_3,v_4\}\{v_7,v_8\}\}$$
     $$\{\{v_1,v_5,v_9\}\{v_2,v_6,v_{10}\}\{v_3,v_4,v_7\}\{v_8\}\} \{\{v_1,v_5,v_9\}\{v_2,v_7\}\{v_3,v_4,v_6\}\{v_8,v_{10}\}\}$$ $$\{\{v_1,v_5,v_9\}\{v_2,v_6\}\{v_3,v_4,v_7\}\{v_8,v_{10}\}\} \{\{v_1,v_5\}\{v_2,v_6,v_10\}\{v_3,v_7,v_8\}\{v_4,v_9\}\}$$
        \begin{center}
        \includegraphics [width=160pt]{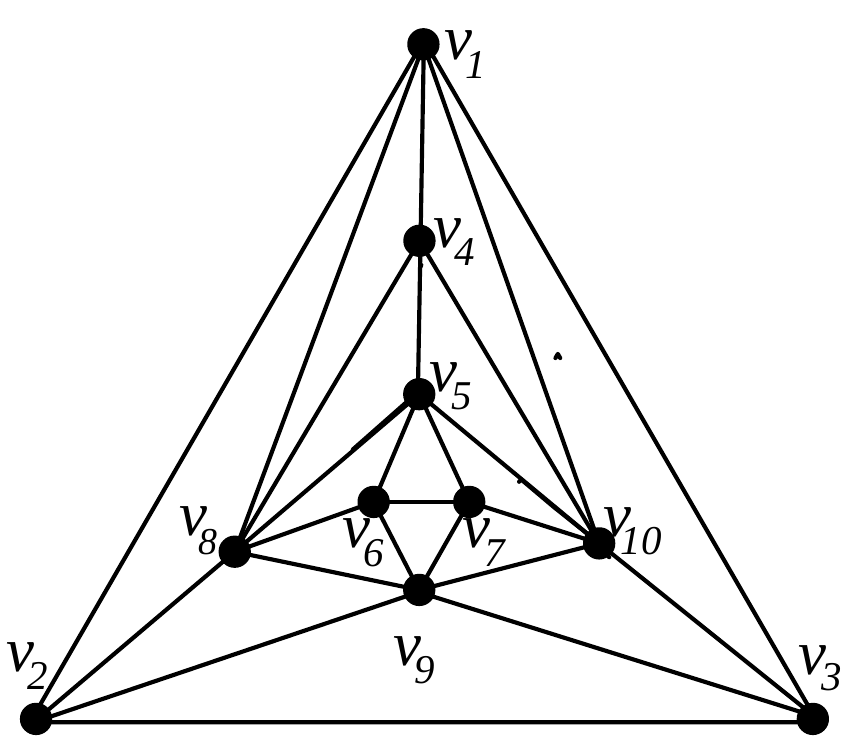}

         \vspace{5mm}
         \includegraphics [width=380pt]{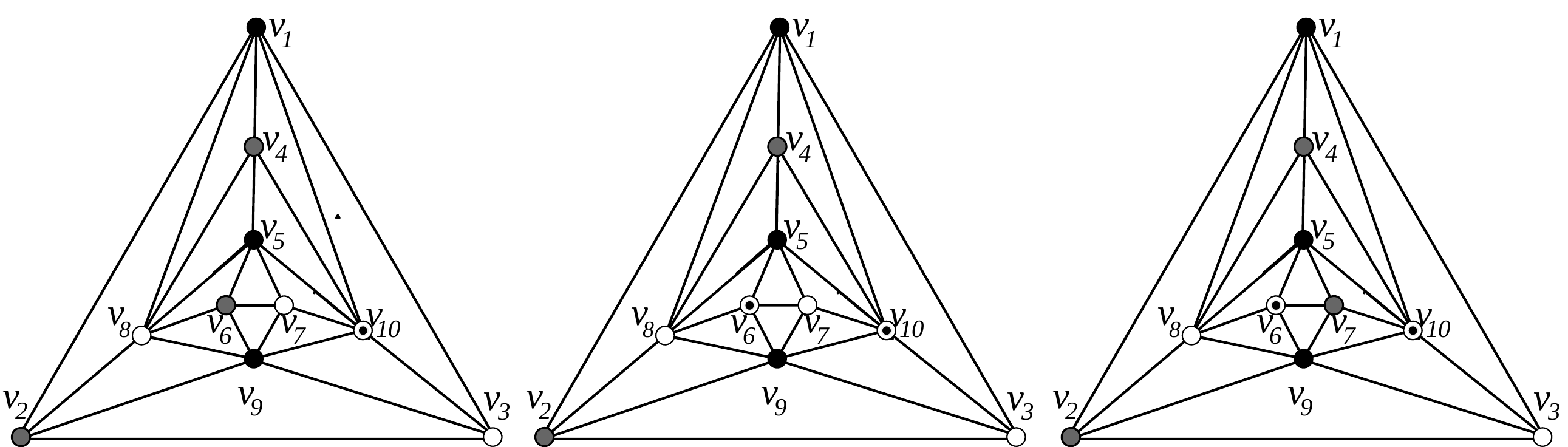}

         \vspace{5mm}
         \includegraphics [width=380pt]{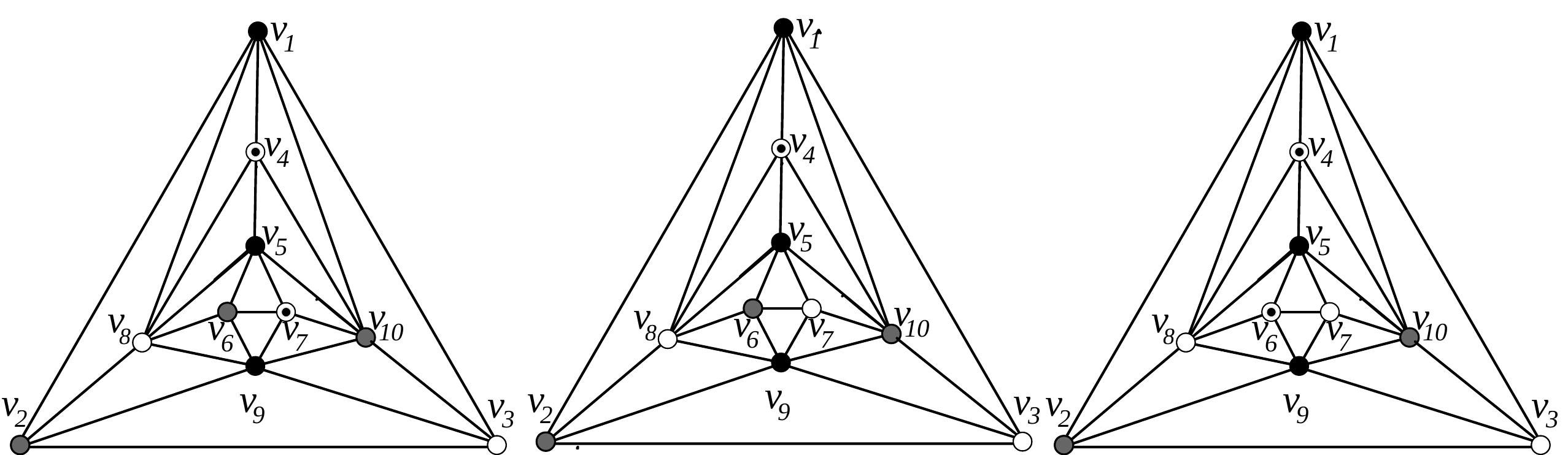}

         \vspace{5mm}
         \includegraphics [width=380pt]{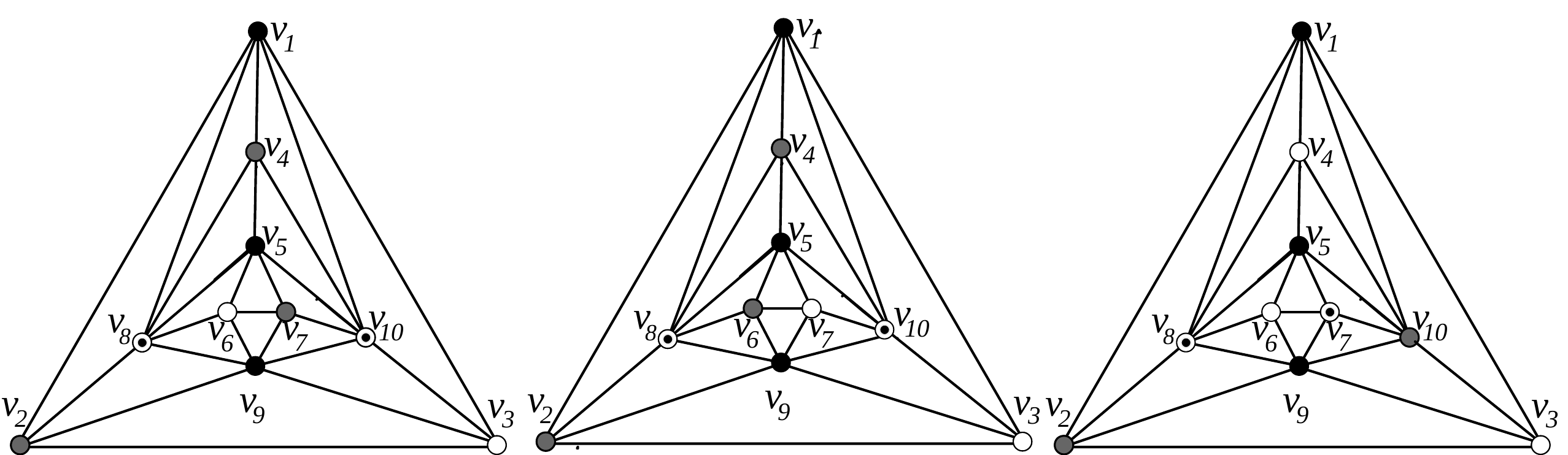}

         \vspace{5mm}
         \includegraphics [width=380pt]{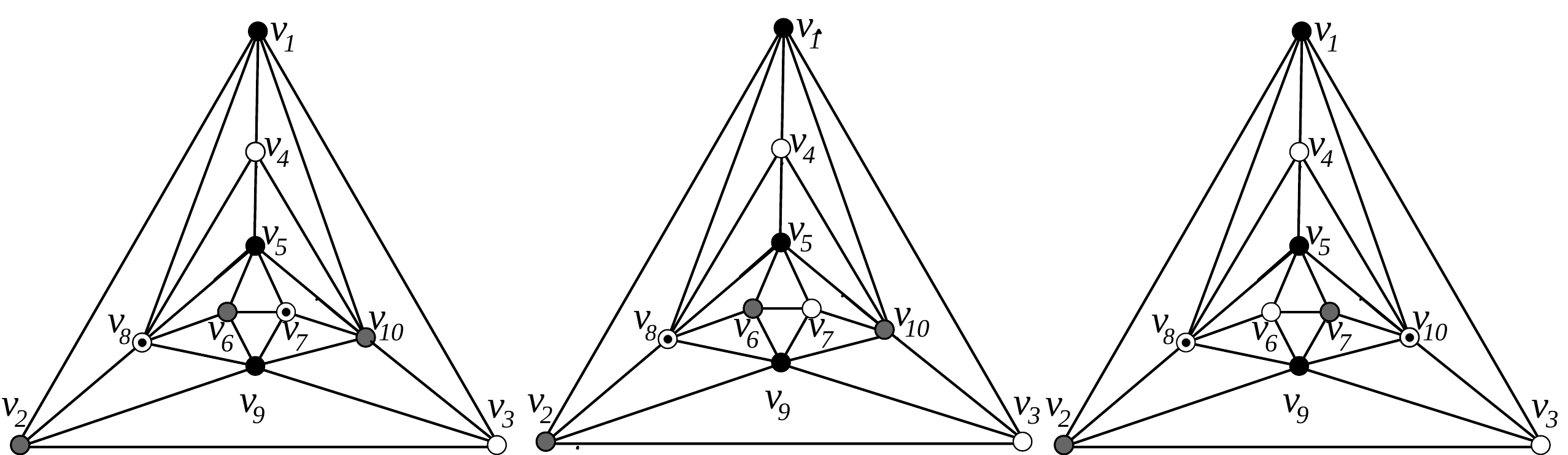}

         \vspace{5mm}
         \includegraphics [width=300pt]{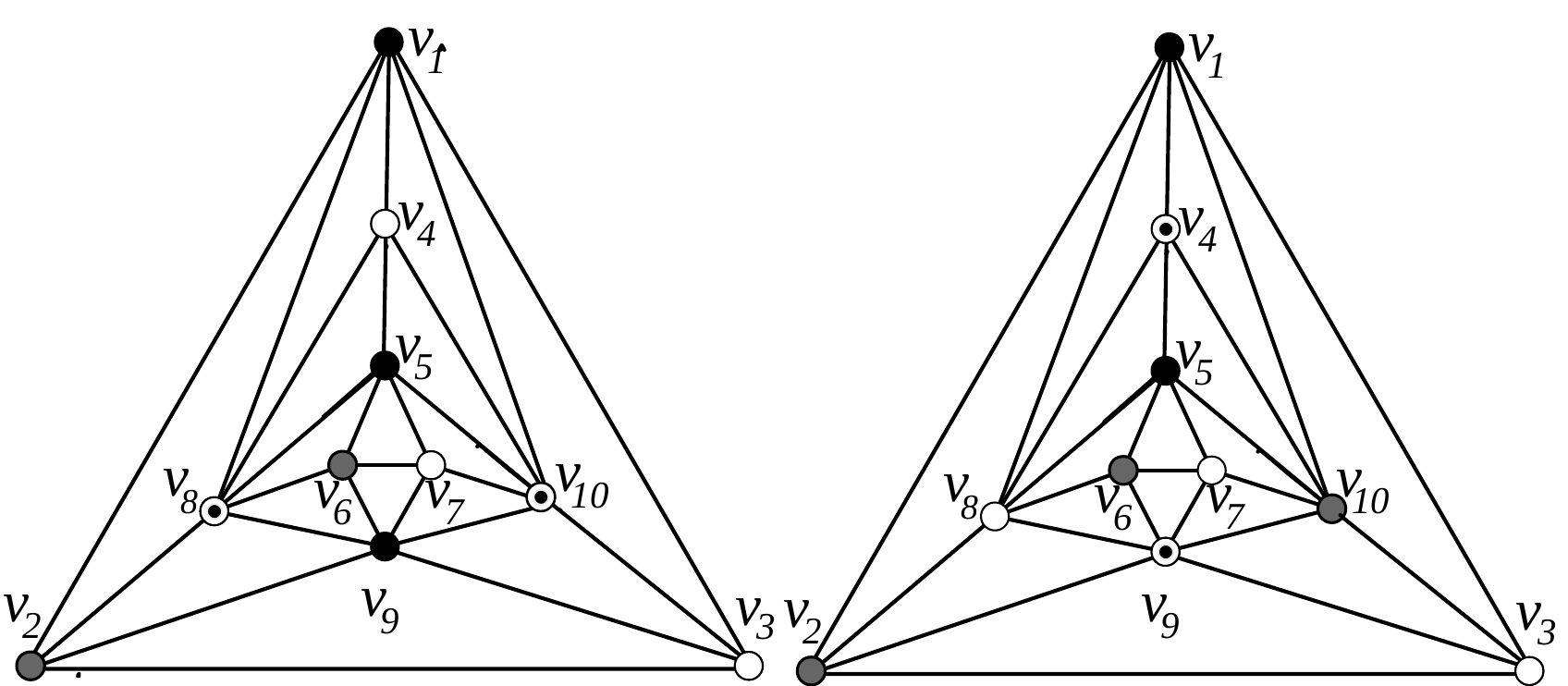}
         \vspace{8mm}
  \end{center}
 5.5 Degree sequence is 4445555666, and it has 6 kinds of different colorings.
     $$\{\{v_1,v_6,v_8\}\{v_2,v_9\}\{v_3,v_5,v_7\}\{v_4,v_{10}\}\} \{\{v_1,v_6,v_8\}\{v_2,v_7,v_9\}\{v_3,v_5\}\{v_4,v_{10}\}\}$$ $$\{\{v_1,v_6,v_8\}\{v_2,v_4\}\{v_3,v_5,v_7\}\{v_9,v_{10}\}\} \{\{v_1,v_6,v_8\}\{v_2,v_7\}\{v_3,v_4\}\{v_5,v_9,v_{10}\}\}$$ $$\{\{v_1,v_6,v_8\}\{v_2,v_7,v_9\}\{v_3,v_4\}\{v_5,v_{10}\}\}  \{\{v_1,v_6,v_8\}\{v_2,v_4\}\{v_3,v_7\}\{v_5,v_9,v_{10}\}\}$$
            \begin{center}
        \includegraphics [width=160pt]{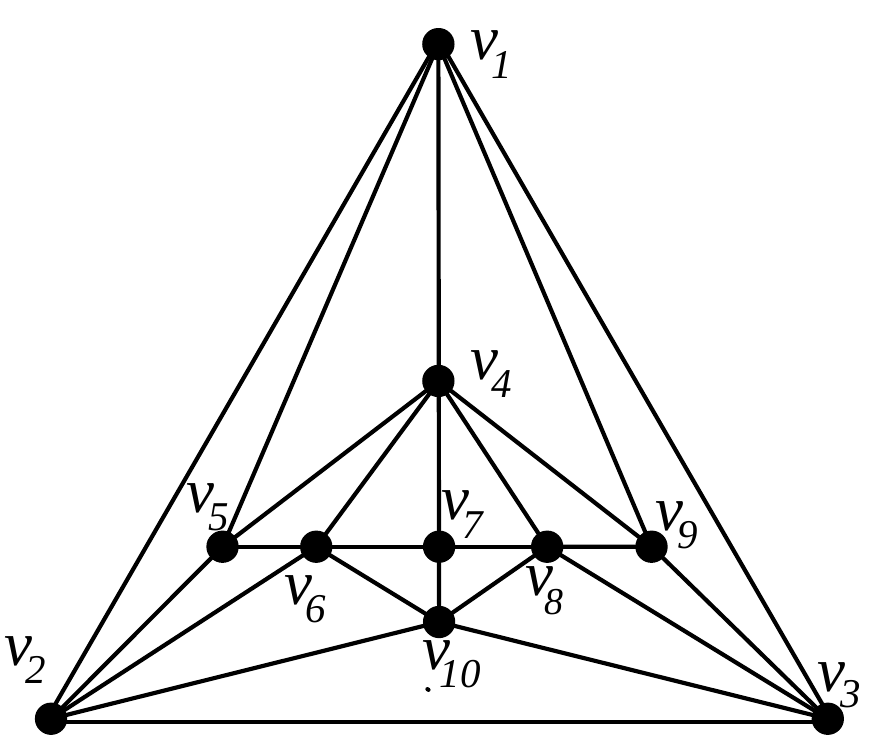}

         \vspace{5mm}
         \includegraphics [width=380pt]{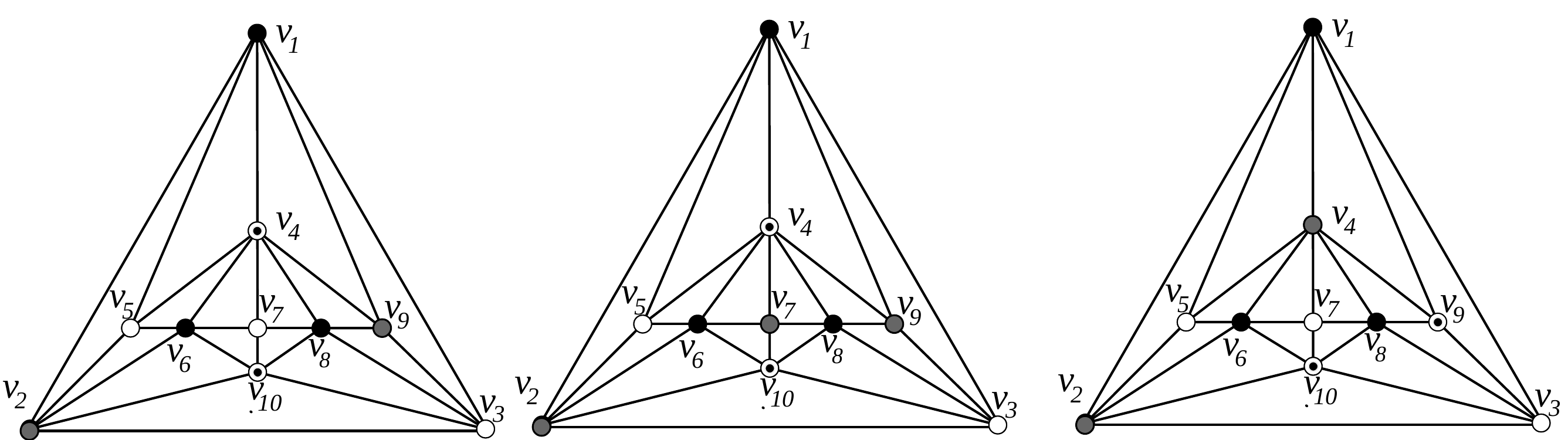}

         \vspace{5mm}
         \includegraphics [width=380pt]{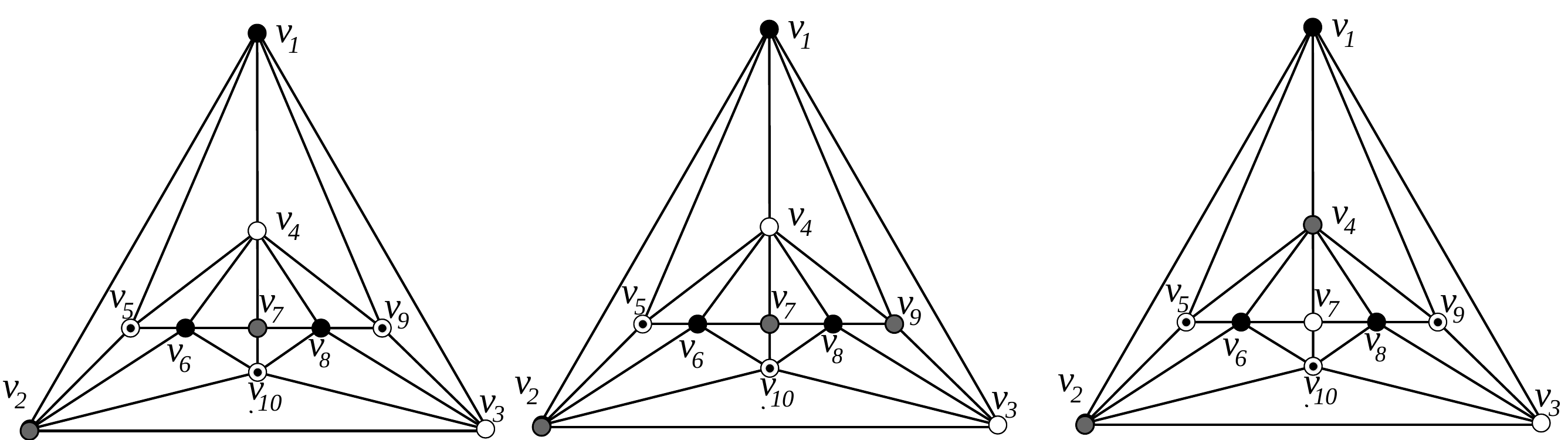}

         \vspace{8mm}
  \end{center}
 5.6 Degree sequence is 4444446666, and it is uniquely 3-colorable and has 28 kinds of different colorings.

 $$ \{\{v_1,v_8\}\{v_2,v_5,v_{10}\}\{v_3,v_4\}\{v_6,v_7,v_9\}\} \{\{v_1\} \{v_2,v_5,v_{10}\}\{v_3,v_4,v_8\}\{v_6,v_7,v_9\}\}$$
 $$ \{\{v_1,v_9\}\{v_2,v_5,v_{10}\}\{v_3,v_4,v_8\}\{v_6,v_7\}\},\{\{v_1,v_{10}\}  \{v_2,v_5\}\{v_3,v_4,v_8\}\{v_6,v_7,v_9\}\}$$
 $$ \{\{v_1,v_7\}\{v_2,v_5,v_{10}\}\{v_3,v_4,v_8\}\{v_6,v_9\}\},\{\{v_1,v_7,v_9\}  \{v_2,v_5,v_{10}\}\{v_3,v_4,v_8\}\{v_6\}\}$$
 $$ \{\{v_1,v_7,v_9\}\{v_2,v_5\}\{v_3,v_4,v_8\}\{v_6,v_{10}\}\}, \{\{v_1,v_6,v_{10}\}  \{v_2,v_7\}\{v_3,v_4,v_8\}\{v_5,v_9\}\}$$
 $$ \{\{v_1,v_6,v_9\}\{v_2,v_7\}\{v_3,v_4,v_8\}\{v_5,v_{10}\}\}, \{\{v_1,v_6,v_7,v_9\}  \{v_2\}\{v_3,v_4,v_8\}\{v_5,v_{10}\}\}$$
 $$ \{\{v_1,v_6,v_7,v_9\}\{v_2,v_{10}\}\{v_3,v_4,v_8\}\{v_5\}\}, \{\{v_1,v_6,v_7\}  \{v_2,v_{10}\}\{v_3,v_4,v_8\}\{v_5,v_9\}\}$$
 $$ \{\{v_1,v_6,v_7,v_9\}\{v_2,v_5,v_{10}\}\{v_3,v_4\}\{v_8\}\}, \{\{v_1,v_6,v_7\}  \{v_2,v_5,v_{10}\}\{v_3,v_4,v_8\}\{v_9\}\}$$
 $$ \{\{v_1,v_6,v_7,v_9\}\{v_2,v_5,v_{10}\}\{v_3,v_4,v_8\}\}\}, \{\{v_1,v_6,v_7,v_9\}  \{v_2,v_5\}\{v_3,v_4,v_8\}\{v_{10}\}\}$$
   $$ \{\{v_1,v_6,v_9\}\{v_2,v_5,v_{10}\}\{v_3,v_4,v_8\}\{v_7\}\}, \{\{v_1,v_6\}  \{v_2,v_5,v_{10}\}\{v_3,v_4,v_8\}\{v_7,v_9\}\}$$
  $$ \{\{v_1,v_6,v_{10}\}\{v_2,v_5\}\{v_3,v_4,v_8\}\{v_7,v_9\}\}, \{\{v_1,v_6,v_7,v_9\}  \{v_2,v_5,v_{10}\}\{v_4,v_8\}\{v_3\}\}$$
   $$ \{\{v_1,v_6,v_7,v_9\}\{v_2,v_5\}\{v_3,v_8\}\{v_4,v_{10}\}\}, \{\{v_1,v_6,v_7,v_9\}  \{v_2,v_5,v_{10}\}\{v_3,v_8\}\{v_4\}\}$$
  $$ \{\{v_1,v_6,v_7\}\{v_2,v_5,v_{10}\}\{v_3,v_8\}\{v_4,v_9\}\}, \{\{v_1,v_6,v_9\} \{v_2,v_5,v_{10}\}\{v_3,v_8\}\{v_4,v_7\}\}$$
   $$ \{\{v_1,v_6\}\{v_2,v_5,v_{10}\}\{v_3,v_8\}\{v_4,v_7,v_9\}\}, \{\{v_1,v_6,v_{10}\} \{v_2,v_5\}\{v_3,v_8\}\{v_4,v_7,v_9\}\}$$
  $$ \{\{v_1,v_8\}\{v_2,v_5,v_{10}\}\{v_3,v_6\}\{v_4,v_7,v_9\}\}, \{\{v_1,v_7,v_9\} \{v_2,v_5,v_{10}\}\{v_3,v_6\}\{v_4,v_8\}\}$$

            \begin{center}
        \includegraphics [width=160pt]{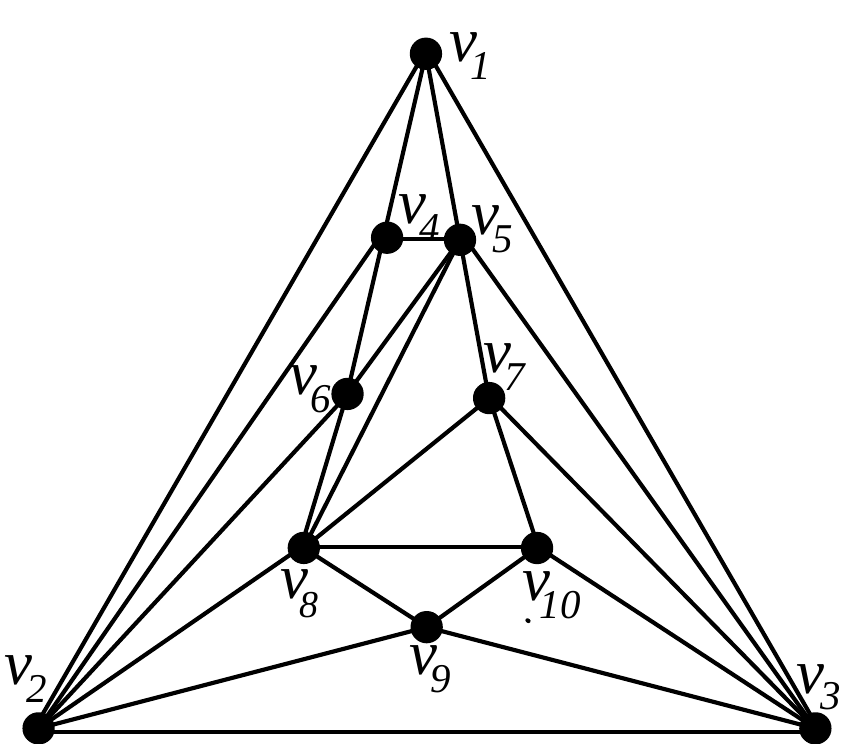}

         \vspace{5mm}
         \includegraphics [width=380pt]{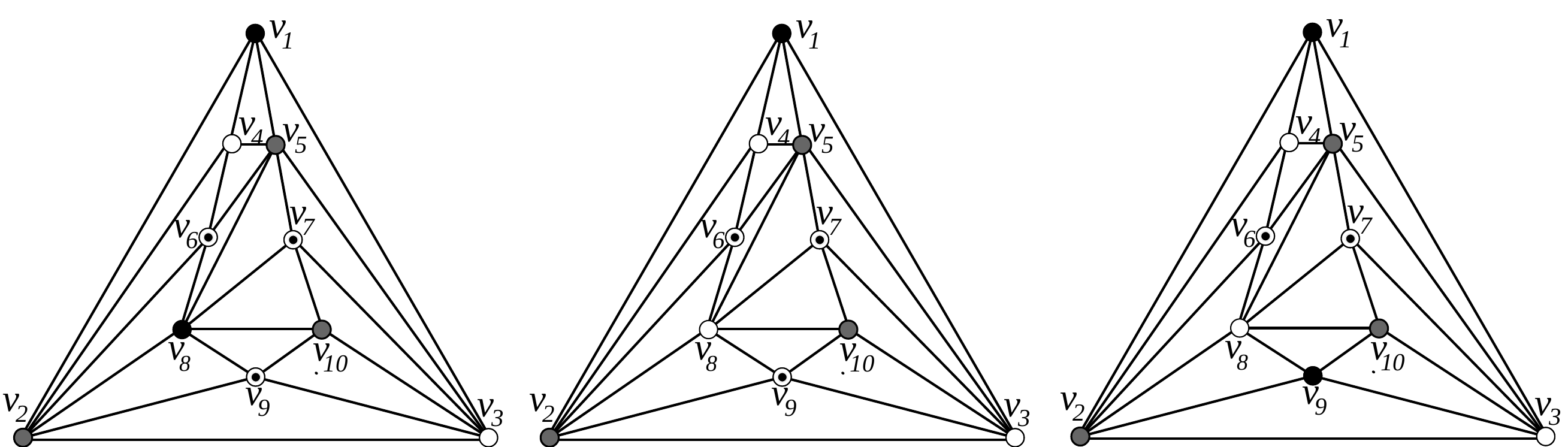}

         \vspace{5mm}
         \includegraphics [width=380pt]{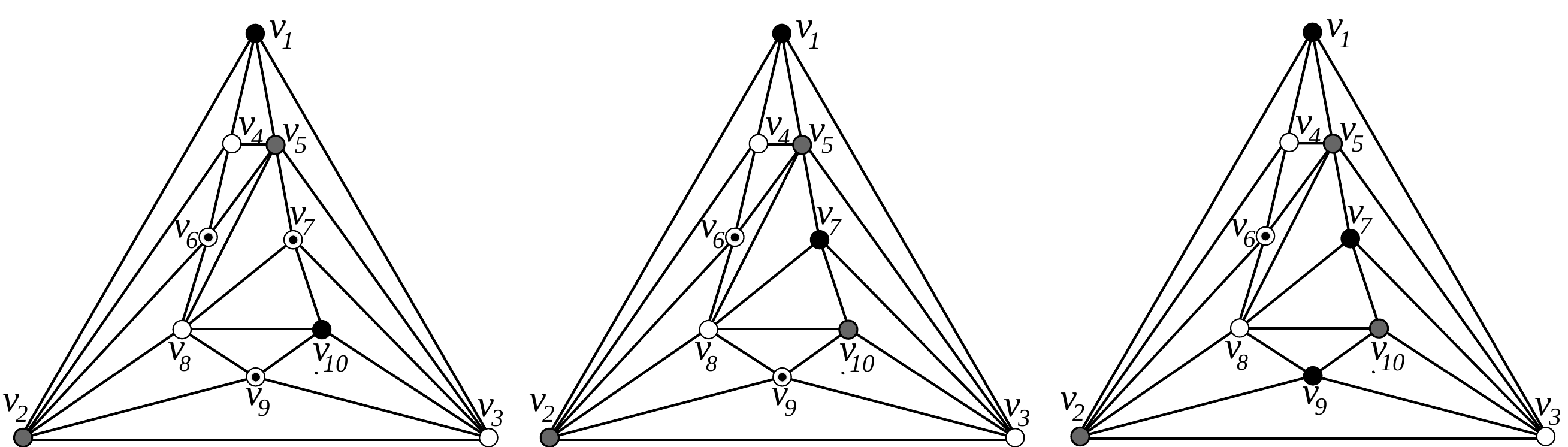}

         \vspace{5mm}
         \includegraphics [width=380pt]{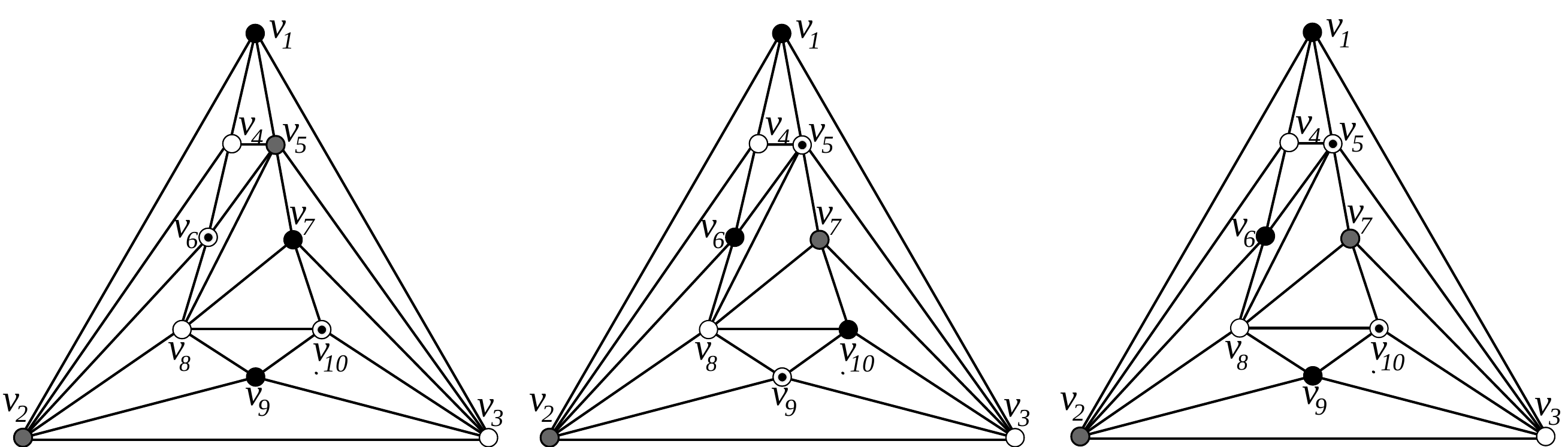}

                  \vspace{5mm}
         \includegraphics [width=380pt]{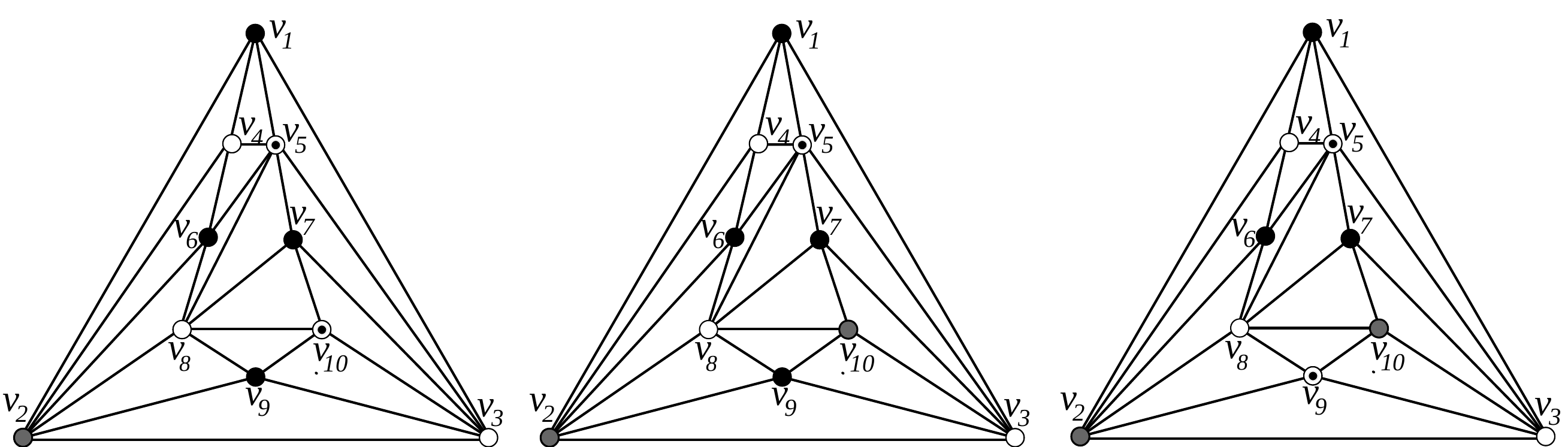}

                  \vspace{5mm}
         \includegraphics [width=380pt]{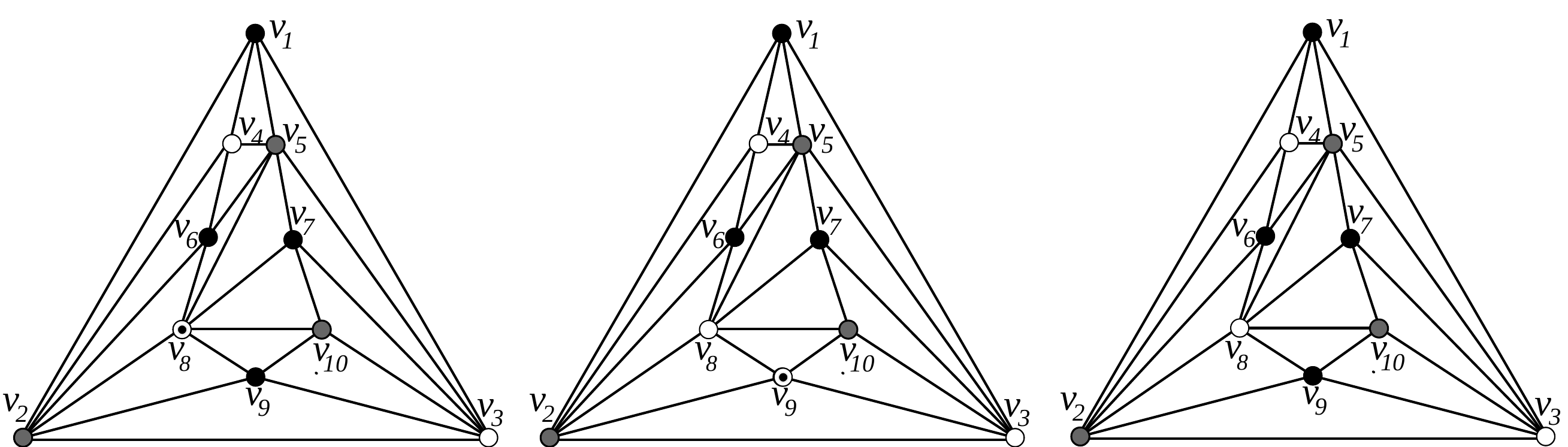}

                  \vspace{5mm}
         \includegraphics [width=380pt]{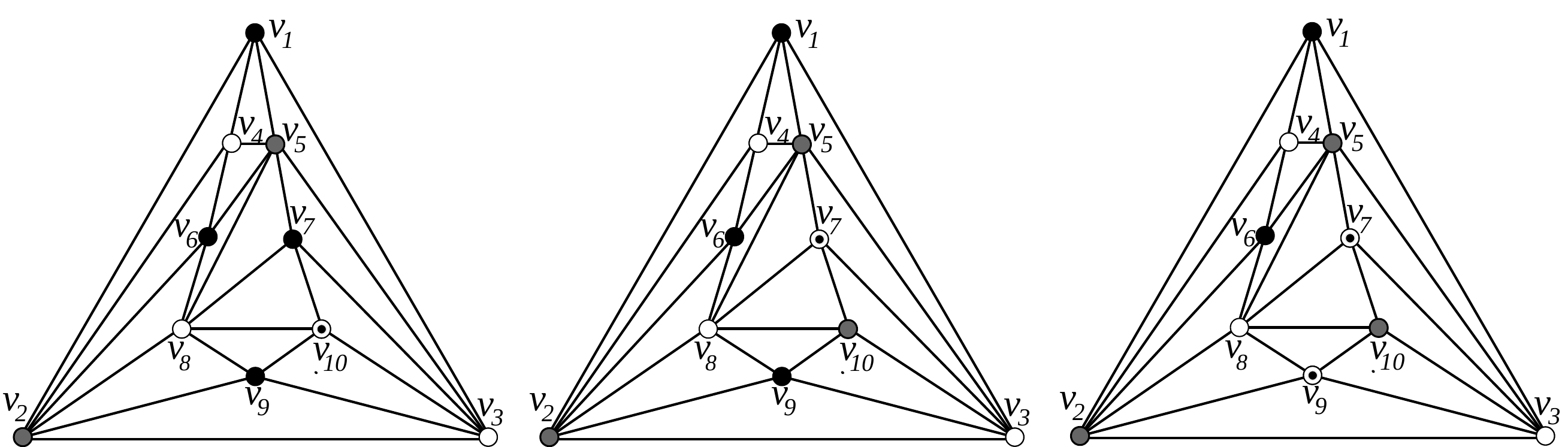}

                  \vspace{5mm}
         \includegraphics [width=380pt]{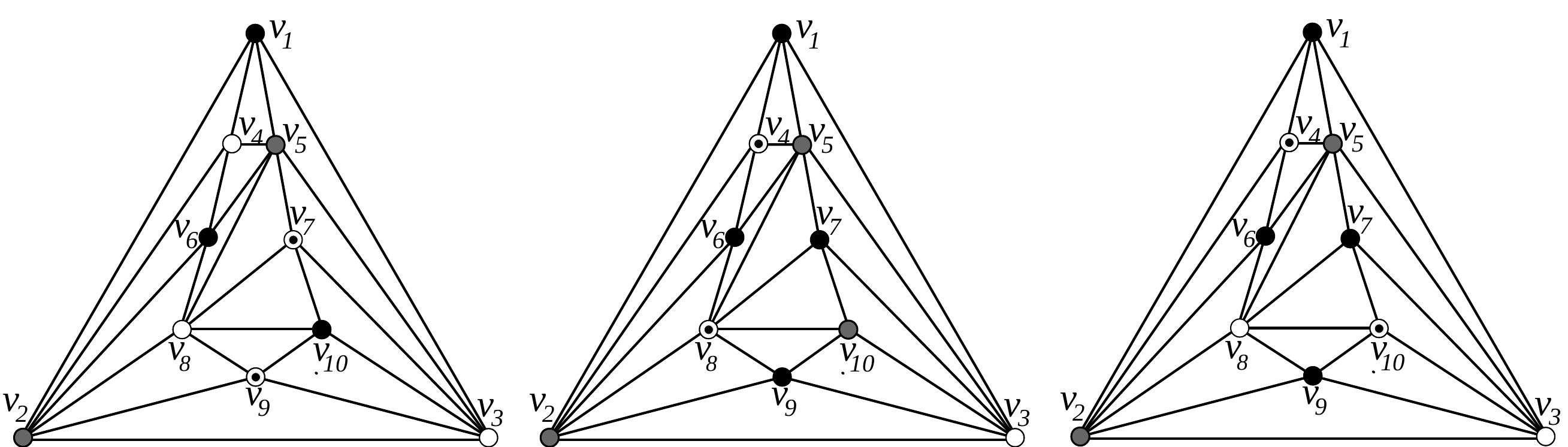}

                  \vspace{5mm}
         \includegraphics [width=380pt]{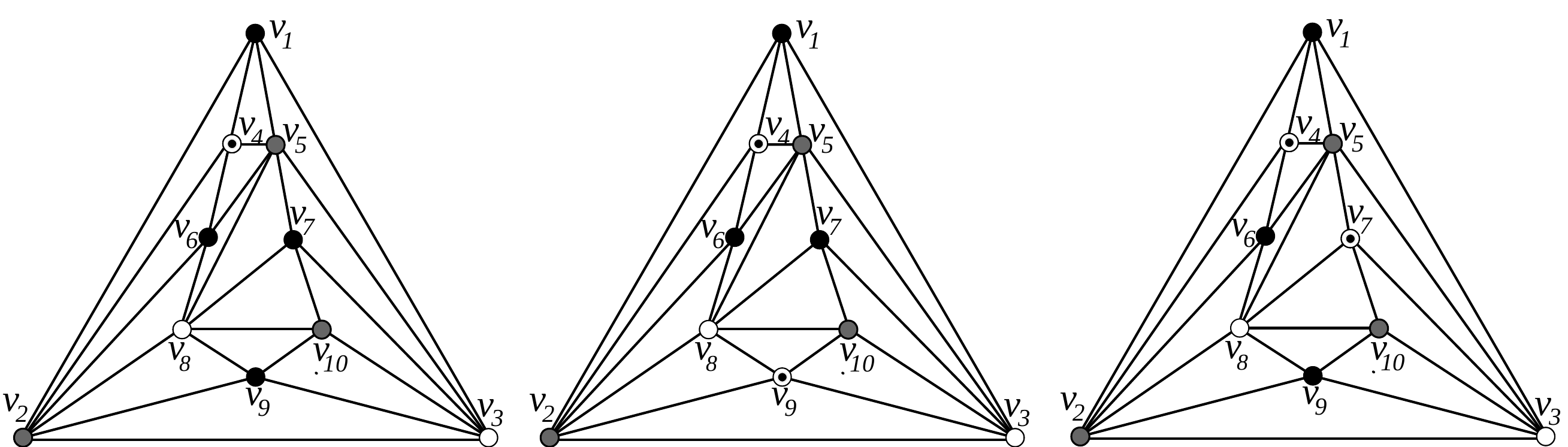}

                  \vspace{5mm}
         \includegraphics [width=380pt]{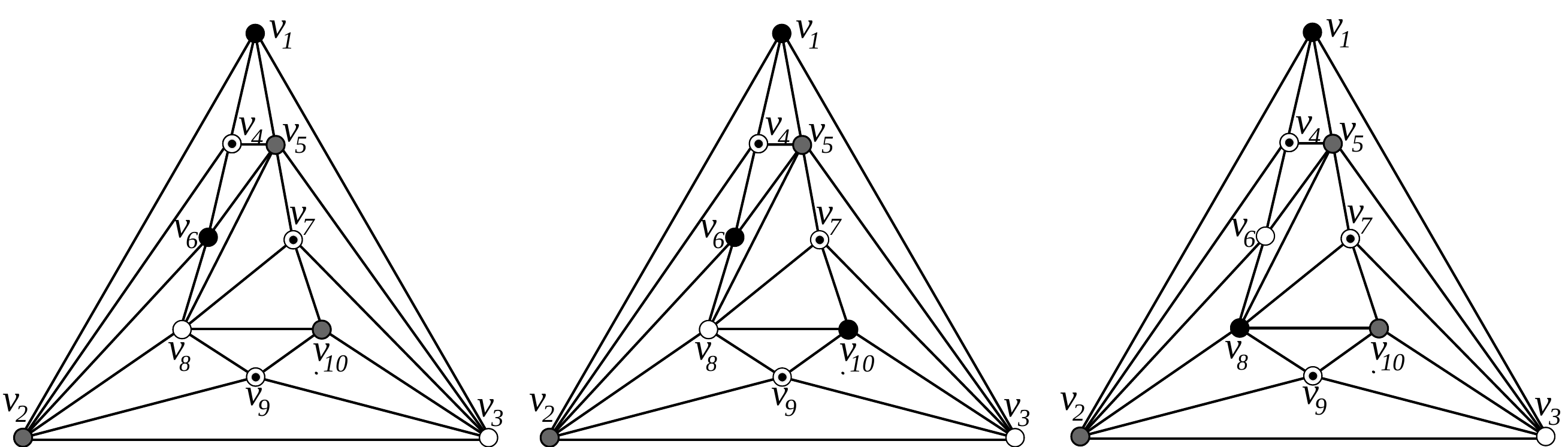}

                  \vspace{5mm}
         \includegraphics [width=160pt]{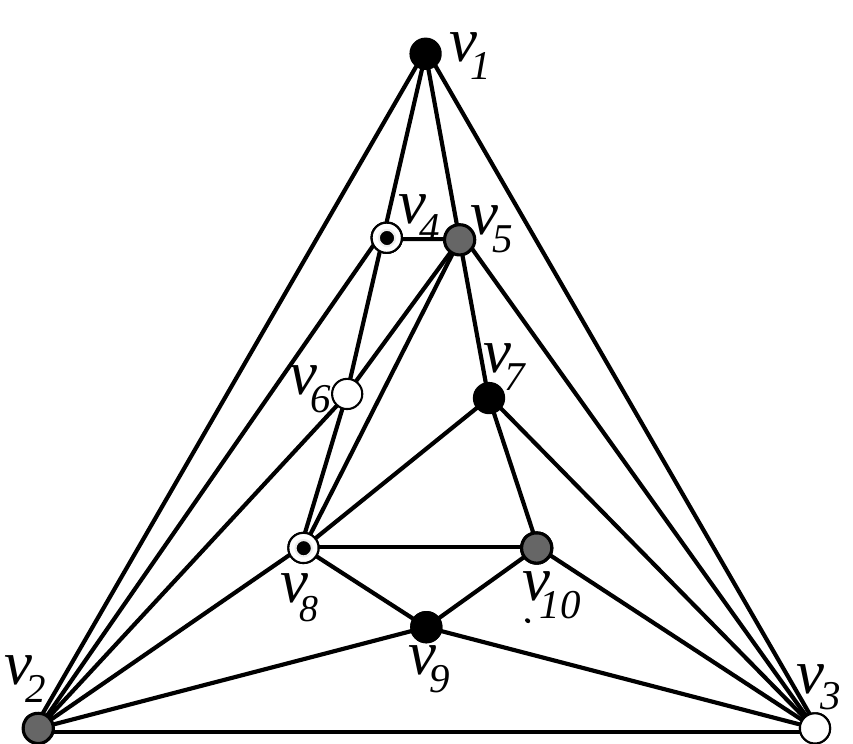}

         \vspace{8mm}
  \end{center}
 5.7 Degree sequence is 4445555556, and it has 4 kinds of different colorings.

  $$ \{\{v_1,v_8,v_{10}\}\{v_2,v_6\}\{v_3,v_4\}\{v_5,v_7,v_9\}\}, \{\{v_1,v_6\} \{v_2,v_4,v_{10}\}\{v_3,v_8\}\{v_5,v_7,v_9\}\}$$
  $$ \{\{v_1,v_{10}\}\{v_2,v_6\}\{v_3,v_4,v_8\}\{v_5,v_7,v_9\}\}, \{\{v_1,v_6\} \{v_2,v_{10}\}\{v_3,v_4,v_8\}\{v_5,v_7,v_9\}\}$$
            \begin{center}
        \includegraphics [width=160pt]{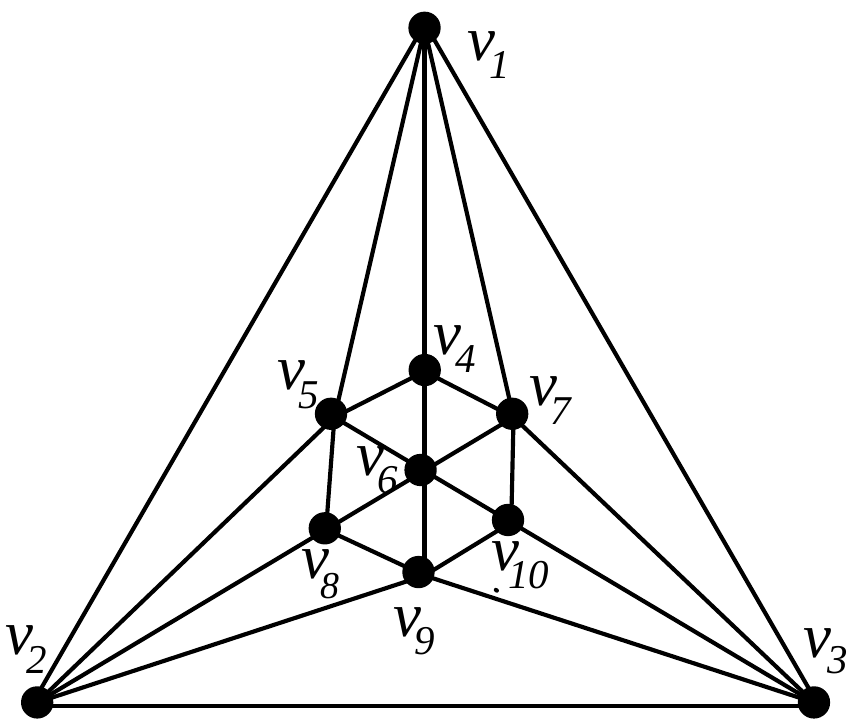}

         \vspace{5mm}
         \includegraphics [width=380pt]{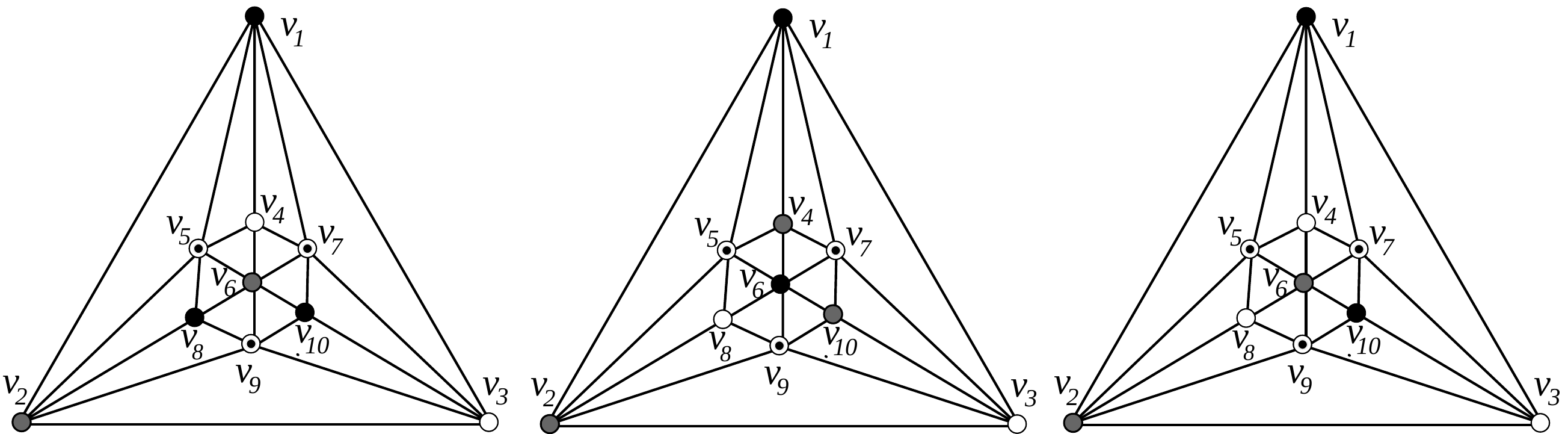}

         \vspace{5mm}
         \includegraphics [width=160pt]{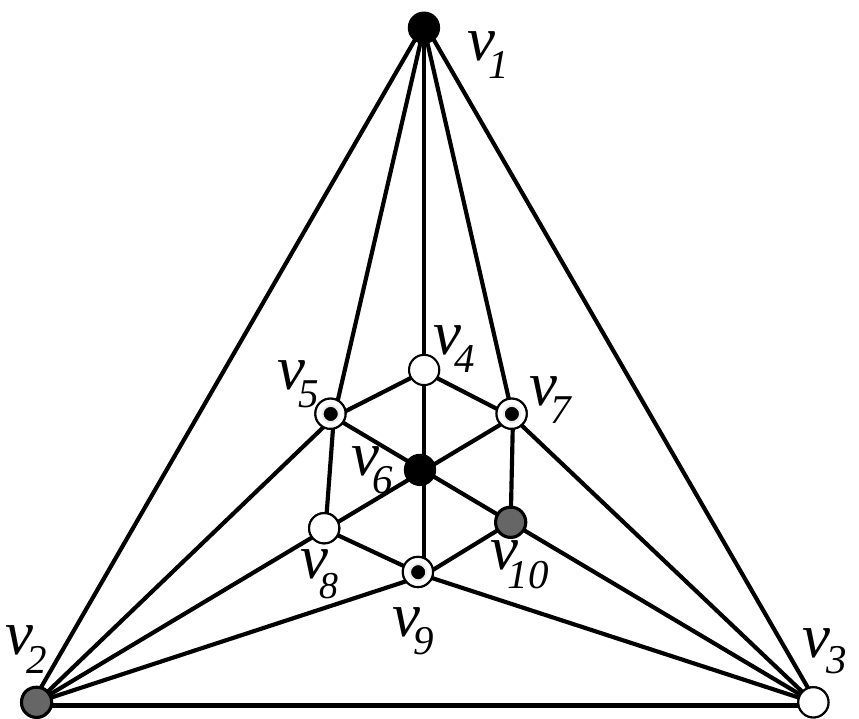}

         \vspace{8mm}
  \end{center}

 5.8 Degree sequence is 4445555555, and it has 8 kinds of different colorings.

  $$ \{\{v_1,v_6,v_9\}\{v_2,v_7\}\{v_3,v_4,v_8\}\{v_5,v_{10}\}\}, \{\{v_1,v_6,v_9\} \{v_2,v_{10}\}\{v_3,v_4,v_8\}\{v_5,v_7\}\}$$
  $$ \{\{v_1,v_6,v_9\}\{v_2,v_4,v_{10}\}\{v_3,v_8\}\{v_5,v_7\}\}, \{\{v_1,v_6,v_9\} \{v_2,v_4,v_{10}\}\{v_3,v_5\}\{v_7,v_8\}\}$$
  $$ \{\{v_1,v_8\}\{v_2,v_4,v_{10}\}\{v_3,v_6\}\{v_5,v_7,v_9\}\}, \{\{v_1,v_6\} \{v_2,v_4,v_{10}\}\{v_3,v_8\}\{v_5,v_7,v_9\}\}$$
  $$ \{\{v_1,v_6\}\{v_2,v_{10}\}\{v_3,v_4,v_8\}\{v_5,v_7,v_9\}\}, \{\{v_1,v_{10}\} \{v_2,v_6\}\{v_3,v_4,v_8\}\{v_5,v_7,v_9\}\}$$
            \begin{center}
        \includegraphics [width=160pt]{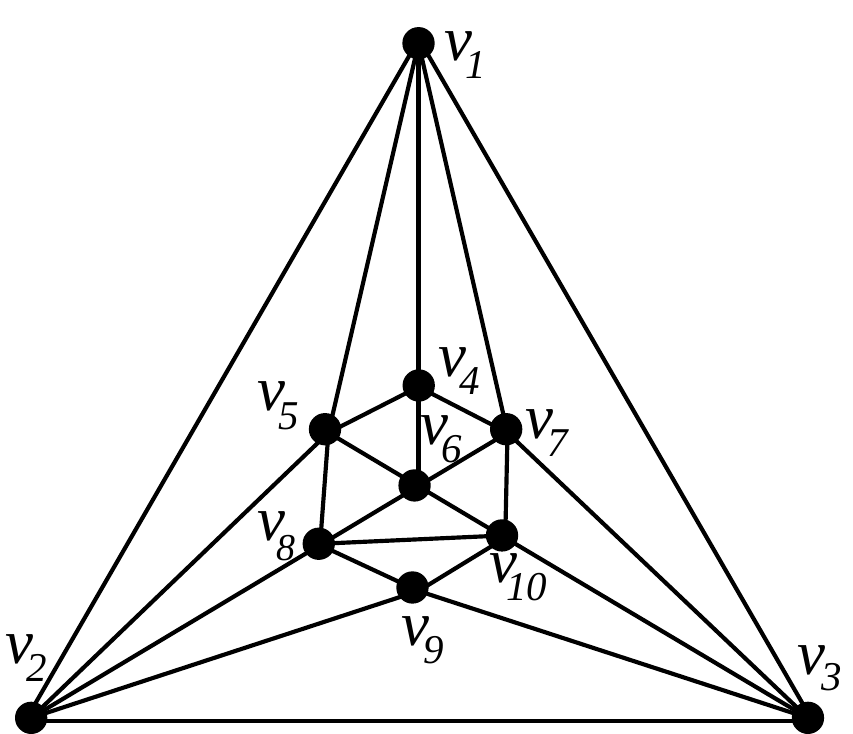}

         \vspace{5mm}
         \includegraphics [width=380pt]{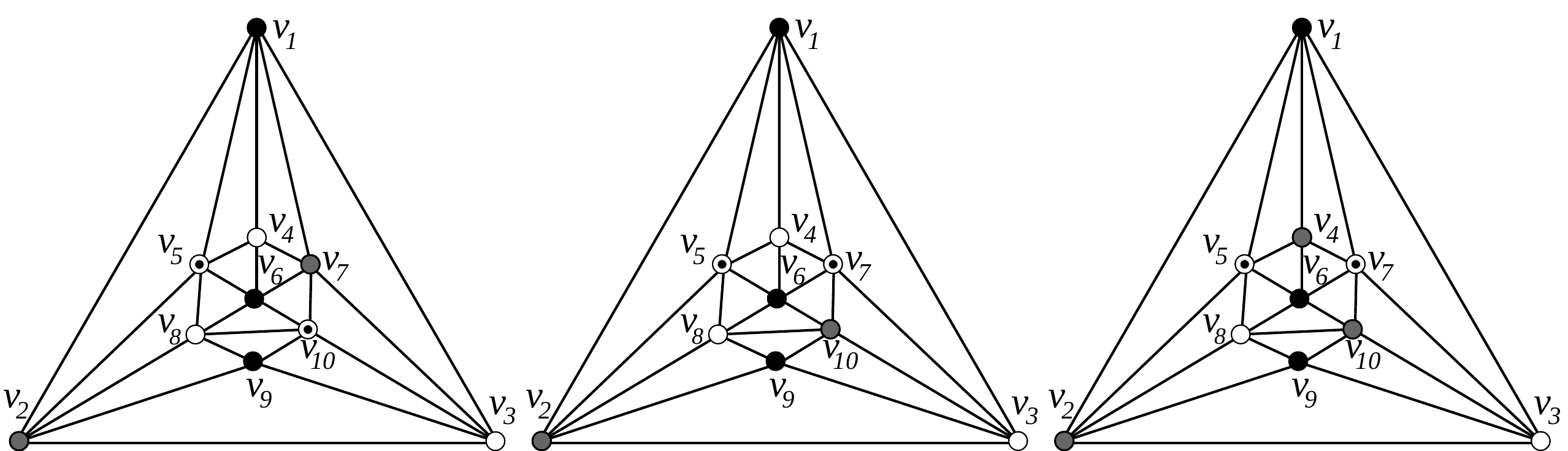}

         \vspace{5mm}
         \includegraphics [width=380pt]{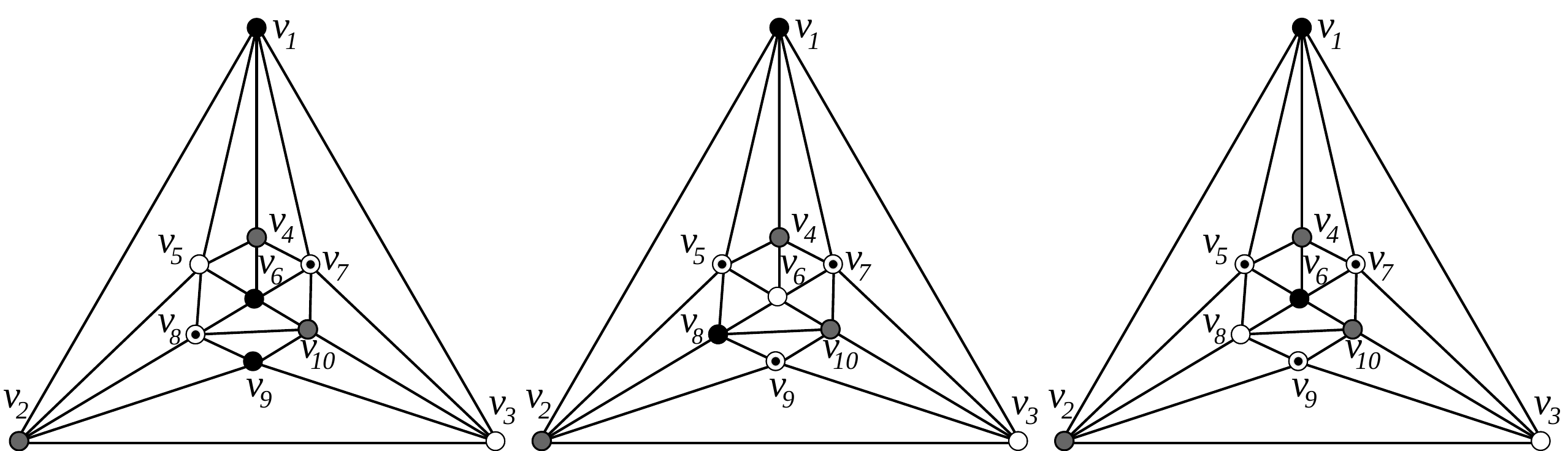}

         \vspace{5mm}
         \includegraphics [width=320pt]{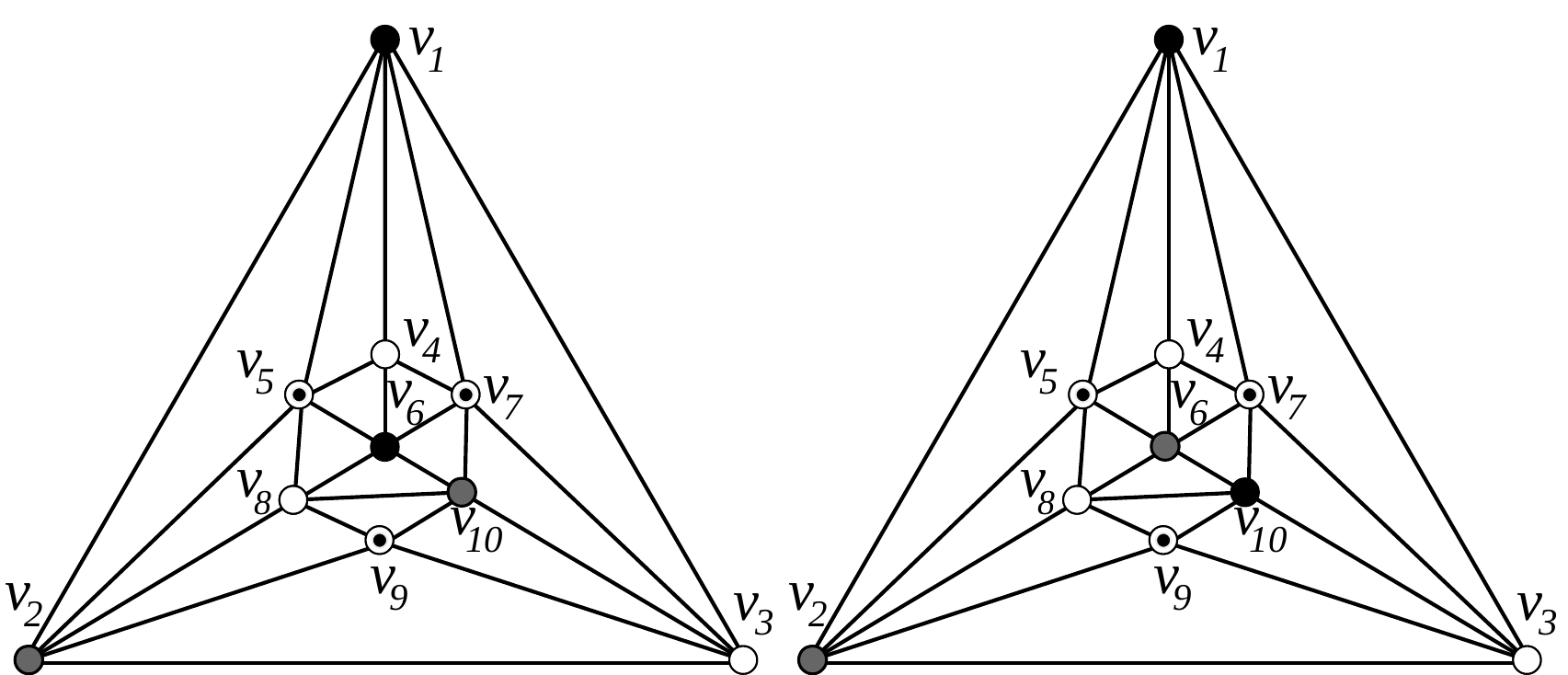}

         \vspace{8mm}
  \end{center}

 5.9 Degree sequence is 4444445667, and it has 20 kinds of different colorings.

  $$ \{\{v_1,v_7,v_{10}\}\{v_2,v_5,v_9\}\{v_3,v_4,v_8\}\{v_6\}\}, \{\{v_1,v_7\} \{v_2,v_5,v_9\}\{v_3,v_4,v_8\}\{v_6,v_{10}\}\}$$
  $$ \{\{v_1,v_7,v_{10}\}\{v_2,v_5\}\{v_3,v_4,v_8\}\{v_6,v_9\}\}, \{\{v_1,v_7,v_{10}\} \{v_2,v_5,v_9\}\{v_3,v_4\}\{v_6,v_8\}\}$$
  $$ \{\{v_1,v_{10}\}\{v_2,v_5,v_9\}\{v_3,v_6,v_7\}\{v_4,v_8\}\}, \{\{v_1,v_6,v_{10}\} \{v_2,v_5,v_9\}\{v_3,v_7\}\{v_4,v_8\}\}$$
  $$ \{\{v_1,v_7,v_{10}\}\{v_2,v_6,v_9\}\{v_3,v_4,v_8\}\{v_5\}\}, \{\{v_1,v_7\} \{v_2,v_6,v_9\}\{v_3,v_4,v_8\}\{v_5,v_{10}\}\}$$
  $$ \{\{v_1,v_7,v_{10}\}\{v_2,v_6\}\{v_3,v_4,v_8\}\{v_5,v_9\}\}, \{\{v_1,v_7,v_{10}\} \{v_2,v_6,v_9\}\{v_3,v_4\}\{v_5,v_8\}\}$$
  $$ \{\{v_1,v_7,v_{10}\}\{v_2,v_5\}\{v_3,v_8\}\{v_4,v_6,v_9\}\}, \{\{v_1,v_7,v_{10}\} \{v_2,v_5,v_9\}\{v_3,v_8\}\{v_4,v_6\}\}$$
  $$ \{\{v_1,v_7,v_{10}\}\{v_2,v_5,v_9\}\{v_3\}\{v_4,v_6,v_8\}\}, \{\{v_1,v_7\} \{v_2,v_5,v_9\}\{v_3,v_8\}\{v_4,v_6,v_{10}\}\}$$
  $$ \{\{v_1,v_7,v_{10}\}\{v_2\}\{v_3,v_5,v_8\}\{v_4,v_6,v_9\}\}, \{\{v_1,v_7,v_{10}\} \{v_2,v_9\}\{v_3,v_5,v_8\}\{v_4,v_6\}\}$$
  $$ \{\{v_1,v_7,v_{10}\}\{v_2,v_9\}\{v_3,v_5\}\{v_4,v_6,v_8\}\}, \{\{v_1,v_7\} \{v_2,v_9\}\{v_3,v_5,v_8\}\{v_4,v_6,v_{10}\}\}$$
  $$ \{\{v_1,v_7,v_{10}\}\{v_2,v_6\}\{v_3,v_5,v_8\}\{v_4,v_9\}\}, \{\{v_1,v_7,v_{10}\} \{v_2,v_6,v_9\}\{v_3,v_5,v_8\}\{v_4\}\}$$
  $$ \{\{v_1,v_7,v_{10}\}\{v_2,v_6,v_9\}\{v_3,v_5\}\{v_4,v_8\}\}, \{\{v_1,v_7\} \{v_2,v_6,v_9\}\{v_3,v_5,v_8\}\{v_4,v_{10}\}\}$$

            \begin{center}
        \includegraphics [width=160pt]{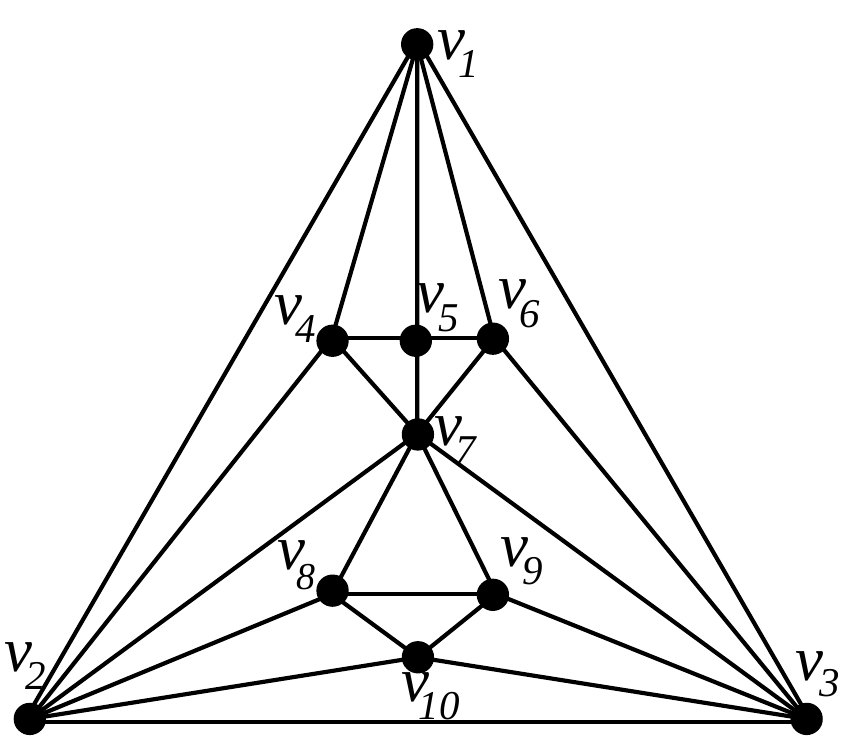}

         \vspace{5mm}
         \includegraphics [width=380pt]{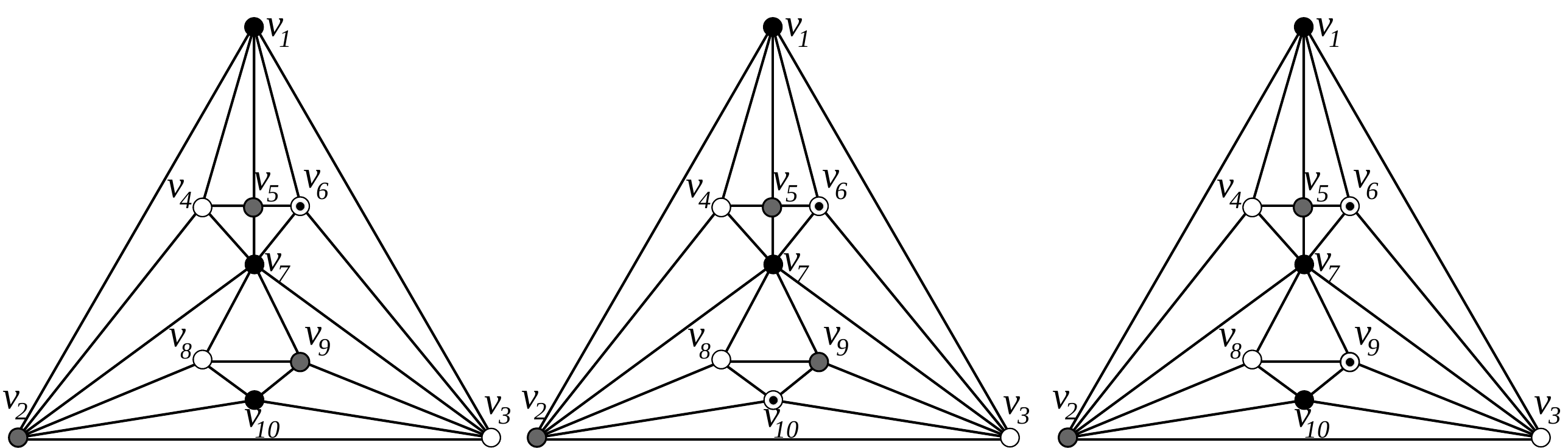}

         \vspace{5mm}
         \includegraphics [width=380pt]{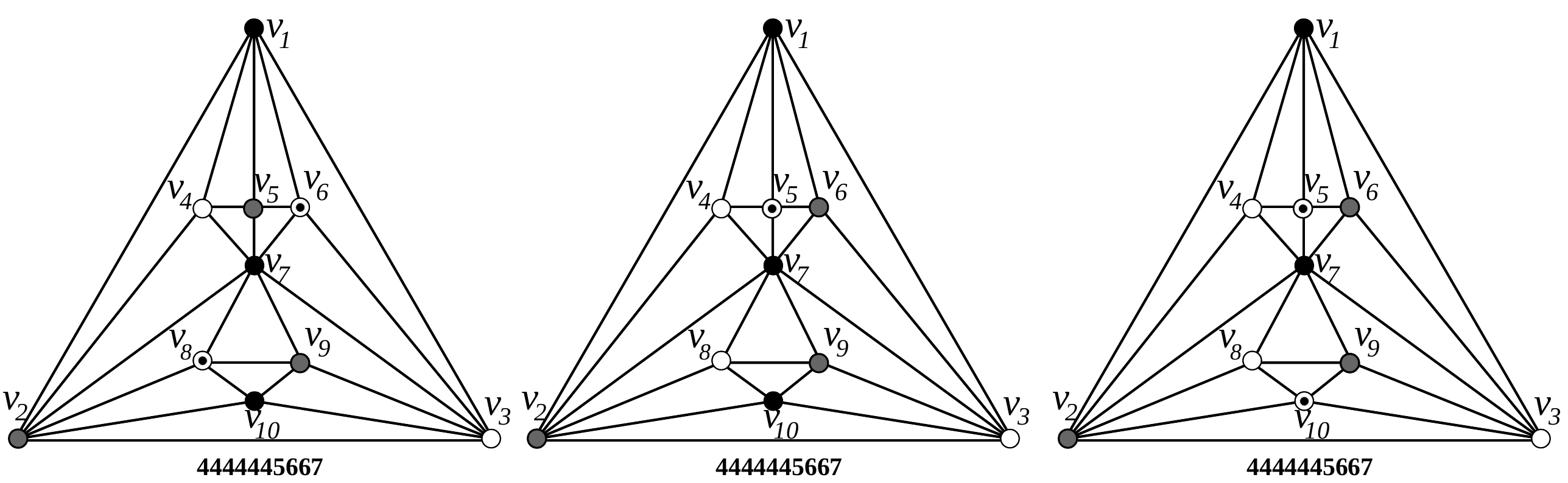}

         \vspace{5mm}
         \includegraphics [width=380pt]{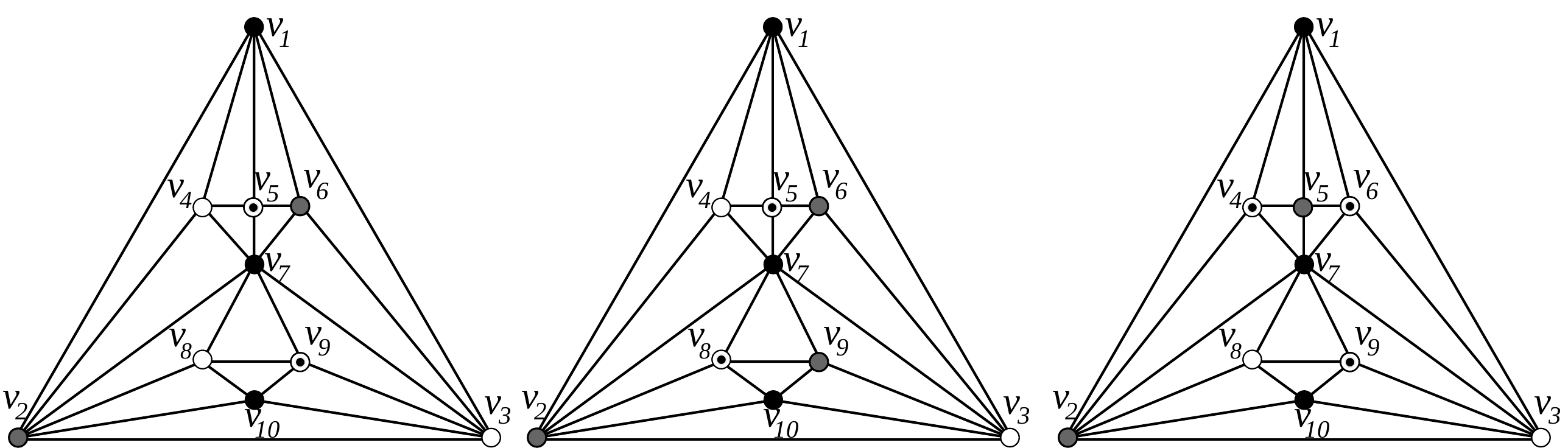}

                  \vspace{5mm}
         \includegraphics [width=380pt]{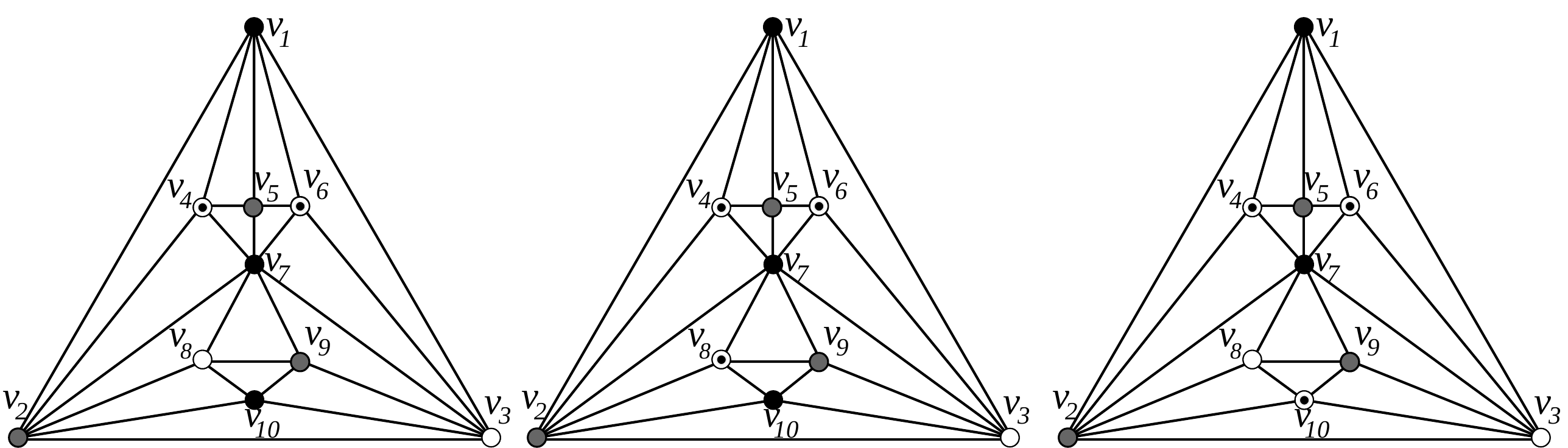}

                  \vspace{5mm}
         \includegraphics [width=380pt]{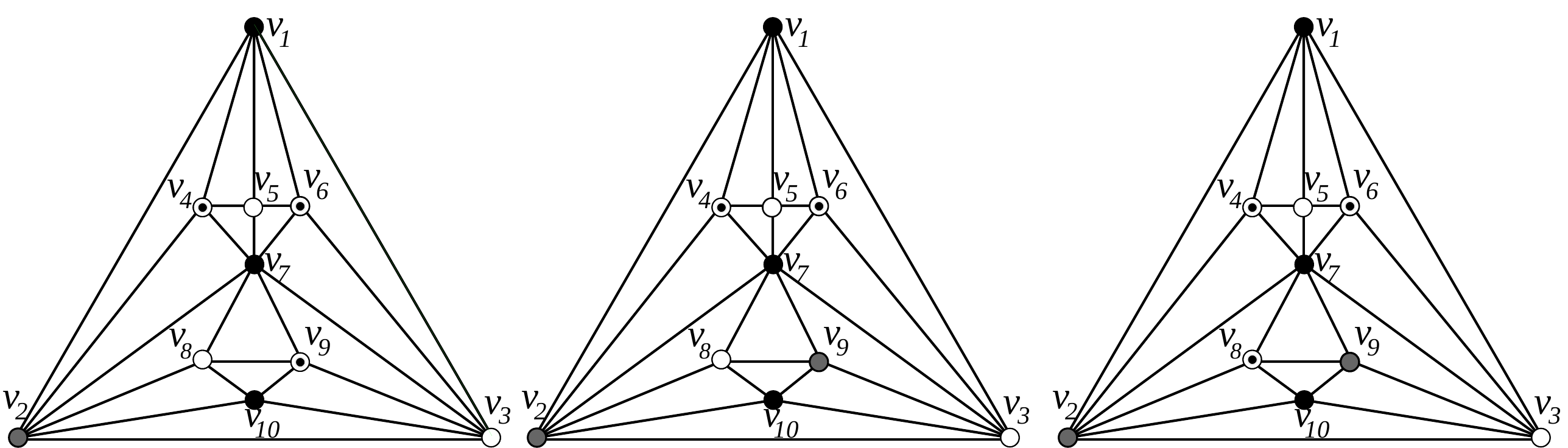}

                  \vspace{5mm}
         \includegraphics [width=380pt]{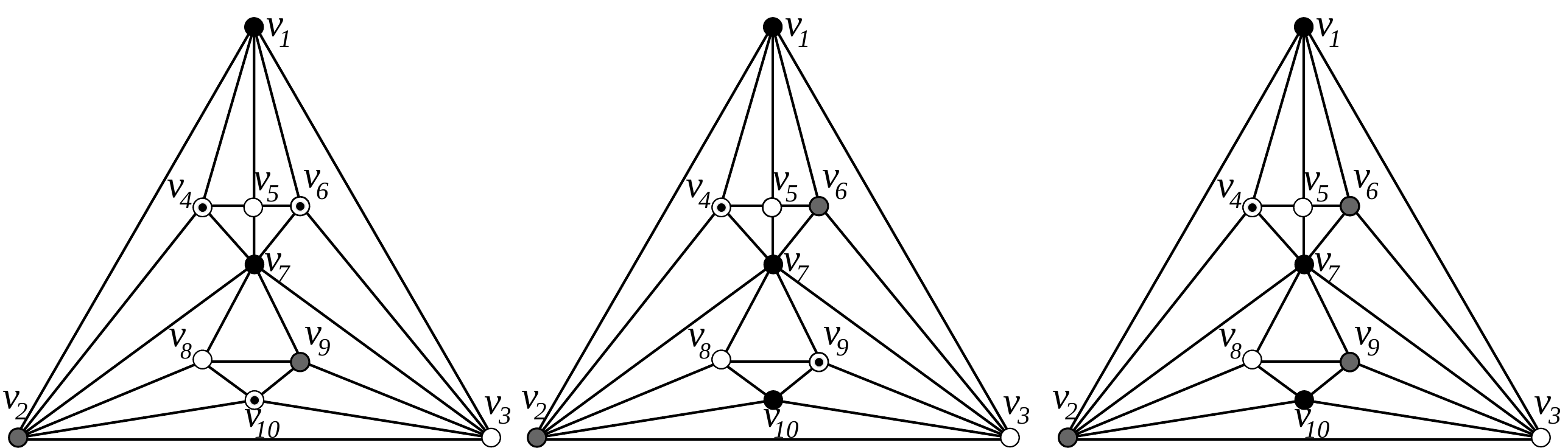}

                  \vspace{5mm}
         \includegraphics [width=320pt]{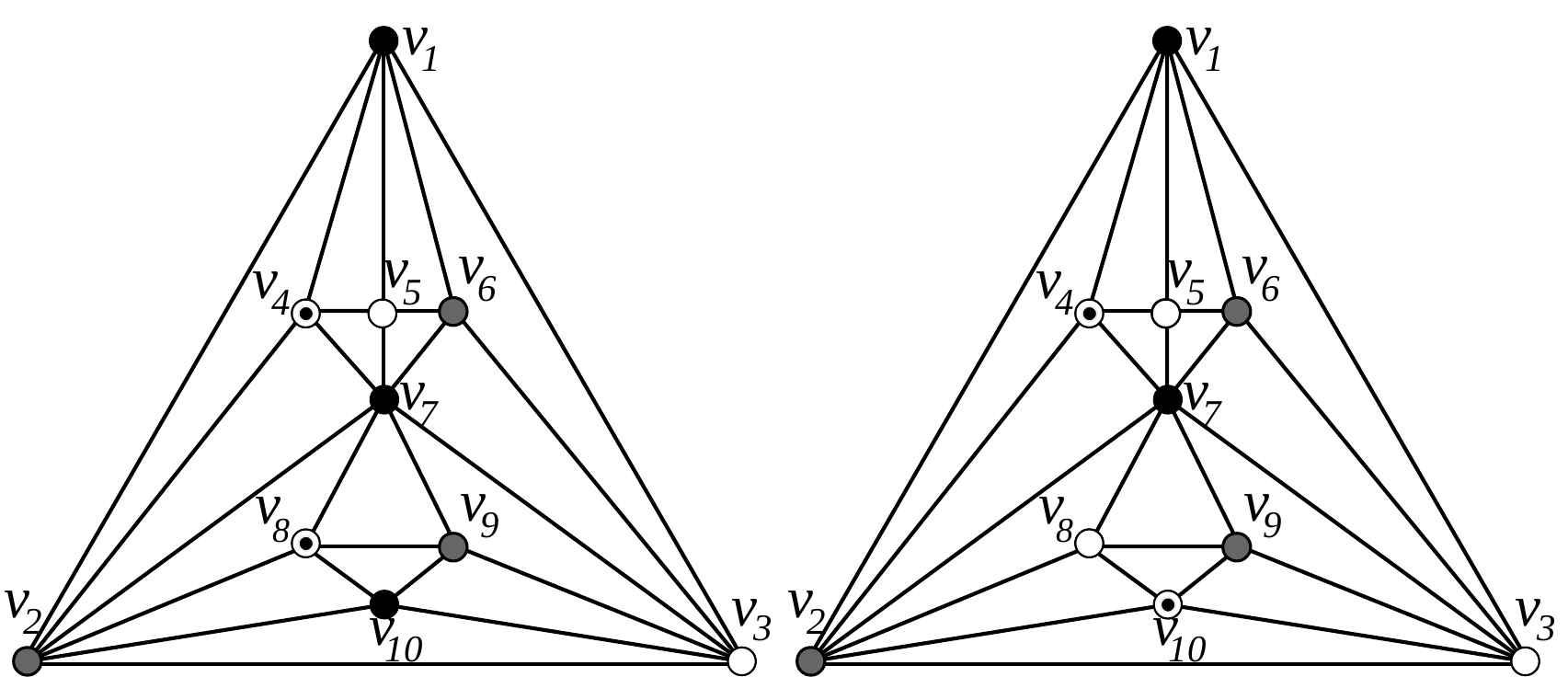}

         \vspace{8mm}
  \end{center}
 5.10 Degree sequence is 4444555566, and it has 14 kinds of different colorings.

  $$ \{\{v_1,v_5,v_{10}\}\{v_2,v_4,v_8\}\{v_3,v_6,v_7\}\{v_9\}\}, \{\{v_1,v_{10}\} \{v_2,v_4,v_8\}\{v_3,v_6,v_7\}\{v_5,v_9\}\}$$
  $$ \{\{v_1,v_6,v_{10}\}\{v_2,v_4,v_8\}\{v_3,v_7\}\{v_5,v_9\}\}, \{\{v_1,v_5,v_{10}\} \{v_2,v_9\}\{v_3,v_6,v_7\}\{v_4,v_8\}\}$$
  $$ \{\{v_1,v_{10}\}\{v_2,v_5,v_9\}\{v_3,v_6,v_7\}\{v_4,v_8\}\}, \{\{v_1,v_6,v_{10}\} \{v_2,v_5,v_9\}\{v_3,v_7\}\{v_4,v_8\}\}$$
  $$ \{\{v_1,v_5,v_{10}\}\{v_2,v_4,v_8\}\{v_3,v_6\}\{v_7,v_9\}\}, \{\{v_1,v_6,v_{10}\} \{v_2,v_4,v_8\}\{v_3,v_5\}\{v_7,v_9\}\}$$
  $$ \{\{v_1,v_5,v_{10}\}\{v_2,v_4,v_8\}\{v_3,v_6\}\{v_7,v_9\}\}, \{\{v_1,v_6,v_{10}\} \{v_2,v_4,v_8\}\{v_3,v_5\}\{v_7,v_9\}\}$$
  $$ \{\{v_1,v_5,v_{10}\}\{v_2,v_9\}\{v_3,v_4,v_8\}\{v_6,v_7\}\}, \{\{v_1,v_{10}\} \{v_2,v_5,v_9\}\{v_3,v_4,v_8\}\{v_6,v_7\}\}$$
  $$ \{\{v_1,v_6,v_{10}\}\{v_2,v_5,v_9\}\{v_3,v_4,v_8\}\{v_7\}\}, \{\{v_1,v_5,v_{10}\} \{v_2,v_6\}\{v_3,v_4,v_8\}\{v_7,v_9\}\}$$
  $$ \{\{v_1,v_6,v_{10}\}\{v_2,v_5\}\{v_3,v_4,v_8\}\{v_7,v_9\}\}, \{\{v_1,v_{8}\} \{v_2,v_5,v_9\}\{v_3,v_6,v_7\}\{v_4,v_{10}\}\}$$
              \begin{center}
        \includegraphics [width=160pt]{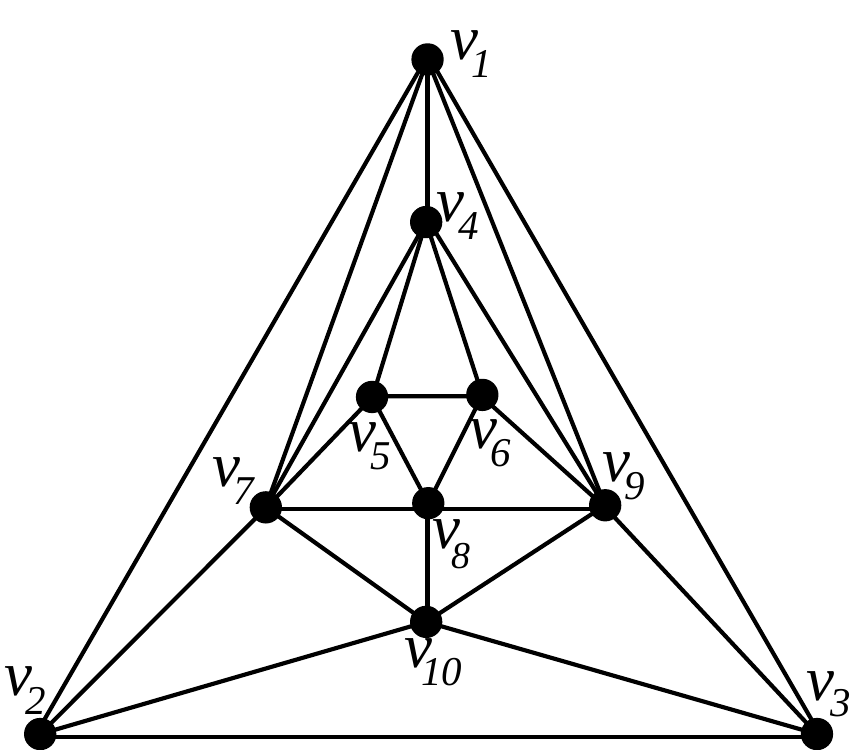}

         \vspace{5mm}
         \includegraphics [width=380pt]{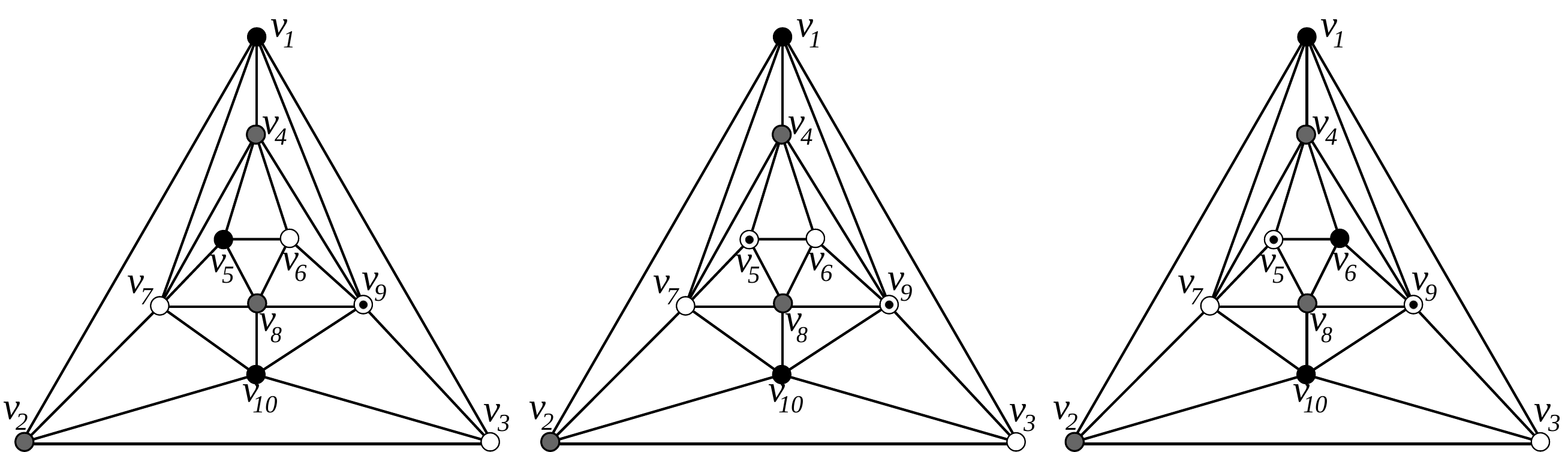}

         \vspace{5mm}
         \includegraphics [width=380pt]{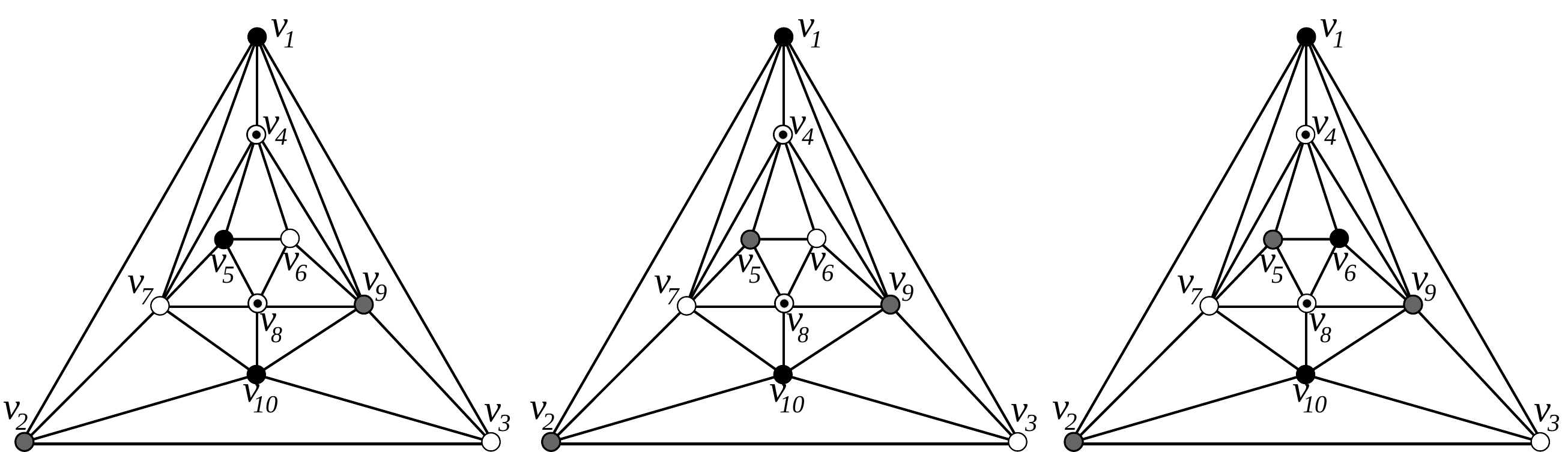}

         \vspace{5mm}
         \includegraphics [width=380pt]{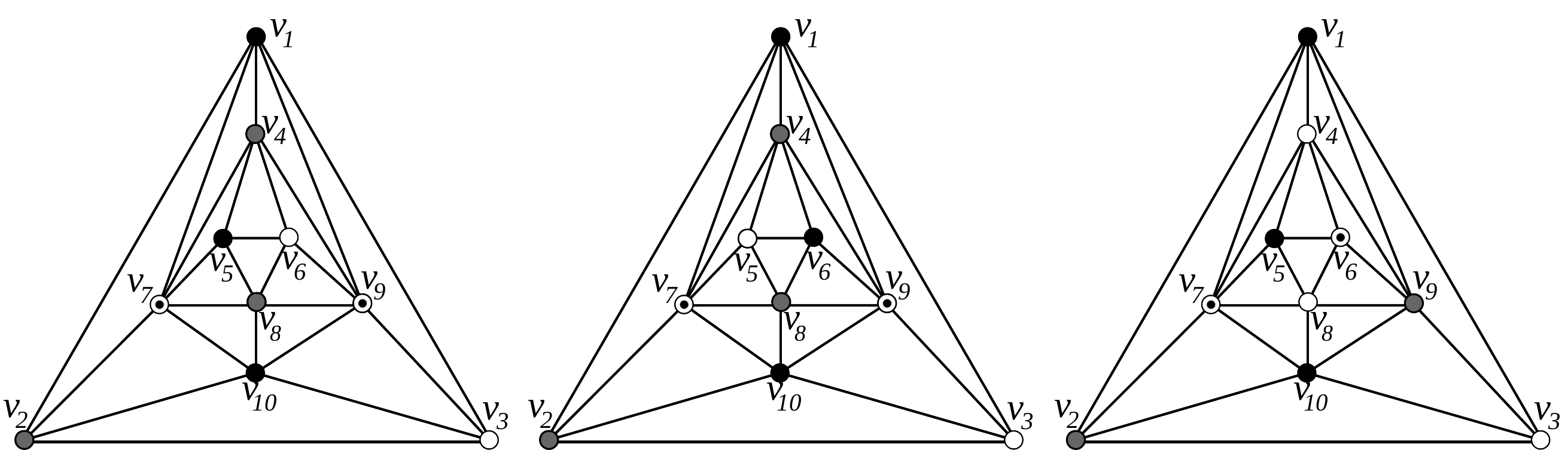}

                  \vspace{5mm}
         \includegraphics [width=380pt]{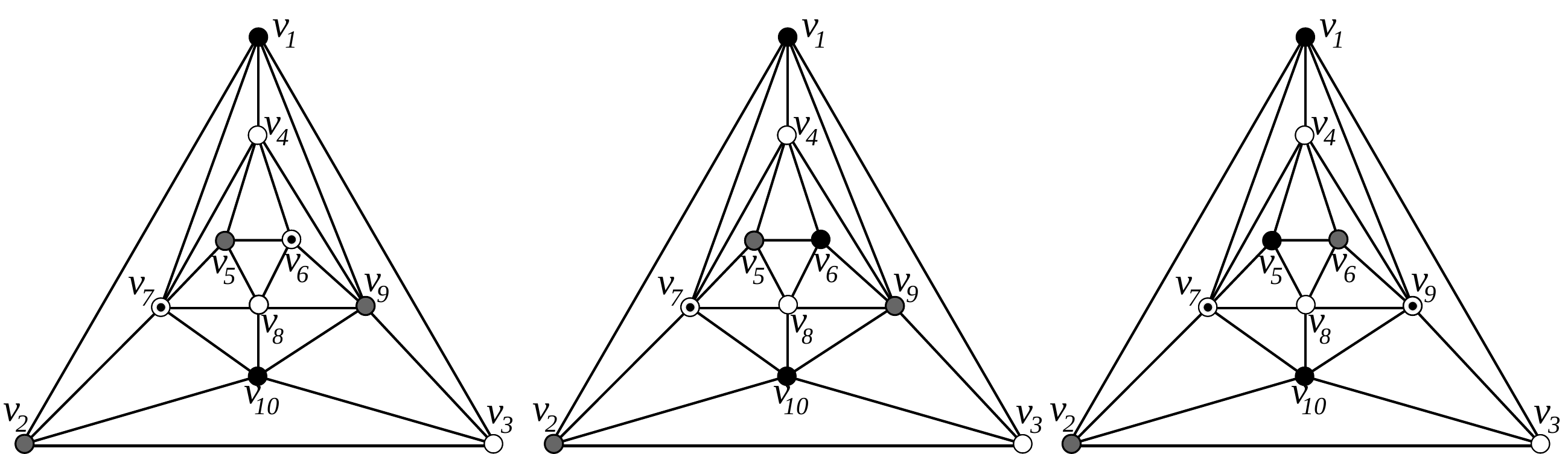}

                  \vspace{5mm}
         \includegraphics [width=320pt]{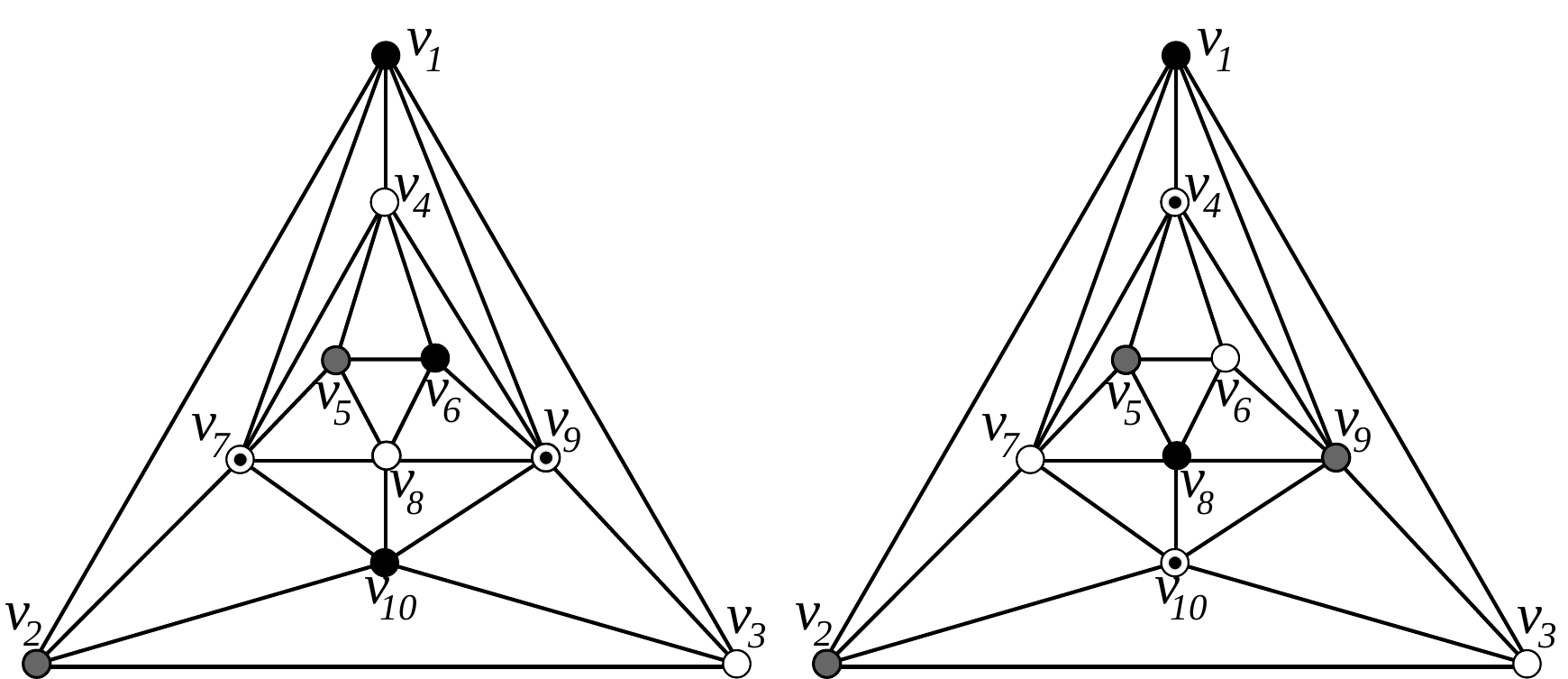}

         \vspace{8mm}
  \end{center}
 5.11 Degree sequence is 4444555566, and it has 13 kinds of different colorings.

  $$ \{\{v_1,v_6,v_{10}\}\{v_2,v_5,v_8\}\{v_3,v_4\}\{v_7,v_9\}\}, \{\{v_1,v_6,v_{10}\} \{v_2,v_5\}\{v_3,v_4,v_8\}\{v_7,v_9\}\}$$
  $$ \{\{v_1,v_6,v_{10}\}\{v_2,v_5,v_8\}\{v_3,v_4,v_9\}\{v_7\}\}, \{\{v_1,v_7\} \{v_2,v_5,v_8\}\{v_3,v_4,v_9\}\{v_6,v_{10}\}\}$$
  $$ \{\{v_1,v_7,v_9\}\{v_2,v_5,v_8\}\{v_3,v_4\}\{v_6,v_{10}\}\}, \{\{v_1,v_7,v_9\} \{v_2,v_5\}\{v_3,v_4,v_8\}\{v_6,v_{10}\}\}$$
  $$ \{\{v_1,v_6,v_{10}\}\{v_2,v_7\}\{v_3,v_4,v_8\}\{v_5,v_9\}\}, \{\{v_1,v_6,v_{10}\} \{v_2,v_7\}\{v_3,v_4,v_9\}\{v_5,v_8\}\}$$
  $$ \{\{v_1,v_6,v_{10}\}\{v_2,v_5\}\{v_3,v_8\}\{v_4,v_7,v_9\}\}, \{\{v_1,v_6,v_{10}\} \{v_2,v_5,v_8\}\{v_3\}\{v_4,v_7,v_9\}\}$$
  $$ \{\{v_1,v_6,v_{10}\}\{v_2,v_5,v_8\}\{v_3,v_9\}\{v_4,v_7\}\}, \{\{v_1,v_7,v_9\} \{v_2,v_5,v_8\}\{v_3,v_6\}\{v_4,v_{10}\}\}$$
  $$ \{\{v_1,v_{10}\}\{v_2,v_5,v_8\}\{v_3,v_6\}\{v_4,v_7,v_9\}\}$$
 \begin{center}
        \includegraphics [width=160pt]{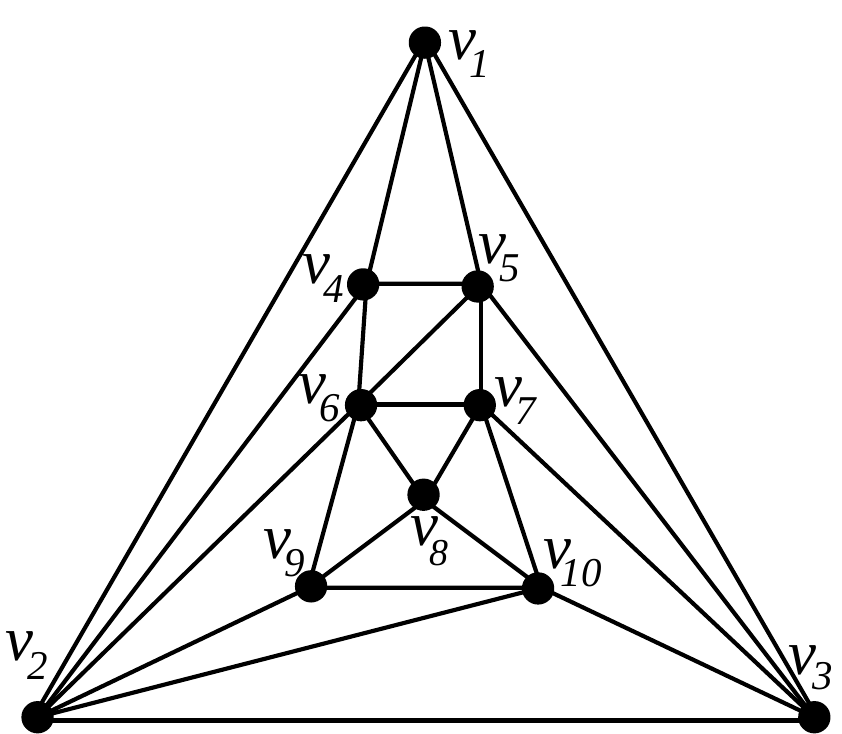}

         \vspace{5mm}
         \includegraphics [width=380pt]{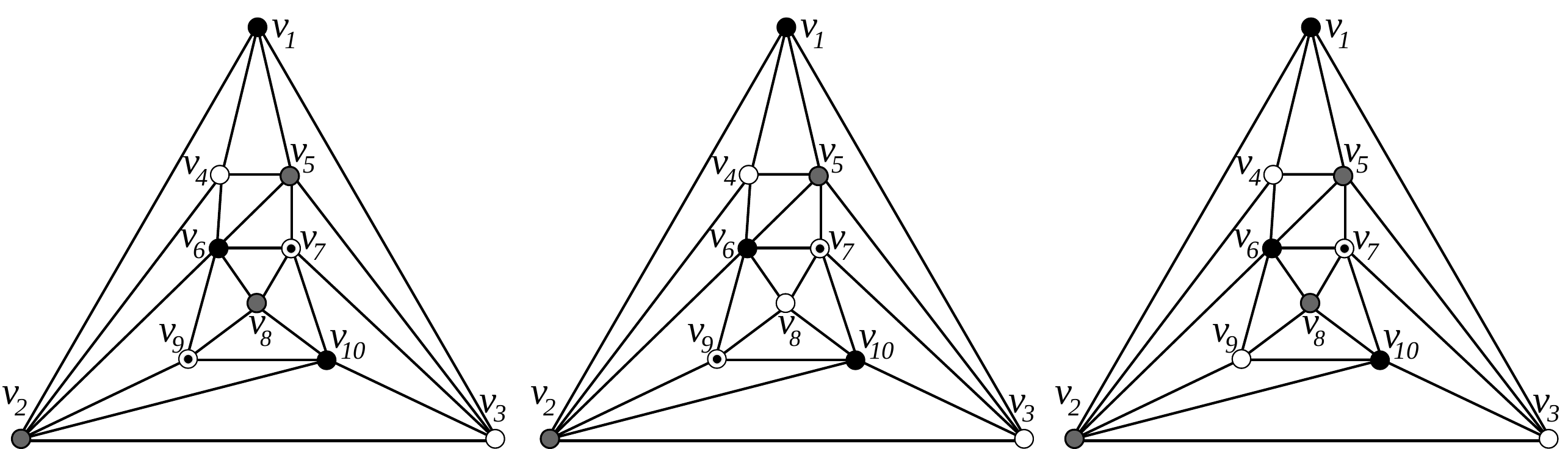}

         \vspace{5mm}
         \includegraphics [width=380pt]{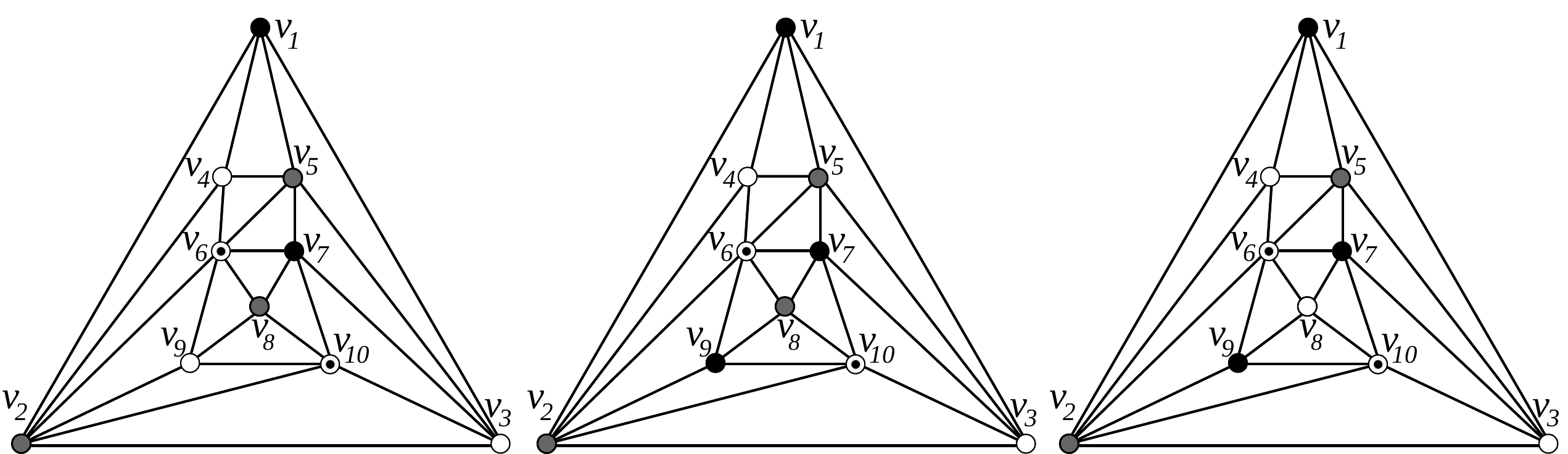}

         \vspace{5mm}
         \includegraphics [width=380pt]{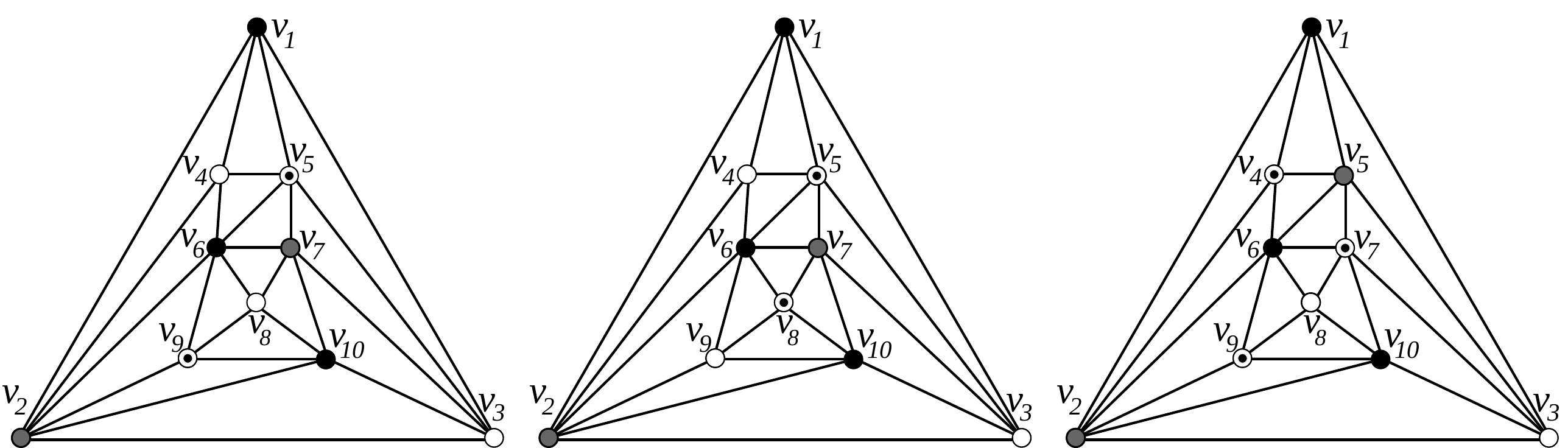}

                  \vspace{5mm}
         \includegraphics [width=380pt]{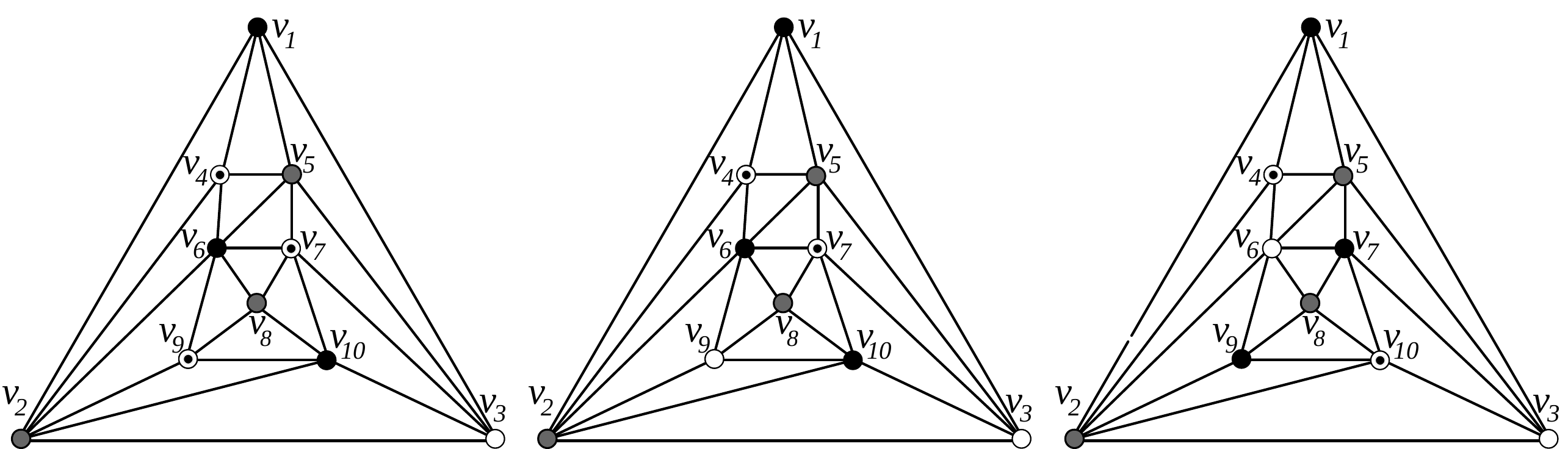}

                  \vspace{5mm}
         \includegraphics [width=160pt]{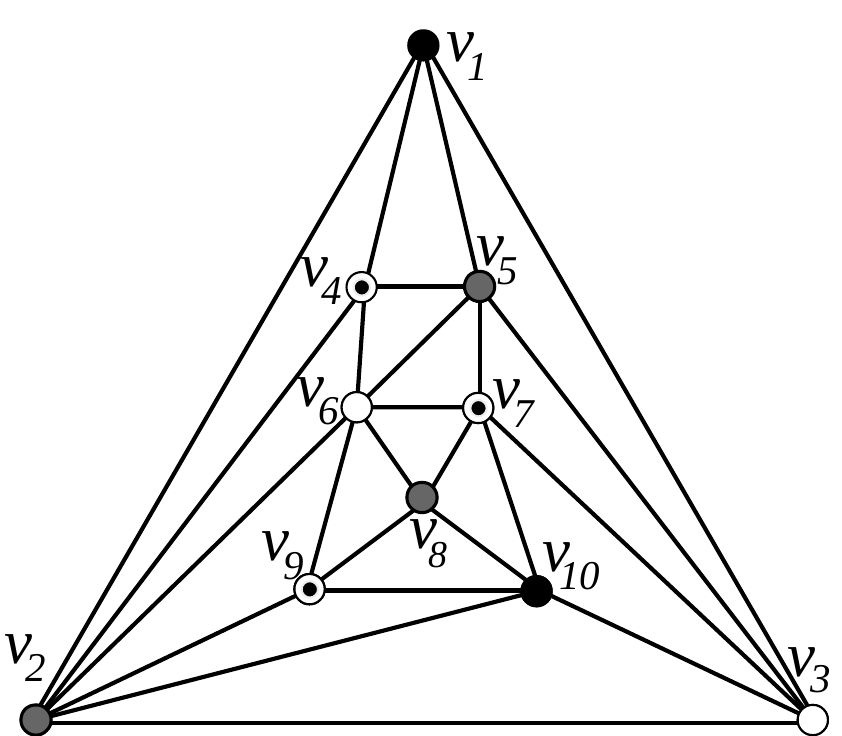}

         \vspace{8mm}
  \end{center}
 5.12 Degree sequence is 4444555566, and it has 5 kinds of different colorings.

  $$ \{\{v_1,v_5,v_{10}\}\{v_2,v_4,v_8\}\{v_3,v_6\}\{v_7,v_9\}\}, \{\{v_1,v_5,v_9\} \{v_2,v_7\}\{v_3,v_4,v_8\}\{v_6,v_{10}\}\}$$
  $$ \{\{v_1,v_5,v_{10}\}\{v_2,v_4,v_8\}\{v_3,v_9\}\{v_6,v_7\}\}, \{\{v_1,v_5,v_{10}\} \{v_2,v_8\}\{v_3,v_4,v_9\}\{v_6,v_7\}\}$$
  $$ \{\{v_1,v_{10}\}\{v_2,v_4,v_8\}\{v_3,v_5,v_9\}\{v_6,v_7\}\}$$

 \begin{center}
        \includegraphics [width=160pt]{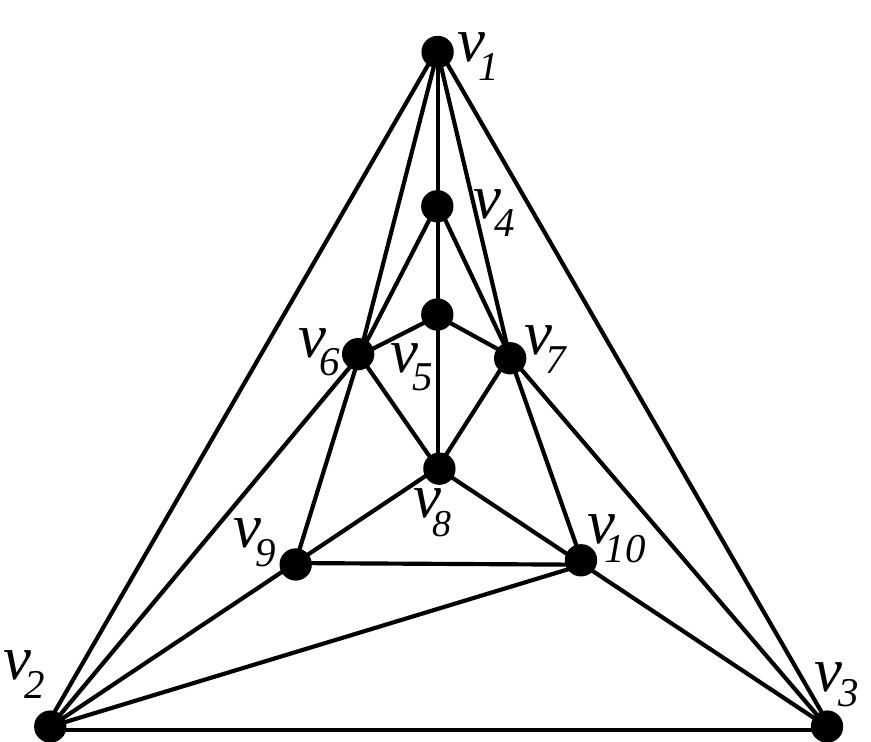}

         \vspace{5mm}
         \includegraphics [width=380pt]{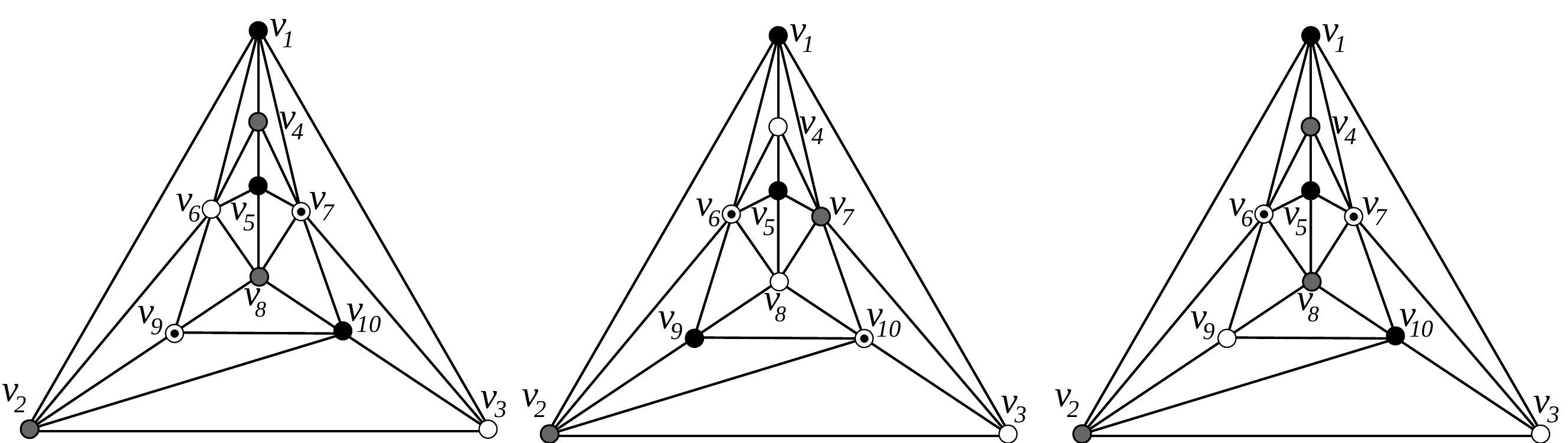}

         \vspace{5mm}
         \includegraphics [width=320pt]{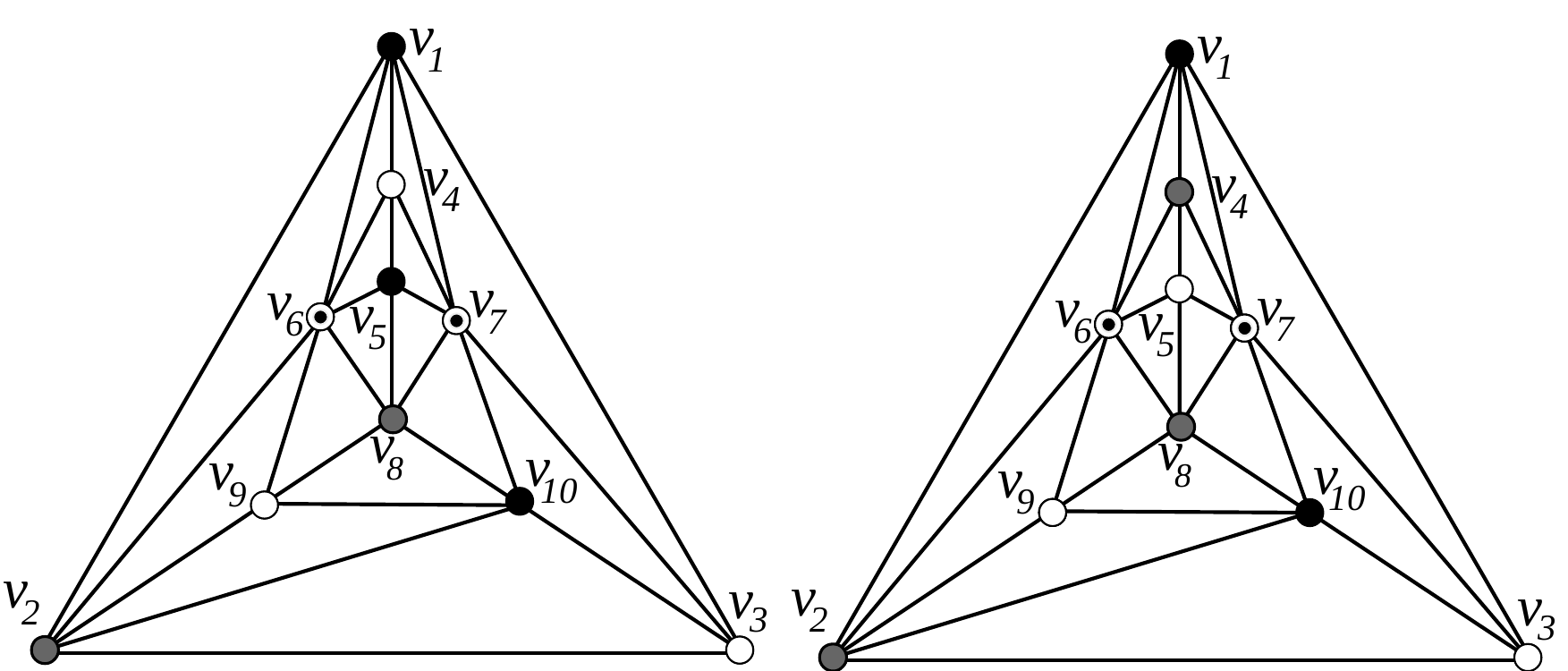}

         \vspace{8mm}
  \end{center}
 5.13 Degree sequence is 4444455577, and it has 15 kinds of different colorings.

  $$ \{\{v_1,v_6,v_8\}\{v_3,v_4\}\{v_2,v_5,v_9\}\{v_7,v_{10}\}\}, \{\{v_1,v_6\} \{v_3,v_4,v_8\}\{v_2,v_5,v_9\}\{v_7,v_{10}\}\}$$
  $$ \{\{v_1,v_6,v_{10}\}\{v_3,v_4,v_8\}\{v_2,v_5,v_9\}\{v_7\}\}, \{\{v_1,v_6,v_9\} \{v_3,v_4,v_8\}\{v_2,v_5\}\{v_7,v_{10}\}\}$$
  $$ \{\{v_1,v_7,v_{10}\}\{v_3,v_4,v_8\}\{v_2,v_5,v_9\}\{v_6\}\}, \{\{v_1,v_7\} \{v_3,v_4,v_8\}\{v_2,v_5,v_9\}\{v_6,v_{10}\}\}$$
  $$ \{\{v_1,v_7,v_{10}\}\{v_3,v_4\}\{v_2,v_5,v_9\}\{v_6,v_8\}\}, \{\{v_1,v_6,v_8\} \{v_3,v_4\}\{v_2,v_9\}\{v_5,v_7,v_{10}\}\}$$
  $$ \{\{v_1,v_6\}\{v_3,v_4,v_8\}\{v_2,v_9\}\{v_5,v_7,v_{10}\}\}, \{\{v_1,v_6,v_{10}\} \{v_3,v_4,v_8\}\{v_2,v_9\}\{v_5,v_7\}\}$$
  $$ \{\{v_1,v_6,v_9\}\{v_3,v_4,v_8\}\{v_2\}\{v_5,v_7,v_{10}\}\}, \{\{v_1,v_6,v_8\} \{v_3\}\{v_2,v_5,v_9\}\{v_4,v_7,v_{10}\}\}$$
  $$ \{\{v_1,v_6\}\{v_3,v_8\}\{v_2,v_5,v_9\}\{v_4,v_7,v_{10}\}\}, \{\{v_1,v_6,v_{10}\} \{v_3,v_8\}\{v_2,v_5,v_9\}\{v_4,v_7\}\}$$
  $$ \{\{v_1,v_6,v_9\}\{v_3,v_8\}\{v_2,v_5\}\{v_4,v_7,v_{10}\}\}$$
 \begin{center}
        \includegraphics [width=160pt]{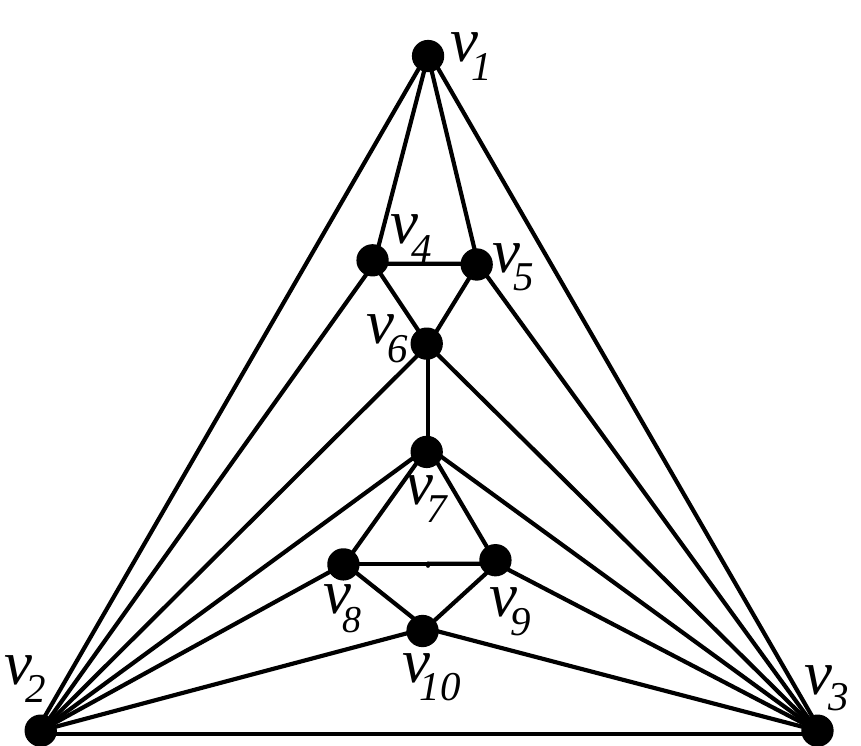}

         \vspace{5mm}
         \includegraphics [width=380pt]{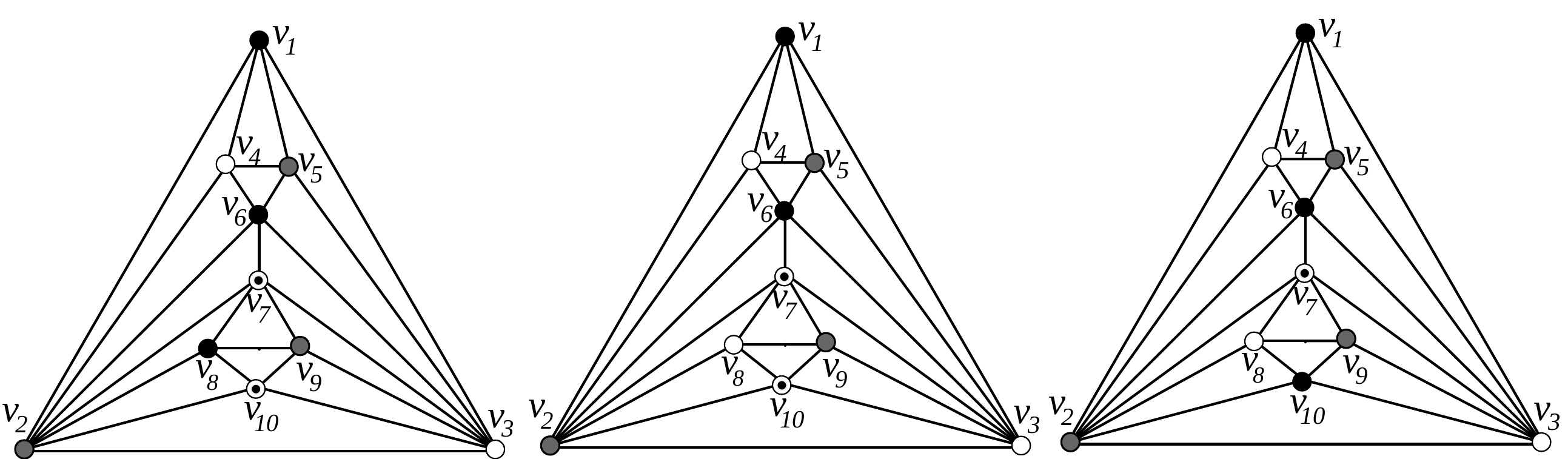}

         \vspace{5mm}
         \includegraphics [width=380pt]{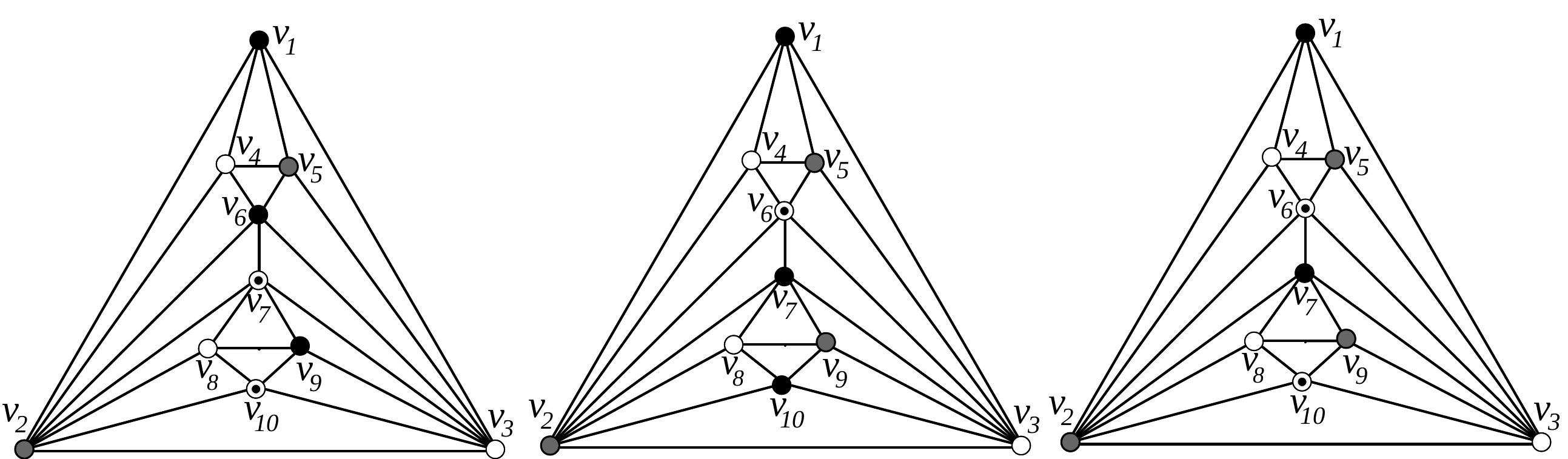}

         \vspace{5mm}
         \includegraphics [width=380pt]{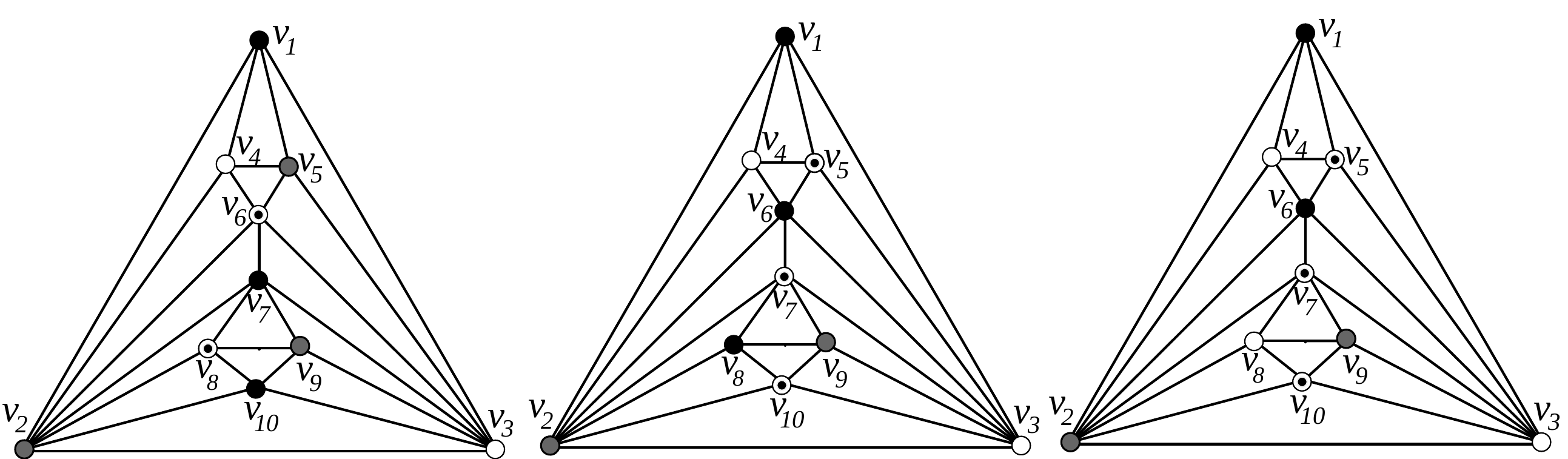}

                  \vspace{5mm}
         \includegraphics [width=380pt]{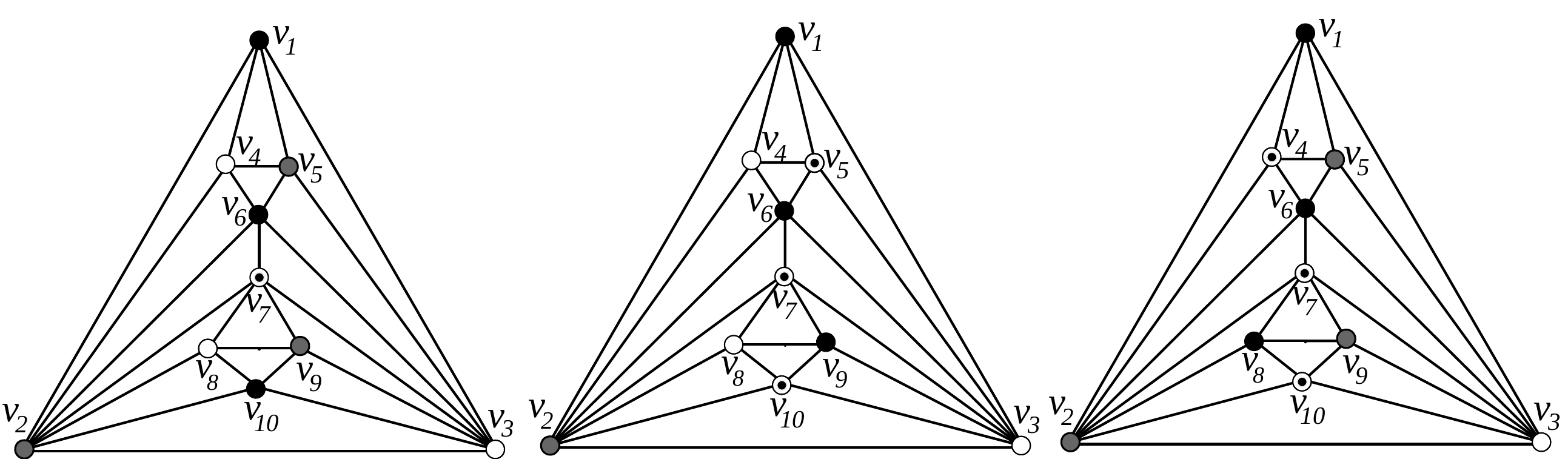}

                  \vspace{5mm}
         \includegraphics [width=380pt]{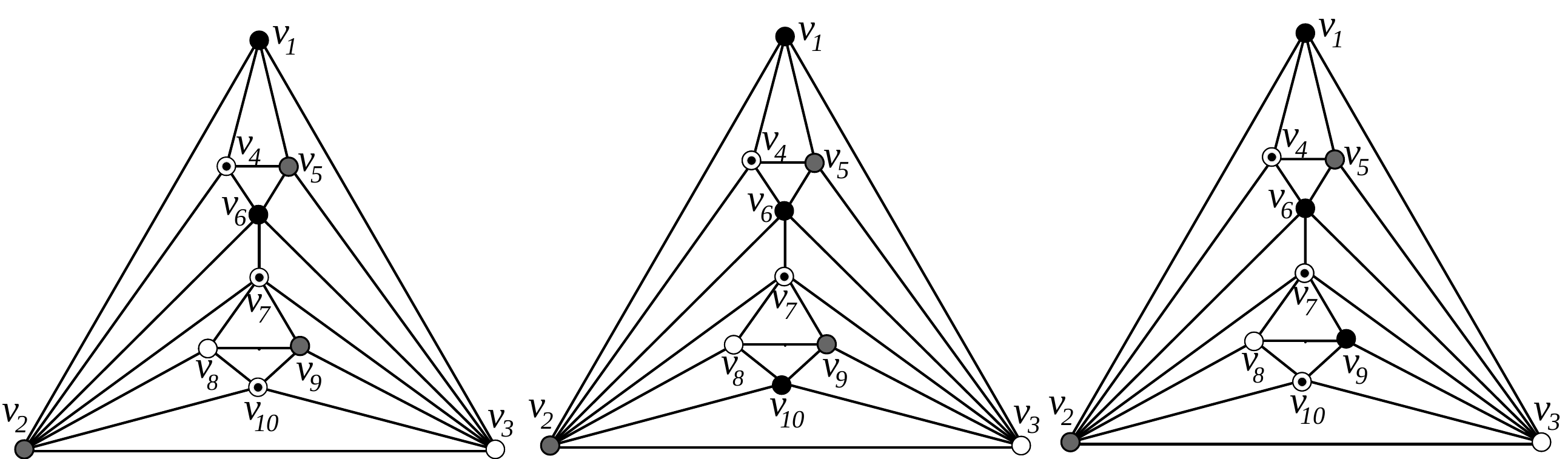}

         \vspace{8mm}
  \end{center}

\newpage


{\footnotesize

\begin{thebibliography}{9}

\setlength{\parskip}{0pt}

\bibitem{AMS2001}
S. Akbari, V. S. Mirrokni and B. S. Sadjad, $K_{r}$-free uniquely
vertex colorable graphs with minimum possible edges,
   \emph{J. Combin. Theory}, Series B, 82 (2001): 316-318.
\bibitem{A1977}
V. A. Aksionov, Chromatically connected vertices in plane graphs,
\emph{Diskret. Analiz}, 31  (1977): 5-16  (in Russian).
\bibitem{AH1977}
K. Appel and W. Haken, The Solution of the Four-Color Map Problem,
\emph{ Sci. Amer.}, 237  (1977): 108-121.
\bibitem{AH1977 (2)}
K. Appel and W. Haken, Every Planar Map is Four-Colorable, II:
Reducibility, \emph{Illinois J. Math.}, 21 (1977): 491-561.
\bibitem{AHK1977}
K. Appel and W. Haken and J. Koch, Every Planar Map is Four
Colorable, I:  Discharging, \emph{Illinois J. Math.}, 21 (1977):
429-490.
\bibitem{BW1978}
L. W. Beineke and R. J. Wilson, Selected Topics in Graph Theory (I),
\emph{Academic Press}, London, New York, San Francisco, 1978.
\bibitem{BLW1976}
N. L. Biggs, E. K. Lloyd and R. J. Wilson, Graph Theory 1736-1936,
\emph{Clarendon Press}, Oxford, 1976.
\bibitem{B1912}
G. D. Birkhoff, A determinantal formula for the number of ways of
coloring a map, \emph{Ann. Of Math}, 14 (1912): 42-46.
\bibitem{B1913}
G. D. Birkhoff, The reducibility of maps, \emph{Amer. J. Math}, 35
(1913): 114-128.
\bibitem{B1946}
G. D. Birkhoff and D. Lewis, Chromatic Polynomials, \emph{Trans.
Amer. Math. Soc.}, 60 (1946):  355-451.
\bibitem{BSV1998}
T. B\"{o}hme, M. Stiebitz, M. Voigt, On uniquely 4-colorable planar
graphs, Preprint No. M 10/98, TU Ilmenau, 1998.
\bibitem{B1978}
B. Bollobas, Uniquely colorable graphs, \emph{J. Comb. Theory},
Series B, 25 (1978): 54-61.
\bibitem{BS1976}
B. Bollobas and N. W. Sauer, Uniquely colourable graphs with large
girth, \emph{Canad. J. Math}, 28 (1976): 1340-1344.
\bibitem{BM2008}
J. A. Bondy and U. S. R. Murty, Graph Theory, \emph{Springer}, 2008.
\bibitem{C1878}
A. Cayley, the solution of a problem which recently achieved some
renown, \emph{Nature}, 18 (1878): 294.
\bibitem{C1879}
A. Cayley, On the colour of maps, \emph{Proc. R. Geogr. Soc}, 1
(1879): 259-261.
\bibitem{CC1993}
C. Chao and Z. Chen, On uniquely 3-colorable graphs, \emph{Discrete
Math.}, 112 (1993): 21-27.
\bibitem{CG1969}
G. Chartrand and D. Geller, Uniquely colorable planar graphs,
\emph{J. Comb. Theory}, 6 (1969): 271-278.
\bibitem{CZ2005}
G. Chartrand and P. Zhang, Introduction to Graph Theory,
\emph{McGraw-Hill Companies Inc}, Singapore, 2005.
\bibitem{CH1968}
F. D. Chartwright and F. Harary, on the coloring of signed graphs,
\emph{Elem. Math.}, 23 (1968): 85-89.
\bibitem{D1980}
D. P. Dailey, Uniqueness of Colorability and Colorability of Planar
4-regular Graphs are NP-complete, \emph{Discrete Mathematics}, 30
(1980): 289-293.
\bibitem{D1980 (2)}
I. G. Dmitriev, Weakly cyclic graphs with integral chromatic number.
\emph{Metody Diskret. Analiz.}, 34 (1980): 3-7.
\bibitem{D1982}
I. G. Dmitriev, Characterization of a class of \emph{k}-trees,
\emph{Metody Diskret. Analiz.}, 38 (1982): 9-18.
\bibitem{DKT2005}
F. M. Dong, K. M. Koh and K. L. Teo, Chromatic Polynomials and
Chromaticity of Graphs, \emph{World Scientific}, Singapore, 2005.
\bibitem{F1975}
S. Fiorini, On the chromatic index of a graph, III:  Uniquely
edge-colorable graphs, \emph{ Quart. J. Math.  (Oxford)},
 (3)26 (1975): 129-140.
\bibitem{FW1977}
S. Fiorini and R. J. Wilson, Edge colouring of graphs,
\emph{Research Notes in Math.}, 16 (1977).
\bibitem{FW1978}
S. Fiorini and R. J. Wilson, Edge colourings of graphs, Selected
Topics in Graph Theory, \emph{Academic Press}, New York  (1978):
103-126.
\bibitem{F1977}
S. Fisk, Geometric coloring theory, \emph{Adv. in Math.}, 24 (1977):
298-340.
\bibitem{F1922}
P. Franklin, The Four Color Problem, \emph{Am. J. Math.}, 44 (1922):
225-236.
\bibitem{F1938}
P. Franklin, Note on the Four Color Problem, \emph{J. Math. and
Phys.}, 16 (1938): 172-184.
\bibitem{F1880}
G. Frederick, Note on the colouring of maps, \emph{Proc. R. Soc.
Edinb}, 10 (1880), 727-728.
\bibitem{FF1998}
R. Fritsch and G. Fritsch , The Four-Color Theorem, \emph{Springer},
1998.
\bibitem{GJS1976}
M. R. Garey, D. S. Johnson and L. J. Stockmeyer, Some simplified
NP-complete graph problems, \emph{Theor. Comput. Sci.}, 1 (1976):
237-267.
\bibitem{G2005}
G. Georges, Formal Proof--The Four Color Theorem, \emph{Notice of
AMS}, 55 (2008):  1382-1393.
\bibitem{GC1967}
T. C. Gleason, Cartwright. D, A note on a matrix criterion for
unique colorability of a signed graph, Psychometrika, 32 (1967):
291-296.

\bibitem{GZ1996}
J. L. Goldwasser and C. Q. Zhang, On the minimal counterexamples to
a conjecture about unique edge-3-coloring, \emph{Congr. Numer.}, 113
(1996): 143-152.  (The result is the same as \textbf{Theorem 3. 1}.
)
\bibitem{GZ2000}
J. L. Goldwasser and C. Q. Zhang, Uniquely edge-colorable graphs and
Snarks, \emph{Graphs and Combinatorics}, 16 (2000): 257-267.
\bibitem{GK1973}
D. Greenwell and H. V. Kronk, Uniquely line-colorable graphs,
\emph{Canad. Math. Bull}, 16  (1973): 525-529.
\bibitem{G1968}
E. J. Grinberg, Plane homogeneous graphs of degree three without
Hamiltonian circuits, \emph{Latvian Math.},  (5)1968: 51-58.
\bibitem{HHR1969}
F. Harary, S. T. Hedetniemi and R. W. Robinson, Uniquely colorable
graphs, \emph{J. Combin. Theory}, 6 (1969): 264-270.
\bibitem{H1969}
F. Haray, Graph Theory, \emph{Addison-Wesley}, Reading, 1969.
\bibitem{H1890}
P. J. Heawood, Map colour theorem, \emph{Q. J. Math. Oxf}, 24
(1890): 332-338.
\bibitem{H1969}
H. Heesch, Untersuchungen zum Vierfarbenproblem,
\emph{Bibliographisches Institut}, Mannheim/Wien/Z\"{u}rich, 1969.
\bibitem{JT1995}
T. R. Jensen and B. Toft, Graph Coloring Problems, \emph{John Wiley
Sons}, New York,  (1995): 48-49.
\bibitem{K1879}
A. B. Kempe, on the geographical problem of the four color,
\emph{Am. J. Math}, 2 (1879): 193-200.
\bibitem{K1879 (2)}
A. B. Kempe, how to colour a map with four colours, \emph{Nature},
21 (1879): 399-400.
\bibitem{K1998}
M. Kriesell, Contractible non-edges in 3-connected graphs, \emph{J
Comb. Theory}, Series B, 4  (1998): 192-201.
\bibitem{M1978}
J. Mayer, Une propri\'{e}t\'{e} des graphes minimaux dans le
probl\`{e}me des quatre couleurs, \emph{Probl\`{e}mes Combinatoires
et Th¨¦orie des Graphes}, Colloques internationaux CNRS No. 260,
Paris, 1978.
\bibitem{MS1977}
L. S. Melnikov and R. Steinberg, one counterexample for two
conjectures on three coloring, \emph{Discrete Mathematics}, 20
(1977): 203-206.
\bibitem{M1974}
V. M\"{u}ller, On colorable critical and uniquely colorable critical
graphs, \emph{Recent Advances in Graph Theory}, Prague,  (1974):
385-386.
\bibitem{M1979}
V. M\"{u}ller, On colorings of graphs without short cycles,
\emph{Discrete Math.}, 26 (1979): 165-176.
\bibitem{N1972}
J. Ne\v{s}et\v{r}il, On critical uniquely colorable graphs,
\emph{Arch. Math.  (Basel)}, 23 (1972): 210-213.
\bibitem{N1973}
J. Ne\v{s}et\v{r}il, On uniquely colorable graphs without short
cycles, \emph{Casopis Pest. Math.}, 98 (1973): 122-125.
\bibitem{O1967}
O. Ore, The four color problem, \emph{Academic Press}, New York,
1967.
\bibitem{OS1970}
O. Ore and J. Stemple, Numerical Calculations on the Four-Color
Problem, J Combin. Theory , 8 (1970): 65-78.
\bibitem{O1974}
L. J. Osterweil, Some classes of uniquely 3-colorable graphs,
\emph{Discrete Math.}, 8 (1974): 59-69.
\bibitem{P1898}
J. Petersen, Sur le th\'{e}or\`{e}me de Tait,
\emph{L'interm\'{e}diaire des Math\'{e}maticiens}, 5 (1898):
225-227.
\bibitem{R1968}
R. C. Read, An introduction to chromatic polynomials, \emph{J.
Combin. Theory}, 4 (1968): 52-71.
\bibitem{R1926}
C. Reynolds, On the problem of coloring maps in four colors,
\emph{Ann. of Math.}, 28 (1926-27): 477-492.
\bibitem{RSST1996}
N. Robertson, D. P. Sanders, P. D. Seymour and R. Thomas, A new
proof of the four colour theorem, \emph{Electron. Res. Announc.
Amer. Math. Soc}, 2 (1996): 17-25.
\bibitem{RSST1997}
N. Robertson, D. P. Sanders, P. D. Seymour and R. Thomas, The four
color theorem, \emph{J. Combin. Theory}, Series B, 70 (1997): 2-44.
\bibitem{T1880}
P. G. Tait, Remarks on the colouring of maps, \emph{Proc. R. Soc.
Edinburgh}, 10 (1880): 501-503.
\bibitem{T1972}
L. S. Thomas, Thirteen colorful variations on Guthrie's four color
conjecture, \emph{America Mathematical monthly}, 79 (1972): 2~43.
\bibitem{TK1986}
L. S. Thomas and P. C. Kainen, The Four-color problem, \emph{Dover
Publications Inc.}, New York, 1986.
\bibitem{T1998}
T. G. Thomas, Unique coloring of planar graphs, Ph. D. thesis,
\emph{Georgia Institute Technology}, 1998.
\bibitem{T1978}
A. G. Thomason, Hamiltonian cycles and uniquely edge colorable
graphs, \emph{Ann. Discrete Math.}, 3 (1978): 259-268.
\bibitem{T1982}
A. G. Thomason, Cubic graphs with three Hamiltonian cycles are not
always uniquely edge Colorable, \emph{Journal of Graph Theory}, 6
(1982): 219-221.
\bibitem{T1946}
W. T. Tutte, On Hamiltonian Circuits, \emph{ J. London Math. Soc.},
21 (1946): 98-101.
\bibitem{T1970}
W. T. Tutte, On chromatic polynomials and the golden ratio, \emph{J.
Combin. Theory}, 9 (1970): 289-296.
\bibitem{T1970 (2)}
W. T. Tutte, More about chromatic polynomials and the golden ratio,
Combinatorial Structures and their Applications  (ed R. K. Guy et
al.), \emph{Gordon and Breach}, New York, 1970: 439-453.
\bibitem{T1974}
W. T. Tutte, Chromatic sums for planar triangulations, V: Special
equations, \emph{Canad. J. Math.}, 26 (1974): 893-907.
\bibitem{T1976}
W. T. Tutte, Hamiltonian circuits, In:  Colloquio Internazionale
sulleteorie Combinatorie, Roma I, \emph{Accademia Nazionale
deiLicei}, 1976: 193-199.
\bibitem{WA1973}
C. C. Wang and E. Artzy, Note on the uniquely colorable graphs,
\emph{Journal of Combin. Theory}, Series B, 15 (1973): 204-206.
\bibitem{W1904}
P. Wernicke, \"{U}ber den kartographischen Vierfarbensatz, \emph{
Math. Ann}, 58 (1904): 479.
\bibitem{W1932}
H. Whitney, On the coloring of graphs, \emph{Ann. of Math.}, 33
(1932): 688-718.
\bibitem{W1940}
C. E. Winn, On the minimum number of polygons in an irreducible map,
American Journal of Mathematics, 62 (1940): 406-416.
\bibitem{X2004}
J. Xu, Recursive formula for calculating the chromatic polynomial of
a graph by vertex deletion, \emph{Acta Mathematica Scientia}, 4
(2004): 577-582.
\bibitem{XL1995}
J. Xu and Z. Liu, The chromatic polynomial between graph and its
complement, \emph{Graph and Combinatorics}, 11 (1995): 337-345.
\bibitem{XW1995}
J. Xu and X. S. Wei, Theorems of uniquely k-colorable graphs,
\emph{Journal of Shaanxi Normal University  (Natural Science
Edition)}, 23 (1995): 59-62.
\bibitem{X1990}
S. J. Xu, The size of uniquely colorable graphs, \emph{J. Combin.
Theory}, Series B, 50 (1990): 319-320.
\bibitem{Z1995}
C. Q. Zhang, Hamiltonian weights and unique edge-3-colorings of
cubic graphs, \emph{Journal of Graph Theory}, 20 (1995): 91-99.
\bibitem{Z1949}
A. A. Zykov, On some properties of linear complexes, \emph{Mat. Sb.
}, 24 (1949): 163-168 (in. Russian); English translation in
\emph{Amer. Math. Soc. Tran.}, 79 (1952).
\end{thebibliography}

}
\end{document}